\documentclass{article} 
\usepackage{iclr2025_conference,times}
\usepackage{hyperref}
\usepackage{url}

\title{ADMM for Nonconvex Optimization under Minimal Continuity Assumption}

\author{Ganzhao Yuan\\
Peng Cheng Laboratory, China\\
\texttt{yuangzh@pcl.ac.cn}}

\iclrfinalcopy 
\usepackage{graphicx}
\usepackage{epstopdf}
\usepackage{enumitem,amsfonts,amssymb,mathtools,amsthm,bm,pifont,bbm,listings,xcolor,threeparttable,xcolor}
\usepackage{ulem,graphicx,subcaption,hyperref,comment,amsmath}
\usepackage[ruled,lined]{algorithm2e}\usepackage{algpseudocode}
\usepackage{multirow}
\usepackage{natbib}
\definecolor{borange}{cmyk}{0.8, 0, 0.8, 0.2} 

\newtheorem{theorem}{Theorem}[section]

\newtheorem{lemma}[theorem]{Lemma}

\newtheorem{definition}[theorem]{Definition}
\newtheorem{assumption}[theorem]{Assumption}
\newtheorem{remark}[theorem]{Remark}

\def\a{\mathbf{a}}\def\b{\mathbf{b}}\def\d{\mathbf{d}}\def\g{\mathbf{g}}\def\p{\mathbf{p}}\def\u{\mathbf{u}}\def\v{\mathbf{v}}\def\w{\mathbf{w}}\def\x{\mathbf{x}}\def\y{\mathbf{y}}\def\z{\mathbf{z}}\def\A{\mathbf{A}}\def\D{\mathbf{D}}\def\G{\mathbf{G}}\def\H{\mathbf{H}}\def\I{\mathbf{I}}\def\O{\mathbf{O}}\def\V{\mathbf{V}}\def\Y{\mathbf{Y}}\def\v{\mathbf{v}}\def\p{\mathbf{p}}

\def\LL{\mathcal{L}}\def\OO{\mathcal{O}}

\def\AAA{\mathbb{A}}
\def\UUU{\mathbb{U}}
\def\BBB{\mathbb{B}}

\def\trans{^\mathsf{T}}\def\Rn{\mathbb{R}}  \def\zero{\mathbf{0}}  
 \def\trans {^\mathsf{T}} \def\fro {\mathsf{F}}

\DeclareMathOperator{\st}{s.t.}
\DeclareMathOperator{\mat}{mat}

\DeclareMathOperator{\Prox}{\operatorname{Prox}}

\newcommand{\bfit}[1]{\textit{\textbf{#1}}}
\newcommand{\beq}{\begin{eqnarray}}
\newcommand{\eeq}{\end{eqnarray}}
\newcommand{\beqq}{\begin{equation}}
\newcommand{\eeqq}{\end{equation}}
\newcommand{\bel}{\begin{align}}
\newcommand{\eel}{\end{align}}
\newcommand{\la}{\langle}
\newcommand{\ra}{\rangle}
\newcommand{\noi}{\noindent}
\newcommand{\nn}{\nonumber}
\def\noi{\noindent}
\def\nn{\nonumber}
\def\la{\langle}
\def\ra{\rangle}
\def\diag{{\rm{diag}}}
\def\Diag{{\rm{Diag}}}
\def\dist{{\rm{dist}}}
\newcommand{\step}[1]{\text{\ding{\numexpr#1+171\relax}}}

\renewcommand{\cite}{\citep}

\def\dist{{\rm{dist}}}

\def\diag{\mathrm{diag}}
\def\Diag{\mathrm{Diag}}
\def\vec{\mathrm{vec}}

\def\trans{^\mathsf{T}}

\def\vec{\mathrm{vec}}

\def\fro{\mathsf{F}}
\def\mat{\mathrm{mat}}

\def\lambdas{\boldsymbol{\lambda}}

\def\thetas{\boldsymbol{\theta}}

\def\ts{\textstyle}

\renewcommand{\frac}[2]{\tfrac{#1}{#2}}

\def\Rn{\mathbb{R}}

\usepackage{dsfont}

\usepackage{listings}\definecolor{backcolour}{rgb}{0.95,0.95,0.92}\definecolor{codegreen}{rgb}{0,0.6,0}\lstdefinestyle{myStyle}{backgroundcolor=\color{backcolour},commentstyle=\color{codegreen},basicstyle=\ttfamily\footnotesize,breakatwhitespace=false,breaklines=true,keepspaces=true,numbers=left,numbersep=5pt,showspaces=false,showstringspaces=false,showtabs=false,tabsize=2,frame = single,numbers=right,}\usepackage{caption} 

\def\lambdas{\boldsymbol{\lambda}}

\def\thetas{\boldsymbol{\theta}}

\def\ts{\textstyle}

\renewcommand{\frac}[2]{\tfrac{#1}{#2}}

\def\Rn{\mathbb{R}}

\usepackage{multirow}

\def\O{\mathcal{O}}

\def\b{{\mathbf{b}}}
\def\p{{\mathbf{p}}}
\def\g{{\mathbf{g}}}

\def\A{{\mathbf{A}}}
\def\x{{\mathbf{x}}}
\def\y{{\mathbf{y}}}
\def\z{{\mathbf{z}}}
\def\d{{\mathbf{d}}}
\def\b{{\mathbf{b}}}

\def\cb{{\mathbf{c}}}
\def\v{{\mathbf{v}}}
\def\V{{\mathbf{V}}}
\def\I{{\mathbf{I}}}
\def\D{{\mathbf{D}}}
\def\u{{\mathbbm{u}}}
\def\H{{\mathbf{Q}}}
\def\w{{\mathbbm{w}}}
\def\ub{{\mathbf{u}}}

\def\oo{L}

\def\gg{{\ddot{\mathbf{g}}}}

\def\LL{{\sf L}}

\def\BIBI{\mathbb{BI}}
\def\SUSU{\mathbb{SU}}

\def\lambdaDown{\underline{\lambda}}
\def\lambdaUp{\overline{\lambda}}

\def\ADown{\underline{\textup{A}}}
\def\AUp{\overline{\textup{A}}}
\def\Crit{{\rm{Crit}}}

\usepackage{booktabs}

\begin{document}

\maketitle

\begin{abstract}
This paper introduces a novel approach to solving multi-block nonconvex composite optimization problems through a proximal linearized Alternating Direction Method of Multipliers (ADMM). This method incorporates an Increasing Penalization and Decreasing Smoothing (IPDS) strategy. Distinguishing itself from existing ADMM-style algorithms, our approach (denoted IPDS-ADMM) imposes a less stringent condition, specifically requiring continuity in just one block of the objective function. IPDS-ADMM requires that the penalty increases and the smoothing parameter decreases, both at a controlled pace. When the associated linear operator is bijective, IPDS-ADMM uses an over-relaxation stepsize for faster convergence; however, when the linear operator is surjective, IPDS-ADMM uses an under-relaxation stepsize for global convergence. We devise a novel potential function to facilitate our convergence analysis and prove an oracle complexity $\mathcal{O}(\epsilon^{-3})$ to achieve an $\epsilon$-approximate critical point. To the best of our knowledge, this is the first complexity result for using ADMM to solve this class of nonsmooth nonconvex problems. Finally, some experiments on the sparse PCA problem are conducted to demonstrate the effectiveness of our approach. \footnote{Future versions of this paper can be found at \url{https://arxiv.org/abs/2405.03233}.}


\end{abstract}

\section{Introduction}

We consider the following multi-block nonconvex nonsmooth composite optimization problem:
\beq\label{eq:main}
\min_{\x_{1},\x_{2},\ldots,\x_{n} }\,\sum_{i=1}^n [ f_i(\x_{i})  + h_i(\x_{i}) ] ,\,s.t.\,[\sum_{i=1}^n \A_{i} \x_{i}] = \b,
\eeq
\noi where $\b\in \Rn^{m\times 1}$, $\A_i \in \Rn^{m\times \mathbf{d}_i}$, $\x_i \in \Rn^{\mathbf{d}_i \times 1}$, and $i\in [n]\triangleq\{1,2,\ldots,n\}$. We assume $f_i(\cdot): \Rn^{\d_i\times 1} \mapsto (-\infty,\infty)$ is differentiable and potentially nonconvex for all $i\in[n]$. The function $h_i(\cdot): \Rn^{\d_i\times 1} \mapsto (-\infty,\infty]$ is assumed to be closed, proper, lower semi-continuous, and potentially nonsmooth. While $h_n(\cdot)$ is convex, we do not require convexity for $h_i(\cdot)$ where $i\in [n-1]$. Additionally, we assume the nonconvex proximal operator of $h_i(\cdot)$ is easy to compute for all $i \in [n]$.

Problem (\ref{eq:main}) has a wide range of applications in machine learning. The function $f_i(\cdot)$ plays a crucial role in handling empirical loss, including neural network activation functions \cite{LiuLC22,ZengLYZ21,WangYCZ19,HuangCH19}. Incorporating multiple nonsmooth regularization terms $h_i(\cdot)$ enables diverse prior information integration, including structured sparsity, low-rank, binary, orthogonality, and non-negativity constraints, enhancing regularization model accuracy. These capabilities extend to various applications such as sparse PCA, overlapping group Lasso, graph-guided fused Lasso, and phase retrieval.




$\blacktriangleright$ \textbf{ADMM Literature}. The Alternating Direction Method of Multipliers (ADMM) is a versatile optimization tool suitable for solving composite constrained problems as in Problem (\ref{eq:main}), which pose challenges for other standard optimization methods, such as the accelerated proximal gradient method \cite{Nesterov03} and the augmented Lagrangian method \cite{zeng2022moreau,lu2012augmented,zhu2020first}. The standard ADMM was initially introduced in \cite{gabay1976dual}, and its complexity analysis for the convex settings was first conducted in \cite{HeY12,monteiro2013iteration}. Since then, numerous papers have explored the iteration complexity of ADMM in diverse settings. These settings include acceleration through multi-step updates \cite{pock2016inertial,li2016majorized,ouyang2015accelerated,shen2017gsos,franca2018admm,hien2021framework,tran2018non}, asynchronous updates \cite{zhang2014asynchronous}, Jacobi updates \cite{deng2017parallel}, non-Euclidean proximal updates \cite{gonccalves2017improved}, and extensions to handle more specific or general functions such as strongly convex functions \cite{nishihara2015general,lin2014convergence,ouyang2015accelerated}, nonlinear constrained functions \cite{Lin2022}, multi-block composite functions \cite{lin2015global,NIPS2017Xu}, finite-sum functions \cite{bian2021stochastic, liu2020accelerated}, and expected-value functions \cite{HuangCH19}. For a comprehensive overview, refer to the book \cite{lin2022alternating}.


$\blacktriangleright$ \textbf{Nonconvex ADMM}. Compared to the classical subgradient methods \cite{li2021weakly,davis2019stochastic} and Smoothing Proximal Gradient Methods (SPGM) \cite{BohmW21,yuan2024smoothing}, designed for general nonconvex optimization, ADMM-type methods potentially offer faster convergence, enhanced parallelization, and greater numerical stability. However, the convergence analysis of the nonconvex ADMM is challenging due to the absence of Fej\'{e}r monotonicity in iterations. In the past decade, significant research has focused on exploring various nonconvex ADMM variants \cite{li2015global,hong2016convergence, yang2017alternating}. \cite{li2015global} establishes the convergence of a class of nonconvex problems when a specific potential function associated with the augmented Lagrangian satisfies the Kurdyka-{\L}ojasiewicz inequality. \cite{yang2017alternating} analyzes ADMM variants for solving low-rank and sparse optimization problems. \cite{hong2016convergence} investigates ADMM variants for nonconvex consensus and sharing problems. Some researchers have examined ADMM variants under weaker conditions, such as restricted weak convexity \cite{wang2019global}, restricted strong convexity \cite{barber2024convergence}, and the Hoffman error bound \cite{zhang2020proximal}. However, existing methods typically assume that at least one block of the objective function is smooth. In contrast, our approach requires only continuity in one block, which is achieved through an Increasing Penalization and Decreasing Smoothing (IPDS) strategy.


$\blacktriangleright$ \textbf{Over-Relaxed and Under-Relaxed ADMM}. Prior studies have analyzed ADMM using either under-relaxation stepsizes $\sigma \in (0,1)$, or over-relaxation stepsizes $\sigma \in [1,2)$, for updating the dual variable. This contrasts with earlier approaches that employed fixed values, such as 1 or the golden ratio $(\sqrt{5}+1)/2$. In nonconvex settings, most existing works require that the associated matrix of the problem be bijective \cite{goncalves2017convergence,yang2017alternating,yashtini2020convergence,yashtini2021multi,boct2020proximal}. However, the work of \cite{Boct2019SIOPT} demonstrates that ADMM can still be applied when the associated matrix is surjective, provided that an under-relaxation stepsize is employed. Inspired by these findings, our work shows that when the associated linear operator is bijective, IPDS-ADMM uses an over-relaxation stepsize for faster convergence. In contrast, when the linear operator is surjective, we employ under-relaxation stepsizes to achieve global convergence.

\begin{table*}[!tb]
  \centering
  \caption{Comparison of existing ADMM approaches for solving the nonconvex problem in Problem (\ref{eq:main}). \texttt{CVX}: convex. \texttt{NC}: nonconvex.  \texttt{LCONT}: Lipschitz continuous. \texttt{WC}: weakly convex. \texttt{RWC}: restricted weakly convex. $\mathbb{F}$: the constraint set is non-empty. $\mathbb{I}$: $\A_n$ is identity. $\mathbb{SU}$: $\A_n$ is surjective with $\lambdas_{\min}(\A_n\A_n\trans)>0$. $\mathbb{IN}$: $\A_n$ is injective with $\lambdas_{\min}(\A_n\trans\A_n)>0$. $\mathbb{BI}$: $\A_n$ is bijective (both surjective and injective). $\mathbb{IM}$: $\rm{Im}([\A_1,\A_2,\ldots,\A_{n-1}]) \subseteq \rm{Im}(\A_n)$ with \rm{Im} being the image of the matrix.}
  \scalebox{0.73}{\begin{tabular}{|p{3.3cm}|p{1.2cm}|p{7.1cm}|p{1.7cm}|p{1.5cm}|p{1.6cm}|}
    \hline
    \multirow{2}{*}{\textbf{\footnotesize Reference}} &
    \multicolumn{3}{c|}{\textbf{\footnotesize Optimization Problems and Main Assumptions}}  & \multirow{2}{*}{\textbf{\footnotesize Complexity}}  & \multirow{2}{*}{\textbf{\footnotesize Parameter $\sigma$}}  \\
    \cline{2-4}
    & \textbf{ Blocks } & \textbf{\footnotesize Functions $f_i(\cdot)$ and $h_i(\cdot)$$^{a}$ } & \textbf{\footnotesize Matrices $\A_i$} &        &   \\
   \hline
   \cite{HeY12} & $n=2$ &  \texttt{CVX}: $f_i,h_i,\forall i \in [2]$ & \text{ $\mathbb{F}$ } & $\O(\epsilon^{-2})$ $^{b}$  &     $\sigma=1$ \\
   \hline
   \cite{li2015global} & $n=2$ &  \texttt{NC}: $h_1,f_2$; $f_1=0$; \textcolor{borange}{$\bm{h}_{\bm{2}}\bm{=0}$} & $\mathbb{SU}$  & $\O(\epsilon^{-2})$  &   $\sigma=1$\\
   \hline
      \cite{yang2017alternating} $^{c}$ & $n=3$ & \texttt{CVX}: $h_1,f_3$; \texttt{NC}: $h_2$; $f_1=f_2=0$; \textcolor{borange}{$\bm{h}_{\bm{3}}\bm{=0}$} & $\mathbb{I}$ & $\O(\epsilon^{-2})$  &   $\sigma\in [1,2)$ \\
   \hline
  \cite{yashtini2020convergence}   & $n=2$ & \texttt{NC}: $f_{[1,2]},h_{[1,2]}$; \textcolor{borange}{$\bm{h}_{\bm{2}}\bm{=0}$}  & $\mathbb{BI}$ & $\O(\epsilon^{-2})$&  $\sigma\in (0,1)$ \\
\hline
   \cite{yashtini2021multi} & $n\geq 2$& \texttt{WC}: $f_{[1,n-1]}$; $h_{[1,n-1]}=0$; \textcolor{borange}{$\bm{h}_{\bm{n}}\bm{=0}$}  & $\mathbb{BI}$, $\mathbb{IM}$  & $\O(\epsilon^{-2})$  &   $\sigma\in (0,1)$  \\
  \hline
  \cite{wang2019global} & $n\geq 2$ & \texttt{RWC}: $h_{[1,n-1]}$; \textcolor{borange}{$\bm{h}_{\bm{n}}\bm{=0}$} &   $\mathbb{IN}$, $\mathbb{IM}$  &  $\O(\epsilon^{-2})$ &  $\sigma=1$ \\
   \hline
\cite{boct2020proximal} & $n=2$ & \texttt{NC}: $h_{1},f_{[1,2]}$; $f_1=0$; \textcolor{borange}{$\bm{h}_{\bm{2}}\bm{=0}$} &  $\mathbb{I}$  & $\O(\epsilon^{-2})$ & $\sigma \in [1,2)$  \\
 \hline
\cite{Boct2019SIOPT} & $n=2$ & \texttt{NC}: $h_{1},f_{[1,2]}$; $f_1=0$; \textcolor{borange}{$\bm{h}_{\bm{2}}\bm{=0}$} &  $\mathbb{SU}$  & $\O(\epsilon^{-2})$ &  $\sigma\in (0,1)$   \\
 \hline
  \cite{HuangCH19}  & $n\geq 2$ & \texttt{NC}: $f_{n}$; $f_{[1,n-1]}=0$; \texttt{CVX}: $h_{[1,n-1]}$; \textcolor{borange}{$\bm{h}_{\bm{n}}\bm{=0}$} & $\mathbb{BI}$ & $\O(\epsilon^{-2})$  & $\sigma=1$    \\ 
 \hline
 \hline
  \cite{li2022riemannian}$^{d}$  & $n=2$ & \texttt{NC}: $f_{1},h_{1}$; \texttt{CVX}: $h_2$;  $f_2=0$; \texttt{LCONT}: \textcolor{borange}{$\bm{h}_{\bm{2}}\bm{\neq0}$}  & $\mathbb{I}$ & $\O(\epsilon^{-4})$  & $\sigma=1$    \\
 \hline
  This work & $n\geq 2$ & \texttt{NC}: $h_{[1,n-1]},f_{[1,n]}$; \texttt{CVX}: $h_n$;  \texttt{LCONT}: $f_n$, \textcolor{borange}{$\bm{h}_{\bm{n}}\bm{\neq 0}$} & $\mathbb{BI}$  &$\O(\epsilon^{-3})$   & $\sigma \in [1,2)$\\
\hline
  This work & $n\geq 2$ & \texttt{NC}: $h_{[1,n-1]},f_{[1,n]}$; \texttt{CVX}: $h_n$;  \texttt{LCONT}: $f_n$, \textcolor{borange}{$\bm{h}_{\bm{n}}\bm{\neq 0}$} & $\mathbb{SU}$  &$\O(\epsilon^{-3})$    &  $\sigma\in (0,1)$ \\
\hline
\end{tabular}}
  \smallskip
\vspace{1pt}
\scalebox{0.95}{\begin{tabular}{@{}p{\linewidth}@{}}
\footnotesize Note $a$: \textcolor{borange}{$\bm{h}_{\bm{n}}\bm{=0}$} denotes that the $n$-th block has no non-smooth part, making the objective function smooth.\\
\footnotesize Note $b$: The iteration complexity relies on the variational inequality of the convex problem. \\
\footnotesize Note $c$: We adapt their application model into our optimization framework in Equation (\ref{eq:main}) with $(L,S,Z)=(\x_1,\x_2,\x_3)$, as their model additionally requires the linear operator for the first two blocks to be injective. \\
\footnotesize Note $d$: This paper studies manifold optimization with a fixed large penalty and small stepsize. \\
\end{tabular}}
\label{tab:ADMM}\vspace{-12pt}
\end{table*}



%


$\blacktriangleright$ \textbf{Existing Challenges.} We consider the linearly constrained optimization problem in Problem (\ref{eq:main}), which involves $(n-1)$ potentially nonsmooth, nonconvex, and non-Lipschitz composite functions $h_i(\cdot)$ for $i\in[n-1]$, and one convex, nonsmooth composite function $h_n(\cdot)$. Existing ADMM-type methods cannot solve this problem, as they require at least one of the composite functions to be smooth (i.e., $h_n(\cdot)=0$). In the special case where $n=2$, $\A_2=\I$, $h_2(\x_2)=\dot{\rho}\|\x_2\|_1$ with $\dot{\rho}>0$, and $h_1(\cdot)$ is the indicator function of orthogonality constraints, the Riemannian ADMM (RADMM) algorithm \cite{li2022riemannian} is the only known method capable of solving Problem (\ref{eq:main}). However, RADMM cannot handle linearly constrained problems, particularly when $\A_n$ is bijective or surjective. Importantly, it results in suboptimal iteration complexity. A comparison of existing nonconvex ADMM approaches is provided in Table \ref{tab:ADMM}.

$\blacktriangleright$ \textbf{Our Contributions.} Our main contributions are summarized as follows. \bfit{(i)} We introduce IPDS-ADMM to solve the nonconvex nonsmooth optimization problem as in Problem (\ref{eq:main}). This approach imposes the least stringent condition, specifically requiring continuity in only one block of the objective function. It employs an Increasing Penalization and Decreasing Smoothing (IPDS) strategy to ensure convergence (See Section \ref{sect:proposed}). \bfit{(ii)} IPDS-ADMM achieves global convergence when the associated matrix is either bijective or surjective. We establish that IPDS-ADMM converges to an $\epsilon$-critical point with a time complexity of $\mathcal{O}(\epsilon^{-3})$ (See Section \ref{sect:global:convergence}). \bfit{(iii)} We have conducted experiments on the sparse PCA problem to demonstrate the effectiveness of our approach. (See Section \ref{sect:exp}). 





$\blacktriangleright$ \textbf{Assumptions.} Through this paper, we impose the following assumptions on Problem (\ref{eq:main}).

\begin{assumption}\label{ass:3}
Each function $f_i(\cdot)$ is $L_{i}$-smooth for all $i\in [n]$ such that $\|\nabla f_i(\x_i) - \nabla f_i(\grave{\x}_i)\| \leq L_{i} \|\x_i-\grave{\x}_i\|$ holds for all $\x_i,\,\grave{\x}_i \in \Rn^{\mathbf{d}_i\times 1}$. This implies that $| f_i(\x_i) - f_i(\grave{\x}_i) - \la \nabla f_i(\grave{\x}_i), \x_i-\grave{\x}_i \ra | \leq \tfrac{L_{i}}{2}\|\x_i-\grave{\x}_i\|_2^2$ (cf. Lemma 1.2.3 in \cite{Nesterov03}).
\end{assumption}

\begin{assumption} \label{ass:1}
The functions $f_n(\cdot)$ and $h_n(\cdot)$ are Lipschitz continuous with some constants $C_f$ and $C_h$, satisfying $\|\nabla f_n(\x_n)\|\leq C_f$ and $\|\partial h_n(\x_n)\|\leq C_h$ for all $\x_n$.
\end{assumption}

\begin{assumption} \label{ass:2}
We define $\lambdaUp\triangleq \lambda_{\text{max}}(\A_n\A_n\trans)$, $\lambdaDown \triangleq \lambda_{\text{min}}(\A_n\A_n\trans)$, $\lambdaDown'=\lambda_{\text{min}}(\A_n\trans\A_n)$. Either of these two conditions holds for matrix $\A_n$:

\begin{enumerate}[label=\textbf{\alph*)}, leftmargin=14pt, itemsep=1pt, topsep=-5pt, parsep=0pt, partopsep=0pt]

\item Condition $\mathbb{BI}$: $\A_n$ is bijective (i.e., $\lambdaDown=\lambdaDown'>0$), and it holds that $\kappa  \triangleq \lambdaUp/\lambdaDown<2$.

\item Condition $\mathbb{SU}$: $\A_n$ is surjective (i.e., $\lambdaDown>0$, and $\lambdaDown'$ could be zero).

\end{enumerate}

\end{assumption}

\begin{assumption} \label{ass:5}
Given any constant $\bar{\beta}\geq 0$, we let $\underline{\Theta}' \triangleq \inf_{\x_1,\x_2,\ldots,\x_n}\,\sum_{i=1}^n [f_i(\x_i)+h_i(\x_i)] + \frac{\bar{\beta}}{2}\|[\sum_{i=1}^n\A_i\x_i]-\b\|_2^2$. We assert that $\underline{\Theta}' > -\infty$.

\end{assumption}

\begin{assumption} \label{ass:proximal:oper}
For all $i\in[n]$, the proximal operator $\Prox_i(\x_i;\mu)\triangleq \min_{\x'_i}\frac{\mu}{2}\|\x'_i - \x_i\|_2^2+h_i(\x'_i)$ can be computed efficiently and exactly for any given $\x_i\in\Rn^{\d_i\times 1}$ and $\mu>0$.
\end{assumption}

\begin{assumption} \label{ass:bound:x}
If $\sum_{i=1}^n [f_i(\x_i)+h_i(\x_i)]<+\infty$, it follows that $\|\x_i\|<+ \infty$ for all $i\in[n]$.

\end{assumption}

\begin{assumption} \label{ass:bound:xxxx}

For any $i\in[n]$, if the vector $\x_i\in\Rn^{\d_i\times 1}$ is bounded, then the set $\Prox_i(\x_i;\mu)$ is also bounded for all $\mu\in(0,\infty)$.
\end{assumption}

\textbf{Remarks}. \bfit{(i)} Assumption \ref{ass:3} is commonly used in the convergence analysis of nonconvex algorithms. \bfit{(ii)} Assumption \ref{ass:1} imposes a continuity assumption only for the last block, allowing other blocks of the function $\{h_i(\x_i)\}_{i=1}^{n-1}$ to be nonsmooth and non-Lipschitz, such as indicator functions of constraint sets. It ensures bounded (sub-)gradients for $f_n(\cdot)$ and $h_n(\cdot)$, a relatively mild requirement that has found use in nonsmooth optimization \cite{li2022riemannian,li2021weakly,HuangCH19,BohmW21}. \bfit{(iii)} Assumption \ref{ass:2} demands a condition on the linear matrix $\A_i$ for the last block ($i=n$), while leaving $\A_i$ unrestricted for $i\in[n-1]$. \bfit{(iv)} Assumption \ref{ass:5} ensures the well-defined nature of the penalty function associated with the problem, as has been used in \cite{goncalves2017convergence}. Furthermore, Assumption \ref{ass:5} can be satisfied if $\sum_{i=1}^n [f_i(\x_i)+h_i(\x_i)]>-\infty$. \bfit{(v)} Assumption \ref{ass:proximal:oper} is frequently employed in nonconvex ADMM frameworks \cite{li2015global,Boct2019SIOPT}. Common examples of functions $h_i(\x_i)$ arising in practical applications include those discussed in \cite{GongZLHY13}, $\ell_0$ regularization, $\ell_{1/2}$ regularization \cite{ZengLWX14}, and indicator functions of cardinality constraints, matrices with orthogonality constraints \cite{lai2014splitting}, and matrices with rank constraints, among others. \bfit{(vi)} Assumptions \ref{ass:bound:x} and \ref{ass:bound:xxxx} are used to guarantee the boundedness of the solution.


$\blacktriangleright$ \textbf{Notations.} We define $[n] \triangleq \{1,2,\ldots,n\}$ and $\x \triangleq \x_{[n]} \triangleq \{\x_1,\x_2,\ldots,\x_n\}$. For any $j\geq i$, we denote $\x_{[i,j]} \triangleq \{\x_{i},\x_{i+1},\ldots.,\x_{j}\}$. We define $\lambda_{\text{min}}(\mathbf{M})$ and $\lambda_{\text{max}}(\mathbf{M})$ as the smallest and largest eigenvalue of the given matrix $\mathbf{M}$, respectively. We denote $\|\A_i\|$ as the spectral norm of the matrix $\A_i$. We denote $\A\x\triangleq \sum_{j=1}^n \A_j\x_j$, and $\|\x^{+} - \x\|_2^2 = \sum_{i=1}^n \|\x_i^{+} - \x_i\|_2^2$. Further notations and technical preliminaries are provided in Appendix \ref{app:sect:preliminaries:additional}.

\section{Motivating Applications} \label{sect:app}

Many machine learning and data science models can be formulated as Problem (\ref{eq:main}). Below, we present two examples, with additional applications provided in Appendix \ref{app:sect:app}.

$\blacktriangleright$ \textbf{Sparse PCA}. Sparse PCA  \cite{chen2016augmented,lu2012augmented} focuses on identifying a subset of informative variables with sparse loadings to enhance interpretability and reduce model complexity. It is formulated as: $\ts\min_{\V \in \Rn^{\dot{d} \times \dot{r}}}  \frac{1}{2\dot{m}}\|\D - \D \V \V\trans \|_{\fro}^2 + \dot{\rho} \|\V\|_1,\,s.t.\, \V\in \mathcal{M}\triangleq\{\V\,|\,\V\trans\V = \mathbf{I}\}$, where $\D\in \Rn^{\dot{m} \times \dot{d}}$ is the data matrix, and $\dot{\rho}\geq 0$. Introducing an additional variable $\Y$, this problem can be formulated as: $\min_{\V,\Y}  \frac{1}{2\dot{m}}\|\D - \D \V \V\trans \|_{\fro}^2 + \dot{\rho} \|\V\|_1 + \iota_{\mathcal{M}}(\Y),\,\st\,-\Y+\V=\zero$. It corresponds to Problem (\ref{eq:main}) with $\x_1 = \vec(\Y)$, $\x_2 = \vec(\V)$, $f_1(\x_1)=0$, $h_1(\x_1)=\iota_{\mathcal{M}}(\Y)$, $f_2(\x_2)=\frac{1}{2\dot{m}}\|\D - \D \V \V\trans \|_{\fro}^2$, $h_2(\x_2) = \dot{\rho} \|\V\|_1$, $\A_1=-\I$, $\A_2=\I$, $\b=\zero$, and Condition $\BIBI$.

$\blacktriangleright$ \textbf{Structured Sparse Phase Retrieval}. Sparse phase retrieval \cite{Duchiiay015} aims to recover a sparse signal from the magnitudes of linear measurements. By incorporating additional linear constraints, recovery accuracy can be further improved. The problem is formulated as: $\min_{\v} \|(\G\v) \odot (\G\v) - \z\|_2^2 + \dot{\rho} \|\v\|_1,\,\st\,\D\v \geq \zero$, where $\dot{\rho}\geq 0$, $\G\in\Rn^{\dot{m}\times \dot{d}}$, $\z\in\Rn^{\dot{m}}$, $\D \in\Rn^{\dot{r}\times \dot{d}}$, with $\D$ being surjective that $\D\D\trans \succ \zero$. Introducing a new variable $\mathbf{y}$, this problem can be formulated as: $\min_{\v,\mathbf{y}} \|(\G\v) \odot (\G\v) - \z\|_2^2 + \dot{\rho} \|\v\|_1 + \iota_{\geq \zero}(\y),\,s.t.\,\mathbf{y}-\D\v=\zero$. This corresponds to Problem (\ref{eq:main}) with $\x_1 = \mathbf{y}$, $\x_2 = \v$, $f_1(\x_1)=0$, $h_1(\x_1)=\iota_{\geq \zero}(\mathbf{y})$, $f_2(\x_2)=\tfrac{1}{2}\|(\G\v) \odot (\G\v) - \b\|_2^2$, $h_2(\x_2) =\dot{\rho}\|\v\|_1$, $\A_1=\I$, $\A_2=-\D$, $\b=\zero$, and Condition $\SUSU$.

\section{The Proposed IPDS-ADMM Algorithm} \label{sect:proposed}

This section describes the proposed IPDS-ADMM algorithm for solving Problem (\ref{eq:main}), featuring with using a new Increasing Penalization and Decreasing Smoothing (IPDS) strategy.

\subsection{Increasing Penalization Strategy}
We employ an increasing penalty update strategy that is crucial to our algorithm. A natural choice for this penalty update rule is to use functions from the $\ell_p$ family. Throughout this paper, we consider the following penalty update rule $\{\beta^t\}_{t=0}^{\infty}$ for any given parameters $\xi,\delta,p\in(0,1)$:
\beq \label{eq:beta:t:beta:t}
\beta^{t} =  \beta^0 (1 +  \xi t^p),~\beta^0 \geq L_n  / (\delta\lambdaUp).
\eeq
\noi Here, $L_n$ and $\lambdaUp$ are defined in Assumption \ref{ass:3} and Assumption \ref{ass:2}, respectively.

We obtain the following useful lemma regarding the penalty update rule.

\begin{lemma} \label{lemma:updating:Sublinear}
(Proof in Appendix \ref{app:lemma:updating:Sublinear}) Given $\xi,\delta,p\in(0,1)$, assume Formulation (\ref{eq:beta:t:beta:t}) is used to choose $\{\beta^t\}_{t=0}^{\infty}$. We have: (\bfit{a}) $\beta^t\leq \beta^{t+1}\leq (1+\xi )\beta^t$, (\bfit{b}) $L_n  \leq \delta \beta^t  \lambdaUp$.

\end{lemma}

\textbf{Remarks} \bfit{(i)} The increasing penalty update strategy is closely coupled with the decreasing smoothing strategy and the diminishing stepsize approach in the literature. These strategies are frequently employed in subgradient methods \cite{li2021weakly}, smoothing gradient methods \cite{BohmW21,sun2021algorithms,YangMOR}, penalty decomposition methods \cite{lu2013sparse}, and stochastic optimization algorithms like ADAM \cite{KingmaB14,ChenJMLR2022}, but are less commonly utilized in ADMM frameworks. We examine this approach within ADMM but limit our discussion to specific form and condition as in Formulation (\ref{eq:beta:t:beta:t}). \bfit{(ii)} The condition $\beta^0 \geq L_n  / (\delta\lambdaUp)$ in Formulation (\ref{eq:beta:t:beta:t}) essentially mandates that the initial penalty value be sufficiently large. This condition can be automatically satisfied since an increasing penalty update is used. \bfit{(iii)} The result $\beta^{t+1} \leq (1 + \xi)\beta^t$ in Lemma \ref{lemma:updating:Sublinear} implies that the penalty parameter grows, but not excessively fast, with a constant $\xi$ to prevent rapid escalation.



\subsection{Decreasing Smoothing Strategy}

IPDS-ADMM is built upon the Moreau envelope smoothing technique \cite{li2022riemannian,zeng2022moreau,sun2021algorithms,BohmW21}. Our algorithm gradually decreases the smoothing parameter at a controlled rate.

Initially, we provide the following useful definition.

\begin{definition}

The Moreau envelope of a proper convex and Lipschitz continuous function $h(\ub): \Rn^{d\times 1} \mapsto\Rn$ with parameter $\mu \in (0,\infty)$ is defined as $h(\ub;\mu) \triangleq \min_{\v \in \Rn^{d \times 1}}  h(\v) + \frac{1}{2 \mu} \| \v - \ub\|_2^2$.

\end{definition}

We present several useful properties of Moreau envelope functions in the following lemmas.

\begin{lemma}\label{lemma:lip:mu:2}

\cite[Chapter 6]{beck2017first} Suppose the function $h(\ub)$ is $C_{h}$-Lipschitz continuous and convex \textit{w.r.t.} $\ub$. We have: (\bfit{a}) The function $h(\ub;\mu)$ is $C_h$-Lipschitz continuous \textit{w.r.t.} $\ub$. (\bfit{b}) The function $h(\ub;\mu)$ is $(1/\mu)$-smooth \textit{w.r.t.} $\ub$, and its gradient can be computed as: $\nabla\,h(\ub;\mu) = \frac{1}{\mu}(\ub -\Prox_{h}(\ub;\mu) )$, where $\Prox_{h}(\ub;\mu) = \arg \min_{\v} h(\v) + \frac{1}{2 \mu} \| \v -\ub\|_2^2$. (\bfit{c}) $0 \leq h(\ub) - h(\ub;\mu) \leq \tfrac{1}{2}\mu C_h^2$.

\end{lemma}

\begin{lemma} \label{lemma:lip:mu0}
(Proof in Appendix \ref{app:lemma:lip:mu0}) Assuming $0<\mu_2<\mu_1$ and fixing $\ub \in \Rn^{d \times 1}$, we have: $0\leq \frac{ h(\ub;\mu_2) - h(\ub;\mu_1)}{\mu_1 - \mu_2}\leq \tfrac{1}{2}C_h^2$.

\end{lemma}

\begin{lemma} \label{lemma:lip:mu}
(Proof in Appendix \ref{app:lemma:lip:mu}) Assuming $0<\mu_2<\mu_1$ and fixing $\ub \in \Rn^{d \times 1}$, we have: $\|\nabla\,h(\ub;\mu_1) -\nabla\,h(\ub;\mu_2)\| \leq (\frac{\mu_1}{\mu_2} - 1) \cdot  C_h$.
\end{lemma}

\begin{lemma} \label{lemma:smoothing:problem}
(Proof in Appendix \ref{app:lemma:smoothing:problem}) Given constants $\{\cb,\mu,\rho\}$, we consider the convex problem: $\bar{\x}_n = \arg \min_{\x_n} h_n(\x_n;\mu) + \frac{\rho}{2} \| \x_n - \cb\|_2^2$. We have: (\bfit{a}) $\bar{\x}_n = \frac{\mu}{1+\mu \rho} ( \frac{1}{\mu}\breve{\x}_n  + \rho \cb )$, where $\breve{\x}_n = \arg \min_{\breve{\x}_n}~h_n(\breve{\x}_n) + \frac{1}{2} \cdot \frac{\rho}{1+\mu\rho} \|\breve{\x}_n-\cb\|_{\fro}^2=\Prox_n(\cb;\mu+1/\rho)$. (\bfit{b}) $\rho (\cb - \bar{\x}_n) \in \partial h(\breve{\x}_n)$. (\bfit{c}) $\|\x_n - \breve{\x}_n\|\leq \mu C_h$.

\end{lemma}

\begin{remark}
(\bfit{i}) Lemmas \ref{lemma:lip:mu0} and \ref{lemma:lip:mu} are derived using standard convex analysis techniques and play a key role in analyzing the proposed IPDS-ADMM algorithm. (\bfit{ii}) Lemma \ref{lemma:smoothing:problem} is essential for establishing the iteration complexity of Algorithm \ref{alg:main} in reaching a critical point. The results of Lemma \ref{lemma:smoothing:problem} are analogous to those of Lemma 1 in \cite{li2022riemannian}.
\end{remark}


\subsection{The Main Algorithm }











This subsection provides the proposed IPDS-ADMM algorithm. Initially, we consider the following alternative optimization problem:
\beq\label{eq:main:1}
\ts \min_{\x_{1},\x_{2},\ldots,\x_{n} }\,h_n(\x_{i};\mu) + [\sum_{i=1}^{n-1} h_i(\x_{i}) ]+  [\sum_{i=1}^{n} f_i(\x_{i})]   ,\,s.t.\,[\sum_{i=1}^n \A_{i} \x_{i}] = \b,
\eeq
\noi where $\mu\rightarrow 0$, and $h_n(\x_{n};\mu) \triangleq \min_{\v \in \Rn^{\d_n \times 1}}  h(\v) + \frac{1}{2 \mu} \| \v - \x_{n}\|_2^2$ is the Moreau envelope of $h_n(\x_n)$ with parameter $\mu$. Lemma \ref{lemma:lip:mu:2} confirms that $h_n(\x_{n},\mu)$ is a $(1/\mu)$-smooth function assuming $h_n(\cdot)$ is convex. We present the augmented Lagrangian function for Problem (\ref{eq:main:1}), as follows:
\beq\label{eq:Lag}
\ts \mathcal{L}(\x,\mathbf{z};\beta,\mu) \triangleq h_{n}(\x_n;\mu) + \{\sum_{i=1}^{n-1} h_i(\x_i)\} + G(\x,\mathbf{z};\beta) ,
\eeq
\noi where $G(\x,\mathbf{z};\beta)$ is differentiable and defined as:
\beq
\ts G(\x,\mathbf{z};\beta) \triangleq   \sum_{i=1}^n f_i(\x_i) + \la [\sum_{i=1}^n\A_i\x_i] - \b,\mathbf{z} \ra + \frac{\beta}{2}\| [\sum_{i=1}^n\A_i\x_i] - \b\|_2^2. ~~~~~~~~~~\,\,\,\nn
\eeq
\noi Here, $\mu \in (0,\infty)$, $\beta\in (0,\infty)$, and $\mathbf{z}\in \Rn^{m\times 1}$ are the smoothing parameter, the penalty parameter, and the dual variable, respectively. We employ an increasing penalty and decreasing smoothing update scheme throughout all iterations $t=\{0,1,\ldots,\infty\}$ with $\beta^t\rightarrow +\infty$ and $\mu^t \propto \tfrac{1}{\beta^t} \rightarrow 0$. Notably, the function $G(\x^t,\mathbf{z}^t;\beta^t)$ is $\LL_i^t$-smooth \textit{w.r.t.} $\x_i$ for all $i\in[m]$, where $\LL_i^t = L_i + \beta^t \|\A_i\|_2^2$. For notation simplicity, for all $i\in[n]$, we denote $\gg_i^t\triangleq \nabla_{\x_i} G(\x^{t+1}_{[1,i-1]},\x^t_i,\x^{t}_{[i+1,n]},\mathbf{z}^t;\beta^t)$ as the gradient of $G(\x,\z^t;\beta^t)$ \textit{w.r.t.} $\x_i$ at the point $\x_i^t$.

In each iteration, we select suitable parameters $\{\beta^t,\mu^t\}$ and sequentially update the variables $(\x_1,\x_2,\ldots,\x_n,\z)$. We employ the proximal linearized method to cyclically update the variables $\{\x_1, \x_2, \ldots, \x_n\}$. Specifically, we update each variable $\x_i$ by solving the following subproblem for all $i\in[n]$: $\x_i^{t+1} \thickapprox \arg \min_{\x_i \in \Rn^{\d_i\times 1}}\mathcal{L}(\x_{[1,i-1]}^t,\x_i,\x^t_{[i+1,n]},\z^t;\beta^t,\mu^t)$. To address the $\x_i$-subproblem, we employ a proximal linearized minimization strategy for all $i\in[n-1]$: $\x_i^{t+1} \in \arg \min_{\x_i} h_i(\x_i)  + \frac{\theta_{1} \LL_i^t}{2}\| \x_i - \x_i^t\|_2^2 + \la \x_i - \x_i^t, \ddot{\g}_i^t,\mathbf{z}^t;\beta^t)\ra$. However, for the final block of the problem, we consider a subtly different proximal linearized minimization strategy: $\x_n^{t+1} =  \arg \min_{\x_n} h_n(\x_n;\mu^t) + \frac{\theta_{2} \LL_n^t}{2}\| \x_n - \x_n^t\|_2^2 + \la \x_n - \x_n^t,\ddot{\g}_n^t \ra$. Importantly, we assign $\theta_1$ to blocks $[1,n-1]$ and $\theta_2$ to block $n$. Our algorithm updates the dual variable $\z^t$ using either an under-relaxation stepsize $\sigma \in(0,1)$ or an over-relaxation stepsize $\sigma \in[1,2)$.

\begin{algorithm} [!h]
\DontPrintSemicolon
\caption{IPDS-ADMM: The Proposed Proximal Linearized ADMM for Problem (\ref{eq:main}). }
\label{alg:main}

Choose suitable parameters $\{p,\xi,\delta\}$ and $\{\sigma,\theta_1,\theta_2\}$ using Formula (\ref{eq:choice:BI}) or Formula (\ref{eq:choice:SU}). \\

Initialize $\{\x^0,\mathbf{z}^0\}$. Choose $\beta^0 \geq L_n  / (\delta\lambdaUp)$.

\For{$t$ from 0 to $T$}{
\textbf{S1}) IPDS Strategy: Set $\beta^t = \beta^0(1+\xi t^p)$, $\mu^t = 1 / (\lambdaUp\delta\beta^t)$.

We define $\gg_i^t\triangleq \nabla_{\x_i} G(\x^{t+1}_{[1,i-1]},\x^t_i,\x^{t}_{[i+1,n]},\mathbf{z}^t;\beta^t)$.

\textbf{S2})  $\x_1^{t+1}  \in \arg \min_{\x_1} h_1(\x_1)  + \la \x_1-\x_1^t,  \gg_1^t \ra + \tfrac{\theta_1 \LL_1^t}{2} \|\x_1 - \x_1^t\|_2^2$

\textbf{S3})  $\x_2^{t+1}  \in \arg\min_{\x_2} h_2(\x_2)  + \la \x_2-\x_2^t,  \gg_2^t \ra+ \tfrac{\theta_1 \LL_2^t}{2} \|\x_2 - \x_2^t\|_2^2$

\quad $\ldots$

\textbf{S4})  $\x_{n-1}^{t+1}  \in \arg\min_{\x_{n-1}} h_{n-1}(\x_{n-1})  + \la \x_{n-1}-\x_{n-1}^t,\gg_{n-1}^t\ra + \tfrac{\theta_1 \LL_{n-1}^t}{2} \|\x_{n-1} - \x_{n-1}^t\|_2^2$

\textbf{S5})  $\x_{n}^{t+1}  \in \arg\min_{\x_{n}} h_{n}(\x_{n};\mu)  + \la \x_{n}-\x_{n}^t,\gg_{n}^t\ra + \tfrac{\theta_2 \LL_{n}^t}{2} \|\x_{n} - \x_{n}^t\|_2^2$. It can be solved using Lemma \ref{lemma:smoothing:problem} as $\x_n^{t+1} =\tfrac{1}{1+\mu\rho} (\breve{\x}_n^{t+1}+\mu\rho \cb)$, where $\breve{\x}_n^{t+1}=\Prox_n(\cb;\mu+1/\rho)$, $\mu=\mu^t$, $\rho\triangleq \theta_2 \LL_n^t$, and $\cb \triangleq \x_n^t - \gg_n^t/\rho$.



\textbf{S6})  $\mathbf{z}^{t+1}  =  \mathbf{z}^t +  \sigma \beta^t ( [\sum_{j=1}^n \A_j \x^{t+1}_j] - \b)$

}

\label{alg:main}
\end{algorithm}

We present IPDS-ADMM in Algorithm \ref{alg:main}, and have the following remarks.

\begin{remark}
(\bfit{i}) Algorithm \ref{alg:main} can be viewed as a generalized cyclic coordinate descent method applied to the augmented Lagrangian function in Equation (\ref{eq:Lag}). (\bfit{ii}) The Moreau envelope smoothing technique has been used in the design of augmented Lagrangian methods \cite{zeng2022moreau}, ADMMs \cite{li2022riemannian,yuanICLR25FADMM}, and minimax optimization \cite{ZhangX0L20}. However, these algorithms typically utilize constant penalties, whereas we adopt an Increasing Penalization and Decreasing Smoothing (IPDS) strategy to improve the iteration complexity of RADMM \cite{li2022riemannian}, reducing it from $\mathcal{O}(\epsilon^{-4})$ to $\mathcal{O}(\epsilon^{-3})$. (\bfit{iii}) Algorithm \ref{alg:main} is a fully splitting algorithm, where each step reduces to computing a proximal operator. For the first $(n-1)$ blocks, we have: $\x_i^{t+1} \in \Prox_{i}(\x_i^t - \gg^t_i/\dot{\rho};1/\dot{\rho})$, where $\dot{\rho}=\theta_1 \LL_{i}^t$. For the last block, Lemma \ref{lemma:smoothing:problem} can be applied to compute the proximal operator of the smoothed function $h_{n}(\x_{n};\mu)$ using the proximal operator of the original function $h_{n}(\x_{n})$. (\bfit{iv}) The point $\breve{\x}_n^{t+1}$ in Step \textbf{S5}) of Algorihtm \ref{alg:main} plays a crucial role. As will be seen later in Theorem \ref{theorem:case:AAAAA:2}, the point $(\x^{t}_1,\x^{t}_2,\ldots,\x^{t}_{n-1},\breve{\x}^{t}_n,\z^t)$, rather than the point $(\x^{t}_1,\x^{t}_2,\ldots,\x^{t}_{n-1},\x^{t}_n,\z^t)$, will serve as an approximate critical point of Problem (\ref{eq:main}) in our complexity results. (\bfit{v}) RADMM \cite{li2022riemannian} uses a fixed large penalty parameter $\mathcal{O}(1/\epsilon)$ and a fixed small smoothing parameter $\mathcal{O}(\epsilon)$ to achieve an $\epsilon$-approximate critical point. However, this leads to overly conservative step sizes for the primal and dual updates, potentially hindering the algorithm's practical performance. (\bfit{vi}) We apply the smoothing strategy only to the last block to bound the dual variables via the primal ones. This leverages the Lipschitz continuity of the smoothed function to estimate $\tfrac{1}{\beta^t}\|\z^{t+1}-\z^t\|_2^2$ and construct a suitable potential function. (\bfit{vii}) Some may be concerned that using an increasing penalty could cause the parameter to grow excessively fast. However, by setting $\xi \ll 1$, we ensure that $\beta^{t+1} \leq (1+\xi)\beta^t$, meaning the penalty grows very slowly in practice.

\end{remark}

\subsection{Choosing Suitable Parameters $\{p,\xi,\delta\}$ and $\{\sigma,\theta_1,\theta_2\}$ } 
\label{sect:sub:para:selection}


Selecting appropriate parameters $\{p,\xi,\delta\}$ and $\{\sigma,\theta_1,\theta_2\}$ is essential to ensuring the global convergence of Algorithm \ref{alg:main}. In our theoretical analysis and empirical experiments, we suggest the following choices for $\{p,\xi,\delta\}$ and $\{\sigma,\theta_1,\theta_2\}$:
\begin{align} 
\BIBI: & ~ p = \tfrac{1}{3},~\xi\in (0,\infty),~\delta \in (0,\frac{1}{3}(\tfrac{2}{\kappa}-1)),\sigma \in [1,2),\theta_1=1.01, \theta_2 = \tfrac{ 1 / \kappa - \delta}{1+\delta} + \tfrac{1}{2\varrho (1+\delta)^2}.\label{eq:choice:BI} \\
\SUSU:&  ~ p=\tfrac{1}{3},~    \xi=\delta  = \sigma = \frac{0.01}{\kappa},~\theta_1=1.01, \theta_2 = 1.5. \label{eq:choice:SU}
\end{align}
\noi Here, $\varrho\triangleq 6 \omega \sigma_1 \kappa$, $\sigma_1 \triangleq \tfrac{\sigma}{(1-|1-\sigma|)^2}$, and $\omega \triangleq 1 + \frac{ \xi}{2 \sigma} + \sigma \xi$. Notably, $\theta_2$ in (\ref{eq:choice:BI}) depends on $(\xi,\delta,\sigma)$.

\begin{remark}
(\bfit{i}) We obverse from (\ref{eq:choice:SU}) that the parameters $\{\xi,\delta,\sigma\}$ is inversely  proportional to the condition number $\kappa$. Such settings are partly consistent with those in \cite{Boct2019SIOPT} (refer to Lemma 5 in \cite{Boct2019SIOPT}). (\bfit{ii}) Introducing the relaxation parameter $\sigma\in(0,2)$ enables handling cases where the matrix is surjective. Specifically, when the matrix is bijective, we can use an over-relaxation step size for faster convergence, whereas for surjective matrices, the algorithm requires conservative step sizes to ensure global convergence.
\end{remark}


%
%

%
%

%
%
%
%
%

\section{Global Convergence} \label{sect:global:convergence}

This section establishes the global convergence of Algorithm \ref{alg:main}.

We begin with a high-level overview of the proof strategy. First, using the Lagrangian function, we derive sufficient decrease conditions for the four parameter sets: primal variables, dual variables, the penalty parameter, and the smoothing parameter. Next, using the first-order optimality conditions and dual update rules, we bound the difference in dual variables using primal by the difference in primal variables. Lastly, we show that the tail error term related to the smoothing parameter is constant, establishing the summability of the sequence linked to a potential function.

We provide the following three useful lemmas.

\begin{lemma} \label{lemma:suf:dec}
(Proof in Appendix \ref{app:lemma:suf:dec}, {\rm{A Sufficient Decrease Property}}) Fix $\varepsilon_3 \triangleq \xi$ and $\varepsilon_1\triangleq \frac{1}{2}\theta_1 - \frac{1}{2}$. Let $\varepsilon_2 \in \Rn$. For all $t\geq1$, we have:
\beq
\mathcal{E}^{t+1} + \Theta^{t+1}_{\oo} - \Theta^{t}_{\oo}  \leq  \ts  (\tfrac{1}{2} - \theta_2 +\varepsilon_2  ) \cdot \LL_n^t \|\x_n^{t+1}-\x_n^t\|_2^2 +   \tfrac{\omega}{\sigma \beta^t} \|\z^{t+1}-\z^t\|_2^2,
\eeq
\noi where 
\beq
&&\mathcal{E}^{t+1} \triangleq  \ts [\varepsilon_1\sum_{i=1}^{n-1} \LL_i^t \| \x_i^{t+1} - \x_i^t \|_2^2] + \varepsilon_2\LL_n^t \| \x_n^{t+1} - \x_n^t \|_2^2 + \frac{\varepsilon_3}{\beta^t}\|\z^{t+1}-\z^t\|_2^2.\nn\\
&&\Theta^{t}_{\oo}  \triangleq  \mathcal{L}(\x^{t},\mathbf{z}^{t};\beta^{t},\mu^{t}) + \tfrac{1}{2} C_h \mu^{t},\,\LL_i^t = L_i + \beta^t \|\A_i\|_2^2,\,\omega \triangleq 1 + \tfrac{ \xi}{2 \sigma} + \sigma \xi.\nn
\eeq



\end{lemma}

\begin{lemma} \label{lemma:first:order}
(Proof in Appendix \ref{app:lemma:first:order}, {\rm{First-Order Optimality Condition}}) Assume $\sigma\in(0,2)$. For all $t\geq1$ and $i\in[n-1]$, we have the following results.

\begin{enumerate}[label=\textbf{(\alph*)}, leftmargin=22pt, itemsep=1pt, topsep=1pt, parsep=0pt, partopsep=0pt]

\item Let $\mathbbm{w}_i^{t+1} \in \partial h_i (\x_i^{t+1} ) + \nabla f_i(\x_i^{t})$, and $\u_i^{t+1} \triangleq \theta_1  \LL_i^t (\x_i^{t+1} - \x_i^{t}) - \beta^t  \A_i \trans [\sum_{j=i}^n \A_j(\x_j^{t+1}-\x_j^t)]$. It holds that: $\zero = \sigma \A_i\trans\mathbf{z}^t +   \A_i\trans(\mathbf{z}^{t+1}-\mathbf{z}^{t}) +  \sigma \mathbbm{w}_i^{t+1} + \sigma \u_i^{t+1}$.

\item Let $\mathbbm{w}_n^{t+1} \triangleq \nabla h_{n} (\x_n^{t+1},\mu^t ) + \nabla f_n(\x_n^{t})$, and $\u_n^{t+1}\triangleq \H^{t} (\x_n^{t+1} - \x_n^{t})$, where $\H^{t} \triangleq \theta_2  \LL_n^{t}\mathbf{I} - \beta^{t}\A_n \trans \A_n$. It holds that: $\zero = \sigma \A_n\trans\mathbf{z}^t +   \A_n\trans(\mathbf{z}^{t+1}-\mathbf{z}^{t}) +  \sigma \mathbbm{w}_n^{t+1} + \sigma \u_n^{t+1} $.

\item We have the following two different identities:
{\small \begin{align}
\BIBI:& \underbrace{\A_n\trans (\mathbf{z}^{t+1} - \mathbf{z}^{t})}_{\triangleq \mathbbm{a}^{t+1}} = (1-\sigma)\underbrace{(\A_n\trans (\mathbf{z}^{t} - \mathbf{z}^{t-1}))}_{\triangleq\mathbbm{a}^t} + \sigma \underbrace{ (\u_n^{t}- \u_n^{t+1} +  \mathbbm{w}_n^{t} - \mathbbm{w}_n^{t+1} ) }_{\mathbbm{c}^{t}}. \label{eq:III} \\
\SUSU:& \underbrace{\A_n\trans (\mathbf{z}^{t+1} - \mathbf{z}^{t}) +  \sigma \u_n^{t+1}}_{ \triangleq \mathbbm{a}^{t+1}}  = (1-\sigma) (\underbrace{\A_n\trans (\mathbf{z}^{t} - \mathbf{z}^{t-1}) +  \sigma \u_n^{t}}_{ \triangleq \mathbbm{a}^{t}})  + \sigma( \underbrace{\sigma \u_n^{t} +   \mathbbm{w}_n^{t} - \mathbbm{w}_n^{t+1}}_{\triangleq\mathbbm{c}^t}).  \label{eq:AAA}
\end{align}
}

\end{enumerate}

\end{lemma}












\begin{lemma} \label{lemma:reg:mu:L}
(Proof in Appendix \ref{app:lemma:reg:mu:L}) For all $t\geq 0$, we have: \bfit{(a)} $\LL_n^t  \leq \beta^t \lambdaUp (1+\delta)$; \bfit{(b)} $\|\H^{t}\|\leq \beta^t \lambdaUp q$, where $q \triangleq   \theta_2 (1+\delta)  -  \lambdaDown' / \lambdaUp$; \bfit{(c)} $\|\u_n^{t+1}\| \leq  q \lambdaUp \beta^t \|\x_n^{t+1} - \x_n^{t}\|$.

\end{lemma}

We provide convergence analysis of Algorithm \ref{alg:main} under two conditions: {Condition} $\BIBI$ using Formulation (\ref{eq:III}), and {Condition} $\SUSU$ using Formulation (\ref{eq:AAA}). We define $\Theta^{t}_{\oo}  \triangleq  \mathcal{L}(\x^{t},\mathbf{z}^{t};\beta^{t},\mu^{t}) + \tfrac{1}{2} C_h \mu^{t}$, and $\omega \triangleq 1 + \frac{ \xi}{2 \sigma} + \sigma \xi$. We define $\sigma_1 \triangleq \tfrac{\sigma}{(1-|1-\sigma|)^2}$, and $\sigma_2\triangleq \tfrac{|1-\sigma|}{\sigma (1-|1-\sigma|)}$, where $\sigma\in(0,2)$. We construct a sequence associated with the potential (or Lyapunov) function for different Conditions $\BIBI$ and $\SUSU$ as follows:
{\small\begin{align}
\ts\hspace{-5pt}\BIBI:~&\ts \Theta^t = \Theta^{t}_{\oo}  + \ts \underbrace{\ts \tfrac{ \omega \sigma_2}{\lambdaDown}}_{\triangleq a}\cdot \underbrace{\ts \tfrac{1}{\beta^t}\|\mathbbm{a}^{t}\|_2^2}_{\triangleq \AAA^t }  + \underbrace{\ts \tfrac{3 \omega \sigma_1}{\lambdaDown}}_{\triangleq b} \cdot \underbrace{\tfrac{1}{\beta^t} (L_{n} \| \x_n^{t} - \x_n^{t-1}\| + \|\u_n^{t}\| )^2}_{\triangleq \BBB^t}. \label{eq:def:III} \\
\ts \hspace{-5pt}\SUSU:~&\ts \Theta^t = \Theta^{t}_{\oo}  + \ts \underbrace{\ts \tfrac{ 2\omega \sigma_2}{\lambdaDown}}_{\triangleq a}\cdot \underbrace{\ts\tfrac{1}{\beta^t}\|\mathbbm{a}^{t}\|_2^2}_{\triangleq \AAA^t }  + \underbrace{\ts \tfrac{6 \omega \sigma_1}{\lambdaDown}}_{\triangleq b} \cdot \underbrace{\tfrac{1}{\beta^t} (L_{n} \| \x_n^{t} - \x_n^{t-1}\| + \sigma \|\u_n^{t}\| )^2}_{\triangleq \BBB^t}.\label{eq:def:AAA}
\end{align}}
\subsection{Analysis for Condition $\BIBI$}

We provide a convergence analysis of Algorithm \ref{alg:main} under Condition $\BIBI$, where $\A_n$ is a bijective matrix. We assume an over-relaxation stepsize is used with $\sigma\in [1,2)$.

The subsequent lemma uses Equation (\ref{eq:III}) to establish an upper bound for the term $\frac{\omega}{\sigma\beta^t}\|\mathbf{z}^{t+1} - \mathbf{z}^t\|_2^2$.

\begin{lemma} \label{lemma:bound:dual:Case:I}

(Proof in Appendix \ref{app:lemma:bound:dual:Case:I}, {\rm{Bounding Dual Using Primal}}) We define $\omega$ as in Lemma \ref{lemma:suf:dec}. For all $t\geq 1$, we have:
\beq \label{eq:bound:dual:I}
\tfrac{\omega}{\sigma \beta^t}\|\mathbf{z}^{t+1}-\mathbf{z}^t\|_2^2 \leq \Theta^t_{+} - \Theta^{t+1}_{+} +  \chi    \LL_n^t \| \x_n^{t+1}  - \x_n^{t} \|_2^2    + \UUU^t ,
\eeq
\noi where $\chi\triangleq \varrho (    \delta + \theta_2+\theta_2\delta-1/\kappa)^2$, $\varrho\triangleq 6 \omega \sigma_1 \kappa$, $\Theta^t_{+} \triangleq a\AAA^t + b\BBB^t $, and $\UUU^t \triangleq  C_{h}^2 \tfrac{b}{\beta^t}\cdot ( \tfrac{\mu^{t-1}}{\mu^t} - 1 )^2$. Here, $\{a,\AAA^t,b,\BBB^t\}$ are defined in Equation (\ref{eq:def:III}).

\end{lemma}

Assume Equation (\ref{eq:choice:BI}) is used to choose $\{p,\xi,\delta,\sigma,\theta_1,\theta_2\}$. We have the following lemma.

\begin{lemma} \label{lemma:rule:I}
(Proof in Appendix \ref{app:lemma:rule:I}) We define $\{\chi,\varrho\}$ in Lemma \ref{lemma:bound:dual:Case:I}. We have the following results:

\begin{enumerate}[label=\textbf{(\alph*)}, leftmargin=22pt, itemsep=1pt, topsep=1pt, parsep=0pt, partopsep=0pt]

\item It holds that $\varepsilon_1 \triangleq  \frac{1}{2} \theta_1- \frac{1}{2}>0$, and $\varepsilon_2 \triangleq  \theta_2-\frac{1}{2} - \chi>0$.
\item For all $t\geq 1$, we have $\mathcal{E}^{t+1}      \leq  \Theta^t - \Theta^{t+1} + \UUU^t$.
\end{enumerate}

\end{lemma}

\subsection{Analysis for Condition $\SUSU$}

We provide a convergence analysis of Algorithm \ref{alg:main} under Condition $\SUSU$, where $\A_n$ is a surjective matrix. We assume an under-relaxation stepsize is used with $\sigma\in (0,1)$.

The following lemma utilizes Equation (\ref{eq:AAA}) to establish an upper bound for the term $\frac{\omega}{\sigma\beta^t}\|\mathbf{z}^{t+1} - \mathbf{z}^t\|_2^2$.

\begin{lemma} \label{lemma:bound:dual:Case:A}
(Proof in Appendix \ref{app:lemma:bound:dual:Case:A}, {\rm{Bounding Dual Using Primal}}) We define $\omega$ as in Lemma \ref{lemma:suf:dec}. For all $t\geq 1$, we have:
\begin{align} \label{eq:bound:dual:A}
 \ts \tfrac{\omega}{\sigma\beta^t} \|\mathbf{z}^{t+1} - \mathbf{z}^{t}\|_2^2 \leq \Theta^t_{+} - \Theta^{t+1}_{+} + \chi \cdot \LL_n^t \|\x_n^{t+1} - \x_n^{t}\|_2^2      + \UUU^t    ,
\end{align}
\noi where $\chi\triangleq \tfrac{2\omega \kappa}{\sigma}\cdot  \{ \sigma^2 q^2   + 3 \delta^2 + 3( \delta + \sigma q  )^2 \}$, $q \triangleq   \theta_2 + \theta_2 \delta $, $\Theta^t_{+} \triangleq a\AAA^t + b\BBB^t $, and $\UUU^t \triangleq  C_{h}^2 \tfrac{ b }{\beta^t}\cdot ( \tfrac{\mu^{t-1}}{\mu^t} - 1 )^2$. Here, $\{a,\AAA^t,b,\BBB^t\}$ are defined in Equation (\ref{eq:def:AAA}).

\end{lemma}

Assume Equation (\ref{eq:choice:SU}) is used to choose $\{p,\xi,\delta,\sigma,\theta_1,\theta_2\}$. We have the following lemma.

\begin{lemma} \label{lemma:rule:A}

(Proof in Appendix \ref{app:lemma:rule:A}) We define $\chi$ in Lemma \ref{lemma:bound:dual:Case:A}. We have the following results:

\begin{enumerate}[label=\textbf{(\alph*)}, leftmargin=22pt, itemsep=1pt, topsep=1pt, parsep=0pt, partopsep=0pt]

\item It holds that $\varepsilon_1 \triangleq  \frac{1}{2} \theta_1- \frac{1}{2}>0$, and $\varepsilon_2 \triangleq \theta_2 - \frac{1}{2} - \chi > 0$.

\item For all $t\geq 1$, we have: $\mathcal{E}^{t+1}      \leq  \Theta^t - \Theta^{t+1} + \UUU^t$.

\end{enumerate}

\end{lemma}

\subsection{Continuing Analysis for Conditions $\BIBI$ and $\SUSU$}

Using Assumption \ref{ass:5}, we show that $\Theta^{t}$ is consistently lower-bounded by the following lemma.

\begin{lemma} \label{lemma:Theta:LB}
(Proof in Appendix \ref{app:lemma:Theta:LB}) For all $t\geq 1$, there exists a constant $\underline{\Theta}$ such that $\Theta^{t}\geq\underline{\Theta}$.
\end{lemma}

The following lemma shows that both $\left(\sum_{t=1}^{\infty} \UUU^t\right)$ and $\left(\sum_{t=1}^{\infty} \mathcal{E}^{t+1}\right)$ are always upper-bounded.

\begin{lemma}\label{lemma:beta:mu}
(Proof in Appendix \ref{app:lemma:beta:mu}) We define $\UUU^t$ as in Lemma \ref{lemma:bound:dual:Case:I} or Lemma \ref{lemma:bound:dual:Case:A}. We define $\overline{\mathcal{E}}$ as in Lemma \ref{lemma:suf:dec}. We have:

\begin{enumerate}[label=\textbf{(\alph*)}, leftmargin=22pt, itemsep=1pt, topsep=1pt, parsep=0pt, partopsep=0pt]

\item There exists a universal positive constant $\overline{\rm{U}}$ such that $\sum_{t=1}^{\infty} \UUU^t \leq \overline{\rm{U}}$.

\item Letting $\overline{\mathcal{E}}\triangleq \Theta^{1} - \underline{\Theta} + \overline{\rm{U}}$, we have: $\sum_{t=1}^{\infty} \mathcal{E}^{t+1} \leq \overline{\mathcal{E}}$.

\end{enumerate}
\end{lemma}

The following lemmas are useful to provide upper bounds for the dual and primal variables.

\begin{lemma} \label{lemma:boundedness:z}
(Proof in Appendix \ref{app:lemma:boundedness:z}) For all $t\geq 1$, there exist constants $\{Z,\ddot{Z}\}$ such that $\tfrac{1}{ \beta^t} \|\z^t\|_2^2 \leq Z$, and $\sum_{t=1}^{\infty} \tfrac{1}{\beta^t} \|\z^{t+1} - \z^t\|_2^2 \leq \ddot{Z}$.
\end{lemma}

\begin{lemma} \label{lemma:boundedness}
(Proof in Appendix \ref{app:lemma:boundedness}) For all $i\in[n]$, we have $\|\x^{t+1}_i\|<+\infty$.
\end{lemma}


Finally, we have the following theorem regrading to the global convergence of IPDS-ADMM.


\begin{theorem} \label{theorem:continuing:analysis}
(Proof in Appendix \ref{app:theorem:continuing:analysis}) We have the following results.

\begin{enumerate}[label=\textbf{(\alph*)}, leftmargin=22pt, itemsep=1pt, topsep=1pt, parsep=0pt, partopsep=0pt]

\item $\sum_{t=1}^T \|\z^{t+1} - \z^t\|_2^2 + \|\beta^t(\x^{t+1} - \x^t)\|_2^2 \leq K\beta^T$, where $K>0$ is some constant. 

\item There exists an index $\bar{t}$ with $\bar{t}\leq T$ such that $\|  \z^{\bar{t}+1}-\z^{\bar{t}}  \|_2^2 + \| \beta^{\bar{t}}(\x^{\bar{t}+1}-\x^{\bar{t}})\|_2^2 \leq \ts   \frac{K\beta^T}{T}$.

\end{enumerate}



\end{theorem}

\begin{remark}
(\bfit{i}) With the choice $\beta^T = \mathcal{O}(T^p)$ with $p\in(0,1)$, we observe $\ddot{e}^{\bar{t}}\triangleq \|  \z^{\bar{t}+1}-\z^{\bar{t}}  \|_2^2 + \| \beta^{\bar{t}}(\x^{\bar{t}+1}-\x^{\bar{t}})\|_2^2 = \mathcal{O} (T^{p-1})$, indicating convergence of $\ddot{e}^{\bar{t}}$ towards 0. (\bfit{ii}) In light of Theorem \ref{theorem:continuing:analysis}, a reasonable stopping criterion for Algorithm \ref{alg:main} is $\|\z^{\bar{t}+1}-\z^{\bar{t}}\| + \| \beta^{\bar{t}}(\x^{\bar{t}+1}-\x^{\bar{t}})\|\leq \epsilon$, where $\epsilon\in(0,1)$ is a user-defined parameter.

\end{remark}
\subsection{Iteration Complexity}

We now establish the iteration complexity of Algorithm \ref{alg:main}. We first restate the following standard definition of approximated critical points.

\begin{definition}
($\epsilon$-\rm{Critical Point}) We define $\Crit(\check{\x},\check{\mathbf{z}}) \triangleq \|\A\check{\x}- \b \|+\sum_{i=1}^n \dist(\zero, \nabla f_i(\check{\x}_i) + \partial h_i(\check{\x}_i)   + \A_i\trans \check{\mathbf{z}})$. A solution $(\check{\x},\check{\mathbf{z}})$ is an $\epsilon$-\rm{critical} point if it holds that: $$\Crit(\check{\x},\check{\mathbf{z}}) \leq \epsilon.$$
\end{definition}

We obtain the following iteration complexity results.

\begin{theorem} \label{theorem:case:AAAAA:2}

(Proof in Appendix \ref{app:theorem:case:AAAAA:2}) We define $\mathbf{q}^{t} \triangleq \{\x^{t}_1,\x^{t}_2,\ldots,\x^{t}_{n-1},\breve{\x}^{t}_n\}$. Let the sequence $\{\mathbf{q}^{t},\z^t\}_{t=0}^T$ be generated by Algorithm \ref{alg:main}. For all $p\in(0,1)$, we have: 
\beq\ts \frac{1}{T} \sum_{t=1}^T \Crit(\mathbf{q}^{t+1},\mathbf{z}^{t+1})  \leq \ts \mathcal{O}(T^{(p-1)/2}) + \mathcal{O}(T^{-p}). \label{eq:ppppp}
\eeq
In particular, with the choice $p=1/3$, we have $\frac{1}{T} \sum_{t=1}^T \Crit(\mathbf{q}^{t+1},\mathbf{z}^{t+1}) \leq \mathcal{O}(T^{-1/3})$. In other words, there exists $\bar{t} \leq T$ such that: $\Crit(\mathbf{q}^{\bar{t}+1},\mathbf{z}^{\bar{t}+1})  \leq \epsilon$, provided that $T\geq\mathcal{O}(\epsilon^{-3})$.




\end{theorem}

\begin{remark}
(\bfit{i}) Minimizing the worst-case complexity of the right-hand side of Inequality (\ref{eq:ppppp}) \textit{w.r.t.} $p$ yields: $\arg \min_{p\in(0,1)}\max(-p,(p-1)/2)=1/3$. Thus, choosing $p=1/3$ achieves the optimal trade-off between the two terms, resulting in the best complexity bounds. (\bfit{ii}) To the best of our knowledge, this represents the first complexity result for using ADMM to solve this class of nonsmooth and nonconvex problems. Remarkably, we observe that it aligns with the iteration bound found in smoothing proximal gradient methods \cite{BohmW21}.

\end{remark}


\subsection{On the Boundedness and Convergence of the Multipliers}

Questions may arise regarding whether the multipliers $\z^t$ in Algorithm \ref{alg:main} are bounded, given that $\|\z^t\|_2^2 \leq Z \beta^t$, as stated in Lemma \ref{lemma:boundedness:z}. We argue that the bounedness of the multipliers is not an issue. We propose the following variable substitution: $\tfrac{\z^t}{\sqrt{\beta^t}} \triangleq \hat{\z}^t$ for all $t$. Consequently, we can implement the following update rule to replace the dual variable update rule of Algorithm \ref{alg:main}: $\ts \hat{\z}^{t+1}   =  \hat{\z}^t  \tfrac{\sqrt{\beta^t}}{\sqrt{\beta^{t+1}}}  + \tfrac{\beta^t}{\sqrt{\beta^{t+1}}} \cdot \sigma  ( \A\x^{t+1} - \b)$. Additionally, $\z^t$ should be replaced with $\sqrt{\beta^t}\cdot \hat{\z}^t$ in the remaining steps of Algorithm \ref{alg:main}. Importantly, such a substitution does not essentially alter the algorithm or our analysis throughout this paper.

We have the following results for the new multipliers $\hat{\z}^t$:

\begin{lemma}
\label{lemma:bound:u}
(Proof in Appendix \ref{app:lemma:bound:u}) We have: (\bfit{a}) $\forall t\geq0,~\|\hat{\z}^t\|_2^2 \leq Z$; (\bfit{b}) $\sum_{t=1}^{\infty} \|\hat{\z}^{t+1} - \hat{\z}^t\|_2^2 \leq 2 (\ddot{Z} +  Z)$. Here, $\{\ddot{Z},Z\}$ are bounded constants defined in Lemma \ref{lemma:boundedness:z}.

\end{lemma}

\begin{remark}
Thanks to the variable substitution, the new multiplier $\|\hat{\z}^t\|$ is bounded and convergent with $\left(\min_{t=1}^T \|\hat{\z}^{t+1} - \hat{\z}^t\|_2^2\right) \leq \frac{1}{T}\sum_{t=1}^{T} \|\hat{\z}^{t+1} - \hat{\z}^t\|_2^2 \leq \mathcal{O}(1/T)$.
\end{remark}

\section{Experiments} \label{sect:exp}

This section assesses the performance of IPDS-ADMM in solving the sparse PCA problem, as shown in Section \ref{sect:app}.

\noi $\blacktriangleright$ \textbf{Compared Methods}. We compare IPDS-ADMM against three state-of-the-art general-purpose algorithms that solve Problem (\ref{eq:main}): (\bfit{i}) the Subgradient method (SubGrad) \cite{li2021weakly,davis2019stochastic}, \bfit{(ii)} the Smoothing Proximal Gradient Method (SPGM) \cite{BohmW21}, (\bfit{iii}) the Riemannian ADMM with fixed and large penalty (RADMM) \cite{li2022riemannian}.

\noi $\blacktriangleright$ \textbf{Experimental Settings}. All methods are implemented in MATLAB on an Intel 2.6 GHz CPU with 64 GB RAM. We incorporate a set of 8 datasets into our experiments, comprising both
randomly generated and publicly available real-world data. Appendix Section \ref{app:sect:exp} describes how to generate the data used in the experiments. For IPDS-ADMM, the relaxation parameter $\sigma$ is set to be around the golden ratio $1.618$, as suggested by \cite{li2016majorized}. Additionally, we set $(\xi,p,\delta,\theta_1,\theta_2)=(1/2,1/3,1/4,1.01,0.60)$. We denote $\dot{\rho}$ as the regularization parameter for sparse PCA model. The penalty parameter for RADMM is set to a reasonably large constant $\beta=100\dot{\rho}$. We fix $\dot{r}=20$ and compare objective values for all methods after running $T'$ seconds, where $T'$ is reasonably large to ensure the proposed method converges. The corresponding MATLAB code is available on the author's research webpage.


\vspace{-8pt}
\begin{figure}[!h]
\centering
\begin{subfigure}{.24\textwidth}\centering\includegraphics[width=1.12\linewidth]{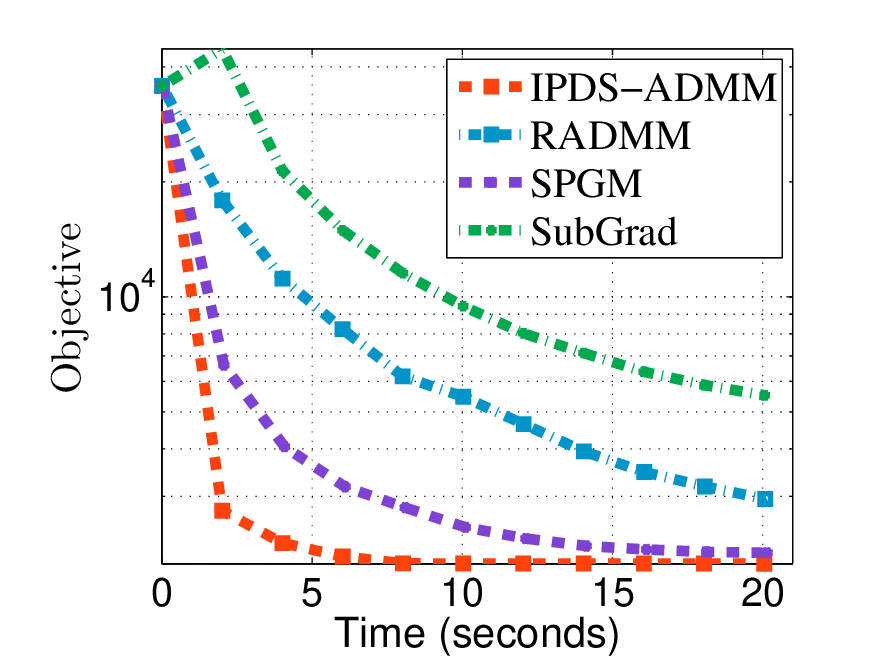}\caption{\scriptsize randn-1500-500}\label{fig:sub1}\end{subfigure}
\begin{subfigure}{.24\textwidth}\centering\includegraphics[width=1.12\linewidth]{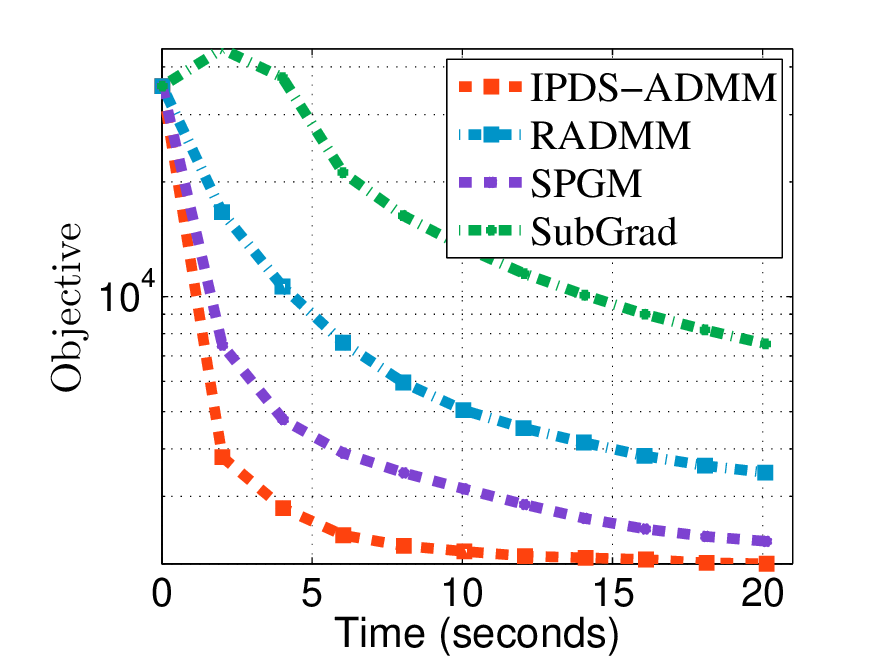}\caption{\scriptsize randn-2500-500}\label{fig:sub2}\end{subfigure}
\begin{subfigure}{.24\textwidth}\centering\includegraphics[width=1.12\linewidth]{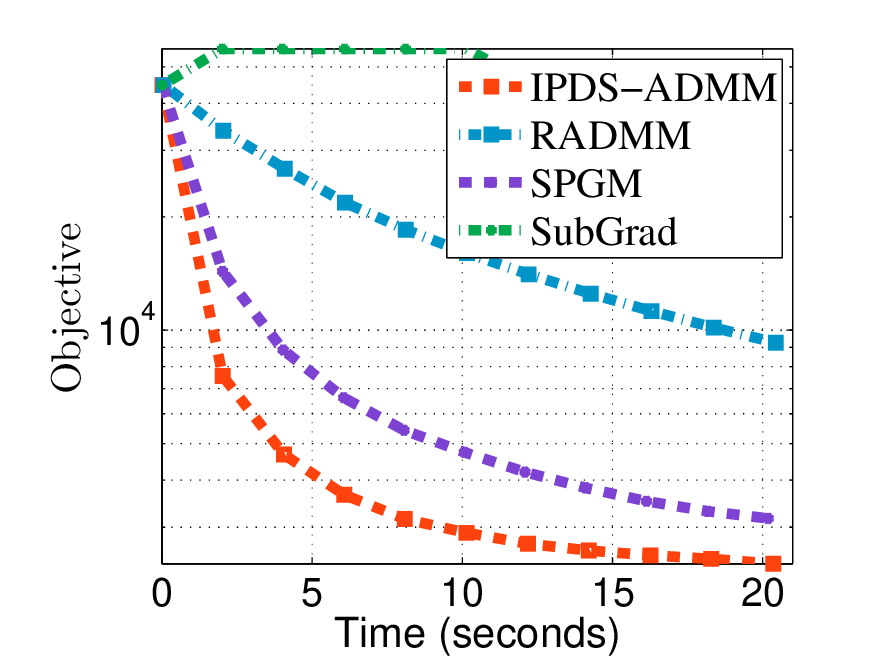}\caption{\scriptsize mnist-1500-780}\label{fig:sub3}\end{subfigure}
\begin{subfigure}{.24\textwidth}\centering\includegraphics[width=1.12\linewidth]{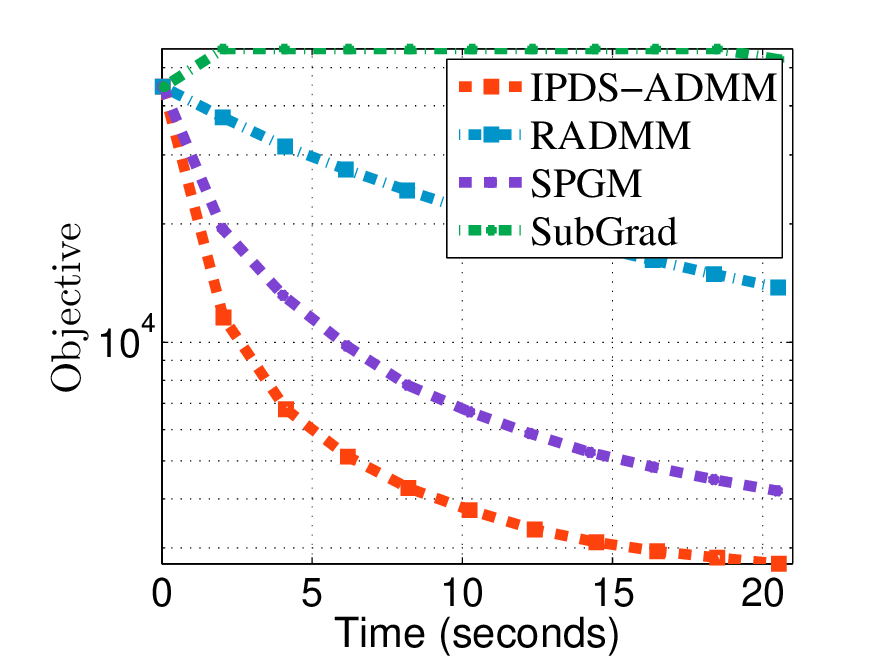}\caption{\scriptsize mnist-2500-780}\label{fig:sub4}\end{subfigure}

\centering
\begin{subfigure}{.24\textwidth}\centering\includegraphics[width=1.12\linewidth]{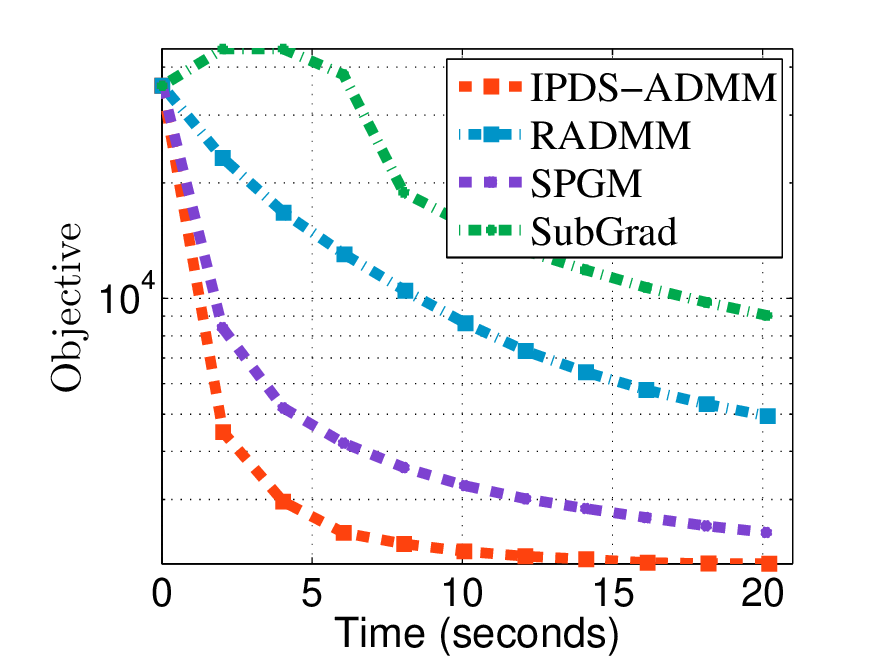}\caption{\scriptsize TDT2-1500-500}\label{fig:sub1}\end{subfigure}
\begin{subfigure}{.24\textwidth}\centering\includegraphics[width=1.12\linewidth]{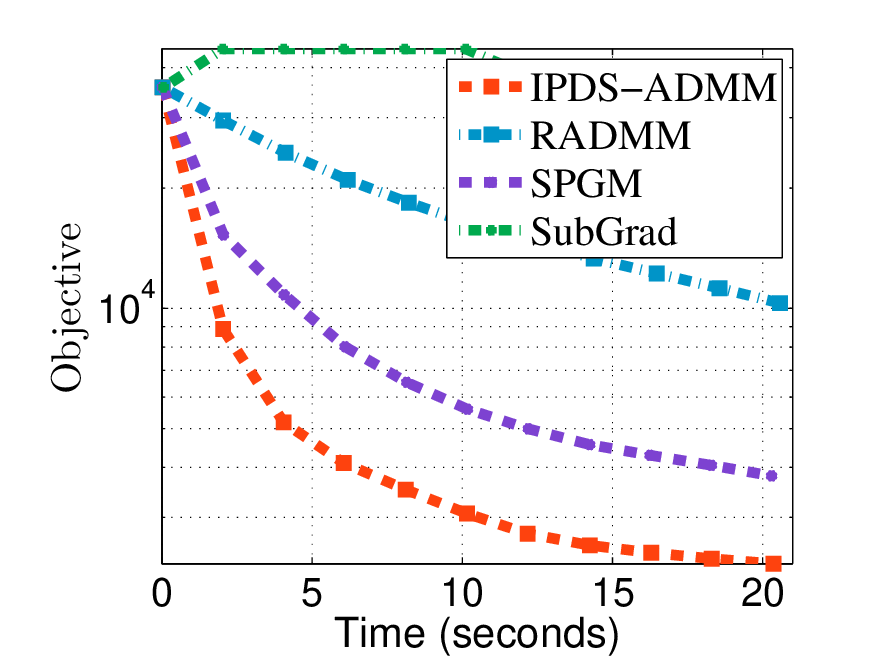}\caption{\scriptsize TDT2-3000-500}\label{fig:sub2}\end{subfigure}
\begin{subfigure}{.24\textwidth}\centering\includegraphics[width=1.12\linewidth]{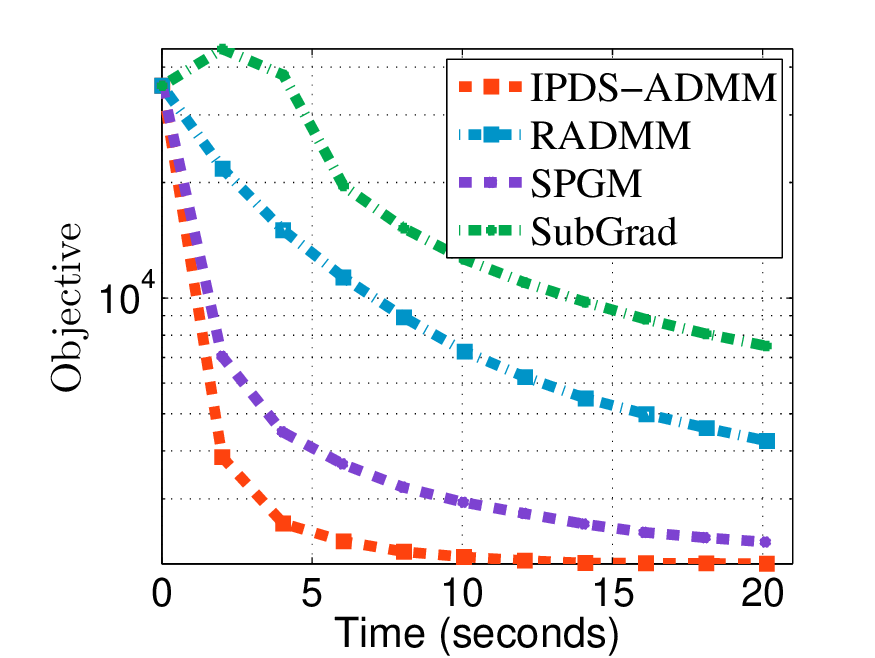}\caption{\scriptsize sector-1500-500}\label{fig:sub3}\end{subfigure}
\begin{subfigure}{.24\textwidth}\centering\includegraphics[width=1.12\linewidth]{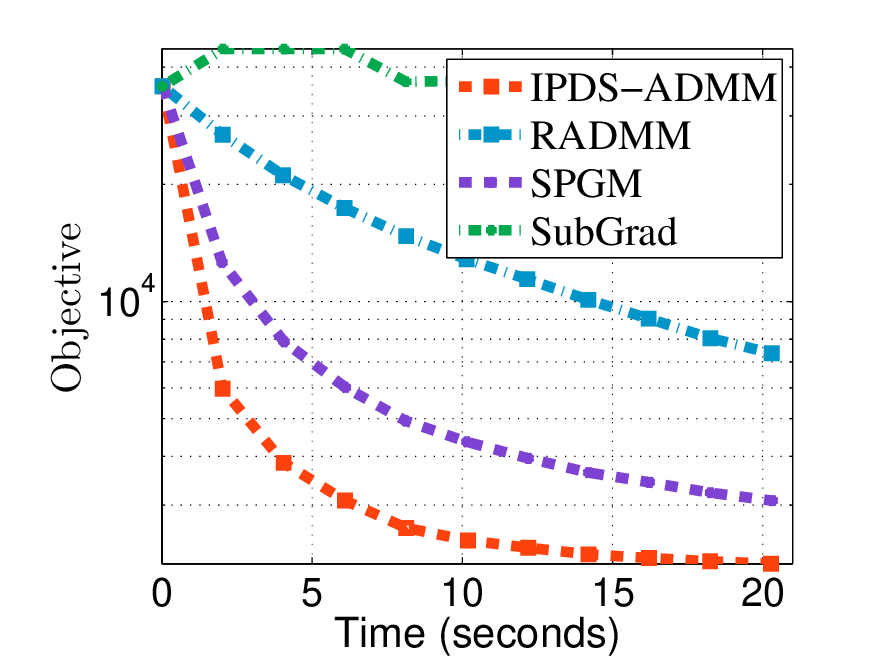}\caption{\scriptsize sector-2500-500}\label{fig:sub4}\end{subfigure}

\caption{Convergence curves of methods for sparse PCA with $\dot{\rho}=100$ and $\beta^0=50\dot{\rho}$.} \label{fig:000}

\end{figure}

\noi $\blacktriangleright$ \textbf{Experiment Results}. We set $\dot{\rho} = 100$ and present the results for $\beta^0=50\dot{\rho}$ in IPDS-ADMM, as shown in Figure \ref{fig:000}. These experimental results provide the following insights: (\bfit{i}) Sub-Grad tends to be less efficient in comparison to other methods. (\bfit{ii}) SPGM, utilizing a variable smoothing strategy, generally demonstrates slower performance than the multiplier-based variable splitting method. This observation corroborates the widely accepted notion that primal-dual methods are typically more robust and quicker than primal-only methods. (\bfit{iii}) The proposed IPDS-ADMM consistently achieves the lowest objective function values among all methods examined. Due to space limitations, additional experimental results are provided in Appendix Section \ref{app:sect:exp}.

%
%
%


\vspace{-4pt}
\section{Conclusions} %
\vspace{-4pt}
In this paper, we introduce IPDS-ADMM, a proximal linearized ADMM that uses an Increasing Penalization and Decreasing Smoothing (IPDS) strategy for solving general multi-block nonconvex composite optimization problems. IPDS-ADMM operates under a relatively relaxed condition, requiring continuity in just one block of the objective function. It incorporates relaxed strategies for dual variable updates when the associated linear operator is either bijective or surjective. We increase the penalty parameter and decrease the smoothing parameter at a controlled pace, and introduce a Lyapunov function for convergence analysis. We also derive the iteration complexity of IPDS-ADMM. Finally, we conduct experiments to demonstrate the effectiveness of our approaches.

\normalem

\section*{Acknowledgments}

This work was supported by NSFC (12271278, 61772570), and Guangdong Natural Science Funds for Distinguished Young Scholar (2018B030306025).

\bibliographystyle{iclr2025_conference}
\bibliography{mybib}

\begin{thebibliography}{64}
\providecommand{\natexlab}[1]{#1}
\providecommand{\url}[1]{\texttt{#1}}
\expandafter\ifx\csname urlstyle\endcsname\relax
  \providecommand{\doi}[1]{doi: #1}\else
  \providecommand{\doi}{doi: \begingroup \urlstyle{rm}\Url}\fi

\bibitem[Barber \& Sidky(2024)Barber and Sidky]{barber2024convergence}
Rina~Foygel Barber and Emil~Y Sidky.
\newblock Convergence for nonconvex admm, with applications to ct imaging.
\newblock \emph{Journal of Machine Learning Research}, 25\penalty0
  (38):\penalty0 1--46, 2024.

\bibitem[Beck(2017)]{beck2017first}
Amir Beck.
\newblock \emph{First-order methods in optimization}.
\newblock SIAM, 2017.

\bibitem[Bertsekas(2015)]{bertsekas2015convex}
Dimitri Bertsekas.
\newblock \emph{Convex optimization algorithms}.
\newblock Athena Scientific, 2015.

\bibitem[Bian et~al.(2021)Bian, Liang, and Zhang]{bian2021stochastic}
Fengmiao Bian, Jingwei Liang, and Xiaoqun Zhang.
\newblock A stochastic alternating direction method of multipliers for
  non-smooth and non-convex optimization.
\newblock \emph{Inverse Problems}, 37\penalty0 (7):\penalty0 075009, 2021.

\bibitem[Bo$\c{t}$ \& Nguyen(2020)Bo$\c{t}$ and Nguyen]{boct2020proximal}
Radu~Ioan Bo$\c{t}$ and Dang-Khoa Nguyen.
\newblock The proximal alternating direction method of multipliers in the
  nonconvex setting: convergence analysis and rates.
\newblock \emph{Mathematics of Operations Research}, 45\penalty0 (2):\penalty0
  682--712, 2020.

\bibitem[Bo$\c{t}$ et~al.(2019)Bo$\c{t}$, Csetnek, and Nguyen]{Boct2019SIOPT}
Radu~Ioan Bo$\c{t}$, Erno~Robert Csetnek, and Dang-Khoa Nguyen.
\newblock A proximal minimization algorithm for structured nonconvex and
  nonsmooth problems.
\newblock \emph{SIAM Journal on Optimization}, 29\penalty0 (2):\penalty0
  1300--1328, 2019.

\bibitem[B{\"{o}}hm \& Wright(2021)B{\"{o}}hm and Wright]{BohmW21}
Axel B{\"{o}}hm and Stephen~J. Wright.
\newblock Variable smoothing for weakly convex composite functions.
\newblock \emph{Journal of Optimization Theory and Applications}, 188\penalty0
  (3):\penalty0 628--649, 2021.

\bibitem[Cand\`{e}s et~al.(2011)Cand\`{e}s, Li, Ma, and Wright]{Candes2011JACM}
E.~Cand\`{e}s, X.~Li, Y.~Ma, and J.~Wright.
\newblock Robust principal component analysis?
\newblock \emph{Journal of the {ACM}}, 58\penalty0 (3), May 2011.

\bibitem[Chen et~al.(2022)Chen, Shen, Zou, and Liu]{ChenJMLR2022}
Congliang Chen, Li~Shen, Fangyu Zou, and Wei Liu.
\newblock Towards practical adam: Non-convexity, convergence theory, and
  mini-batch acceleration.
\newblock \emph{Journal of Machine Learning Research}, 23\penalty0
  (229):\penalty0 1--47, 2022.

\bibitem[Chen et~al.(2016)Chen, Ji, and You]{chen2016augmented}
Weiqiang Chen, Hui Ji, and Yanfei You.
\newblock An augmented lagrangian method for $\ell_1$-regularized optimization
  problems with orthogonality constraints.
\newblock \emph{SIAM Journal on Scientific Computing}, 38\penalty0
  (4):\penalty0 B570--B592, 2016.

\bibitem[Davis \& Drusvyatskiy(2019)Davis and
  Drusvyatskiy]{davis2019stochastic}
Damek Davis and Dmitriy Drusvyatskiy.
\newblock Stochastic model-based minimization of weakly convex functions.
\newblock \emph{SIAM Journal on Optimization}, 29\penalty0 (1):\penalty0
  207--239, 2019.

\bibitem[Deng et~al.(2017)Deng, Lai, Peng, and Yin]{deng2017parallel}
Wei Deng, Ming-Jun Lai, Zhimin Peng, and Wotao Yin.
\newblock Parallel multi-block admm with o (1/k) convergence.
\newblock \emph{Journal of Scientific Computing}, 71:\penalty0 712--736, 2017.

\bibitem[Duchi \& Ruan(2018)Duchi and Ruan]{Duchiiay015}
John~C Duchi and Feng Ruan.
\newblock {Solving (most) of a set of quadratic equalities: composite
  optimization for robust phase retrieval}.
\newblock \emph{Information and Inference: A Journal of the IMA}, 8\penalty0
  (3):\penalty0 471--529, 2018.

\bibitem[Franca et~al.(2018)Franca, Robinson, and Vidal]{franca2018admm}
Guilherme Franca, Daniel Robinson, and Rene Vidal.
\newblock Admm and accelerated admm as continuous dynamical systems.
\newblock In \emph{International conference on machine learning}, pp.\
  1559--1567, 2018.

\bibitem[Gabay \& Mercier(1976)Gabay and Mercier]{gabay1976dual}
Daniel Gabay and Bertrand Mercier.
\newblock A dual algorithm for the solution of nonlinear variational problems
  via finite element approximation.
\newblock \emph{Computers \& Mathematics with Applications}, 2\penalty0
  (1):\penalty0 17--40, 1976.

\bibitem[Gon{\c{c}}alves et~al.(2017{\natexlab{a}})Gon{\c{c}}alves, Melo, and
  Monteiro]{goncalves2017convergence}
Max~LN Gon{\c{c}}alves, Jefferson~G Melo, and Renato~DC Monteiro.
\newblock Convergence rate bounds for a proximal admm with over-relaxation
  stepsize parameter for solving nonconvex linearly constrained problems.
\newblock \emph{arXiv:1702.01850}, 2017{\natexlab{a}}.

\bibitem[Gon{\c{c}}alves et~al.(2017{\natexlab{b}})Gon{\c{c}}alves, Melo, and
  Monteiro]{gonccalves2017improved}
Max~LN Gon{\c{c}}alves, Jefferson~G Melo, and Renato~DC Monteiro.
\newblock Improved pointwise iteration-complexity of a regularized admm and of
  a regularized non-euclidean hpe framework.
\newblock \emph{SIAM Journal on Optimization}, 27\penalty0 (1):\penalty0
  379--407, 2017{\natexlab{b}}.

\bibitem[Gong et~al.(2013)Gong, Zhang, Lu, Huang, and Ye]{GongZLHY13}
Pinghua Gong, Changshui Zhang, Zhaosong Lu, Jianhua Huang, and Jieping Ye.
\newblock A general iterative shrinkage and thresholding algorithm for
  non-convex regularized optimization problems.
\newblock In \emph{International Conference on Machine Learning (ICML)},
  volume~28, pp.\  37--45, 2013.

\bibitem[He \& Yuan(2012)He and Yuan]{HeY12}
Bingsheng He and Xiaoming Yuan.
\newblock On the $\mathcal{O}(1/n)$ convergence rate of the douglas-rachford
  alternating direction method.
\newblock \emph{SIAM Journal on Numerical Analysis}, 50\penalty0 (2):\penalty0
  700--709, 2012.

\bibitem[Hien et~al.(2022)Hien, Phan, and Gillis]{hien2021framework}
Le~Thi~Khanh Hien, Duy~Nhat Phan, and Nicolas Gillis.
\newblock Inertial alternating direction method of multipliers for non-convex
  non-smooth optimization.
\newblock \emph{Computational Optimization and Applications}, 83\penalty0
  (1):\penalty0 247--285, 2022.

\bibitem[Hong et~al.(2016)Hong, Luo, and Razaviyayn]{hong2016convergence}
Mingyi Hong, Zhi-Quan Luo, and Meisam Razaviyayn.
\newblock Convergence analysis of alternating direction method of multipliers
  for a family of nonconvex problems.
\newblock \emph{SIAM Journal on Optimization}, 26\penalty0 (1):\penalty0
  337--364, 2016.

\bibitem[Huang et~al.(2019)Huang, Chen, and Huang]{HuangCH19}
Feihu Huang, Songcan Chen, and Heng Huang.
\newblock Faster stochastic alternating direction method of multipliers for
  nonconvex optimization.
\newblock In \emph{International Conference on Machine Learning (ICML)},
  volume~97, pp.\  2839--2848, 2019.

\bibitem[Kingma \& Ba(2015)Kingma and Ba]{KingmaB14}
Diederik~P. Kingma and Jimmy Ba.
\newblock Adam: {A} method for stochastic optimization.
\newblock In Yoshua Bengio and Yann LeCun (eds.), \emph{International
  Conference on Learning Representations (ICLR)}, 2015.

\bibitem[Lai \& Osher(2014)Lai and Osher]{lai2014splitting}
Rongjie Lai and Stanley Osher.
\newblock A splitting method for orthogonality constrained problems.
\newblock \emph{Journal of Scientific Computing}, 58\penalty0 (2):\penalty0
  431--449, 2014.

\bibitem[Lei~Yang(2021)]{YangMOR}
Shuhuang~Xiang Lei~Yang, Xiaojun~Chen.
\newblock Sparse solutions of a class of constrained optimization problems.
\newblock \emph{Mathematics of Operations Research}, 2021.

\bibitem[Li \& Pong(2015)Li and Pong]{li2015global}
Guoyin Li and Ting~Kei Pong.
\newblock Global convergence of splitting methods for nonconvex composite
  optimization.
\newblock \emph{SIAM Journal on Optimization}, 25\penalty0 (4):\penalty0
  2434--2460, 2015.

\bibitem[Li et~al.(2022)Li, Ma, and Srivastava]{li2022riemannian}
Jiaxiang Li, Shiqian Ma, and Tejes Srivastava.
\newblock A riemannian admm.
\newblock \emph{arXiv preprint arXiv:2211.02163}, 2022.

\bibitem[Li et~al.(2016)Li, Sun, and Toh]{li2016majorized}
Min Li, Defeng Sun, and Kim-Chuan Toh.
\newblock A majorized admm with indefinite proximal terms for linearly
  constrained convex composite optimization.
\newblock \emph{SIAM Journal on Optimization}, 26\penalty0 (2):\penalty0
  922--950, 2016.

\bibitem[Li et~al.(2021)Li, Chen, Deng, Qu, Zhu, and Man-Cho~So]{li2021weakly}
Xiao Li, Shixiang Chen, Zengde Deng, Qing Qu, Zhihui Zhu, and Anthony
  Man-Cho~So.
\newblock Weakly convex optimization over stiefel manifold using riemannian
  subgradient-type methods.
\newblock \emph{SIAM Journal on Optimization}, 31\penalty0 (3):\penalty0
  1605--1634, 2021.

\bibitem[Lin et~al.(2022{\natexlab{a}})Lin, Ma, and Xu]{Lin2022}
Qihang Lin, Runchao Ma, and Yangyang Xu.
\newblock Complexity of an inexact proximal-point penalty method for
  constrained smooth non-convex optimization.
\newblock \emph{Computational Optimization and Applications}, 82\penalty0
  (1):\penalty0 175--224, 2022{\natexlab{a}}.

\bibitem[Lin et~al.(2015{\natexlab{a}})Lin, Ma, and Zhang]{lin2014convergence}
Tian-Yi Lin, Shi-Qian Ma, and Shu-Zhong Zhang.
\newblock On the sublinear convergence rate of multi-block admm.
\newblock \emph{Journal of the Operations Research Society of China},
  3:\penalty0 251--274, 2015{\natexlab{a}}.

\bibitem[Lin et~al.(2015{\natexlab{b}})Lin, Ma, and Zhang]{lin2015global}
Tianyi Lin, Shiqian Ma, and Shuzhong Zhang.
\newblock On the global linear convergence of the admm with multiblock
  variables.
\newblock \emph{SIAM Journal on Optimization}, 25\penalty0 (3):\penalty0
  1478--1497, 2015{\natexlab{b}}.

\bibitem[Lin et~al.(2022{\natexlab{b}})Lin, Li, and Fang]{lin2022alternating}
Zhouchen Lin, Huan Li, and Cong Fang.
\newblock \emph{Alternating direction method of multipliers for machine
  learning}.
\newblock Springer, 2022{\natexlab{b}}.

\bibitem[Liu et~al.(2019)Liu, Li, and Shen]{liu2019one}
Dekai Liu, Song Li, and Yi~Shen.
\newblock One-bit compressive sensing with projected subgradient method under
  sparsity constraints.
\newblock \emph{IEEE Transactions on Information Theory}, 65\penalty0
  (10):\penalty0 6650--6663, 2019.

\bibitem[Liu et~al.(2022)Liu, Liu, and Chen]{LiuLC22}
Wei Liu, Xin Liu, and Xiaojun Chen.
\newblock Linearly constrained nonsmooth optimization for training
  autoencoders.
\newblock \emph{SIAM Journal on Optimization}, 32\penalty0 (3):\penalty0
  1931--1957, 2022.

\bibitem[Liu et~al.(2020)Liu, Shang, Liu, Kong, Jiao, and
  Lin]{liu2020accelerated}
Yuanyuan Liu, Fanhua Shang, Hongying Liu, Lin Kong, Licheng Jiao, and Zhouchen
  Lin.
\newblock Accelerated variance reduction stochastic admm for large-scale
  machine learning.
\newblock \emph{IEEE Transactions on Pattern Analysis and Machine
  Intelligence}, 43\penalty0 (12):\penalty0 4242--4255, 2020.

\bibitem[Lu \& Zhang(2012)Lu and Zhang]{lu2012augmented}
Zhaosong Lu and Yong Zhang.
\newblock An augmented lagrangian approach for sparse principal component
  analysis.
\newblock \emph{Mathematical Programming}, 135:\penalty0 149--193, 2012.

\bibitem[Lu \& Zhang(2013)Lu and Zhang]{lu2013sparse}
Zhaosong Lu and Yong Zhang.
\newblock Sparse approximation via penalty decomposition methods.
\newblock \emph{SIAM Journal on Optimization}, 23\penalty0 (4):\penalty0
  2448--2478, 2013.

\bibitem[Monteiro \& Svaiter(2013)Monteiro and Svaiter]{monteiro2013iteration}
Renato~DC Monteiro and Benar~F Svaiter.
\newblock Iteration-complexity of block-decomposition algorithms and the
  alternating direction method of multipliers.
\newblock \emph{SIAM Journal on Optimization}, 23\penalty0 (1):\penalty0
  475--507, 2013.

\bibitem[Mordukhovich(2006)]{Mordukhovich2006}
Boris~S. Mordukhovich.
\newblock Variational analysis and generalized differentiation i: Basic theory.
\newblock \emph{Berlin Springer}, 330, 2006.

\bibitem[Nesterov(2003)]{Nesterov03}
Y.~E. Nesterov.
\newblock \emph{Introductory lectures on convex optimization: a basic course},
  volume~87 of \emph{Applied Optimization}.
\newblock Kluwer Academic Publishers, 2003.

\bibitem[Nishihara et~al.(2015)Nishihara, Lessard, Recht, Packard, and
  Jordan]{nishihara2015general}
Robert Nishihara, Laurent Lessard, Ben Recht, Andrew Packard, and Michael
  Jordan.
\newblock A general analysis of the convergence of admm.
\newblock In \emph{International Conference on Machine Learning (ICML)}, pp.\
  343--352, 2015.

\bibitem[Ouyang et~al.(2015)Ouyang, Chen, Lan, and
  Pasiliao~Jr]{ouyang2015accelerated}
Yuyuan Ouyang, Yunmei Chen, Guanghui Lan, and Eduardo Pasiliao~Jr.
\newblock An accelerated linearized alternating direction method of
  multipliers.
\newblock \emph{SIAM Journal on Imaging Sciences}, 8\penalty0 (1):\penalty0
  644--681, 2015.

\bibitem[Pock \& Sabach(2016)Pock and Sabach]{pock2016inertial}
Thomas Pock and Shoham Sabach.
\newblock Inertial proximal alternating linearized minimization (ipalm) for
  nonconvex and nonsmooth problems.
\newblock \emph{SIAM Journal on Imaging Sciences}, 9\penalty0 (4):\penalty0
  1756--1787, 2016.

\bibitem[Rockafellar \& Wets.(2009)Rockafellar and Wets.]{Rockafellar2009}
R.~Tyrrell Rockafellar and Roger J-B. Wets.
\newblock Variational analysis.
\newblock \emph{Springer Science \& Business Media}, 317, 2009.

\bibitem[Shen et~al.(2017)Shen, Liu, Yuan, and Ma]{shen2017gsos}
Li~Shen, Wei Liu, Ganzhao Yuan, and Shiqian Ma.
\newblock Gsos: Gauss-seidel operator splitting algorithm for multi-term
  nonsmooth convex composite optimization.
\newblock In \emph{International Conference on Machine Learning (ICML)}, pp.\
  3125--3134, 2017.

\bibitem[Sun \& Sun(2023)Sun and Sun]{sun2021algorithms}
Kaizhao Sun and Xu~Andy Sun.
\newblock Algorithms for difference-of-convex programs based on
  difference-of-moreau-envelopes smoothing.
\newblock \emph{INFORMS Journal on Optimization}, 5\penalty0 (4):\penalty0
  321--339, 2023.

\bibitem[Tran~Dinh(2018)]{tran2018non}
Quoc Tran~Dinh.
\newblock Non-ergodic alternating proximal augmented lagrangian algorithms with
  optimal rates.
\newblock \emph{Advances in Neural Information Processing Systems (NeurIPS)},
  31, 2018.

\bibitem[Tsakiris \& Vidal(2018)Tsakiris and Vidal]{TsakirisV18}
Manolis~C. Tsakiris and Ren{\'{e}} Vidal.
\newblock Dual principal component pursuit.
\newblock \emph{Journal of Machine Learning Research}, 19:\penalty0
  18:1--18:50, 2018.

\bibitem[Wang et~al.(2019{\natexlab{a}})Wang, Yu, Chen, and Zhao]{WangYCZ19}
Junxiang Wang, Fuxun Yu, Xiang Chen, and Liang Zhao.
\newblock {ADMM} for efficient deep learning with global convergence.
\newblock In \emph{ACM International Conference on Knowledge Discovery {\&}
  Data Mining (SIGKDD)}, pp.\  111--119, 2019{\natexlab{a}}.

\bibitem[Wang et~al.(2019{\natexlab{b}})Wang, Yin, and Zeng]{wang2019global}
Yu~Wang, Wotao Yin, and Jinshan Zeng.
\newblock Global convergence of admm in nonconvex nonsmooth optimization.
\newblock \emph{Journal of Scientific Computing}, 78\penalty0 (1):\penalty0
  29--63, 2019{\natexlab{b}}.

\bibitem[Xu et~al.(2017)Xu, Liu, Lin, and Yang]{NIPS2017Xu}
Yi~Xu, Mingrui Liu, Qihang Lin, and Tianbao Yang.
\newblock Admm without a fixed penalty parameter: Faster convergence with new
  adaptive penalization.
\newblock In \emph{Advances in Neural Information Processing Systems
  (NeurIPS)}, volume~30, 2017.

\bibitem[Yang et~al.(2017)Yang, Pong, and Chen]{yang2017alternating}
Lei Yang, Ting~Kei Pong, and Xiaojun Chen.
\newblock Alternating direction method of multipliers for a class of nonconvex
  and nonsmooth problems with applications to background/foreground extraction.
\newblock \emph{SIAM Journal on Imaging Sciences}, 10\penalty0 (1):\penalty0
  74--110, 2017.

\bibitem[Yashtini(2021)]{yashtini2021multi}
Maryam Yashtini.
\newblock Multi-block nonconvex nonsmooth proximal admm: Convergence and rates
  under kurdyka--$\l$ojasiewicz property.
\newblock \emph{Journal of Optimization Theory and Applications}, 190\penalty0
  (3):\penalty0 966--998, 2021.

\bibitem[Yashtini(2022)]{yashtini2020convergence}
Maryam Yashtini.
\newblock Convergence and rate analysis of a proximal linearized {ADMM} for
  nonconvex nonsmooth optimization.
\newblock \emph{Journal of Global Optimization}, 84\penalty0 (4):\penalty0
  913--939, 2022.

\bibitem[Yuan(2024)]{yuan2024smoothing}
Ganzhao Yuan.
\newblock Smoothing proximal gradient methods for nonsmooth sparsity
  constrained optimization: Optimality conditions and global convergence.
\newblock In \emph{International Conference on Machine Learning (ICML)}, 2024.

\bibitem[Yuan(2025)]{yuanICLR25FADMM}
Ganzhao Yuan.
\newblock Admm for structured fractional minimization.
\newblock In \emph{International Conference on Learning Representations
  (ICLR)}, 2025.

\bibitem[Zeng et~al.(2014)Zeng, Lin, Wang, and Xu]{ZengLWX14}
Jinshan Zeng, Shaobo Lin, Yao Wang, and Zongben Xu.
\newblock $l_{1/2}$ regularization: Convergence of iterative half thresholding
  algorithm.
\newblock \emph{IEEE Transactions on Signal Processing}, 62\penalty0
  (9):\penalty0 2317--2329, 2014.

\bibitem[Zeng et~al.(2021)Zeng, Lin, Yao, and Zhou]{ZengLYZ21}
Jinshan Zeng, Shao{-}Bo Lin, Yuan Yao, and Ding{-}Xuan Zhou.
\newblock On {ADMM} in deep learning: Convergence and saturation-avoidance.
\newblock \emph{Journal of Machine Learning Research}, 22:\penalty0
  199:1--199:67, 2021.

\bibitem[Zeng et~al.(2022)Zeng, Yin, and Zhou]{zeng2022moreau}
Jinshan Zeng, Wotao Yin, and Ding-Xuan Zhou.
\newblock Moreau envelope augmented lagrangian method for nonconvex
  optimization with linear constraints.
\newblock \emph{Journal of Scientific Computing}, 91\penalty0 (2):\penalty0 61,
  2022.

\bibitem[Zhang \& Luo(2020)Zhang and Luo]{zhang2020proximal}
Jiawei Zhang and Zhi-Quan Luo.
\newblock A proximal alternating direction method of multiplier for linearly
  constrained nonconvex minimization.
\newblock \emph{SIAM Journal on Optimization}, 30\penalty0 (3):\penalty0
  2272--2302, 2020.

\bibitem[Zhang et~al.(2020)Zhang, Xiao, Sun, and Luo]{ZhangX0L20}
Jiawei Zhang, Peijun Xiao, Ruoyu Sun, and Zhi{-}Quan Luo.
\newblock A single-loop smoothed gradient descent-ascent algorithm for
  nonconvex-concave min-max problems.
\newblock In \emph{Advances in Neural Information Processing Systems
  (NeurIPS)}, 2020.

\bibitem[Zhang \& Kwok(2014)Zhang and Kwok]{zhang2014asynchronous}
Ruiliang Zhang and James Kwok.
\newblock Asynchronous distributed admm for consensus optimization.
\newblock In \emph{International Conference on Machine Learning (ICML)}, pp.\
  1701--1709, 2014.

\bibitem[Zhu et~al.(2023)Zhu, Zhao, and Zhang]{zhu2020first}
Daoli Zhu, Lei Zhao, and Shuzhong Zhang.
\newblock A first-order primal-dual method for nonconvex constrained
  optimization based on the augmented lagrangian.
\newblock \emph{Mathematics of Operations Research}, 2023.

\end{thebibliography}

\onecolumn

\appendix

{\huge Appendix}

The organization of the appendix is as follows:

Appendix \ref{app:sect:preliminaries:additional} covers notations, technical preliminaries, and relevant lemmas.

Appendix \ref{app:sect:app} provides additional motivating applications.

Appendix \ref{app:sect:proposed} contains proofs related to Section \ref{sect:proposed}.

Appendix \ref{app:sect:global:convergence} offers proofs related to Section \ref{sect:global:convergence}.

Appendix \ref{app:sect:exp} includes additional experiments details and results.






\section{Notations, Technical Preliminaries, and Relevant Lemmas} \label{app:sect:preliminaries:additional}

\subsection{Notations}

We use the following notations in this paper.

\begin{itemize}[leftmargin=12pt,itemsep=0.2ex]

\item $[n]$: $\{1,2,...,n\}$.

\item $\x$: $\x  \triangleq \{\x_1,\x_2,\ldots,\x_n\}= \x_{[n]}$.

\item $\x_{[i,j]}$: $\x_{[i,j]} \triangleq \{\x_{i},\x_{i+1},\x_{i+2},\ldots.,\x_{j}\}$, where $j\geq i$.


\item $\LL_i^t$: $\LL_i^t = L_i + \beta^t \|\A_i\|_2^2$. Note that the function $G(\x,\mathbf{z}^t;\beta^t)$ is $\LL_i^t$-smooth \textit{w.r.t.} $\x_i$.

\item $\sigma_1$: $\sigma_1 \triangleq \tfrac{\sigma}{(1-|1-\sigma|)^2} \in \mathbb{R}$, where $\sigma\in(0,2)$. Refer to Lemma \ref{useful:lemma:1}.

\item $\sigma_2$: $\sigma_2\triangleq \tfrac{|1-\sigma|}{\sigma (1-|1-\sigma|)}\in \mathbb{R}$, where $\sigma\in(0,2)$. Refer to Lemma \ref{useful:lemma:1}.

\item $\|\mathbf{x}\|$: Euclidean norm: $\|\mathbf{x}\|=\|\mathbf{x}\|_2 = \sqrt{\la \mathbf{x},\mathbf{x}\ra}$.

\item $\la \mathbf{a},\mathbf{b}\ra$ : Euclidean inner product, i.e., $\la \mathbf{a},\mathbf{b}\ra =\sum_{i}{\mathbf{a}_{i}\mathbf{b}_{i}}$.

\item $\A \trans$ : the transpose of the matrix $\A$.

\item $\x_i$: the $i$-th block of the vector $\x\in \Rn^{(\d_1+\d_2+\ldots+\d_n)\times 1}$ with $\x_i\in\Rn^{\d_i\times 1}$.

\item $\lambdaUp$: the largest eigenvalue of the matrix $\mathbf{A}_n\mathbf{A}_n\trans$.

\item $\lambdaDown$: the smallest eigenvalue of the matrix $\mathbf{A}_n\mathbf{A}_n\trans$.

\item $\lambdaDown'$: the smallest eigenvalue of the matrix $\mathbf{A}_n\trans\mathbf{A}_n$.

\item $\|\A\|$: the spectral norm of the matrix $\A$. 


\item $\I_r$  : $\I_r \in \mathbb{R}^{r\times r}$, the identity matrix; the subscript is omitted at times.



\item $\iota_{\Omega}(\mathbf{\x})$ : Indicator function of a set $\Omega$ with $\iota_{\Omega}(\mathbf{\x})=0$ if $\mathbf{\x} \in\Omega$ and otherwise $+\infty$.

\item $\vec(\mathbf{V})$ : Vector formed by stacking the column vectors of $\mathbf{V}$ with $\vec(\mathbf{V}) \in \mathbb{R}^{d'\times r'}$.

\item $\mat(\mathbf{x})$ : Convert $\mathbf{x} \in \mathbb{R}^{(d'\cdot r') \times 1}$ into a matrix with $\mat(\vec(\mathbf{V}))=\V$ with $\mat(\mathbf{x}) \in \mathbb{R}^{d'\times r'}$.



\item  $\dist(\Omega,\Omega')$ : distance between two sets with $\dist(\Omega,\Omega') \triangleq \inf_{\mathbf{w}\in \Omega,\mathbf{w}'\in \Omega'} \|\mathbf{w}-\mathbf{w}'\|$.


\end{itemize}

\subsection{Technical Preliminaries}

We present some tools in non-smooth analysis including Fr\'{e}chet subdifferential, and limiting (Fr\'{e}chet) subdifferential \cite{Mordukhovich2006,Rockafellar2009,bertsekas2015convex}. For any extended real-valued (not necessarily convex) function $F: \Rn^n \rightarrow (-\infty,+\infty]$, its domain is defined by $$\text{dom}(F)\triangleq \{\x\in\Rn^n: -\infty< F(\x) <+\infty\}.$$ 

The Fr\'{e}chet subdifferential of $F$ at $\x\in\text{dom}(F)$, denoted as $\hat{\partial}F(\x)$, is defined as $$\hat{\partial}{F}(\x) \triangleq \{\mathbf{v}\in\Rn^n: \lim_{\mathbf{z} \rightarrow \x} \inf_{\mathbf{z}\neq \x} \frac{ {F}(\mathbf{z}) - {F}(\x) - \la \mathbf{v},\mathbf{z}-\x \ra  }{\|\mathbf{z}-\x\|}\geq 0\}.$$ The limiting subdifferential of ${F}(\x)$ at $\x\in\text{dom}({F})$ is defined as: $$\partial{F}(\x)\triangleq \{\mathbf{v}\in \Rn^n: \exists \x^k \rightarrow \x, {F}(\x^k)  \rightarrow {F}(\x), \mathbf{v}^k \in\hat{\partial}{F}(\x^k) \rightarrow \mathbf{v},\forall k\}.$$ 

These subdifferentials satisfy the following key properties:

\begin{enumerate}[label=\textbf{(\alph*)}, leftmargin=22pt, itemsep=1pt, topsep=1pt, parsep=0pt, partopsep=0pt]

\item It holds that $\hat{\partial}{F}(\x) \subseteq \partial{F}(\x)$. 

\item If $F(\cdot)$ is differentiable at $\x$, then $\hat{\partial}{F}(\x) = \partial{F}(\x) = \{\nabla F(\x)\}$ with $\nabla F(\x)$ being the gradient of $F(\cdot)$ at $\x$.
     
\item When $F(\cdot)$ is convex, both $\hat{\partial}{F}(\x)$ and $\partial{F}(\x)$ reduce to the classical subdifferential for convex functions, i.e., $\hat{\partial}{F}(\x) = \partial{F}(\x) = \{\mathbf{v}: F(\mathbf{z})-F(\x)-\la\mathbf{v},\mathbf{z}-\x \ra\geq 0,\forall \mathbf{z}\}$. 

\end{enumerate}

\subsection{Relevant Lemmas}


We present several useful lemmas, each independent of context and specific methodology.






\begin{lemma}\label{lemma:relation}
\rm{(Pythagoras Relation)} For any vectors $\mathbf{a}\in\Rn^n$, $\b\in\Rn^n$, $\mathbf{c}\in\Rn^n$, we have:
\beq
\frac{1}{2}\|\mathbf{a}-\b\|_{2}^2  - \frac{1}{2} \|\mathbf{c}-\b\|_{2}^2 &=&  \frac{1}{2}\|\mathbf{a}-\mathbf{c}\|_{2}^2 + \la\b-\mathbf{c},\mathbf{c}-\mathbf{a}\ra. \nn \\
 \frac{1}{2}\|\b\|_2^2 - \frac{1}{2}\|\mathbf{c}-\b\|_2^2 &=& \frac{1}{2}\|\mathbf{c}\|_2^2  +\la \b-\mathbf{c}, \mathbf{c}\ra. \nn
\eeq
\end{lemma}

\begin{lemma}\label{useful:lemma:1}
Assume $\sigma \in (0,2)$. Let $\b^+ = \sigma\mathbf{a} + (1-\sigma)\b$, where $\b^+ \in \Rn^{n}$, $\b \in \Rn^{n}$, and $\mathbf{a} \in \Rn^{n}$. We have:
\beq
\ts\tfrac{1}{\sigma}\|\b^+\|_2^2 \leq  \sigma_1 \|\mathbf{a}\|_2^2 + \sigma_2( \|\b\|_2^2 - \|\b^+\|_2^2 ),\nn
\eeq
\noi where $\sigma_1 \triangleq \frac{\sigma}{(1-|1-\sigma|)^2}$, and $\sigma_2\triangleq \frac{|1-\sigma|}{\sigma (1-|1-\sigma|)}$.

\begin{proof}
\textbf{Part (a).} When $\sigma=1$, we have $\sigma_1=1$, $\sigma_2=0$, and $\b^+ = \mathbf{a}$. The conclusion of this lemma clearly holds.

\noi \textbf{Part (b).} We now focus on the case when $\sigma\neq1$. Noticing $|1-\sigma| \neq 0$ and $1-|1-\sigma|\neq 0$, we rewrite $\b^+ = (1-\sigma)\b + \sigma\mathbf{a}$ into the following equivalent equality:
\beq
\ts \b^+  =  (1-|1-\sigma|)\cdot    \tfrac{\sigma \mathbf{a}}{1-|1-\sigma|}  +  |1-\sigma| \cdot    \tfrac{(1-\sigma)\b}{|1-\sigma|}   . \nn
\eeq
\noi Using the fact that the function $\|\cdot\|_2^2$ is convex and $|1-\sigma|\in(0,1)$, we derive the following results:
\beq
\ts \|\b^+\|_2^2&\leq& \ts   (1-|1-\sigma|)  \cdot \|  \tfrac{\sigma \mathbf{a}}{1-|1-\sigma|}   \|_2^2 + |1-\sigma| \cdot \|\tfrac{(1-\sigma)\b}{|1-\sigma|} \|_2^2 \nn\\
&=&  \ts \tfrac{\sigma^2}{ 1 - |1-\sigma| }   \cdot \|\mathbf{a}\|_2^2 + |1-\sigma| \cdot \|\b\|_2^2.\nn
\eeq
\noi Subtracting $(|1-\sigma| \cdot \|\b^+\|_2^2 )$ from both sides of the above inequality, we have:
\beq
\ts (1- |1-\sigma|)\|\b^+\|_2^2\leq      \tfrac{\sigma^2}{ 1 - |1-\sigma| }  \cdot \|\mathbf{a}\|_2^2 + |1-\sigma| (\|\b\|_2^2 - \|\b^+\|_2^2).\nn
\eeq
\noi Dividing both sides by $\sigma(1 - |1-\sigma|)$, we have:
\beq
\ts \tfrac{1}{\sigma}\|\b^+\|_2^2\leq \tfrac{\sigma}{ ( 1 - |1-\sigma|)^2 }  \|\mathbf{a}\|_2^2 + \tfrac{|1-\sigma|}{\sigma (1- |1-\sigma|)} (\|\b\|_2^2 - \|\b^+\|_2^2) . \nn
\eeq
\noi Using the definition of $\sigma_1$ and $\sigma_2$, we finish the proof of this lemma.

\end{proof}

\end{lemma}

\begin{lemma} \label{lemma:qtqt}

We let $t\geq 1$, and $q\in(0,1)$. We have: $\tfrac{1}{q}(t+1)^q-\tfrac{1}{q} \geq \tfrac{1}{2} t^q$.

\begin{proof}

We let $h(t) \triangleq (t+1)^q-1-\tfrac{q}{2} t^q$.

Initially, we prove that $f(q) \triangleq 2^q - \tfrac{q}{2} -1 \geq 0$ for all $q\geq0$. Given $\nabla f(q) = 2^q \log(2) - \tfrac{1}{2} \geq 2^0 \log(2) - \tfrac{1}{2} = 0.1931 >0$, the function $f(q)$ is increasing for all $q\geq 0$. Combining with the fact that $f(0)=0$, we have: $f(q)\geq 0$ for all $q\geq 0$.

We derive the following inequalities:
 \beq
 \nabla h(t)  =  q t^{q-1} \cdot \{  (\frac{t+1}{t})^{q-1} - \tfrac{q}{2}\}  \overset{\step{1}}{\geq}  q t^{p-1} \cdot \{  2^{q-1} - \tfrac{q}{2} \}  \overset{\step{2}}{\geq}   q t^{q-1} \cdot \{  \tfrac{q/2 + 1}{2}  - \tfrac{q}{2} \} \overset{\step{3}}{\geq}   0, \nn
 \eeq
 \noi where step \step{1} uses $\frac{t+1}{t}\leq 2$ and $q-1\leq 0$; step \step{2} uses $2^q \geq \tfrac{q}{2} +1$ for all $q\geq 0$; step \step{3} uses $1-q\geq0$. Therefore, $h(t)$ is an increasing function.

Finally, noticing that $h(1) = 2^q-1-\tfrac{q}{2} \geq 0$, we conclude that $h(t)\geq 0$ for all $t\geq 1$.

\end{proof}

\end{lemma}


\begin{lemma} \label{eq:two:two:before}

We let $p\in(0,1)$ and $t\geq 1$. We have: $(t+1)^p - t^p \leq   p  t^{p-1}$.

\begin{proof}

We notice that $h(t) \triangleq t^p$ is concave for all $t\geq 1$ and $p\in(0,1)$ since $\nabla h(t) = p t^{p-1}$ and $\nabla^2 h(t) = p (p-1) t^{p-2} <0$. It follows that: $\forall x,y\geq 1, h(y) - h(x) \leq \la y-x,\nabla h(x)\ra$. Letting $x=t$ and $y=t+1$, for all $t\geq 1$ and $p\in(0,1)$, we have: $(t+1)^p - t^p \leq   p  t^{p-1}$.

\end{proof}
\end{lemma}

\begin{lemma} \label{eq:two:two}

We let $p\in(0,1)$. We have: $\sum_{t=1}^{\infty} (\frac{(t+1)^p -  t^p}{t^{p}}     )^2 \leq 2$.

\begin{proof}
We have:
\beq
\ts  \sum_{t=1}^{\infty} ( \frac{(t+1)^p -  t^p}{t^{p}}     )^2   \overset{\step{1}}{\leq}  \ts \sum_{t=1}^{\infty} \frac{1}{t^{2p}} t^{2p-2}  = \sum_{t=1}^{\infty}  t^{-2}  \overset{\step{2}}{\leq} \ts  2, \nn
\eeq
\noi where step \step{1} uses Lemma \ref{eq:two:two:before} and $p\leq1$; step \step{2} uses $\sum_{t=1}^{\infty} \frac{1}{t^2} \leq \sum_{t=1}^{\infty} \frac{1}{t^2} = \frac{\pi^2}{6} <2$.
\end{proof}
\end{lemma}

\begin{lemma} \label{eq:lp:function:sum:upperbound}

We let $p\in(0,1)$. We have: $\tfrac{1}{2} T^{1-p} \leq \sum_{t=1}^T t^{-p}\leq \tfrac{T^{(1-p)}}{1-p}$.

\begin{proof}

We define $h(x) = x^{-p}$ and $g(x) = \frac{1}{1-p} x^{1-p}$. Clearly, we have: $\nabla g(x) = h(x)$.

By the integral test for convergence \footnote{\url{https://en.wikipedia.org/wiki/Integral_test_for_convergence}}, we obtain: $\int_{1}^{T+1} h(x)dx \leq \sum_{t=1}^T h(t) \leq h(1) + \int_{1}^{T} h(x)dx$.

\textbf{Part (a).} We have: $\ts \sum_{t=1}^T t^{-p} \overset{}{\geq}  \ts  \int_{1}^{T+1} x^{-p} dx  \overset{\step{1}}{=}  \ts g(T+1) - g(1) =  \tfrac{1}{1-p} (T+1)^{1-p}-\tfrac{1}{1-p}  \overset{\step{2}}{\geq}  \tfrac{1}{2} T^{1-p}$, where step \step{1} uses $\nabla g(x)= h(x) = x^{-p}$; step \step{2} uses Lemma \ref{lemma:qtqt} with $q=1-p$ and $t=T$.

\textbf{Part (b).} We have: $\ts \sum_{t=1}^T t^{-p}  \overset{}{\leq}  \ts h(1) + \int_{1}^{T} x^{-p} dx \overset{\step{1}}{=} \ts  1 + g(T) - g(1) =  1 + \tfrac{1}{1-p} (T)^{1-p}-\tfrac{1}{1-p}\overset{}{=}  \tfrac{T^{(1-p)} - p }{1-p} < \tfrac{T^{(1-p)}}{1-p}$, where step \step{1} uses $h(1)=1$, and $\nabla g(x)= h(x) = x^{-p}$.

\end{proof}

\end{lemma}

\begin{lemma} \label{lemma:et1:Psi}

Let $\sigma\in(0,2)$, and $e^{t+1} - |1-\sigma| e^{t} \leq \sigma  p^t$ for all $t\geq 1$. We have: $e^{t} \leq e^1 + \sigma_3 \max_{i=1}^{t-1}  p^i$, where $\sigma_3 = \frac{\sigma}{1 - |1-\sigma|} \in [1,\infty) $.

\begin{proof}

Given $\sigma\in(0,2)$, we define $\sigma_{\star}\triangleq |1-\sigma| \in [0,1)$.

We derive the following results:
\beq
 t= 1,~~~  e^{2}  & \leq& \sigma_{\star} e^{1} +\sigma  p^1 \nn\\
 t= 2,~~~  e^{3} & \leq& \sigma_{\star} e^{2} + \sigma p^2 \leq  \sigma_{\star}^2 e^{1} + \sigma_{\star} \sigma  p^1 + \sigma  p^2 \nn\\
 t= 3,~~~  e^{4} & \leq& \sigma_{\star} e^{3} + \sigma p^3 \leq  \sigma_{\star}^3 e^{1} + \sigma_{\star}^2 \sigma  p^1 + \sigma_{\star} \sigma  p^2 + \sigma  p^3 \nn\\
 &\ldots& \nn\\
 t= T,  ~~~ \ts e^{T+1} & \leq& \ts \sigma_{\star} e^{T} + \sigma  p^T \leq  \sigma_{\star}^{T} e^{1} + \sigma \sum_{i=1}^T \sigma_{\star}^{T-i}  p^i .\nn
\eeq

Therefore, we have:
\beq
\ts e^{T+1} &\leq& \ts \sigma_{\star}^{T} e^1 + \sigma \sum_{i=1}^T \sigma_{\star}^{T-i}  p^i \nn\\
&\overset{\step{1}}{\leq}& \ts e^1  + \sigma \{\max_{i=1}^{T}  p^i\} \{\sum_{i=1}^T \sigma_{\star}^{T-i}\} \nn\\
&\overset{\step{2}}{\leq}& \ts e^1 + \sigma  \{\max_{i=1}^{T}  p^i\} \frac{ 1}{1 - \sigma_{\star} },\nn
\eeq
\noi where step \step{1} uses $\sigma_{\star}^{T}\leq 1$; step \step{2} uses the fact that:
\beq
\ts \sum_{i=1}^T \sigma_{\star}^{T-i} =  \sigma_{\star}^{T-1} + \ldots + \sigma_{\star}^{1}  +   \sigma_{\star}^{0} = \frac{ 1- \sigma_{\star}^{T}}{1 - \sigma_{\star} } \leq \frac{ 1}{1 - \sigma_{\star} }. \nn
\eeq
\end{proof}
\end{lemma}

\begin{lemma} \label{lemma:ftheta}
Assume $\kappa \in [1,2)$, $\delta \in (0,\tfrac{1}{3} (\tfrac{2}{\kappa}-1))$. For any $\varrho>0$, we define
\beq
f(\theta) = \theta - \tfrac{1}{2} - \varrho (\delta + \theta + \theta\delta - 1/\kappa)^2.\nn
\eeq
\noi We have $f(\bar{\theta})\geq \tfrac{1}{ 8 \varrho}$, where $\bar{\theta} =\tfrac{1}{2 \varrho (1+\delta)^2 } + \tfrac{1/\kappa - \delta}{\delta+1}$.

\begin{proof}

Initially, given $\kappa \in [1,2)$, $\delta \in (0,\tfrac{1}{3} (\tfrac{2}{\kappa}-1))$, we have:
\beq \label{eq:first:fff}
\tfrac{1}{2} < \tfrac{1/\kappa - \delta}{1+\delta} < 1.
\eeq

Setting the gradient of $f(\theta)$ \textit{w.r.t.} $\theta$ yields: $1 - 2 \varrho (\delta + \theta + \delta \theta - 1/\kappa) (1+\delta)=0$. It follows that the solution $\bar{\theta} = \tfrac{1}{2 \varrho (1+\delta)^2 } + \tfrac{1/\kappa - \delta}{\delta+1}$ is the maximizer of the concave function $f(\theta)$. We have:
\beq
f(\bar{\theta}) & \overset{\step{1}}{=}& \bar{\theta} - \tfrac{1}{2} - \varrho (   \delta + \theta_2 + \delta \theta_2 - 1/\kappa      )^2\nn\\
& \overset{}{=}&   \frac{1}{4 (1+\delta)^2 \varrho}  + \frac{1/\kappa - \delta}{\delta+1} - \tfrac{1}{2}        \nn\\
& \overset{\step{2}}{\geq}& \frac{1}{4 (1+\delta)^2 \varrho}  + 0 \nn\\
& \overset{\step{3}}{\geq}& \frac{1}{4 (1+1/3)^2 \varrho}   \nn\\
& \overset{\step{4}}{\geq}& \frac{1}{ 8 \varrho} , \nn
\eeq
\noi where step \step{1} uses the definitions of $f(\theta)$ and $\bar{\theta}$; step \step{2} uses the first Inequality in (\ref{eq:first:fff}); step \step{3} uses the fact that $\delta\leq \frac{1}{3}$; step \step{4} uses $4 \times (1+1/3)^2< 8$.

\end{proof}

\end{lemma}

\section{Additional Motivating Applications}
\label{app:sect:app}

$\blacktriangleright$ \textbf{Robust Sparse Regression}. Robust sparse regression \cite{liu2019one} utilizes the $\ell_1$-norm of the residuals to ensure robustness against outliers while enforcing sparsity via $\ell_0$-norm constraints to identify key variables. The problem is formulated as: $\min_{\v} \|\G\v - \z\|_1, \, \st \, \v \in\Omega \triangleq\{\v\,|\,\|\v\|_0 \leq \dot{s}\}$, where $\dot{s}\geq 0$ is an integer, $\G\in\Rn^{\dot{m}\times \dot{d}}$, and $\z\in\Rn^{\dot{m}}$. By introducing a new variable $\mathbf{y}$, this problem can be formulated as: $\min_{\v,\mathbf{y}}  \iota_{\Omega}(\v)+\|\y\|_1,\,s.t.\, -\G\v+\y=-\z$. It corresponds to Problem (\ref{eq:main}) with $\x_1 = \mathbf{v}$, $\x_2=\y$, $f_1(\x_1)=f_2(\x_2)=0$, $h_1(\x_1) =\iota_{\Omega}(\mathbf{v})$, $h_2(\x_2) = \|\y\|_1$, and $\A_1=-\G$, $\A_2=\I$, $\b=-\z$, and Condition $\BIBI$.

$\blacktriangleright$ \textbf{Dual Principal Component Pursuit}. Dual principal component pursuit \cite{TsakirisV18} is used primarily in subspace clustering and outlier detection, aiming to robustly represent data structures across different subspaces in the presence of noise and outliers. The problem is formulated as: $\min_{\V} \|\G\V\|_{2,1},\,\st\,\V \in \Omega \triangleq \{\V\,|\,\V\trans\V=\I\}$, where $\G\in\Rn^{\dot{m}\times \dot{d}}$, and $\|\Y\|_{2,1}\triangleq \sum_{i}\|\Y(i,:)\|$. By introducing a new variable $\Y$, this problem can be formulated as: $\min_{\V,\Y}  \iota_{\Omega}(\V)+\|\Y\|_{2,1},\,\st\,-\G\V+\Y=\zero$. It corresponds to Problem (\ref{eq:main}) with $\x_1 = \vec(\mathbf{V})$, $\x_2=\vec(\Y)$, $f_1(\x_1)=f_2(\x_1)=0$, $h_1(\x_1) =\iota_{\Omega}(\mathbf{V})$, $h_2(\x_2) = \|\Y\|_{2,1}$, and $\A_1=-\G$, $\A_2=\I$, $\b=\zero$, and Condition $\BIBI$.

$\blacktriangleright$ \textbf{Robust Low-Rank Approximation}. Robust low-rank approximation \cite{Candes2011JACM} uses the $\ell_1$-norm of the residuals to ensure robustness against outliers while imposing a low-rank constraint on the solution matrix The problem is formulated as: $\min_{\V} \|\G(\V) - \z\|_1,\,\st\,\V\triangleq\{\V\,|\,\text{rank}(\V) \leq \dot{s}\}$, where $\dot{s}\geq 0$ is an integer, $\G(\cdot): \Rn^{\dot{d}\times \dot{r}}\mapsto \Rn^{\dot{m}}$, and $\z\in\Rn^{\dot{m}}$. By introducing a new variable $\mathbf{y}$, this problem can be formulated as: $\min_{\V,\mathbf{y}} \iota_{\Omega}(\V)+\|\y\|_1 ,\,s.t.\, -\G(\V)+\y=-\z$. It corresponds to Problem (\ref{eq:main}) with $\x_1 = \vec(\mathbf{V})$, $\x_2=\y$, $f_1(\x_1)=f_2(\x_1)=0$, $h_1(\x_1) =\iota_{\Omega}(\mathbf{V})$, $h_2(\x_2) = \|\y\|_1$, $\A_1\x_1 = -\G(\mathbf{V})$, $\A_2=\I$, $\b=-\z$, and Condition $\BIBI$.

\section{Proofs for Section \ref{sect:proposed}} \label{app:sect:proposed}

\subsection{Proof of Lemma \ref{lemma:updating:Sublinear}}\label{app:lemma:updating:Sublinear}

\begin{proof}

Consider the update rule $\beta^{t} =  \beta^0 +  \beta^0 \xi t^p$, where $p\in(0,1)$.

\textbf{Part (a).}We have:
\beq
\ts \beta^{t+1} -  \beta^t   - \xi\beta^t   \overset{\step{1}}{=}   \ts  \beta^0 \xi  ( (t+1)^p - t^p ) - \xi \beta^0 \overset{\step{2}}{\leq}    \beta^0 \xi - \beta^0 \xi =   0, \nn
\eeq
\noi where step \step{1} uses the update rule $\beta^{t} =  \beta^0 +  \beta^0 \xi t^p$; step \step{2} uses the fact that the function $h(t)\triangleq (t+1)^p - t^p$ is monotonically decreasing \textit{w.r.t.} $t$ that: $h(t) \leq h(0) = 1$.

\textbf{Part (b).}We derive: $L_n \leq \beta^0 \delta \lambdaUp  \overset{\step{1}}{\leq} \beta^t \delta \lambdaUp$, where step \step{1} uses $\beta^t\geq\beta^0$.


\end{proof}

\subsection{Proof of Lemma \ref{lemma:lip:mu0}}\label{app:lemma:lip:mu0}

\begin{proof}

We let $\ub \in \Rn^d$ be a fixed constant vector. We assume $0<\mu_2<\mu_1$.

We define: $h(\ub;\mu_1) \triangleq \min_{\v}  h(\v) + \tfrac{1}{2\mu_1} \| \v - \ub\|_2^2$, and $h(\ub;\mu_2) \triangleq \min_{\v}  h(\v) + \tfrac{1}{2 \mu_2} \| \v - \ub\|_2^2$.
 
We define $\p_1 \triangleq\Prox_{h}(\ub;\mu_1) \triangleq \arg \min_{\v} h(\v) + \frac{1}{2 \mu_1} \| \v -\ub\|_2^2$.

We define $\p_2 \triangleq\Prox_{h}(\ub;\mu_2) \triangleq \arg \min_{\v} h(\v) + \frac{1}{2 \mu_2} \| \v -\ub\|_2^2$.

We let $\g_1 \in \partial h(\Prox_{h}(\ub;\mu_1))$, and $\g_2\in \partial h(\Prox_{h}(\ub;\mu_2))$.

\noi Initially, by the optimality of $\p_1\triangleq\Prox_{h}(\ub;\mu_1)$ and $\p_2\triangleq\Prox_{h}(\ub;\mu_2)$, we obtain:
\beq
&& \ub - \p_1  \in \mu_1 \partial h(\Prox_{h}(\ub;\mu_1)) = \mu_1 \g_1, \label{eq:optimality:mu1}\\
&&\ub - \p_2 \in \mu_2  \partial h(\Prox_{h}(\ub;\mu_2))= \mu_2 \g_2 .\label{eq:optimality:mu2}
\eeq

\noi \textbf{Part (a).} We now prove that $0\leq \frac{ h(\ub;\mu_2) - h(\ub;\mu_1)}{\mu_1 - \mu_2}$. We have:
\beq
h(\ub;\mu_1)  - h(\ub;\mu_2)&\overset{\step{1}}{=}&  \tfrac{1}{2\mu_1} \|\ub - \p_1\|_2^2 - \tfrac{1}{2\mu_2} \|\ub - \p_2\|_2^2 + h(\p_1)  - h(\p_2)  \nn\\
&\overset{\step{2}}{\leq}& \tfrac{1}{2\mu_1}\| \ub - \p_1 \|_2^2 - \tfrac{1}{2\mu_2}\| \ub - \p_2\|_2^2 + \la \p_1-\p_2,\g_1\ra  \nn\\
&\overset{\step{3}}{=}& \tfrac{\mu_1}{2}\|  \g_1\|_2^2 - \tfrac{\mu_2}{2}\| \g_2\|_2^2+ \la \mu_2 \g_2 -\mu_1 \g_1,\g_1\ra  \nn\\
&=&   - \tfrac{\mu_1}{2}\|  \g_1\|_2^2 - \tfrac{\mu_2}{2}\| \g_2\|_2^2  + \mu_2 \la  \g_2,\g_1\ra \nn\\
&\overset{\step{4}}{\leq}&  - \tfrac{\mu_2}{2}\|  \g_1\|_2^2 - \tfrac{\mu_2}{2}\| \g_2\|_2^2 + \mu_2 \la \g_2 ,\g_1\ra  \nn\\
&\overset{}{=}&   -\tfrac{\mu_2}{2} \|\g_2 - \g_1\|_2^2 \leq 0, \nn
\eeq
\noi where step \step{1} uses the definition of $h(\ub;\mu)$; step \step{2} uses the convexity of $h(\cdot)$; step \step{3} uses the optimality of $\p_1\triangleq\Prox_{h}(\ub;\mu_1)$ and $\p_2\triangleq\Prox_{h}(\ub;\mu_2)$ as in (\ref{eq:optimality:mu1}) and (\ref{eq:optimality:mu2}); step \step{4} uses $\mu_2 < \mu_1$.

\noi \textbf{Part (b).} We now prove that $\frac{ h(\ub;\mu_2) - h(\ub;\mu_1)}{\mu_1 - \mu_2}\leq \tfrac{1}{2}C_h^2$. We have:
\beq
h(\ub;\mu_2) - h(\ub;\mu_1)&\overset{\step{1}}{=}& \tfrac{1}{2\mu_2} \|\ub - \p_2\|_2^2 - \tfrac{1}{2\mu_1} \|\ub - \p_1\|_2^2 + h(\p_2)   -h(\p_1)  \nn\\
&\overset{\step{2}}{\leq}&  \tfrac{1}{2\mu_2} \|\ub - \p_2\|_2^2 - \tfrac{1}{2\mu_1} \|\ub - \p_1 \|_2^2+ \la \p_2 -\p_1 ,  \g_2\ra    \nn\\
&\overset{\step{3}}{=}&  \tfrac{\mu_2}{2} \| \g_2 \|_2^2 - \tfrac{\mu_1}{2} \|\g_1\|_2^2+ \la \mu_1 \g_1- \mu_2 \g_2,  \g_2\ra    \nn\\
&=&  - \tfrac{\mu_2}{2} \|\g_2\|_2^2 - \tfrac{\mu_1}{2} \| \g_1\|_2^2  + \mu_1\la \g_2, \g_1\ra  \nn\\
&\overset{\step{4}}{\leq}&  \tfrac{\mu_1}{2}\|\g_2\|_2^2  - \tfrac{\mu_2}{2} \|\g_2\|_2^2   \nn\\
&\overset{\step{5}}{\leq}&   \tfrac{ \mu_1 - \mu_2   }{2} \cdot C_h^2,\nn
\eeq
\noi where step \step{1} uses the definition of $h(\ub;\mu)$; step \step{2} uses the convexity of $h(\cdot)$; step \step{3} uses the optimality of $\p_1\triangleq\Prox_{h}(\ub;\mu_1)$ and $\p_2\triangleq\Prox_{h}(\ub;\mu_2)$ as in (\ref{eq:optimality:mu1}) and (\ref{eq:optimality:mu2}); step \step{4} uses the inequality that: $-\tfrac{1}{2} \| \g_1\|_2^2  +   \la \g_1, \g_2\ra\leq \tfrac{1}{2} \|\g_2\|_2^2$ for all $\g_1,\g_2\in \mathbb{R}^{d\times 1}$; step \step{5} uses $\|\g_2\|\leq C_h$.
 \end{proof}

\subsection{Proof of Lemma \ref{lemma:lip:mu}}\label{app:lemma:lip:mu}

\begin{proof}

We let $\ub$ be a fixed constant vector. We assume $0<\mu_2<\mu_1$.

We define: $h(\ub;\mu) \triangleq \min_{\v \in \Rn^{d \times 1}}  h(\v) + \frac{1}{2 \mu} \| \v - \ub\|_2^2$.

We define: $\Prox_h(\ub;\mu) \triangleq \arg \min_{\v \in \Rn^{d \times 1}}  h(\v) + \frac{1}{2 \mu} \| \v - \ub\|_2^2$.

Using Claim (\bfit{b}) of Lemma \ref{lemma:lip:mu:2}, we establish that $h(\ub;\mu)$ is smooth \textit{w.r.t.} $\ub$, and its gradient can be computed as:
\beq
\nabla h(\ub;\mu) = \mu^{-1}(\ub-\Prox_{h}(\ub;\mu)).\nn
\eeq

We examine the following mapping $\mathcal{H}(\upsilon)\triangleq \upsilon (\ub-\Prox_{h}(\ub;\tfrac{1}{\upsilon}))$ with $\mathcal{H}(\upsilon):\Rn \mapsto \Rn^n$. We derive:
     \beq
     \ts     \lim_{\delta \rightarrow 0} \frac{ \mathcal{H}(\upsilon+\delta) - \mathcal{H}(\upsilon) }{\delta} &=& \ts \lim_{\delta \rightarrow 0} \frac{   (\upsilon +\delta)  (\ub-\Prox_{h}(\ub; \tfrac{1}{\upsilon +\delta}))   -  \upsilon (\ub-\Prox_{h}(\ub;\tfrac{1}{\upsilon}))  }{\delta} \nn\\
     &=&\ts\lim_{\delta \rightarrow 0} \frac{\delta\ub - (\upsilon+\delta) \Prox_{h}(\ub;\tfrac{1}{\upsilon}) + \upsilon \Prox_{h}(\ub;\tfrac{1}{\upsilon}) }{\delta} = \ub-\Prox_{h}(\ub;\tfrac{1}{\upsilon}).\nn
     \eeq
     \noi
     Therefore, the first-order derivative of the mapping $\mathcal{H}(\upsilon)$ \textit{w.r.t.} $\upsilon$ always exists and can be computed as $\nabla_{\upsilon} \mathcal{H}(\upsilon) =\ub-\Prox_{h}(\ub;\tfrac{1}{\upsilon})$, leading to:
     \beq
     \forall \upsilon,\upsilon'>0,\,\frac{\|\mathcal{H}(\upsilon) - \mathcal{H}(\upsilon')\|}{|\upsilon-\upsilon'|} \leq \|\ub-\Prox_{h}(\ub;\tfrac{1}{\upsilon})\|.\nn
     \eeq
     Letting $\upsilon=1/\mu_1$ and $\upsilon'=1/\mu_2$, we derive:
     \beq
      \frac{\|\nabla h(\ub;\mu_1) - \nabla h(\ub;\mu_2) \|}{|1/\mu_1-1/\mu_2|} \overset{}{\leq} \|\ub-\Prox_{h}(\ub;\mu_1)\| \overset{\step{1}}{=}  \mu_1 \| \partial h(\Prox_{h}(\ub;\mu_1)) \| \overset{\step{2}}{\leq}   \mu_1 C_h,\nn
     \eeq
     \noi where step \step{1} uses the optimality of $\Prox_{h}(\ub;\mu)$ that $\zero \in \partial h( \Prox_{h}(\ub;\mu) ) + \frac{1}{\mu} ( \Prox_{h}(\ub;\mu) - \ub)$ for all $\mu$; step \step{2} uses the Lipschitz continuity of $h(\cdot)$. We further obtain:
     \beq
     \|\nabla\,h(\ub;\mu_1) -\nabla\,h(\ub;\mu_2)\|\leq |1/\mu_1-1/\mu_2| \cdot \mu_1 C_h = (\mu_1/\mu_2-1) \cdot C_h.\nn
     \eeq

\end{proof}

\subsection{Proof of Lemma \ref{lemma:smoothing:problem}}\label{app:lemma:smoothing:problem}

\begin{proof}

The proof of this lemma is similar to that of Lemma 1 in \cite{li2022riemannian}. For completeness, we include the proof here.

We consider the following strongly convex problems:
\beq
\bar{\x}_n &=&  \arg \min_{\x_n} h_n(\x_n;\mu) + \tfrac{\rho}{2} \|\x_n - \cb\|_2^2 \nn \\
\Leftrightarrow~(\bar{\x}_n,\breve{\x}_n) &=&  \arg \min_{\x_n,\breve{\x}_n} h_n(\breve{\x}_n) + \frac{1}{2\mu}\| \x_n-\breve{\x}_n\|_2^2 + \tfrac{\rho}{2} \|\x_n - \cb\|_2^2.\nn
\eeq
We have the following first-order optimality conditions:
\beq
\zero &=& \tfrac{1}{\mu} (\bar{\x}_n - \breve{\x}_n) + \rho (\bar{\x}_n - \cb) \label{eq:gUY:0}\\
\zero &\in& \partial h_n(\breve{\x}_n) + \tfrac{1}{\mu} ( \breve{\x}_n - \bar{\x}_n). \label{eq:gUY}
\eeq
\noi \textbf{Part (a).} Using (\ref{eq:gUY:0}), we obtain: $\bar{\x}_n = \frac{1}{1/\mu+\rho} ( \frac{1}{\mu} \breve{\x}_n  + \rho \cb )$. Plugging this equation into (\ref{eq:gUY}) yields:
\beq
\zero &\in&   \partial h_n(\breve{\x}_n) + \tfrac{1}{\mu} ( \breve{\x}_n -  \frac{1}{ 1/\mu+\rho} ( \frac{1}{\mu} \breve{\x}_n  + \rho \cb ) ) \nn\\
& = &  \partial h_n(\breve{\x}_n) +     \frac{\rho}{1+\mu\rho} (\breve{\x}_n - \cb). \nn
\eeq
\noi The inclusion above implies that:
\beq
\breve{\x}_n = \arg \min_{\breve{\x}_n} h_n(\breve{\x}_n) + \frac{1}{2 } \cdot \frac{\rho}{1+\mu\rho} \|\breve{\x}_n-\cb\|_{2}^2.\nn
\eeq

\noi \textbf{Part (b).} We derive:
\beq
- \rho (\bar{\x}_n - \cb) \overset{\step{1}}{=}  \tfrac{1}{\mu} (\bar{\x}_n - \breve{\x}_n)   \overset{\step{2}}{\in}   \partial h_n(\breve{\x}_n),\nn
\eeq
\noi where step \step{1} uses (\ref{eq:gUY:0}); step \step{2} uses (\ref{eq:gUY}).

\noi \textbf{Part (c).} Using (\ref{eq:gUY}), we have: $\breve{\x}_n - \bar{\x}_n = -\mu \partial h_n(\breve{\x}_n)$. This leads to $\|\breve{\x}_n - \bar{\x}_n\| \leq \mu C_h$.

\end{proof}

\section{Proofs for Section \ref{sect:global:convergence}} \label{app:sect:global:convergence}

\subsection{Proof of Lemma \ref{lemma:suf:dec}} \label{app:lemma:suf:dec}

\begin{proof}

\textbf{Part (a).} We now focus on sufficient decrease for variables $\{\x_1,\x_2,\ldots,\x_{n-1}\}$. We define $\Phi_i^t = G(\x^{t+1}_{[1,i-1]},\x^{t+1}_{i},\x^t_{[i+1,n]},\mathbf{z}^t;\beta^t) - G(\x^{t+1}_{[1,i-1]},\x^t_{i},\x^t_{[i+1,n]},\mathbf{z}^t;\beta^t) + h_i(\x_i^{t+1}) - h_i(\x_i^{t})$, where $i\in[n-1]$.

Noticing the function $G(\x^{t+1}_{[1,i-1]},\x_{i},\x^t_{[i+1,n]},\mathbf{z}^t;\beta^t)$ is $\LL_i^t$-smooth \textit{w.r.t.} $\x_i$ for the $t$-th iteration, we have:
\beq
&& G(\x^{t+1}_{[1,i-1]},\x^{t+1}_{i},\x^t_{[i+1,n]},\mathbf{z}^t;\beta^t) - G(\x^{t+1}_{[1,i-1]},\x^t_{i},\x^t_{[i+1,n]},\mathbf{z}^t;\beta^t)  \nn\\
&\leq&  \la \x_i^{t+1}-\x_i^t,  \nabla_{\x_i}G(\x^{t+1}_{[1,n-1]},\x_{i}^t,\x_{[i+1,n]}^t,\mathbf{z}^t;\beta^t)\ra + \tfrac{\LL_i^t}{2}\|\x_i^{t+1} - \x_i^t\|_2^2. \label{eq:smooth:S:x}
\eeq
\noi Given $\x_i^{t+1}$ is the minimizer of the following optimization problem:
\beq
\ts \x_i^{t+1}  \in  \arg \min_{\x_i}\, h_i(\x_i) + \la\x_i - \x_i^{t}, \nabla_{\x_i}G(\x^{t+1}_{[1,n-1]},\x_{i}^t,\x_{[i+1,n]}^t,\mathbf{z}^t;\beta^t)\ra + \tfrac{\theta_1 \LL_i^t}{2}\|\x_i-\x_i^t\|_2^2.\nn
\eeq
\noi The optimality of $\x_i^{t+1}$ leads to:
\beq \label{eq:optimality:x}
 \ts h_i(\x_i^{t+1}) - h_i(\x_i^{t}) + \la\x_i^{t+1} - \x_i^{t}, \nabla_{\x_i}G(\x^{t+1}_{[1,n-1]},\x_{i}^t,\x_{[i+1,n]}^t,\mathbf{z}^t;\beta^t) \ra \overset{\text{}}{\leq}  \ts  - \tfrac{\theta_1 \LL_i^t}{2}\|\x_i^{t+1}-\x_i^t\|_2^2.
\eeq
\noi Combining equations (\ref{eq:smooth:S:x}) and (\ref{eq:optimality:x}), we derive the following expressions:
\beq
\Phi_i^t  \leq  ( \tfrac{1}{2} -\tfrac{\theta_1}{2} ) \cdot \LL_i^t\|\x_i^{t+1} - \x_i^t\|_2^2. \nn
\eeq
\noi Telescoping the above inequality over $i$ from $1$ to $(n-1)$ leads to:
\beq
\ts \sum_{i=1}^{n-1} \Phi_i^t \leq  \sum_{i=1}^{n-1} \{ ( \tfrac{1}{2} -\tfrac{\theta_1}{2} ) \cdot \LL_i^t\|\x_i^{t+1} - \x_i^t\|_2^2\}.\nn
\eeq
Therefore, we obtain:
\beq
\ts \mathcal{L}(\x_{[1,n-1]}^{t+1},\x_n^{t},\mathbf{z}^t;\beta^t,\mu^t)-\mathcal{L}(\x^{t},\mathbf{z}^t;\beta^t,\mu^t) \leq \sum_{i=1}^{n-1} \{ ( \tfrac{1}{2} -\tfrac{\theta_1}{2} ) \cdot \LL_i^t\|\x_i^{t+1} - \x_i^t\|_2^2\}. \label{eq:dec:x:1}
\eeq

\textbf{Part (b).} We now focus on sufficient decrease for variable $\{\x_n\}$. Noticing the function $G(\x^{t+1}_{[1,n-1]},\x_{n},\mathbf{z}^t;\beta^t)$ is $\LL_n^t$-smooth \textit{w.r.t.} $\x_n$ for the $t$-th iteration, we have:
\beq
&& G(\x^{t+1}_{[1,n-1]},\x^{t+1}_{n},\mathbf{z}^t;\beta^t) - G(\x^{t+1}_{[1,n-1]},\x^t_{n},\mathbf{z}^t;\beta^t) \nn\\
&\leq&  \la \x_n^{t+1}-\x_n^t,  \nabla_{\x_n}G(\x^{t+1}_{[1,n-1]},\x_{n}^t,\mathbf{z}^t;\beta^t)\ra + \tfrac{\LL_n^t}{2}\|\x_n^{t+1} - \x_n^t\|_2^2. \label{eq:smooth:S:x:22}
\eeq

Since $h_{n}(\x_n;\mu^t)$ is convex, we have:
\beq \label{eq:optimality:xn:addd:2}
&& h_{n}( \x_n^{t+1};\mu^t) - h_{n}( \x_n^{t};\mu^t ) \nn\\
& \leq& \la  \x_n^{t+1}-  \x_n^{t+1} ,\nabla h_{n}( \x_i^{t+1};\mu^t ) \ra \nn\\
&\overset{\step{1}}{=}&  \la  \x_n^{t+1}-  \x_n^{t+1} , -\nabla_{\x_n}G(\x^{t+1}_{[1,n-1]},\x_{n}^t,\mathbf{z}^t;\beta^t)\ra - \theta_2 \LL_n^t (\x_n^{t+1}-\x_n^t)\ra ,
\eeq
\noi where step \step{1} uses the the first-order optimality condition of $\x_n^{t+1}$ that:
\beq
\ts  \zero = \nabla h_{n}( \x_n^{t+1};\mu^t ) + \nabla_{\x_n}G(\x^{t+1}_{[1,n-1]},\x_{n}^t,\mathbf{z}^t;\beta^t)\ra + \theta_2 \LL_n^t (\x_n^{t+1}-\x_n^t).\nn
\eeq
\noi Adding Inequalities (\ref{eq:smooth:S:x:22}) and (\ref{eq:optimality:xn:addd:2}) together, we have:
\beq \label{eq:adding:dec:x:i:2}
&& h_{n}(\x_n^{t+1};\mu^t)   - h_{n}(\x_n^{t};\mu^t)  + G(\x^{t+1}_{[1,n-1]},\x^{t+1}_{n},\mathbf{z}^t;\beta^t)  - G(\x^{t+1}_{[1,n-1]},\x^t_{n},\mathbf{z}^t;\beta^t)  \nn\\
&\leq &  \tfrac{\LL_n^t}{2}\|\x_n^{t+1} - \x_n^t\|_2^2 - \theta_2 \LL_n^t \|\x_n^{t+1}-\x_n^t\|_2^2\nn\\
&=&  (\tfrac{1}{2} - \theta_2 ) \cdot \LL_n^t \|\x_n^{t+1}-\x_n^t\|_2^2.\nn
\eeq
\noi This results in the following inequality:
\beq
\mathcal{L}(\x^{t+1},\mathbf{z}^t;\beta^t,\mu^t)-\mathcal{L}(\x_{[1,n-1]}^{t+1},\x_n^{t},\mathbf{z}^t;\beta^t,\mu^t) \leq (\tfrac{1}{2} - \theta_2 ) \cdot \LL_n^t \|\x_n^{t+1}-\x_n^t\|_2^2. \label{eq:dec:x:2}
\eeq

\textbf{Part (c).} We now focus on sufficient decrease for variable $\{\z\}$. We have:
\beq
&&\mathcal{L}(\x^{t+1},\mathbf{z}^{t+1};\beta^t,\mu^t)-\mathcal{L}(\x^{t+1},\mathbf{z}^{t};\beta^t,\mu^t)\nn\\
& = & \la \A\x^{t+1} - \b,\z^{t+1} - \z^t \ra \nn\\
&\overset{\step{1}}{=}& \la \frac{1}{\sigma \beta^t} (\z^{t+1}-\z^t),\z^{t+1} - \z^t \ra\nn\\
 &=& \frac{1}{\sigma \beta^t} \|\z^{t+1}-\z^t\|_2^2,\label{eq:dec:z}
\eeq
\noi where step \step{1} uses $\z^{t+1}= \z^t + \sigma \beta^t (\A \x^{t+1} - \b  ) $ with $\A \x^{t+1}\triangleq \sum_{j=1}^n \A_j \x_j^{t+1}$.

\textbf{Part (d).} We now focus on sufficient decrease for variable $\{\beta\}$. We have:
\beq
&& \ts \mathcal{L}(\x^{t+1},\mathbf{z}^{t+1};\beta^{t+1},\mu^t)-\mathcal{L}(\x^{t+1},\mathbf{z}^{t+1};\beta^t,\mu^t) \nn\\
&=&\ts (\tfrac{\beta^{t+1}}{2} - \tfrac{\beta^{t}}{2}) \| \A\x^{t+1} - \b\|_2^2 \nn\\
&\overset{\step{1}}{=}& (\tfrac{\beta^{t+1}}{2} - \tfrac{\beta^{t}}{2}) \| \frac{1}{\sigma \beta^t} (\z^{t+1}-\z^t)\|_2^2 \nn\\
&\overset{\step{2}}{\leq}& \ts (\tfrac{ (1+\xi)\beta^{t}}{2} - \tfrac{\beta^{t}}{2}) \| \frac{1}{\sigma \beta^t} (\z^{t+1}-\z^t)\|_2^2 \nn\\
&=&\ts \tfrac{ \xi}{2 \sigma}  \cdot \tfrac{1}{\sigma \beta^t} \|\z^{t+1}-\z^t\|_2^2,\label{eq:dec:beta}
\eeq
\noi where step \step{1} uses $\z^{t+1}= \z^t + \sigma \beta^t (\A \x^{t+1} - \b  ) $; step \step{2} uses Lemma \ref{lemma:updating:Sublinear} that $\beta^{t+1} \leq \beta^t (1+\xi)$.

\textbf{Part (e).} We now focus on sufficient decrease for variable $\{\mu\}$. We have:
\beq
&&\mathcal{L}(\x^{t+1},\mathbf{z}^{t+1};\beta^{t+1},\mu^{t+1})-\mathcal{L}(\x^{t+1},\mathbf{z}^{t+1};\beta^{t+1},\mu^t) \nn\\
& = & h_{n}(\x^{t+1}_n;\mu^{t+1}) -  h_{n}(\x^{t+1}_n;\mu^{t}) \nn\\
&\overset{\step{1}}{\leq}& \tfrac{1}{2} C_h (\mu^{t} - \mu^{t+1}) ,  \label{eq:dec:mu}
\eeq
\noi where step \step{1} uses Lemma \ref{lemma:lip:mu0}.

Combining Inequalities (\ref{eq:dec:x:1}), (\ref{eq:dec:x:2}), (\ref{eq:dec:z}), (\ref{eq:dec:beta}), and (\ref{eq:dec:mu}), we have:
\beq
&&\mathcal{L}(\x^{t+1},\mathbf{z}^{t+1};\beta^{t+1},\mu^{t+1}) - \mathcal{L}(\x^{t},\mathbf{z}^{t};\beta^{t},\mu^{t}) \nn\\
&\leq& \ts  [\sum_{i=1}^{n-1} \{ ( \tfrac{1}{2} -\tfrac{\theta_1}{2} ) \cdot \LL_i^t\|\x_i^{t+1} - \x_i^t\|_2^2\}] + (\tfrac{1}{2} - \theta_2 ) \cdot \LL_n^t \|\x_n^{t+1}-\x_n^t\|_2^2 \nn\\
&&   + (1 + \frac{ \xi}{2 \sigma}  ) \cdot \frac{1}{\sigma \beta^t} \|\z^{t+1}-\z^t\|_2^2+ \tfrac{1}{2} C_h (\mu^{t} - \mu^{t+1}).
\eeq

We define $\Theta^{t}_{\oo}  \triangleq  \mathcal{L}(\x^{t},\mathbf{z}^{t};\beta^{t},\mu^{t}) + \tfrac{1}{2} C_h \mu^{t}$, $\varepsilon_3 \triangleq \xi$, $\varepsilon_1\triangleq \frac{1}{2}\theta_1 - \frac{1}{2}$, and $\mathcal{E}^{t+1} \triangleq  \ts  \frac{\varepsilon_3}{\beta^t}\|\z^{t+1}-\z^t\|_2^2 + \varepsilon_2 \LL_n^t \| \x_n^{t+1} - \x_n^t \|_2^2 + \varepsilon_1\sum_{i=1}^{n-1} \LL_i^t \| \x_i^{t+1} - \x_i^t \|_2^2$. We have:
\beq
&& \mathcal{E}^{t+1} + \Theta^{t+1}_{\oo} - \Theta^{t}_{\oo} \nn\\
& \leq & \ts  (\tfrac{1}{2} - \theta_2 + \varepsilon_2 ) \cdot \LL_n^t \|\x_n^{t+1}-\x_n^t\|_2^2 +   (1 + \frac{ \xi}{2 \sigma} + \sigma \xi  ) \cdot \frac{1}{\sigma \beta^t} \|\z^{t+1}-\z^t\|_2^2. \nn
\eeq
\end{proof}

\subsection{Proof of Lemma \ref{lemma:first:order}} \label{app:lemma:first:order}

\begin{proof}

For any $i\in[n]$, we define $\u_i^{t+1} \triangleq \thetas_i  \LL_i^t [\x_i^{t+1} - \x_i^{t}] - \beta^t  \A_i \trans [\sum_{j=i}^n \A_j(\x_j^{t+1}-\x_j^t)]$, and let $\mathbbm{w}_i^{t+1} \in \partial h_i (\x_i^{t+1})+\nabla f_i(\x_i^{t})$.

We notice that $\x_i^{t+1}$ is the minimizer of the following problem:
\beq
\ts  \x_i^{t+1}  \in  \arg \min_{\x_i}\, \tfrac{\theta \LL_i^t}{2}\|\x_i-\x_i^t\|_2^2 + h_i(\x_i) + \la\x_i - \x_i^{t}, \nabla_{\x_i}G(\x_{[1,i-1]}^{t+1},\x_{[i,n]}^{t},\mathbf{z}^t;\beta^t)\ra. \nn
\eeq
\noi Using the necessary first-order optimality condition of the solution $\x_i^{t+1}$, we have:
\beq \label{eq:opt:x:first}
\nabla_{\x_i}G(\x_{[1,i-1]}^{t+1},\x_{[i,n]}^{t},\mathbf{z}^t;\beta^t)   \in  -\partial h_i(\x_i^{t+1} ) - \theta  \LL_i^t (\x_i^{t+1} - \x_i^t).
\eeq

Using the definition of the function $G(\x,\mathbf{z};\beta) \triangleq   \la [\sum_{j=1}^n\A_j\x_j]-\b,\mathbf{z} \ra  + \frac{\beta}{2}\|[\sum_{j=1}^n\A_j\x_j]-\b\|_2^2 + \sum_{j=1}^n f_j(\x_j)$, we have:
\beq \label{eq:grad:xn}
\ts && \nabla_{\x_i}G(\x_{[1,i-1]}^{t+1},\x_{[i,n]}^{t},\mathbf{z}^t;\beta^t)\nn\\
& =& \ts  \nabla f_i(\x_i^t) + \A_i \trans\mathbf{z}^t +   \beta^t \A_i \trans \{  [\sum_{j=1}^{i-1} \A_j \x_j^{t+1}] + [\sum_{j=i}^{n} \A_j \x_j^{t}] - \b \} \nn\\
&=& \ts  \nabla f_i(\x_i^t) + \A_i \trans\mathbf{z}^t  + \beta^t \A_i \trans \{ \A\x^{t+1} - \b  +    [\sum_{j=i}^n \A_j(\x_j^t-\x_j^{t+1})] \} \nn\\
&\overset{\step{1}}{=}& \ts  \nabla f_i(\x_i^t) + \A_i \trans\mathbf{z}^t  +  \frac{1}{\sigma}\A_i\trans (\mathbf{z}^{t+1}-\mathbf{z}^t) +   \beta^t \A_i \trans \{\sum_{j=i}^n \A_j(\x_j^t-\x_j^{t+1})\},
\eeq
\noi where step \step{1} uses the update rule of $\mathbf{z}^{t+1}$ that $\mathbf{z}^{t+1}-\mathbf{z}^t =     \sigma \beta^t  (\sum_{i=1}^n \A_i\x_i^{t+1} - \b)$. Combining the Equalities (\ref{eq:opt:x:first}) and (\ref{eq:grad:xn}), we obtain the following result:
\beq
 \ts \zero  &\in&  \partial h_i(\x_i^{t+1} ) + \thetas_i  \LL_i^t [ \x_i^{t+1} - \x_i^{t}] +  \nabla f_i(\x_i^t)\nn\\
 && \ts + \A_i \trans\mathbf{z}^t +   \beta^t  \A_i \trans [\sum_{j=i}^n \A_j(\x_j^t-\x_j^{t+1})] +  \frac{1}{\sigma}\A_i\trans (\mathbf{z}^{t+1}-\mathbf{z}^t). \nn
\eeq
\noi Using the definition of $\mathbbm{w}_i^{t+1}$ and $\u_i^{t+1}$ for all $i\in[n]$, we have: $\zero = \mathbbm{w}_i^{t+1}  + \u_i^{t+1} + \A_i \trans\mathbf{z}^t +    \frac{1}{\sigma}\A_i\trans (\mathbf{z}^{t+1}-\mathbf{z}^t)$. Multiplying both sides by $\sigma \in (0,2)$, for all $t\geq0$, we have:
\beq
\ts  \zero = \sigma \mathbbm{w}_i^{t+1} + \sigma \A_i\trans\mathbf{z}^t +   \A_i\trans(\mathbf{z}^{t+1}-\mathbf{z}^{t})+\sigma \u_i^{t+1}  . \label{eq:combine:1}
\eeq
\noi Given that $t$ can take on any integer value, for all $t\geq1$, we derive:
\beq
\ts  \zero = \sigma \mathbbm{w}_i^{t} + \sigma \A_i\trans\mathbf{z}^{t-1} +   \A_i\trans(\mathbf{z}^{t}-\mathbf{z}^{t-1})+\sigma \u_i^{t}  .  \label{eq:combine:2}
\eeq
\noi Combining Equality (\ref{eq:combine:1}) and Equality (\ref{eq:combine:2}), for all $t\geq1$, we have:
\beq\label{eq:combine:3}
\A_i\trans(\mathbf{z}^{t+1}-\mathbf{z}^{t}) = (1-\sigma) \A_i\trans (\mathbf{z}^t -\mathbf{z}^{t-1} ) - \sigma (\mathbbm{w}_i^{t+1} - \mathbbm{w}_i^{t}) -    \sigma(\u_i^{t+1} - \u_i^{t} ).
\eeq
\noi In view of (\ref{eq:combine:3}), we let $i=n$ and arrive at the following two distinct identities:
\begin{align}
\BIBI:& \underbrace{\A_n\trans (\mathbf{z}^{t+1} - \mathbf{z}^{t})}_{\triangleq \mathbbm{a}^{t+1}} = (1-\sigma)\underbrace{(\A_n\trans (\mathbf{z}^{t} - \mathbf{z}^{t-1}))}_{\triangleq\mathbbm{a}^t} + \sigma \underbrace{ (\u_n^{t}- \u_n^{t+1} +  \mathbbm{w}_n^{t} - \mathbbm{w}_n^{t+1} ) }_{\mathbbm{c}^{t}}. \nn\\
\SUSU:& \underbrace{\A_n\trans (\mathbf{z}^{t+1} - \mathbf{z}^{t}) +  \sigma \u_n^{t+1}}_{ \triangleq \mathbbm{a}^{t+1}}  = (1-\sigma) (\underbrace{\A_n\trans (\mathbf{z}^{t} - \mathbf{z}^{t-1}) +  \sigma \u_n^{t}}_{ \triangleq \mathbbm{a}^{t}})  + \sigma( \underbrace{\sigma \u_n^{t} +   \mathbbm{w}_n^{t} - \mathbbm{w}_n^{t+1}}_{\triangleq\mathbbm{c}^t}).  \nn
\end{align}

\end{proof}

\subsection{Proof of Lemma \ref{lemma:reg:mu:L}}\label{app:lemma:reg:mu:L}

\begin{proof}



We denote $\H^t\triangleq \theta_2 \LL_n^t  \mathbf{I} -  \beta^t  \A_n\trans\A_n \in \Rn^{\mathbf{d}_i \times \mathbf{d}_i}$.

We assume $\A_n\trans \A_n $ has the singular value decomposition $\A_n\trans \A_n = \tilde{\mathbf{U}}\trans \diag(\lambdas)\tilde{\mathbf{U}}$, where $\tilde{\mathbf{U}}\in \Rn^{\mathbf{d}_i\times \mathbf{d}_i}$, $\lambdas\in\Rn^{\mathbf{d}_i\times 1}$, and $\tilde{\mathbf{U}}\trans\tilde{\mathbf{U}}=\tilde{\mathbf{U}}\tilde{\mathbf{U}}\trans=\mathbf{I}_{\mathbf{d}_i}$. Here, $\diag(\lambdas)$ denotes a diagonal matrix with $\lambdas$ as the main diagonal entries.

\textbf{Part (a).} We derive:
\beq
\ts  \LL_n^t  \triangleq  L_n + \beta^t\lambdaUp   \overset{\step{1}}{\leq}  \beta^t \lambdaUp (\delta+1) , \label{eq:bound:Lnt}
\eeq
\noi where step \step{1} uses Lemma \ref{lemma:updating:Sublinear} that $L_n\leq \delta \beta^t \lambdaUp$.

\textbf{Part (b).} We have:
\beq 
\|\H^t\|  \overset{\step{1}}{=}  \|\theta_2 \LL_n^t  - \beta^t \lambdas \|_{\infty}  \overset{\step{2}}{=}  \theta_2\LL_n^t-\min(\beta^t \lambdas)  \overset{\step{3}}{\leq}   \lambdaUp \beta^t \cdot (  \underbrace{{ \theta_2 (1+\delta)}     -  \lambdaDown' / \lambdaUp}_{\triangleq q } )  ,\nn
\eeq
\noi where step \step{1} uses $\|\theta_2 \LL_n^t  \mathbf{I} -  \beta^t  \A_n\trans\A_n\| = \|  \tilde{\mathbf{U}}\trans \diag( \theta_2 \LL_n^t - \beta^t \lambdas)\tilde{\mathbf{U}} \| = \|\theta_2 \LL_n^t  - \beta^t \lambdas \|_{\infty}$; step \step{2} uses the fact that $\|\rho - \x\|_{\infty} = \max(\rho - \x) = \rho - \min(\x)$ whenever $\rho \geq \max(\x)$ for all $\rho$ and $\x$; step \step{3} uses Inequality (\ref{eq:bound:Lnt}).

\textbf{Part (c).} Given $\u_n^{t+1}\triangleq \H^{t} (\x_n^{t+1} - \x_n^{t})$ as presented in Lemma \ref{lemma:first:order}, we have:
$\|\u_n^{t+1}\| \leq \|\H^{t}\|\cdot \|\x_n^{t+1} - \x_n^{t}\|\leq q  \lambdaUp \beta^t \|\x_n^{t+1} - \x_n^{t}\|$.

\end{proof}

\subsection{Proof of Lemma \ref{lemma:bound:dual:Case:I}} \label{app:lemma:bound:dual:Case:I}

\begin{proof}

For any $\sigma \in[1,2)$, we define $\sigma_1 \triangleq \frac{\sigma}{(1-|1-\sigma|)^2}$, and $\sigma_2\triangleq \frac{|1-\sigma|}{\sigma (1-|1-\sigma|)}$.

We define $\mathbbm{w}_n^{t+1} = \nabla h_{n} (\x_n^{t+1};\mu^t ) + \nabla f_n(\x_n^{t})$.

We define $\mathbbm{a}^{t+1} \triangleq  \A_n\trans (\z^{t+1}-\z^{t})$, and $\mathbbm{c}^{t} \triangleq  \u_n^{t} - \u_n^{t+1} +  \mathbbm{w}_n^{t} - \mathbbm{w}_n^{t+1}$.

We define $a = \tfrac{ \omega \sigma_2}{\lambdaDown}$, and $\AAA^t\triangleq \tfrac{1}{\beta^t}\|\mathbbm{a}^{t}\|_2^2$.

We define $b = \tfrac{3 \omega \sigma_1}{\lambdaDown}$, and $\BBB^t \triangleq \tfrac{1}{\beta^t} (L_{n} \| \x_n^{t} - \x_n^{t-1}\| + \|\u_n^{t}\| )^2$.

We define $\UUU^t \triangleq \frac{C_{h}^2 b}{\beta^t}\cdot(\frac{\mu^{t-1}}{\mu^t} -  1 )^2$.

First, we bound the term $\| \mathbbm{c}^{t} \|$. For all $t\geq 1$, we have:
\beq
&& \ts \| \mathbbm{c}^{t} \| =  \|\mathbbm{w}_n^{t} - \mathbbm{w}_n^{t+1} + \u_n^{t} - \u_n^{t+1}   \| \nn\\
&\overset{\step{1}}{ \leq }&  \ts\|  \nabla h_{n} (\x_n^{t+1};\mu^t ) - \nabla h_{n} (\x_n^{t} ;\mu^{t-1}) \| + \| \nabla f_n(\x_n^{t}) - \nabla f_n(\x_n^{t-1}) \| + \|\u_n^{t} - \u_n^{t+1}   \| \nn\\
&\overset{\step{2}}{ \leq }& \ts \|  \nabla h_{n} (\x_n^{t+1};\mu^t ) - \nabla h_{n} (\x_n^{t};\mu^{t-1} ) \| + L_{n} \| \x_n^{t} - \x_n^{t-1}\| + \|\u_n^{t} - \u_n^{t+1}   \| \nn\\
&\overset{}{ = }&  \ts\|  \nabla h_{n} (\x_n^{t+1};\mu^t ) - \nabla h_{n} (\x_n^{t};\mu^t ) + \nabla h_{n} (\x_n^{t};\mu^t ) - \nabla h_{n} (\x_n^{t};\mu^{t-1} ) \| \nn\\
&& \ts + L_{n} \| \x_n^{t} - \x_n^{t-1}\| + \|\u_n^{t} - \u_n^{t+1}   \| \nn\\
&\overset{\step{3}}{ \leq }& \ts \frac{1}{\mu^t}\| \x_n^{t+1}  - \x_n^{t} \| + ( \frac{\mu^{t-1}}{\mu^t} - 1 ) C_h + L_{n} \| \x_n^{t} - \x_n^{t-1}\| + \|\u_n^{t}\| + \| \u_n^{t+1}   \|,
\label{eq:c:bound}
\eeq
\noi where step \step{1} uses the triangle inequality; step \step{2} uses the fact that $f_n(\x)$ is $L_n$-smooth; step \step{3} uses Lemma \ref{lemma:lip:mu} and Lemma \ref{lemma:lip:mu:2}.

Second, we bound the term $\frac{\omega \sigma_1}{\lambdaDown\beta^t}\| \mathbbm{c}^{t} \|_2^2$. For all $t\geq 1$, we have:
\beq
&&\frac{\omega \sigma_1}{\lambdaDown\beta^t}\| \mathbbm{c}^{t} \|_2^2 \nn\\
&\overset{\step{1}}{ \leq } &  \frac{3 \omega \sigma_1}{\lambdaDown\beta^t} ( \frac{1}{\mu^t} \| \x_n^{t+1}  - \x_n^{t} \| + \| \u_n^{t+1}   \|)^2 + \underbrace{\frac{3 \omega \sigma_1}{\lambdaDown\beta^t}C_{h}^2( \tfrac{\mu^{t-1}}{\mu^t} - 1 )^2  }_{\triangleq \UUU^t} + \underbrace{ \frac{3 \omega \sigma_1}{\lambdaDown\beta^t} (L_{n} \| \x_n^{t} - \x_n^{t-1}\| + \|\u_n^{t}\| )^2 }_{\triangleq b\BBB^t }    \nn\\
&\overset{\step{2}}{ = } & \frac{ 3 \omega \sigma_1 }{\lambdaDown \beta^t} \{ ( \tfrac{1}{\mu^t} \| \x_n^{t+1}  - \x_n^{t} \| + \| \u_n^{t+1}   \|)^2 + (L_{n} \| \x_n^{t+1} - \x_n^{t}\| + \|\u_n^{t+1}\| )^2\}    + \UUU^t + b(\BBB^t - \BBB^{t+1}) \nn\\
&\overset{\step{3}}{ \leq } &  \frac{ 3 \omega \sigma_1 }{\lambdaDown \beta^t} \cdot 2(  ( \delta  + q  )\lambdaUp \beta^t \| \x_n^{t+1}  - \x_n^{t} \| )^2 + \UUU^t + b(\BBB^t - \BBB^{t+1}) \nn\\
&\overset{}{ = } & \underbrace{ 6 \omega \sigma_1 \kappa  (   \delta  + q    )^2}_{\triangleq \chi } \cdot \lambdaUp \beta^t \cdot \| \x_n^{t+1}  - \x_n^{t} \|_2^2    + \UUU^t + b(\BBB^t - \BBB^{t+1}) \nn\\
&\overset{\step{4}}{ \leq } & \chi \LL_n^t \| \x_n^{t+1}  - \x_n^{t} \|_2^2    + \UUU^t + b(\BBB^t - \BBB^{t+1}), \label{eq:c:squae:bound}
\eeq
\noi where step \step{1} uses Inequality \ref{eq:c:bound} and the fact that $(a+b+c)^2\leq 3a^2+3b^2+3c^2$ for all $a, b, c\in\Rn$; step \step{2} uses the definitions of $\{b, \BBB^t, \UUU^t\}$; step \step{3} uses Lemma \ref{lemma:reg:mu:L} that: $\tfrac{1}{\mu^t} \leq \delta \lambdaUp \beta^t $, $L_n \leq \delta\lambdaUp \beta^t$, and $\|\u_n^{t+1}\| \leq \|\H^{t}\|\cdot \|\x_n^{t+1} - \x_n^{t}\|\leq q  \lambdaUp \beta^t \|\x_n^{t+1} - \x_n^{t}\|$; step \step{4} uses $\beta^t \lambdaUp \leq \LL_n^t \triangleq \beta^t \lambdaUp + L_{n}$.

Finally, we derive the following inequalities for all $t\geq 1$:
\beq
\frac{\omega}{\sigma \beta^t}\|\mathbf{z}^{t+1}-\mathbf{z}^t\|_2^2 &\overset{\step{1}}{\leq}&  \frac{\omega}{\lambdaDown \sigma \beta^t} \|\A_n\trans (\mathbf{z}^{t+1}-\mathbf{z}^t)\|_2^2 = \tfrac{\omega}{\sigma \lambdaDown \beta^t}\|\mathbbm{a}^{t}\|_2^2 \nn\\
& \overset{\step{2}}{\leq}&   \tfrac{ \sigma_2 \omega}{\lambdaDown}   ( \frac{1}{\beta^t} \|\mathbbm{a}^{t}\|_2^2 - \frac{1}{\beta^{t}}\|\mathbbm{a}^{t+1}\|_2^2) + \frac{ \omega\sigma_1 }{\lambdaDown \beta^t}\|\mathbbm{c}^t\|_2^2 \nn\\
& \overset{\step{3}}{\leq}&   \underbrace{\frac{\sigma_2 \omega}{\lambdaDown}\cdot \frac{1}{\beta^t} \|\mathbbm{a}^{t}\|_2^2}_{\triangleq a\AAA^t} - \frac{\sigma_2 \omega}{\lambdaDown}\cdot \frac{1}{\beta^{t+1}}\|\mathbbm{a}^{t+1}\|_2^2  +  \frac{ \omega\sigma_1 }{ \lambdaDown}\cdot \tfrac{1}{\beta^t}\|\mathbbm{c}^t\|_2^2  \nn\\
& \overset{\step{4}}{\leq}&  a (\AAA^t - \AAA^{t+1}) + \chi \LL_n^t \| \x_n^{t+1}  - \x_n^{t} \|_2^2    + \UUU^t + b (\BBB^t - \BBB^{t+1}), \nn
\eeq
\noi where step \ding{172} uses $\lambdaDown \|\z\|_2^2 \leq \|\A_n\trans \z\|_2^2$ for all $\z$; step \step{2} uses Lemma \ref{useful:lemma:1} with $\b=\mathbbm{a}^{t}$, $\b^+ = \mathbbm{a}^{t+1}$, and $\mathbf{a}=\mathbbm{c}^{t}$ that:
\beq \label{eq:bound:at1:222}
\frac{1}{\sigma\beta^t} \|\mathbbm{a}^{t+1}\|_2^2 \leq   \frac{\sigma_2}{\beta^t}(   \|\mathbbm{a}^{t}\|_2^2 - \|\mathbbm{a}^{t+1}\|_2^2) + \frac{\sigma_1}{\beta^t}\|\mathbbm{c}^t\|_2^2; \nn
\eeq
\noi step \step{3} uses $-\frac{1}{\beta^t}\leq-\frac{1}{\beta^{t+1}}$; step \step{4} uses Inequality (\ref{eq:c:squae:bound}).


\end{proof}

\subsection{Proof of Lemma \ref{lemma:rule:I}}\label{app:lemma:rule:I}
\begin{proof}


We define $\varepsilon_1 \triangleq \tfrac{1}{2}\theta_1 - \tfrac{1}{2}$, and $\varepsilon_2 \triangleq  \theta_2-\frac{1}{2} - \chi$.

We define $f(\theta_2) \triangleq \left(\theta_2 -\frac{1}{2}\right) -  \varrho (  \delta + \theta_2 + \delta \theta_2 - 1/\kappa )^2$, and $\chi\triangleq \varrho (    \delta + \theta_2+\theta_2\delta-1/\kappa)^2$.

We define $\Theta^t \triangleq \Theta^{t}_{\oo} + \Theta^t_{+}$, where $\Theta^t_{+} \triangleq a\AAA^t + b\BBB^t$.

We define $\mathcal{E}^{t+1} \triangleq  \ts [\varepsilon_1\sum_{i=1}^{n-1} \LL_i^t \| \x_i^{t+1} - \x_i^t \|_2^2] + \varepsilon_2\LL_n^t \| \x_n^{t+1} - \x_n^t \|_2^2 + \frac{\varepsilon_3}{\beta^t}\|\z^{t+1}-\z^t\|_2^2$.

\textbf{Part (a).} With the choice $\theta_1=1.01$, it clearly holds that $\varepsilon_1 \triangleq \tfrac{1}{2}\theta_1 - \tfrac{1}{2} >0$.

Using Lemma \ref{lemma:ftheta} with $\theta = \theta_2$, we have $f(\theta_2)\geq \tfrac{1}{8\varrho}>0$, leading to $\varepsilon_2 \triangleq  \left(\theta_2 -\frac{1}{2}\right) - \chi>0$.

\textbf{Part (b).} Using Lemmas \ref{lemma:suf:dec} and \ref{lemma:bound:dual:Case:I}, we derive the following two respective inequalities:
\beq
 \ts \mathcal{E}^{t+1} + \Theta^{t+1}_{\oo} - \Theta^{t}_{\oo} & \leq & \ts  (\tfrac{1}{2} - \theta_2 +\varepsilon_2  ) \cdot \LL_n^t \|\x_n^{t+1}-\x_n^t\|_2^2 +   \tfrac{\omega}{\sigma \beta^t} \|\z^{t+1}-\z^t\|_2^2\label{eq:di:I:add1} \\
\ts \tfrac{\omega}{\sigma \beta^t}\|\mathbf{z}^{t+1}-\mathbf{z}^t\|_2^2 &\leq &\Theta^t_{+} - \Theta^{t+1}_{+} +  \chi \LL_n^t \| \x_n^{t+1}  - \x_n^{t} \|_2^2    + \UUU^t. \label{eq:di:I:add2}
\eeq
\noi Adding Inequalities (\ref{eq:di:I:add1}) and (\ref{eq:di:I:add2}) together, we have:
\beq
 \ts \mathcal{E}^{t+1} + \Theta^{t+1}- \Theta^t -  \UUU^t
 \leq  \ts \LL_n^t \|\x_n^{t+1}  - \x_n^{t}\|_2^2 \cdot \{     \tfrac{1}{2} - \theta_2 + \varepsilon_2 + \chi \}      \overset{\step{1}}{=}  0, \nn
\eeq
\noi where step \step{1} uses the definition of $\varepsilon_2 \triangleq  \theta_2-\frac{1}{2} - \chi$.

\end{proof}

\subsection{Proof of Lemma \ref{lemma:bound:dual:Case:A}} \label{app:lemma:bound:dual:Case:A}

\begin{proof}

For any $\sigma \in (0,1)$, we define $\sigma_1 \triangleq \frac{\sigma}{(1-|1-\sigma|)^2}$, and $\sigma_2\triangleq \frac{|1-\sigma|}{\sigma (1-|1-\sigma|)}$.

We define $\mathbbm{w}_n^{t+1} = \nabla h_{n} (\x_n^{t+1};\mu^t ) + \nabla f_n(\x_n^{t})$.

We define $\mathbbm{a}^{t+1} \triangleq  \A_n\trans (\z^{t+1}-\z^{t}) + \sigma \u_n^{t}$, and $\mathbbm{c}^{t} \triangleq \sigma \u_n^{t} +  \mathbbm{w}_n^{t} - \mathbbm{w}_n^{t+1}$.

We define $a \triangleq \tfrac{2\omega \sigma_2}{\lambdaDown}$, and $\AAA^t\triangleq \tfrac{1}{\beta^t}\|\mathbbm{a}^{t}\|_2^2$.

We define $b \triangleq \tfrac{6 \omega \sigma_1}{ \lambdaDown }$, and $\BBB^t \triangleq \tfrac{1}{\beta^t} (L_{n} \| \x_n^{t} - \x_n^{t-1}\| + \sigma \|\u_n^{t}\| )^2$.

We define $\UUU^t \triangleq \tfrac{ C_{h}^2 b }{  \beta^t}\cdot ( \frac{\mu^{t-1}}{\mu^t} - 1 )^2$.

First, we bound the term $\| \mathbbm{c}^{t} \|$. For all $t\geq 1$, we have:
\beq
&& \| \mathbbm{c}^{t} \| =  \|\mathbbm{w}_n^{t} - \mathbbm{w}_n^{t+1} + \sigma\u_n^{t}   \| \nn\\
&\overset{\step{1}}{ \leq }& \|  \nabla h_{n} (\x_n^{t+1};\mu^t ) - \nabla h_{n} (\x_n^{t};\mu^{t-1} ) \| + \| \nabla f_n(\x_n^{t}) - \nabla f_n(\x_n^{t-1}) \| +\sigma \|  \u_n^{t}  \| \nn\\
&\overset{\step{2}}{ \leq }& \|  \nabla h_{n} (\x_n^{t+1};\mu^t ) - \nabla h_{n} (\x_n^{t};\mu^{t-1} ) \| + L_{n} \| \x_n^{t} - \x_n^{t-1}\| + \sigma \|\u_n^{t} \| \nn\\
&\overset{}{ = }& \|  \nabla h_{n} (\x_n^{t+1};\mu^t ) - \nabla h_{n} (\x_n^{t};\mu^t ) + \nabla h_{n} (\x_n^{t};\mu^t) - \nabla h_{n} (\x_n^{t};\mu^{t-1} ) \| + L_{n} \| \x_n^{t} - \x_n^{t-1}\| + \sigma \|\u_n^{t}  \| \nn\\
&\overset{\step{3}}{ \leq }& \frac{1}{\mu^t}\| \x_n^{t+1}  - \x_n^{t} \| + ( \frac{\mu^{t-1}}{\mu^t} - 1 ) C_h + L_{n} \| \x_n^{t} - \x_n^{t-1}\| + \sigma \|\u_n^{t}\| , \label{eq:c:bound}
\eeq
\noi where step \step{1} uses the triangle inequality; step \step{2} uses the fact that $f_n(\x)$ is $L_n$-smooth; step \step{3} uses Lemma \ref{lemma:lip:mu:2} and Lemma \ref{lemma:lip:mu}.

Second, we bound the term $\frac{2 \omega \sigma}{\lambdaDown\beta^t}\|\u_n^{t} \|_2^2 + \frac{2 \omega }{\sigma \lambdaDown\beta^t}\| \mathbbm{c}^{t} \|_2^2$. For all $t\geq 1$, we have:
\beq
&& \frac{2 \omega \sigma}{\lambdaDown\beta^t}\|\u_n^{t+1} \|_2^2 + \frac{2 \omega\sigma_1 }{ \lambdaDown\beta^t}\| \mathbbm{c}^{t} \|_2^2\nn\\
&\overset{\step{1}}{ \leq } &   \frac{2 \omega \sigma}{\lambdaDown\beta^t}\|\u_n^{t+1} \|_2^2  +  \frac{6 \omega \sigma_1 }{ \lambdaDown\beta^t}  ( \tfrac{1}{\mu^t} \| \x_n^{t+1}  - \x_n^{t} \|   )^2  + \underbrace {\frac{6 \omega \sigma_1 }{  \lambdaDown\beta^t} ( \tfrac{\mu^{t-1}}{\mu^t} - 1 )^2 C_{h}^2}_{\UUU^t} + \underbrace {\frac{6 \omega \sigma_1 }{ \lambdaDown\beta^t} (L_{n} \| \x_n^{t} - \x_n^{t-1}\| + \sigma \|\u_n^{t}\| )^2}_{ \triangleq b \BBB^t  } \nn\\
&\overset{\step{2}}{ = } & \frac{2 \omega \sigma_1}{ \beta^t \lambdaDown} \cdot \{ \frac{\sigma}{\sigma_1}\|\u_n^{t+1} \|_2^2  + 3 ( \tfrac{1}{\mu^t} \| \x_n^{t+1}  - \x_n^{t} \|   )^2 +   3(L_{n} \| \x_n^{t+1} - \x_n^{t}\| + \sigma \|\u_n^{t+1}\| )^2 \}      + \UUU^t + b (\BBB^t - \BBB^{t+1}) \nn\\
&\overset{\step{3}}{\leq} &  \frac{2 \omega \sigma_1}{ \beta^t \lambdaDown}\cdot \lambdaUp^2 (\beta^t)^2 \cdot \{ \frac{ \sigma}{\sigma_1}   q^2   + 3 \delta^2 + 3( \delta + \sigma q  )^2 \}  \|\x_n^{t+1} - \x_n^{t}\|_2^2    + \UUU^t + b (\BBB^t - \BBB^{t+1}) \nn\\
&\overset{\step{4}}{\leq} & \underbrace{\frac{2 \omega \kappa}{\sigma} \cdot  \{ \sigma^2 q^2   + 3 \delta^2 + 3( \delta + \sigma q  )^2 \} }_{\triangleq \chi}\cdot \lambdaUp \beta^t \|\x_n^{t+1} - \x_n^{t}\|_2^2    + \UUU^t + b (\BBB^t - \BBB^{t+1})\nn\\
&\overset{\step{5}}{\leq} & \chi \cdot \LL_n^t \|\x_n^{t+1} - \x_n^{t}\|_2^2    + \UUU^t + b (\BBB^t - \BBB^{t+1}), \label{eq:c:squae:bound:2:2}
\eeq
\noi where step \step{1} uses Inequality (\ref{eq:c:bound}) and the fact that $(a+b+c)^2\leq 3a^2+3b^2+3c^2$ for all $a, b, c\in\Rn$; step \step{2} uses the definitions of $\{b,\BBB^t,\UUU^t\}$; step \step{3} uses $\|\u_n^{t+1}\| \leq \|\H^{t}\| \|\x_n^{t+1} - \x_n^{t}\|\leq \beta^t \lambdaUp q \|\x_n^{t+1} - \x_n^{t}\|$ and $L_n \leq \lambdaUp \beta^t \delta$, as has been shown respectively in Lemma \ref{lemma:reg:mu:L} and Lemma \ref{lemma:updating:Sublinear}, as well as the fact that $\frac{1}{\mu^t} = \beta^t \lambdaUp \delta$; step \step{4} uses $\kappa=\lambdaUp/\lambdaDown$, and the fact that $\sigma_1 = \frac{1}{\sigma}$ when $\sigma\in (0,1)$; step \step{5} uses $\beta^t \lambdaUp \leq \LL_n^t \triangleq \beta^t \lambdaUp + L_{n}$.

Finally, for all $t\geq 1$, we derive:
\beq
&&\tfrac{\omega}{\sigma\beta^t} \|\mathbf{z}^{t+1} - \mathbf{z}^{t}\|_2^2\nn\\
&\overset{\step{1}}{\leq}&\tfrac{\omega}{\sigma\beta^t\cdot \lambdaDown } \|\A_n\trans (\mathbf{z}^{t+1} - \mathbf{z}^{t})\|_2^2 \nn\\
&\overset{\step{2}}{=}&  \tfrac{\omega}{ \lambdaDown }  \cdot  \tfrac{1}{\sigma \beta^t} \| \mathbbm{a}^{t+1}  - \sigma \u_n^{t+1} \|_2^2 \nn\\
&\overset{\step{3}}{\leq}&  \tfrac{2\omega}{\lambdaDown} \cdot \{ \tfrac{1}{\sigma \beta^t } \| \mathbbm{a}^{t+1}  \|_2^2 +  \tfrac{\sigma  }{\beta^t}   \|  \u_n^{t+1} \|_2^2\} \nn\\
&\overset{\step{4}}{\leq}&  \tfrac{2\omega}{ \lambdaDown }  \cdot \{ \tfrac{\sigma_2}{\beta^t} \|\mathbbm{a}^{t}\|_2^2 - \tfrac{\sigma_2}{\beta^{t}}\|\mathbbm{a}^{t+1}\|_2^2  + \tfrac{\sigma_1}{  \beta^t}\|\mathbbm{c}^t\|_2^2 \}  +  \tfrac{2\omega \sigma  }{\beta^t \lambdaDown }   \|  \u_n^{t+1} \|_2^2 \nn\\
&\overset{\step{5}}{\leq}&  \underbrace{\ts \tfrac{2\omega}{ \lambdaDown }  \tfrac{\sigma_2}{\beta^t} \|\mathbbm{a}^{t}\|_2^2}_{\ts \triangleq a\AAA^t} - \underbrace{\ts \tfrac{2\omega}{ \lambdaDown } \tfrac{\sigma_2}{\beta^{t+1}}\|\mathbbm{a}^{t+1}\|_2^2}_{\ts \triangleq a\AAA^{t+1} }  +\tfrac{2\omega}{ \lambdaDown } \tfrac{\sigma_1}{  \beta^t}\|  \mathbbm{c}^{t} \|_2^2    +  \tfrac{2\omega \sigma  }{\beta^t \lambdaDown }   \|  \u_n^{t+1} \|_2^2 \nn\\
&\overset{\step{6}}{\leq}&  a\AAA^t -a\AAA^{t+1} + \chi \LL_n^t \|\x_n^{t+1}-\x_n^{t} \|_2^2     + \UUU^t + b\BBB^t - b\BBB^{t+1}, \nn
\eeq
\noi where step \step{1} uses the fact that $\lambdaDown \|\x\|_2^2 \leq \|\A_n\trans \x\|_2^2$ for all $\x$; step \step{2} uses the definition of $\mathbbm{a}^{t+1}$; step \step{3} uses the inequality $\| \mathbf{a} + \b \|_2^2 \leq 2\|\mathbf{a}\|_2^2 + 2 \|\b\|_2^2$ for all $\mathbf{a}$ and $\b$; step \step{4} uses Lemma \ref{useful:lemma:1} with $\b=\mathbbm{a}^t$, $\b^+=\mathbbm{a}^{t+1}$, and $\mathbf{a} = \mathbbm{c}^t$ that
\beq
 \frac{1}{\sigma \beta^t} \|\mathbbm{a}^{t+1}\|_2^2 &\leq & \frac{\sigma_1}{\beta^t}\|\mathbbm{c}^t\|_2^2 + \frac{\sigma_2}{\beta^t}  (   \|\mathbbm{a}^{t}\|_2^2 - \|\mathbbm{a}^{t+1}\|_2^2);\nn
\eeq
\noi step \step{5} uses $-\frac{1}{\beta^t}\leq - \frac{1}{\beta^{t+1}}$ and $\sigma_1=\frac{1}{\sigma}$ when $\sigma\in(0,1)$; step \step{6} uses Inequality (\ref{eq:c:squae:bound:2:2}).

\end{proof}

\subsection{Proof of Lemma \ref{lemma:rule:A}}\label{app:lemma:rule:A}
\begin{proof}

We define $\mathcal{E}^{t+1} \triangleq  \ts [\varepsilon_1\sum_{i=1}^{n-1} \LL_i^t \| \x_i^{t+1} - \x_i^t \|_2^2] + \varepsilon_2\LL_n^t \| \x_n^{t+1} - \x_n^t \|_2^2 + \frac{\varepsilon_3}{\beta^t}\|\z^{t+1}-\z^t\|_2^2$.

We define $\Theta^t \triangleq \Theta^{t}_{\oo} + \Theta^t_{+}$, where $\Theta^t_{+} \triangleq a\AAA^t + b\BBB^t$.

\textbf{Part (a).} We assume $\xi=\delta  = \sigma = \frac{c}{\kappa}$, where $c\in(0,1)$. We have:
\beq
\omega & \triangleq & 1 + \frac{\xi}{\sigma} = 2 \label{eq:delta:h:1} \\
q  &\triangleq&   \theta_2 + \theta_2 \delta   \overset{\step{1}}{\leq} \theta_2 + \theta_2 c. \label{eq:delta:h}
\eeq
\noi where step \step{1} uses $\delta = c/\kappa \leq c$ since $\kappa\geq 1$. We further obtain:
\beq
\varepsilon_2 & \triangleq & \theta_2 - \frac{1}{2} - \frac{6 \omega  \kappa}{\sigma }   \{   \tfrac{ 1    }{ 3 } \sigma^2  q^2 +   ( \delta    + \sigma   q )^2  + \delta^2 \}\nn\\
& \overset{\step{1}}{\geq} &  \theta_2 - \frac{1}{2} - \frac{12  }{c }   \{   \tfrac{ 1    }{ 3 } c^2  q^2 +   ( c    + c   q )^2  + c^2 \}\nn\\
& \overset{}{=} &  \theta_2 - \frac{1}{2} - 12 c   \{   \tfrac{ 1    }{ 3 }   q^2 +   ( 1    +    q )^2  + 1 \}\nn\\
& \overset{\step{2}}{\geq} &  \theta_2 - \frac{1}{2} - 12 c   \{   \tfrac{ (\theta_2 + \theta_2 c)^2    }{ 3 }    +   ( 1    +    \theta_2 + \theta_2 c )^2  + 1 \}\nn\\
& \overset{\step{3}}{>} &   0.02, \nn
\eeq
\noi where step \step{1} uses (\ref{eq:delta:h:1}), $\sigma\leq c$, $\delta\leq c$; step \step{2} uses (\ref{eq:delta:h}); step \step{3} uses the choice $c = 0.01$ and $\theta_2=1.5$.

 \textbf{Part (b).} Using Lemmas \ref{lemma:suf:dec} and \ref{lemma:bound:dual:Case:A}, we derive the following two respective inequalities:
\beq
\ts  \mathcal{E}^{t+1}  + \Theta_{\oo}^{t+1}  - \Theta_{\oo}^{t}  &  \leq & \ts  (\tfrac{1}{2} - \theta_2 +\varepsilon_2  ) \cdot \LL_n^t \|\x_n^{t+1}-\x_n^t\|_2^2 +   \tfrac{\omega}{\sigma \beta^t} \|\z^{t+1}-\z^t\|_2^2, \nn\\
\ts  \tfrac{\omega}{\sigma\beta^t} \|\mathbf{z}^{t+1} - \mathbf{z}^{t}\|_2^2 + \Theta^{t+1}_{+} - \Theta^{t}_{+} &\leq& \chi \LL_n^t \|\x_n^{t+1} - \x_n^{t}\|_2^2      + \UUU^t.\nn
\eeq
\noi Adding the two inequalities above together leads to:
\beq
 \ts \mathcal{E}^{t+1} + \Theta^{t+1}- \Theta^t
- \UUU^t  \leq \ts \LL_n^t \|\x_n^{t+1} - \x_n^{t}\|_2^2 \cdot \{  \tfrac{1}{2} - \theta_2 +\varepsilon_2 + \chi \}   \overset{\step{1}}{=}   0, \nn
\eeq
\noi where step \step{1} uses the definition of $\varepsilon_2 \triangleq \theta_2 - \frac{1}{2} - \chi$.

\end{proof}

\subsection{Proof of Lemma \ref{lemma:Theta:LB}}
\label{app:lemma:Theta:LB}

\begin{proof}

The proof of this lemma closely resembles that of Theorem 6 in \cite{Boct2019SIOPT}.

We denote $\underline{\Theta} \triangleq \underline{\Theta}'-  \mu^0 C_h^2$, where $\underline{\Theta}'$ is defined in Assumption \ref{ass:5}

Initially, for all $t\geq 1$, we have:
\beq \label{eq:key:inequalty}
\Theta^t &\overset{\step{1}}{=}& \ts \mathcal{L}(\x^{t},\mathbf{z}^{t};\beta^{t},\mu^t) +  \tfrac{1}{2}C_h \mu^t + a\AAA^t + b\BBB^t \nn\\
&\overset{\step{2}}{\geq}&  \ts \mathcal{L}(\x^{t},\mathbf{z}^{t};\beta^{t},\mu^t)   \nn\\
&\overset{\step{3}}{=}&  \ts h_{n}(\x^t_n;\mu^t) + \{\sum_{i=1}^{n-1} h_i(\x^t_i)\} +\sum_{i=1}^n f_i(\x^t_i) + \la \A\x^t - \b,\mathbf{z} \ra + \frac{\beta^t}{2}\| \A\x^t - \b\|_2^2  \nn\\
&\overset{\step{4}}{\geq}&  \ts - \mu^0 C_h^2 + \{\sum_{i=1}^{n} h_i(\x^t_i)\} +\{\sum_{i=1}^n f_i(\x^t_i)\} + \la \A\x^t - \b,\mathbf{z} \ra  + \frac{\beta^t}{2}\| \A\x^t - \b\|_2^2  \nn\\
&\overset{\step{5}}{\geq}&  \ts  \la \A\x^t - \b,\mathbf{z}^t\ra    +  \underbrace{\ts \underline{\Theta}'- \mu^0 C_h^2}_{\ts \triangleq \underline{\Theta}},
\eeq
\noi where step \step{1} uses the definition of $\Theta^t$; step \step{2} uses $\tfrac{1}{2}C_h \mu^t + a\AAA^t + b\BBB^t\geq 0$; step \step{3} uses the definition of $\mathcal{L}(\x^{t},\mathbf{z}^{t};\beta^{t},\mu^t)$ in Equation (\ref{eq:Lag}); step \step{4} uses $0 \leq h_n(\ub) - h_n(\ub;\mu) \leq \mu C_h^2$ as shown in Lemma \ref{lemma:lip:mu:2}, and the fact that $\mu^t \leq \mu^0$; step \step{5} uses Assumption \ref{ass:5}.

We now conclude the proof of this lemma through contradiction. Suppose that there exists $t_0\geq 1$ such that $\Theta^{t_0}<\underline{\Theta}$. We derive the following inequalities:
\beq
\ts \sum_{t=1}^{T} (\Theta^t -  \underline{\Theta} ) &=& \ts  [\sum_{t=1}^{t_0-1}( \Theta^t - \underline{\Theta} ) ] + [\sum_{t=t_0}^{T}(\Theta^t - \underline{\Theta})] \nn\\
&\leq & \ts  [ \sum_{t=1}^{t_0-1}(\Theta^t - \underline{\Theta}) ] + (T + 1 - t_0)\cdot \max_{t=t_0}^{T} (\Theta^{t} - \underline{\Theta}) \nn\\
&\overset{\step{1}}{\leq}& \ts  [ \sum_{t=1}^{t_0-1}(\Theta^t - \underline{\Theta}) ] + (T + 1 - t_0)\cdot (\Theta^{t_0} - \underline{\Theta}), \label{eq:close:look}
\eeq
\noi where step \step{1} uses $\Theta^{t}  \leq \Theta^{t_0}$ for all $t\geq t_0$. We closely examine Inequality (\ref{eq:close:look}). As $t_0$ is finite, the sum $\sum_{t=1}^{t_0-1}(\Theta^t - \underline{\Theta})$ is upper bounded. Considering the negativity of the term $(\Theta^{t_0}- \underline{\Theta})$, we deduce from Inequality (\ref{eq:close:look}):
\beq
\ts\lim_{T \rightarrow \infty}\,\sum_{t=1}^{T} (\Theta^t - \underline{\Theta})  = - \infty. \label{eq:go:to:inf}
\eeq
\noi Meanwhile, for all $t \geq 1$, the following inequalities hold:
\beq
\ts \Theta^t - \underline{\Theta} &\overset{\step{1}}{ \geq }&  \tfrac{1}{\sigma\beta^{t-1}}\la \mathbf{z}^{t}-\mathbf{z}^{t-1},\mathbf{z}^t \ra \nn\\
&\overset{\step{2}}{=}& \tfrac{1}{2\sigma}  \{\tfrac{1}{\beta^{t-1}}\|\mathbf{z}^t\|_2^2   - \tfrac{1}{\beta^{t-1}}\|\mathbf{z}^{t-1}\|_2^2+ \tfrac{1}{\beta^{t-1}}\|\mathbf{z}^t-\mathbf{z}^{t-1}\|_2^2 \}  \nn\\
&\overset{\step{3}}{\geq}& \tfrac{1}{2\sigma} \{ \tfrac{1}{\beta^{t}}\|\mathbf{z}^t\|_2^2   - \tfrac{1}{\beta^{t-1}}\|\mathbf{z}^{t-1}\|_2^2+ 0 \} , \label{eq:Qt:Q:low}
\eeq
\noi where step \ding{172} uses Inequality (\ref{eq:key:inequalty}) and $\mathbf{z}^{t+1} = \mathbf{z}^t + \sigma \beta^t (\A\x^{t+1} - \b)$; step \step{2} uses the Pythagoras relation in Lemma \ref{lemma:relation}; step \step{3} uses $\tfrac{1}{\beta^{t-1}}\geq\tfrac{1}{\beta^t}$.

\noi Telescoping Inequality (\ref{eq:Qt:Q:low}) over $t$ from $1$ to $T$, we have:
\beq
\ts \sum_{t=1}^{T} ( \Theta^t - \underline{\Theta}) \geq \frac{1}{2\sigma}\cdot \{\frac{1}{\beta^{T}}\|\mathbf{z}^{T}\|_2^2   - \frac{1}{\beta^{0}}\|\mathbf{z}^{0}\|_2^2\}  \geq - \frac{1}{2\sigma\beta^{0}}\|\mathbf{z}^{0}\|_2^2. \label{eq:df:right}
\eeq
\noi The finiteness of the right-hand-side in (\ref{eq:df:right}) contradicts with (\ref{eq:go:to:inf}).

Therefore, we conclude that $\Theta^t\geq \underline{\Theta}$ for all $t\geq 1$.

\end{proof}

\subsection{Proof of Lemma \ref{lemma:beta:mu}}\label{app:lemma:beta:mu}


\begin{proof}

We define $\overline{\rm{U}} \triangleq 3C_{h}^2\tfrac{b}{\beta^0}$.

We define $\UUU^t \triangleq  C_{h}^2 \tfrac{ b }{\beta^t}\cdot ( \frac{\mu^{t-1}}{\mu^t} - 1  )^2$, where $\beta^t  = \beta^0 (1 +\xi t^p)$, $\mu^t \propto \frac{1}{\beta^t}$.

\textbf{Part (a).} Letting $T\in [1,\infty)$, we obtain:
\beq\label{eq:mmmu:mu:sum}
 \ts \sum_{t=1}^T ( \frac{\mu^{t-1}}{\mu^t} - 1)^2 & \overset{\step{1}}{=}& \ts \sum_{t=1}^T ( \tfrac{\beta^t}{\beta^{t-1}} - 1)^2 =  \ts ( \tfrac{\beta^1}{\beta^{0}}-1)^2 + \sum_{t=2}^T ( \tfrac{\beta^t}{\beta^{t-1}} - 1)^2 \nn\\
& \overset{\step{2}}{=}& \ts (1+\xi 1^p-1)^2  + \sum_{t=1}^{T-1}  ( \tfrac{\beta^{t+1}}{\beta^{t}} - 1)^2 \nn\\
& \overset{\step{3}}{\leq }& \ts  1  + \sum_{t=1}^{\infty} \tfrac{(   \xi (t+1)^p - \xi t^p )^2}{( 1 +  \xi t^p)^2}  \nn\\
& \overset{\step{4}}{\leq }& \ts  1  + \sum_{t=1}^{\infty} (  \tfrac{ (t+1)^p - t^p }{t^p} )^2 \nn\\
& \overset{\step{5}}{\leq }& \ts  1  + 2,
\eeq
\noi where step \step{1} uses $\mu^t \propto \frac{1}{\beta^t}$; step \step{2} uses $\beta^1 = \beta^0 (1+\xi 1^p)$; step \step{3} uses the definition of $\beta^t  = \beta^0 +  \beta^0 \xi t^p$; step \step{4} uses $\tfrac{1}{( 1 +  \xi t^p)^2}\leq \tfrac{1}{(\xi t^p)^2}$; step \step{5} uses Lemma \ref{eq:two:two}.

We further obtain:
\beq
\ts \sum_{t=1}^{\infty} \UUU^t   \overset{\step{1}}{\leq}  \ts C_{h}^2 \tfrac{ b }{\beta^0}\cdot\{  \sum_{t=1}^{\infty}  ( \frac{\mu^{t-1}}{\mu^t} - 1 )^2 \}   \overset{\step{2}}{\leq} \ts  3 C_{h}^2 \tfrac{ b }{\beta^0} \triangleq \overline{\rm{U}} ,\nn
\eeq
\noi where step \step{1} uses $\beta^t\geq \beta^0$; step \step{2} uses Inequality (\ref{eq:mmmu:mu:sum}).

\textbf{Part (b).} For both conditions $\BIBI$ and $\SUSU$, we have from Lemmas (\ref{lemma:rule:I}) and 
(\ref{lemma:rule:A}):
\beq
\mathcal{E}^{t+1}      \leq  \Theta^t - \Theta^{t+1} + \UUU^t. \nn
\eeq

\noi Telescoping this inequality over $t$ from $1$ to $T$, we have:
\beq
\ts \sum_{t=1}^{T} \mathcal{E}^{t+1}  \overset{}{\leq}  \ts\Theta^{1} - \Theta^{T+1} + \sum_{t=1}^{T} \UUU^t  \overset{\step{1}}{\leq}  \ts \Theta^{1} - \underline{\Theta} + \overline{\rm{U}}  \triangleq  \overline{\mathcal{E}}, \label{eq:sum:of:square}
\eeq
\noi where step \step{1} uses Lemma \ref{lemma:Theta:LB} that $\Theta^{t}\geq\underline{\Theta}$ for all $t$, and Lemma \ref{lemma:beta:mu}.

\end{proof}



\subsection{ Proof of Lemma \ref{lemma:boundedness:z}}
\label{app:lemma:boundedness:z}

\begin{proof}

Given $\sigma \in (0,2)$, we define $\sigma_3 \triangleq \frac{ \sigma}{1 - |1-\sigma|} \in [1,\infty)$.

We define $\mathbbm{w}_n^{t+1} \triangleq \nabla h_{n} (\x_n^{t+1},\mu^t ) + \nabla f_n(\x_n^{t})$.

We define $\u_n^{t+1}\triangleq \H^{t} (\x_n^{t+1} - \x_n^{t})$, where $\H^{t} \triangleq \theta_2  \LL_n^{t}\mathbf{I} - \beta^{t}\A_n \trans \A_n$.

We define $Z \triangleq \tfrac{3}{\lambdaDown} \cdot (\tfrac{1}{\beta^0}\lambdaUp \|\z^1\|_2^2 + 2 \sigma_3 C_h^2 + 2 \sigma_3 C_f^2 + \sigma_3 q^2 \lambdaUp  \frac{\overline{\mathcal{E}}}{\varepsilon_2})$. We define $\ddot{Z} \triangleq \overline{\mathcal{E}}/\varepsilon_3$.

First, we have:
\beq
\ts \max_{i=1}^{\infty} \{ \|\w_{n}^{i+1} \|_2^2\} & = & \ts \max_{i=1}^{\infty}  \{  \| \nabla h_{n} (\x_n^{i+1},\mu^i ) + \nabla f_n(\x_n^{i}) \|_2^2 \}  \nn\\
&\overset{\step{1}}{\leq}& 2 C^2_h   + 2 C^2_f \label{eq:bound:zzzz:1}
\eeq
\noi where step \step{1} uses $\|\mathbf{a}+\mathbf{b}\|_2^2\leq2 \|\mathbf{a}\|_2^2+2\|\mathbf{b}\|_2^2$, and Assumption \ref{ass:1}.

Second, we have:
\beq\label{eq:bound:zzzz:2}
\ts \max_{i=1}^{\infty} \{  \frac{1}{\beta^i}\|\u_{n}^{i+1}\|_2^2  \} & = &\ts  \max_{i=1}^{\infty}  \{ \frac{1}{\beta^i}\|\H^{i} (\x_n^{i+1} - \x_n^{i})\|_2^2  \}  \nn\\
&\overset{\step{1}}{\leq}& \ts \max_{i=1}^{\infty}  \{ \frac{1}{\beta^i}   (q \lambdaUp \beta^i)^2 \|\x_{n}^{i+1}-\x_{n}^{i}\|_2^2  \}   \nn\\
&\overset{\step{2}}{\leq}& \ts q^2 \lambdaUp \cdot \sum_{i=1}^{\infty}  \{    \LL_n^i \|\x_{n}^{i+1}-\x_{n}^{i}\|_2^2  \}  \nn\\
&\overset{\step{3}}{\leq}& \ts q^2 \lambdaUp \cdot  \tfrac{\overline{\mathcal{E}}}{\varepsilon_2} ,
\eeq
\noi where step \step{1} uses $\|\H^{t}\|\leq  \beta^t \lambdaUp q$ for all $t\geq 0$, as shown in Lemma \ref{lemma:reg:mu:L}(b); step \step{2} uses $\beta^i\lambdaUp \leq \LL_n^i \triangleq \beta^i\lambdaUp + L_n$; step \step{3} uses $ \sum_{t=1}^{\infty}  \varepsilon_2\LL_n^t \| \x_n^{t+1} - \x_n^t \|_2^2 \leq \sum_{t=1}^{\infty} \mathcal{E}^{t+1} \leq \overline{\mathcal{E}}$.

\textbf{Part (a).} Using Lemma \ref{lemma:first:order}(b), we have:
\beq
\A_n\trans \z^{t+1} = (1-\sigma) \cdot \A_n\trans \z^t - \sigma \{\w_{n}^{t+1}+ \u_{n}^{t+1}\}.\nn
\eeq
\noi For all $t\geq 1$, we have:
\beq
\|\A_n\trans \z^{t+1}\|  \leq |1-\sigma| \cdot \|\A_n\trans \z^t \| + \sigma \{\|\w_{n}^{t+1} \| +  \|\u_{n}^{t+1}\| \}. \nn
\eeq
\noi Applying Lemma \ref{lemma:et1:Psi} with $e^{t}\triangleq \|\A_n\trans \z^{t}\|$ and $ p^t \triangleq \|\w_{n}^{t+1} \| +  \|\u_{n}^{t+1}\|$, for all $t\geq 1$, we obtain:
\beq\label{eq:def:Psi}
\ts \|\A_n\trans \z^{t}\|_2^2 &\leq&  \ts ( \|\A_n\trans \z^1\| + \sigma_3 \max_{i=1}^{t-1} \{ \|\w_{n}^{i+1} \| +  \|\u_{n}^{i+1}\| \})^2\nn\\
&\overset{\step{1}}{\leq}& \ts 3 \cdot \{ \lambdaUp \|\z^1\|_2^2 + \sigma_3 \max_{i=1}^{t-1} \|\w_{n}^{i+1} \|_2^2 + \sigma_3 \max_{i=1}^{t-1} \|\u_{n}^{i+1}\|_2^2 \} \nn\\
&\overset{}{=}& \ts \beta^t \cdot 3 \cdot \{ \lambdaUp \tfrac{1}{\beta^t}\|\z^1\|_2^2 + \sigma_3 \max_{i=1}^{t-1}  \tfrac{1}{\beta^t} \|\w_{n}^{i+1} \|_2^2 + \sigma_3 \max_{i=1}^{t-1} \tfrac{1}{\beta^t} \|\u_{n}^{i+1}\|_2^2 \} \nn\\
&\overset{\step{2}}{\leq}& \ts \beta^t \cdot 3 \cdot \{ \tfrac{1}{\beta^0}\lambdaUp \|\z^1\|_2^2 + \sigma_3 \max_{i=1}^{\infty} \tfrac{1}{\beta^i}\|\w_{n}^{i+1} \|_2^2 + \sigma_3 \max_{i=1}^{\infty} \tfrac{1}{\beta^i}\|\u_{n}^{i+1}\|_2^2 \} \nn\\
&\overset{\step{3}}{\leq}& \ts \beta^t \cdot \underbrace{ \ts 3 \cdot \{ \tfrac{1}{\beta^0}\lambdaUp \|\z^1\|_2^2 + 2 \sigma_3 C_h^2 + 2 \sigma_3 C_f^2 + \sigma_3 q^2 \lambdaUp  \frac{\overline{\mathcal{E}}}{\varepsilon_2} \} }_{\triangleq Z \lambdaDown}, \nn
\eeq
\noi where step \step{1} use $(a+b+c)^2\leq 3 (a^2+b^2+c^2)$, Assumption \ref{ass:2} that $\|\A_n\|_2^2\leq \lambdaUp$; step \step{2} uses $\beta^i\leq \beta^t$ for all $i\leq t$; \step{3} uses Inequalities (\ref{eq:bound:zzzz:1}) and (\ref{eq:bound:zzzz:2}). This further leads to $$\|\z^t\|_2^2\leq \tfrac{1}{\lambdaDown}\|\A_n\trans \z^{t}\|_2^2= \tfrac{1}{\lambdaDown}\cdot \lambdaDown Z\beta^t.$$

\textbf{Part (b).} We have:
\beq
\ts \tfrac{1}{\varepsilon_3}\sum_{t=1}^{\infty} \tfrac{\varepsilon_3}{\beta^t}\|\z^{t+1}-\z^t\|_2^2 \overset{\step{2}}{\leq} \tfrac{1}{\varepsilon_3}\ts\sum_{t=1}^{\infty} \mathcal{E}^{t+1}   \overset{\step{1}}{\leq}\overline{\mathcal{E}}/\varepsilon_3 \triangleq \ddot{Z}\ts, \nn
\eeq
\noi where step \step{1} uses the definition of $\mathcal{E}^{t+1} \triangleq  \ts [\varepsilon_1\sum_{i=1}^{n-1} \LL_i^t \| \x_i^{t+1} - \x_i^t \|_2^2] + \varepsilon_2\LL_n^t \| \x_n^{t+1} - \x_n^t \|_2^2 + \frac{\varepsilon_3}{\beta^t}\|\z^{t+1}-\z^t\|_2^2$ in Lemma \ref{lemma:suf:dec}; step \step{2} uses Lemma \ref{lemma:beta:mu}(b).

\end{proof}

\subsection{ Proof of Lemma \ref{lemma:boundedness}}
\label{app:lemma:boundedness}

\begin{proof}

We let $\sigma \in (0,2)$.

First, we derive the following inequalities:
\beq \label{eq:bound:AxbZ}
\ts  \la  \A\x^{t+1} - \b,\mathbf{z}^{t+1}\ra &\overset{}{=}& \tfrac{1}{\sigma\beta^{t}}\la \mathbf{z}^{t+1}-\mathbf{z}^{t},\mathbf{z}^{t+1} \ra \nn\\
&\overset{\step{1}}{=}& \tfrac{1}{2\sigma}  \{\tfrac{1}{\beta^{t}}\|\mathbf{z}^{t+1}\|_2^2   - \tfrac{1}{\beta^{t}}\|\mathbf{z}^{t}\|_2^2+ \tfrac{1}{\beta^{t}}\|\mathbf{z}^{t+1}-\mathbf{z}^{t}\|_2^2 \}  \nn\\
&\overset{}{\geq }& -\tfrac{1}{2\sigma \beta^{t}}\|\mathbf{z}^{t}\|_2^2,
\eeq
\noi where step \step{1} uses the Pythagoras relation in Fact \ref{lemma:relation}.

In view of Lemmas \ref{lemma:rule:I} and \ref{lemma:rule:A}, given $\mathcal{E}^{i+1}\geq0$ for all $i\geq 1$, we have:
\beq
 0   \leq  \Theta^i - \Theta^{i+1} + \UUU^i. \nn
\eeq
\noi Telescoping this inequality over $i$ from $1$ to $t$, we have:
\beq
0 \leq\ts \Theta^1 - \Theta^{t+1} + \sum_{i=1}^t  \UUU^i  \overset{\step{1}}{\leq}  \ts \Theta^1 - \Theta^{t+1} + \overline{\rm{U}}, \nn
\eeq
\noi where step \step{1} uses Lemma \ref{lemma:beta:mu}(b). For all $t\geq1$, we derive the following results:
\beq
\Theta^1 + \overline{\rm{U}} &\geq & \Theta^{t+1} \nn\\
&\overset{\step{1}}{=}& \Theta^{t+1}_{\oo} + a\AAA^{t+1} + b\BBB^{t+1}\nn\\
&\overset{\step{2}}{=}& \ts \mathcal{L}(\x^{t+1},\mathbf{z}^{t+1};\beta^{t+1},\mu^{t+1}) + \tfrac{1}{2} C_h \mu^{t+1} + a\AAA^{t+1} + b\BBB^{t+1} \nn\\
&\overset{\step{3}}{=}& \ts  \sum_{i=1}^n f_i(\x^{t+1}_i) + \la \A\x^{t+1} - \b,\mathbf{z}^{t+1} \ra + \frac{\beta^{t+1}}{2}\| \A \x^{t+1} - \b\|_2^2 \nn\\
 &&\ts   +  \{\sum_{i=1}^{n-1} h_i(\x^{t+1}_i)\} + h_{n}(\x^{t+1}_n;\mu^{t+1}) + \tfrac{1}{2} C_h \mu^{t+1} + a\AAA^{t+1} + b\BBB^{t+1}\nn\\
&\overset{\step{4}}{\geq}& \ts  \sum_{i=1}^n [ f_i(\x^{t+1}_i)+ h_i(\x^{t+1}_i)] + \la \A\x^{t+1}  - \b,\mathbf{z}^{t+1} \ra     -  \tfrac{1}{2}\mu^{t+1} C_h^2 \nn\\
&\overset{\step{5}}{\geq}&  \ts  \sum_{i=1}^n [ f_i(\x^{t+1}_i)+ h_i(\x^{t+1}_i)] - \tfrac{1}{2\sigma \beta^{t}}\|\mathbf{z}^{t}\|_2^2 -  \tfrac{1}{2}\mu^{t+1} C_h^2 \nn\\
&\overset{\step{6}}{\geq}&  \ts  \sum_{i=1}^n [ f_i(\x^{t+1}_i)+ h_i(\x^{t+1}_i)] - \tfrac{1}{2\sigma \beta^{t}}\|\mathbf{z}^{t}\|_2^2 - \tfrac{1}{2} \mu^{0} C_h^2, \nn
\eeq
\noi where step \step{1} uses the definition of $\Theta^{t+1}$; step \step{2} uses uses the definition of $\Theta^{t+1}_{\oo}$ in Lemma \ref{lemma:suf:dec}; step \step{3} uses the definition of $\mathcal{L}(\x^{t+1},\mathbf{z}^{t+1};\beta^{t+1},\mu^{t+1})$ in (\ref{eq:Lag}); step \step{4} uses $\frac{\beta^{t+1}}{2}\| \A \x^{t+1} - \b\|_2^2\geq 0$, $\tfrac{1}{2} C_h \mu^{t+1} + a\AAA^{t+1} + b\BBB^{t+1} \geq 0$, and the fact that $h_{n}(\x^{t+1}_n;\mu^{t+1}) \geq h_{n}(\x^{t+1}_n)  - \tfrac{1}{2}\mu^{t+1} C_h^2$; step \step{5} uses Inequality (\ref{eq:bound:AxbZ}); step \step{6} uses $\mu^{t}\leq \mu^0$ for all $t$.

We further obtain:
\beq
\ts\sum_{i=1}^n [ f_i(\x^{t+1}_i)+ h_i(\x^{t+1}_i)] &\leq&\ts \Theta^1 + \overline{\rm{U}} + \tfrac{1}{2\sigma \beta^{t}}\|\mathbf{z}^{t}\|_2^2 +  \tfrac{1}{2}\mu^{0} C_h^2\nn\\
&\overset{\step{1}}{<}& + \infty,\nn
\eeq
\noi where step \step{1} uses the boundedness of $\tfrac{1}{\beta^{t}}\|\mathbf{z}^{t}\|_2^2$ for all $t\geq0$, as shown in Lemma \ref{lemma:boundedness:z}. According to Assumption \ref{ass:5}, we have $\|\x^{t+1}_i\|<+\infty$ for all $i\in[n]$.

\end{proof}

\subsection{Proof of Theorem \ref{theorem:continuing:analysis}}
\label{app:theorem:continuing:analysis}
\begin{proof}

We define $K\triangleq  \tfrac{\overline{\mathcal{E}}}{K'}$, where $K' \triangleq \min \{ \min(\varepsilon_1,\varepsilon_2) \ADown^2,\varepsilon_3\}$, and $\ADown\triangleq \min_{i=1}^n \|\A_i\|$.

We define $\mathcal{E}^{t+1} \triangleq  \ts [\varepsilon_1\sum_{i=1}^{n-1} \LL_i^t \| \x_i^{t+1} - \x_i^t \|_2^2] + \varepsilon_2\LL_n^t \| \x_n^{t+1} - \x_n^t \|_2^2 + \frac{\varepsilon_3}{\beta^t}\|\z^{t+1}-\z^t\|_2^2$.

\textbf{Part (a).} We have:
\beq
\ts \overline{\mathcal{E}} &\overset{\step{1}}{\geq}& \ts \sum_{t=1}^{T} \mathcal{E}^{t+1} \nn\\
&\overset{\step{2}}{ = }&  \ts \sum_{t=1}^{T} \{\varepsilon_1\sum_{i=1}^{n-1} \LL_i^t \| \x_i^{t+1} - \x_i^t \|_2^2] + \varepsilon_2\LL_n^t \| \x_n^{t+1} - \x_n^t \|_2^2 + \frac{\varepsilon_3}{\beta^t}\|\z^{t+1}-\z^t\|_2^2 \} \nn\\
&\overset{\step{3}}{\geq}& \ts \tfrac{1}{\beta^T} \sum_{t=1}^{T} \{  [\varepsilon_1\sum_{i=1}^{n-1}  \frac{\LL_i^t}{\beta^t} \| \beta^t(\x_i^{t+1} - \x_i^t) \|_2^2] + \varepsilon_2 \frac{\LL_n^t}{\beta^t}  \|\beta^t( \x_n^{t+1} - \x_n^t) \|_2^2 + \varepsilon_3\|\z^{t+1}-\z^t\|_2^2 \} \nn\\
&\overset{\step{4}}{\geq}& \ts \tfrac{1}{\beta^T} \sum_{t=1}^{T} \{  [\varepsilon_1\sum_{i=1}^{n-1}  \ADown^2 \| \beta^t(\x_i^{t+1} - \x_i^t) \|_2^2] + \varepsilon_2 \ADown^2  \|\beta^t( \x_n^{t+1} - \x_n^t) \|_2^2 + \varepsilon_3\|\z^{t+1}-\z^t\|_2^2 \} \nn\\
&\overset{\step{5}}{\geq}& \ts \tfrac{1}{\beta^T}
 \cdot K' \cdot \sum_{t=1}^{T} \{  \sum_{i=1}^n \| \beta^t(\x_i^{t+1} - \x_i^t) \|_2^2 + \|\z^{t+1}-\z^t\|_2^2 \} \nn\\
&\overset{\step{6}}{=}& \ts \tfrac{1}{\beta^T}
 \cdot K' \cdot \sum_{t=1}^{T} \{  \| \beta^t(\x^{t+1} - \x^t) \|_2^2 + \|\z^{t+1}-\z^t\|_2^2 \},\nn
\eeq
\noi where step \step{1} uses Lemma (\ref{lemma:beta:mu})(b); step \step{2} uses the definition of $\mathcal{E}^{t+1}$; step \step{3} uses $\beta^T \geq \beta^t$ for all $t\leq T$; step \step{4} uses $\tfrac{\LL_i^t}{\beta^t} = \tfrac{L_i + \beta^t \|\A_i\|_2^2}{\beta^t}\geq\|\A_i\|_2^2\geq \ADown^2$; step \step{5} uses the definition of $K' \triangleq \min \{ \min(\varepsilon_1,\varepsilon_2) \ADown^2,\varepsilon_3\}$; step \step{6} uses $\sum_{i=1}^n \|\x_i^{t+1} - \x_i^t\|_2^2=\|\x^{t+1}-\x^t\|_2^2$. Therefore, we obtain:
\beq
\ts \sum_{t=1}^{T} \{  \| \beta^t(\x^{t+1} - \x^t) \|_2^2 + \|\z^{t+1}-\z^t\|_2^2 \} \leq \tfrac{\overline{\mathcal{E}}}{K'} \beta^T = K \beta^T. \nn
\eeq

\textbf{Part (b).} By dividing both sides of the above inequality by $T$, we obtain:
\beq
\ts\frac{K \beta^T}{T}&  \geq&  \ts \frac{1}{T}\sum_{t=1}^{T} \{  \| \beta^t(\x^{t+1} - \x^t) \|_2^2 + \|\z^{t+1}-\z^t\|_2^2 \} \nn\\
&\geq & \ts \min_{t=1}^T\{  \| \beta^t(\x^{t+1} - \x^t) \|_2^2 + \|\z^{t+1}-\z^t\|_2^2 \}.\nn
\eeq
\noi We conclude that there exists an index $\bar{t}$ with $\bar{t}\leq T$ such that $\|  \z^{\bar{t}+1}-\z^{\bar{t}}  \|_2^2 + \| \beta^{\bar{t}}(\x^{\bar{t}+1}-\x^{\bar{t}})\|_2^2 \leq \ts   \frac{K\beta^T}{T}$.

\end{proof}

\subsection{Proof of Theorem \ref{theorem:case:AAAAA:2}} \label{app:theorem:case:AAAAA:2}

To prove this theorem, we first provide the following lemma.

\begin{lemma} \label{eq:critical:123}
We define $\mathbf{q}^{t} \triangleq \{\x^{t}_1,\x^{t}_2,\ldots,\x^{t}_{n-1},\breve{\x}^{t}_n\}$. We have:
\begin{enumerate}[label=\textbf{(\alph*)}, leftmargin=16pt, itemsep=1pt, topsep=1pt, parsep=0pt, partopsep=0pt]

\item $\|\A \mathbf{q}^{t+1} - \b\| \leq c_{1} \|\z^{t+1}-\z^t\| + c_2 (\beta^t)^{-1}$.

\item $\dist(\zero, \partial h_n(\breve{\x}^{t+1}_n)  + \nabla_{\x_n} f_n(\breve{\x}^{t+1}_n) +  \A_n\trans \z^{t+1}) \leq  c_3 \|\z^{t+1}-\z^t\| + c_4 \|\beta^t(\x^{t+1} - \x^t)\| + c_5 (\beta^t)^{-1}$.

\item $\sum_{i=1}^{n-1} \dist (\zero,\partial h_i (\x^{t+1}_i) + \nabla_{\x_i} f_i(\x^{t+1}_i) +  \A_i\trans \z^{t+1}) \leq  c_6 \|\z^{t+1}-\z^t\| + c_7\| \beta^t(\x^{t+1} - \x^t)\|$.

\end{enumerate}

Here, $c_1 =\tfrac{1}{\sigma \beta^0}$, $c_2 = \AUp \tfrac{C_h}{ \delta \lambdaUp}$, $c_3 = (1 - \tfrac{1}{\sigma}) \AUp$, $c_4 = q \lambdaUp +  \tfrac{L_n}{\beta^0}$, $c_5 = \tfrac{L_n C_h}{\delta \lambdaUp}  $, $c_6 = (1 -   \tfrac{1}{\sigma})\AUp (n-1)$, and $c_7 =  \tfrac{ \overline{L} \sqrt{n-1} }{ \beta^0} + \theta_1   (\tfrac{\overline{L}}{\beta^0} +  \AUp^2 ) + \AUp^2 (n-1)$. Furthermore, $\AUp\triangleq \max_{i=1}^n \|\A_i\|$, and $\overline{L}\triangleq \max_{i=1}^n L_i$.

\begin{proof}

We define $\AUp\triangleq \max_{i=1}^n \|\A_i\|$, and $\overline{L}\triangleq \max_{i=1}^n L_i$.

We define $\u_i^{t+1} = \theta_1  \LL_i^t (\x_i^{t+1} - \x_i^{t}) - \beta^t  \A_i \trans [\sum_{j=i}^n \A_j(\x_j^{t+1}-\x_j^t)]$ with $i\in[n-1]$.

We define $\u_n^{t+1}\triangleq \H^{t} (\x_n^{t+1} - \x_n^{t})$ with $\H^{t} \triangleq \theta_2  \LL_n^{t}\mathbf{I} - \beta^{t}\A_n \trans \A_n$.

\textbf{Part (a).} We have:
\beq
 \ts \|\A \mathbf{q}^{t+1} - \b\|& = &\ts \| [ \sum_{i=1}^{n} \A_i \x_i^{t+1} ]  - \A_n \x^{t+1}_n + \A_n \breve{\x}^{t+1}_n - \b\|   \nn\\
&\overset{}{ \leq }&\ts \|\sum_{i=1}^{n} \A_i \x_i^{t+1}  - \b\| + \| \A_n(\x^{t+1}_n - \breve{\x}^{t+1}_n) \|   \nn\\
&\overset{\step{1}}{ \leq }&\ts \|\A \x^{t+1}  - \b\| + \AUp \mu^t C_h     \nn\\
&\overset{\step{2}}{ = }&\ts \| \tfrac{1}{\sigma \beta^t} (\z^{t+1}-\z^{t}) \| + \AUp  \tfrac{C_h}{ \delta \lambdaUp \beta^t},   \nn\\
&\overset{\step{3}}{\leq }&\ts \underbrace{\tfrac{1}{\sigma \beta^0}}_{\triangleq c_1} \|\z^{t+1}-\z^{t}\|+ \underbrace{\AUp ( \tfrac{C_h}{ \delta \lambdaUp } ) }_{\triangleq c_2} \cdot (\beta^t)^{-1} ,   \nn
\eeq
\noi where step \step{1} uses $\|\A_n\|\leq \AUp$ and Lemma \ref{lemma:smoothing:problem}(c); step \step{2} uses $\z^{t+1}=\z^t + \beta^t \sigma (\A\x^{t+1}-\b)$, and $\mu^t = \frac{1}{\delta \lambdaUp \beta^t}$; step \step{3} uses $\beta^0 \leq \beta^t$.

\textbf{Part (b).} We first have the following inequalities:
\beq
\|\nabla f_n(\x_n^t) - \nabla  f_n(\breve{\x}^{t+1}_n)  \| &=&\|\nabla f_n(\x_n^t)-\nabla f_n(\x_n^{t+1})+\nabla f_n(\x_n^{t+1})-\nabla  f_n(\breve{\x}^{t+1}_n)\|\nn\\
&\overset{}{ \leq }& \|\nabla f_n(\x_n^t)-\nabla f_n(\x_n^{t+1})\| + \|\nabla f_n(\x_n^{t+1})-\nabla  f_n(\breve{\x}^{t+1}_n)\|  \nn\\
&\overset{\step{1}}{ \leq }& L_n \| \x_n^{t+1} - \x_n^t\| + L_n \| \breve{\x}^{t+1}_n - \x_n^{t+1}\|  \nn\\
&\overset{\step{2}}{ \leq }& L_n \|\x_n^{t+1} - \x_n^t \| + L_n \mu^t C_h \nn\\
&\overset{\step{3}}{ \leq }&  L_n\| \x_n^{t+1} - \x_n^t\| + \underbrace{L_n \tfrac{1}{ (\delta \lambdaUp ) } C_h }_{\triangleq c_5} \cdot \tfrac{1}{\beta^t}, \label{eq:fxnt1:nt}
\eeq
\noi where step \step{1} uses the fact that $f_n(\x_n)$ is $L_n$-smooth; step \step{2} uses Lemma \ref{lemma:smoothing:problem}(c) that: $\| \breve{\x}^{t+1}_n - \x_n^{t+1}\|\leq \mu^t C_h$; step \step{3} uses $\mu^t = \frac{1}{\delta \lambdaUp \beta^t}$.

We further obtain:
\beq \label{eq:inp:imply:C}
&&  \ts \dist(\zero, \partial h_n (\breve{\x}^{t+1}_n) + \nabla f_i(\breve{\x}^{t+1}_i) +  \A_i\trans \z^{t+1} ) \nn\\
&\overset{\step{1}}{ = }& \ts  \|\theta_2 \LL_n^t (\cb^t- \x_n^{t+1}) + \nabla f_i(\breve{\x}^{t+1}_i) +  \A_i\trans \z^{t+1} \|  \nn\\
&\overset{\step{2}}{ = }& \ts  \|  \theta_2 \LL_n^t ( \x_n^{t} - \frac{1}{\theta_2 \LL_n^t} \ddot{\g}_n^t - \x_n^{t+1}   ) + \nabla  f_i(\breve{\x}^{t+1}_i) +  \A_i\trans \z^{t+1} \| \nn\\
&\overset{\step{3}}{ = }& \ts  \|  (\theta_2 \LL_n^t   -   \beta^t \A_n \trans  \A_n)(   \x_n^{t} -  \x_n^{t+1})  + \nabla  f_n(\breve{\x}^{t+1}_n) - \nabla f_n(\x_n^t)  + (1 - \frac{1}{\sigma})\A_n\trans (\mathbf{z}^{t+1}-\mathbf{z}^t) \| \nn\\
&\overset{\step{4}}{ \leq }& \ts  \|  \H(   \x_n^{t} -  \x_n^{t+1}) + (1 - \frac{1}{\sigma})\A_n\trans (\mathbf{z}^{t+1}-\mathbf{z}^t) \| + \|\nabla  f_n(\breve{\x}^{t+1}_n) - \nabla f_n(\x_n^t)   \| \nn\\
&\overset{\step{5}}{ \leq }& \ts  (1 - \frac{1}{\sigma}) \AUp \|\mathbf{z}^{t+1}-\mathbf{z}^t\| + \|\H(   \x_n^{t} -  \x_n^{t+1})\| + \|\nabla  f_n(\breve{\x}^{t+1}_n) - \nabla f_n(\x_n^t)\| \nn\\
&\overset{\step{6}}{ \leq }& \ts  \underbrace{ (1 - \frac{1}{\sigma}) \AUp}_{\triangleq c_3 } \|\mathbf{z}^{t+1}-\mathbf{z}^t\| + \{ \underbrace{  q \lambdaUp + L_n \tfrac{1}{\beta^0}}_{\triangleq c_4} \}  \cdot \|  \beta^t( \x_n^{t} -  \x_n^{t+1}) \| +  \tfrac{c_5}{ \beta^t },
\eeq
\noi where step \step{1} uses the optimality condition as shown in Lemma \ref{lemma:smoothing:problem}(b) that:
\beq
\rho (\cb^t - \x^{t+1}_n) \in \partial h_n(\breve{\x}^{t+1}_n),~\text{with}~\rho = \theta_2 \LL_n^t\nn;
\eeq
\noi step \step{2} uses $\cb^t =  \x_n^t - \ddot{\g}_n^t/\rho$ as shown in Algorithm \ref{alg:main}; step \step{3} uses the fact that:
\beq
\ddot{\g}_n^t = \ts  \nabla f_n(\x_n^t) + \A_n \trans\mathbf{z}^t  +  \frac{1}{\sigma}\A_n\trans (\mathbf{z}^{t+1}-\mathbf{z}^t) +   \beta^t \A_n \trans \A_n(\x_n^t-\x_n^{t+1}),\nn
\eeq
\noi step \step{4} uses the definition of $\H \triangleq \theta_2  \LL_n^{t}\mathbf{I} - \beta^{t}\A_n \trans \A_n$, as in Lemma \ref{lemma:first:order}; step \step{5} uses $\|\A_n\|\leq \AUp$; step \step{6} uses $\|\H^{t}\|\leq \beta^t \lambdaUp q$ as shown in Lemma \ref{lemma:reg:mu:L}, Inequality (\ref{eq:fxnt1:nt}), and the fact that $\beta^0\leq \beta^t$.

\textbf{Part (c).} We first have the following inequalities:
\beq
\ts \sum_{i=1}^{n-1} \| \u_i^{t+1}\| &\overset{\step{1}}{ = }& \ts \| \sum_{i=1}^{n-1} \{ \theta_1  \LL_i^t (\x_i^{t+1} - \x_i^{t}) - \beta^t  \A_i \trans [\sum_{j=i}^n \A_j(\x_j^{t+1}-\x_j^t)]\}\| \nn\\
&\overset{}{ \leq }& \ts \theta_1 \| \sum_{i=1}^{n-1}    \LL_i^t (\x_i^{t+1} - \x_i^{t})\| + \beta^t\| \sum_{i=1}^{n-1}    \A_i \trans [\sum_{j=i}^n \A_j(\x_j^{t+1}-\x_j^t)]\|\nn\\
&\overset{\step{2}}{ \leq }& \ts \theta_1  ( \max_{i=1}^n \LL_i^t) \sum_{i=1}^{n-1} \| \x_i^{t+1} - \x_i^{t}  \|  + \beta^t \AUp^2 (n-1) \|\x^{t+1}-\x^t\|\nn\\
&\overset{\step{3}}{ \leq }& \ts  \theta_1  (\max_{i=1}^{n-1} \LL_i^t) \sqrt{n-1} \| \x^{t+1} - \x^{t}  \| + \AUp^2  (n-1) \|   \beta^t( \x^{t+1}-\x^t ) \| \nn\\
&\overset{\step{4}}{ \leq }& \ts \theta_1  ( \tfrac{\overline{L}}{\beta^0} +  \AUp^2) \cdot  \| \beta^t(\x^{t+1} - \x^{t})  \|+ \AUp^2 (n-1) \| \beta^t(\x^{t+1}-\x^t)  \| \nn\\
&\overset{ }{ = }& \ts \{ \underbrace{ \ts  \theta_1   (\tfrac{\overline{L}}{\beta^0} +  \AUp^2 ) + \AUp^2 (n-1)}_{\triangleq c'_7} \} \cdot \| \beta^t(\x^{t+1} - \x^{t})  \|,  \label{eq:uuuuu}
\eeq
\noi where step \step{1} uses the definition of $\u_i^{t+1}$ for all $i\in[n-1]$; step \step{2} uses $\AUp\triangleq \max_{i=1}^n \|\A_i\|$; step \step{3} uses $\sum_{j=1}^{n-1} \|    \x_j^{t+1}-\x_j^t  \| \leq \sqrt{n-1} \|\x^{t+1}-\x^t\|$; step \step{4} uses $\LL_i^t = L_i + \beta^t \|\A_i\|_2^2 \leq \frac{\beta^t L_i}{\beta^0} + \beta^t \AUp^2 \leq \tfrac{\beta^t \overline{L}}{\beta^0} + \beta^t \AUp^2$ for all $i$.

We have the following results:
\beq \label{eq:inp:imply:C}
&&  \ts \sum_{i=1}^{n-1} \dist(\partial h_i (\x^{t+1}_i) + \nabla  f_i(\x^{t+1}_i) +  \A_i\trans \z^{t+1} )  \nn\\
&\overset{\step{1}}{ = }&  \ts \sum_{i=1}^{n-1} \|(1 -   \frac{1}{\sigma} ) \A_i\trans(\mathbf{z}^{t+1}-\mathbf{z}^{t}) - \nabla f_i(\x_i^{t}) -  \u_i^{t+1} + \nabla  f_i(\x^{t+1}_i) \| \nn\\
&\overset{}{ \leq }&  \ts \sum_{i=1}^{n-1} \|  (1 -   \frac{1}{\sigma} ) \A_i\trans(\mathbf{z}^{t+1}-\mathbf{z}^{t}) \|  +  \sum_{i=1}^{n-1} \| \nabla f_i(\x_i^{t}) -  \nabla  f_i(\x^{t+1}_i) \|+\sum_{i=1}^{n-1} \|   \u_i^{t+1} \|   \nn\\
&\overset{\step{2}}{ \leq }&  \ts \underbrace{(1 -   \frac{1}{\sigma} ) \AUp (n-1)}_{\triangleq c_6} \| \mathbf{z}^{t+1}-\mathbf{z}^{t} \|_2^2 +  \tfrac{\overline{L}}{\beta^0} \sqrt{n-1}\|\beta^t( \x^{t} - \x^{t+1})\|+\sum_{i=1}^{n-1} \|\u_i^{t+1} \| \nn\\
&\overset{\step{3}}{ = }& c_6 \| \z^{t+1}-\z^{t} \|_2^2 + \underbrace{ \ts (\tfrac{ \overline{L} \sqrt{n-1} }{ \beta^0} + c'_7)  }_{\triangleq c_7} \cdot  \| \beta^t( \x^{t} - \x^{t+1})  \|_2^2,   \nn
\eeq
\noi where step \step{1} uses Part (\bfit{a}) in Lemma \ref{lemma:first:order} that:
\begin{align}
i\in[n-1],~\partial h_i (\x_i^{t+1} ) &\ni -\u_i^{t+1} -\A_i\trans\mathbf{z}^t -   \frac{1}{\sigma}\A_i\trans(\mathbf{z}^{t+1}-\mathbf{z}^{t}) - \nabla f_i(\x_i^{t}) ;\nn
\end{align}
\noi step \step{2} uses $\|\A_i\| \leq \AUp$, $f_i(\x_i)$ is $L_i$-smooth, $L_i\leq \overline{L}$, and $\beta^0\leq \beta^t$; step \step{3} uses Inequality (\ref{eq:uuuuu}).

\end{proof}
\end{lemma}

Now, we proceed to prove the theorem.

\begin{proof}

We define $\Crit({\x},{\mathbf{z}}) \triangleq \|\A{\x}- \b \| + \sum_{i=1}^n \dist(\zero, \nabla f_i({\x}_i) + \partial h_i({\x}_i)   + \A_i\trans {\mathbf{z}})$.

We define $\mathbf{q}^{t} \triangleq \{\x^{t}_1,\x^{t}_2,\ldots,\x^{t}_{n-1},\breve{\x}^{t}_n\}$.

First, we deduce from Theorem \ref{theorem:continuing:analysis}(a) that
\beq
\ts K\beta^T &\geq& \ts  \sum_{t=1}^T \|\z^{t+1} - \z^t\|_2^2 + \|\beta^t(\x^{t+1} - \x^t)\|_2^2\nn\\
&\overset{\step{1}}{ \geq }& \ts \tfrac{1}{2T} ( \sum_{t=1}^T \|\z^{t+1} - \z^t\| + \|\beta^t(\x^{t+1} - \x^t)\|)^{2} \nn\
\eeq
\noi where step \step{1} uses the inequality $\|\a\|_2^2\geq \tfrac{1}{2T}(\|\a\|_1)^2$ for all $\a\in\Rn^{2T}$. This leads to:
\beq
\ts \sum_{t=1}^T \{ \|\z^{t+1} - \z^t\| + \|\beta^t(\x^{t+1} - \x^t)\|\} \leq \ts \sqrt{K\beta^T \cdot 2T} =  \ts \OO(T^{(p+1)/2}). \label{eq:sum:l1:norm:upperbound}
\eeq

Second, using lemma \ref{eq:critical:123}, for all $t\geq 0$, we have:
\beq
&& \Crit(\mathbf{q}^{t+1},\mathbf{z}^{t+1}) \nn\\
& \leq  &  \underbrace{ (c_{1} + c_{3} + c_{6})}_{\triangleq d_1} \|\z^{t+1}-\z^t\| + \underbrace{  (c_4+c_7) }_{\triangleq d_2} \|\beta^t(\x^{t+1} - \x^t)\| + \underbrace{(c_2 + c_5)}_{ \triangleq d_3 }  (\beta^t)^{-1}. \label{eq:c:p:bound}
\eeq
We further derive:
\beq
&& \ts \tfrac{1}{T} \sum_{t=1}^T \Crit(\mathbf{q}^{t+1},\mathbf{z}^{t+1}) \nn\\
&\overset{\step{1}}{\leq }& \ts  \tfrac{1}{T} \max(d_1,d_2) \sum_{t=0}^T\{\|\z^{t+1}-\z^t\| + \|\beta^t(\x^{t+1} - \x^t)\| \} +   \frac{d_3}{T} \sum_{t=0}^T (\beta^t)^{-1}  \nn\\
&\overset{\step{2}}{\leq }& \ts  \OO(T^{(p-1)/2}) +  \frac{d_3}{T} \sum_{t=0}^T (\beta^t)^{-1}  \nn\\
&\overset{\step{3}}{\leq }& \ts  \OO(T^{(p-1)/2}) + \OO(T^{1-p-1}),\nn
\eeq
\noi Here, step \step{1} uses Inequality \ref{eq:c:p:bound}; step \step{2} uses Inequality (\ref{eq:sum:l1:norm:upperbound}); step \step{3} uses $\beta^t=\OO(t^p)$, and the fact that $\sum_{t=1}^T t^{-p}\leq \tfrac{T^{(1-p)}}{1-p}$ if $p\in(0,1)$ (refer to Lemma \ref{eq:lp:function:sum:upperbound}).

In particular, with the choice $p=1/3$, we have: $\frac{1}{T} \sum_{t=1}^T \Crit(\mathbf{q}^{t+1},\mathbf{z}^{t+1})\leq \mathcal{O}(T^{-1/3})$.

\end{proof}

\subsection{Proof of Lemma \ref{lemma:bound:u}} \label{app:lemma:bound:u}

\begin{proof}

We let $\tfrac{\z^t}{\sqrt{\beta^t}} \triangleq \hat{\z}^t$ for all $t$.

Initially, we derive:
\beq
\ts \sum_{t=1}^{\infty} (1 - \sqrt{\frac{ \beta^t}{\beta^{t+1}}} )^2 &\overset{\step{1}}{= }& \ts \sum_{t=1}^{\infty} (1 - \sqrt{\frac{ 1 +  \xi t^p  }{ 1 +  \xi (t+1)^p}} )^2 \nn\\
&\overset{\step{2}}{\leq}& \ts \sum_{t=1}^{\infty} (1 - \sqrt{\tfrac{ t^p }{  (t+1)^p}} )^2 \nn\\
&\overset{}{=}& \ts \sum_{t=1}^{\infty} \frac{ \{ (t+1)^{p/2} - t^{p/2} \}^2 }{ (t+1)^p }  \nn\\
&\overset{\step{3}}{\leq}&\ts    \sum_{t=1}^{\infty} \frac{ \{  \frac{p}{2} \cdot t^{(p/2-1)} \}^2 }{ t^p } \nn \\
&\overset{\step{4}}{\leq}&\ts  \frac{1}{4} \sum_{t=1}^{\infty} \frac{ t^{(p-2)}  }{ t^p } \nn \\
&\overset{\step{5}}{\leq}&\ts 1/2 ,\label{eq:beta:beta:sqrt}
\eeq
\noi where step \step{1} uses $\beta^t =\beta^0 (1 +\xi t^p ) $ for all $t\geq0$; step \step{2} uses $\frac{ 1 + \xi t^p  }{1 +  \xi (t+1)^p} \leq \frac{ \xi t^p  }{ \xi (t+1)^p}$; step \step{3} uses Lemma \ref{eq:two:two:before} that $(t+1)^{p/2} - t^{p/2} \leq \frac{p}{2} t^{(p/2-1)}$ for all $t\geq 1$ and $\frac{p}{2} \in (0,1)$; step \step{4} uses $p\leq 1$ and $\frac{1}{t+1}\leq \frac{1}{t}$; step \step{5} uses $\sum_{t=1}^{\infty}\frac{1}{t^2}=\frac{\pi^2}{6}<2$.

\textbf{Part (a).} We have: $\|\hat{\z}^{t}\|_2^2  =  \|\frac{\z^t}{\sqrt{\beta^t}} \|_2^2  =  \frac{1}{\beta^t} \| \z^t \|_2^2 \overset{\step{1}}{\leq} Z < +\infty$, where step \step{1} uses Lemma \ref{lemma:boundedness:z}.

\textbf{Part (b).} We derive the following inequalities:
\beq
\ts \sum_{t=1}^{\infty} \|\hat{\z}^{t+1} - \hat{\z}^t\|_2^2 &\overset{\step{1}}{=}& \ts  \sum_{t=1}^{\infty} \| \frac{\z^{t+1}}{\sqrt{\beta^{t+1}}} - \frac{\z^t}{\sqrt{\beta^t}}\|_2^2\nn \\
& = & \ts \sum_{t=1}^{\infty} \| \tfrac{\z^{t+1}-\z^{t}}{\sqrt{\beta^{t+1}}} - \z^t (\frac{1}{\sqrt{\beta^t}} - \tfrac{1}{\sqrt{\beta^{t+1}}} ) \|_2^2\nn\\
&\overset{\step{2}}{\leq }& \ts 2\sum_{t=1}^{\infty} \| \frac{\z^{t+1}-\z^{t}}{\sqrt{\beta^{t+1}}} \|_2^2 + 2\sum_{t=1}^{\infty} \| \z^t (\frac{1}{\sqrt{\beta^t}} - \frac{1}{\sqrt{\beta^{t+1}}} ) \|_2^2\nn\\
&\overset{\step{3}}{\leq }& \ts 2\sum_{t=1}^{\infty} \frac{1}{\beta^t} \|\z^{t+1}-\z^{t}\|_2^2 + 2\sum_{t=1}^{\infty} \tfrac{1}{\beta^t}\| (1 - \sqrt{\tfrac{\beta^t}{\beta^{t+1}}} ) \cdot \z^t  \|_2^2\nn\\
&\overset{\step{4}}{\leq }& \ts 2 \ddot{Z} +  2 Z  \cdot \sum_{t=1}^{\infty} (1 - \sqrt{\tfrac{\beta^t}{\beta^{t+1}}} )^2\nn\\
&\overset{\step{5}}{\leq }& \ts 2 \ddot{Z} + 2  Z \cdot \tfrac{1}{2} \nn\\
&<& \ts 2 (\ddot{Z} +  Z),\nn
\eeq
\noi where step \step{1} uses the definition $\tfrac{\z^t}{\sqrt{\beta^t}} \triangleq \hat{\z}^t$ for all $t$; step \step{2} uses $\|\mathbf{a}-\mathbf{b}\|_2^2\leq 2 \|\mathbf{a}\|_2^2+ 2 \|\mathbf{b}\|_2^2$; step \step{3} uses $\frac{1}{\beta^{t+1}}\leq \frac{1}{\beta^{t}}$; step \step{4} uses $\sum_{t=1}^{\infty} \frac{1}{\beta^t} \|\z^{t+1}-\z^{t}\|_2^2\leq \ddot{Z}$ and $\tfrac{1}{\beta^t} \|\z^{t}\|_2^2\leq Z$, as shown in Lemma \ref{lemma:boundedness:z}; step \step{5} uses Inequality (\ref{eq:beta:beta:sqrt}).

\end{proof}

\section{Additional Experiment Details and Results}
\label{app:sect:exp}

We offer further experimental details in Sections \ref{app:sect:data:set} and \ref{app:sect:O:proj}, and include additional results in Section \ref{app:sect:add:exp}.

\subsection{Datasets} \label{app:sect:data:set}

We incorporate four datasets in our experiments, including both randomly generated data and publicly available real-world data. These datasets serve as our data matrices $\mathbf{D}\in\mathbb{R}^{\dot{m}\times \dot{d}}$. The dataset names are as follows: `TDT2-$\dot{m}$-$\dot{d}$', `sector-$\dot{m}$-$\dot{d}$', `mnist-$\dot{m}$-$\dot{d}$', and `randn-$\dot{m}$-$\dot{d}$'. Here, ${\text{randn(}m,n)}$ refers to a function that generates a standard Gaussian random matrix with dimensions $m\times n$. The matrix $\mathbf{D}\in\mathbb{R}^{\dot{m}\times \dot{d}}$ is constructed by randomly selecting $\dot{m}$ examples and $\dot{d}$ dimensions from the original real-world dataset (\url{http://www.cad.zju.edu.cn/home/dengcai/Data/TextData.html},\url{https://www.csie.ntu.edu.tw/~cjlin/libsvm/}). We normalize each column of $\mathbf{D}$ to have a unit norm and center the data by subtracting the mean.


\subsection{Projection on Orthogonality Constraints} \label{app:sect:O:proj}

When $h(\x) = \iota_{\mathcal{M}}(\mat(\x))$ with $\Omega\triangleq \{\V\,|\,\V\trans\V=\I\}$, computing the proximal operator reduces to the following optimization problem:
\beq
\textstyle \bar{\x} \in \arg\min_{\x} \frac{\mu}{2}\|\x-\x'\|_2^2,\,s.t.\,\mat(\x)\in\mathcal{M} \triangleq \{\mathbf{V}\,|\,\mathbf{V}\trans \mathbf{V} = \mathbf{I}\}.\nn
\eeq
\noi This is the nearest orthogonality matrix problem, and the optimal solution can be computed as $\bar{\x}=\vec(\hat{\mathbf{U}}\hat{\mathbf{V}}\trans)$, where $\mat(\x') = \hat{\mathbf{U}}\Diag(\mathbf{s})\hat{\mathbf{U}}\trans$ is the singular value decomposition of the matrix $\mat(\x')$. Please refer to \cite{lai2014splitting}.

\subsection{Additional Experiment Results}
\label{app:sect:add:exp}

\begin{figure}[!t]
\centering
\begin{subfigure}{.24\textwidth}\centering\includegraphics[width=1.12\linewidth]{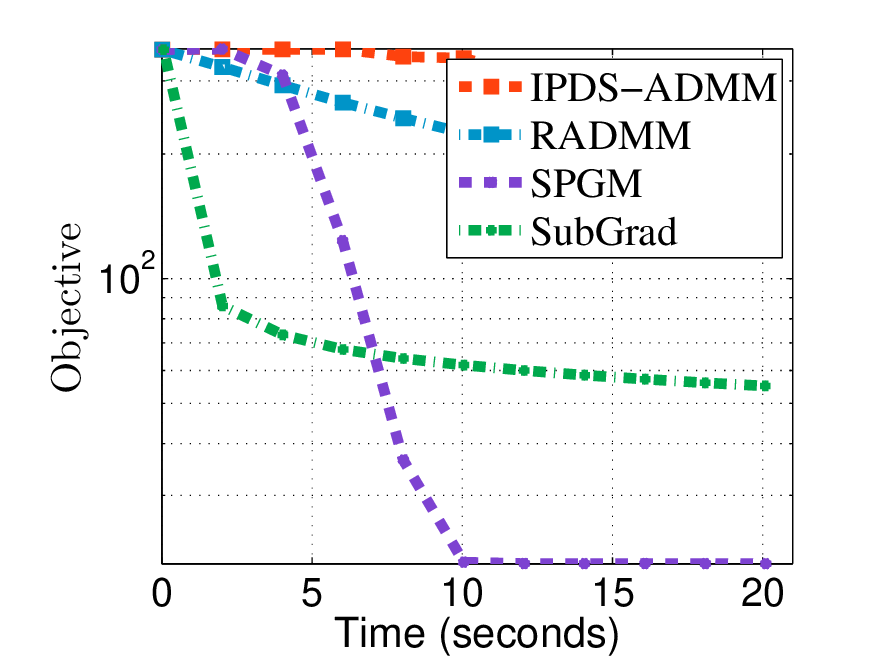}\caption{\scriptsize randn-1500-500}\label{fig:sub1}\end{subfigure}
\begin{subfigure}{.24\textwidth}\centering\includegraphics[width=1.12\linewidth]{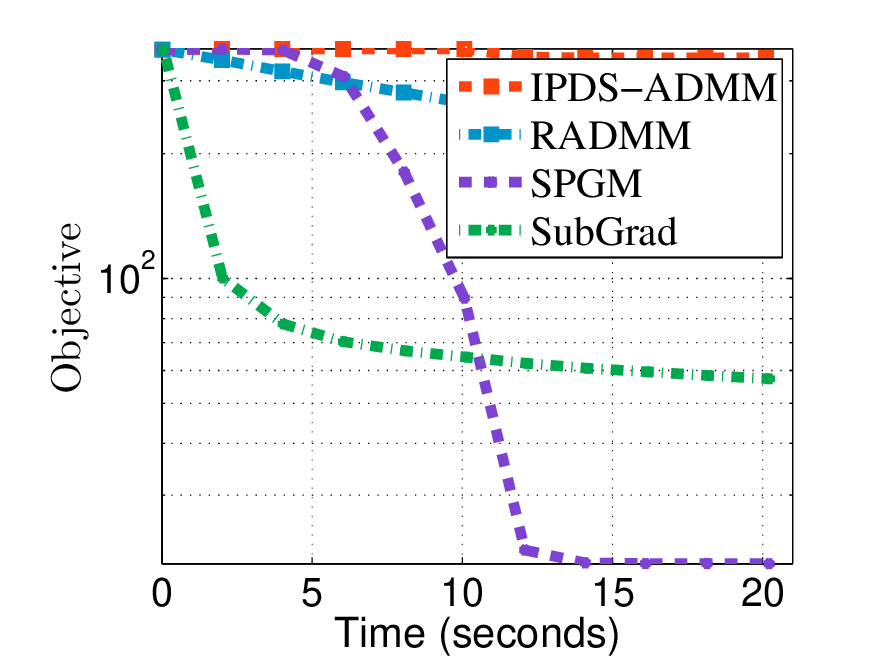}\caption{\scriptsize randn-2500-500}\label{fig:sub2}\end{subfigure}
\begin{subfigure}{.24\textwidth}\centering\includegraphics[width=1.12\linewidth]{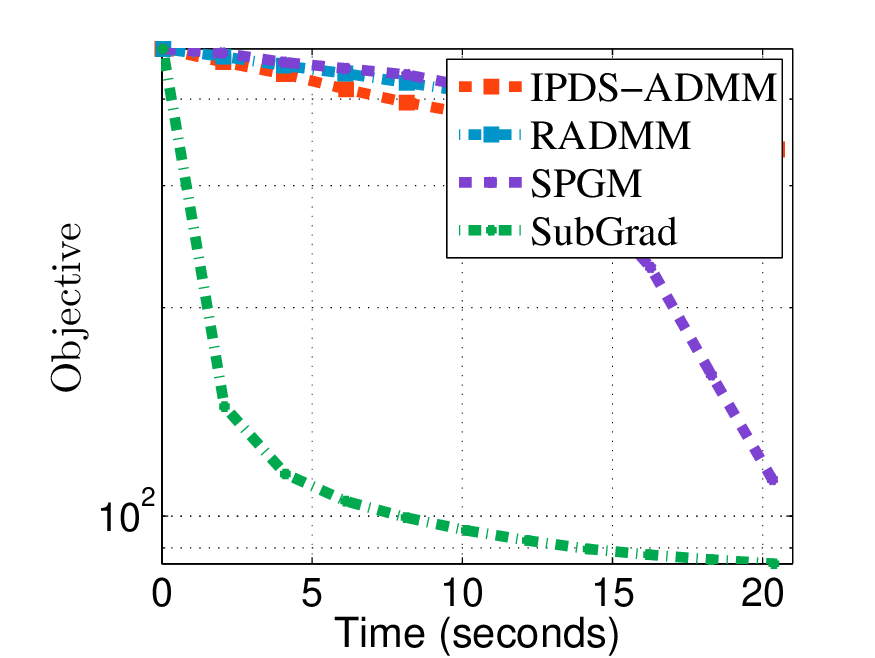}\caption{\scriptsize mnist-1500-780}\label{fig:sub3}\end{subfigure}
\begin{subfigure}{.24\textwidth}\centering\includegraphics[width=1.12\linewidth]{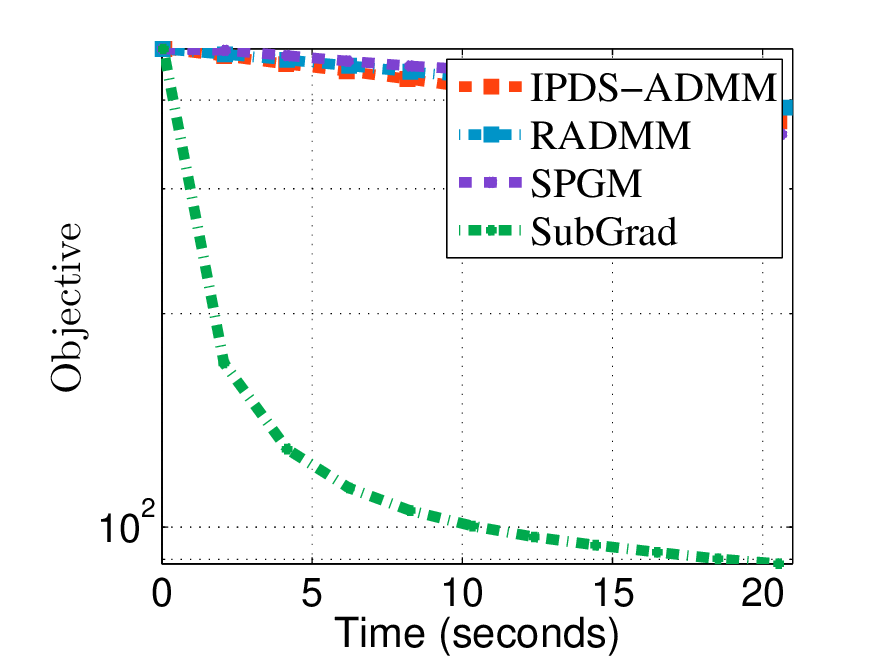}\caption{\scriptsize mnist-2500-780}\label{fig:sub4}\end{subfigure}

\centering
\begin{subfigure}{.24\textwidth}\centering\includegraphics[width=1.12\linewidth]{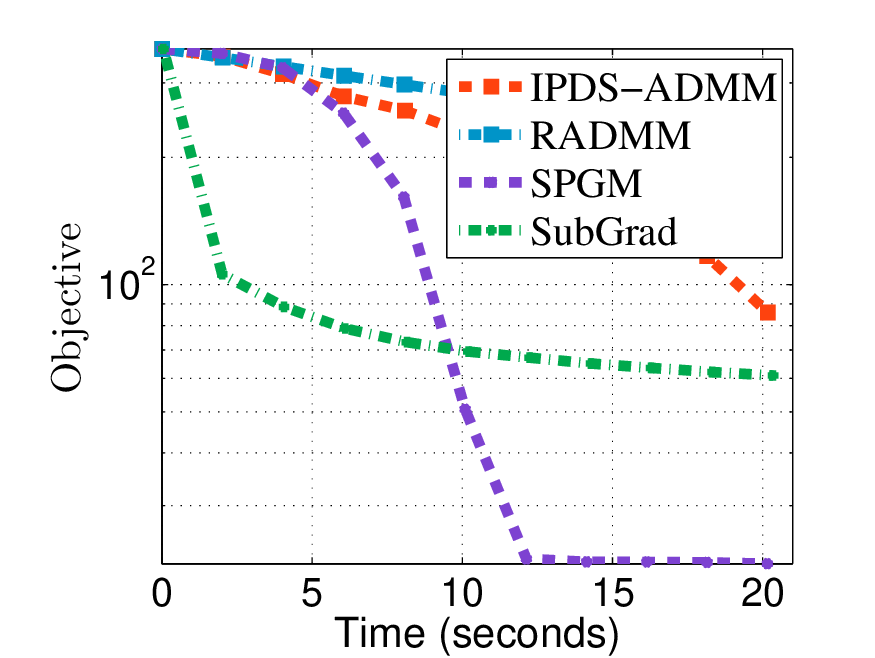}\caption{\scriptsize TDT2-1500-500}\label{fig:sub1}\end{subfigure}
\begin{subfigure}{.24\textwidth}\centering\includegraphics[width=1.12\linewidth]{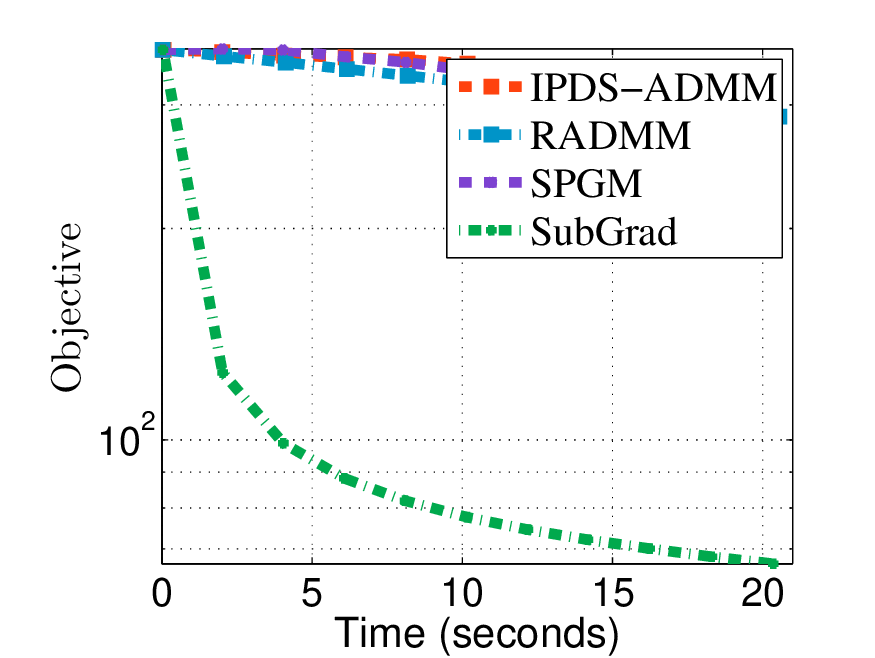}\caption{\scriptsize TDT2-3000-500}\label{fig:sub2}\end{subfigure}
\begin{subfigure}{.24\textwidth}\centering\includegraphics[width=1.12\linewidth]{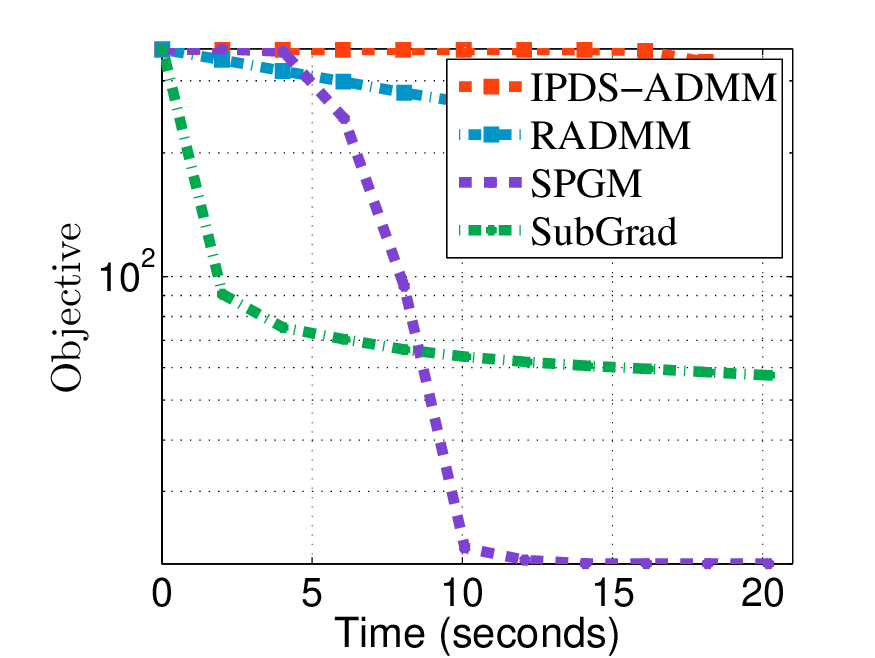}\caption{\scriptsize sector-1500-500}\label{fig:sub3}\end{subfigure}
\begin{subfigure}{.24\textwidth}\centering\includegraphics[width=1.12\linewidth]{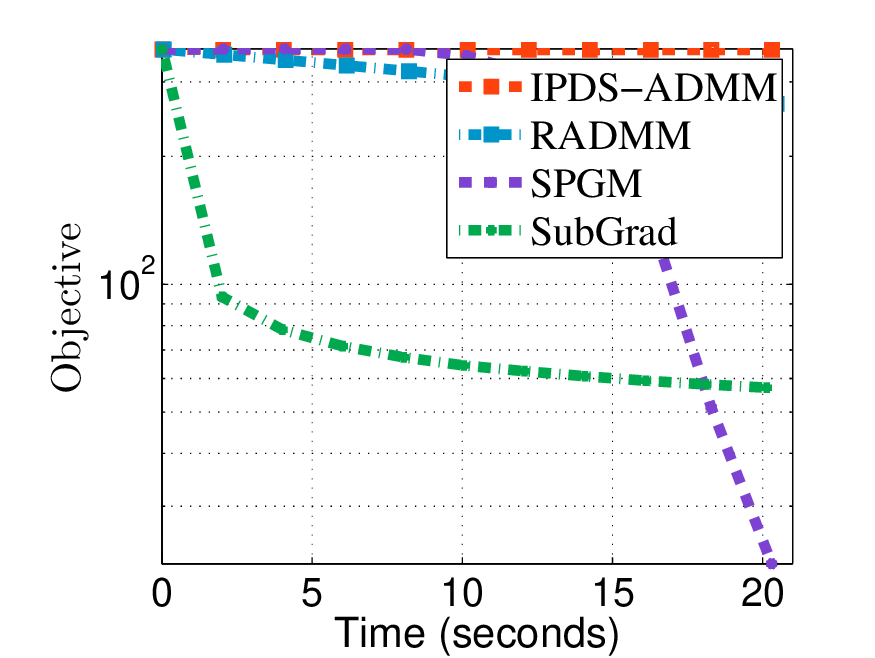}\caption{\scriptsize sector-2500-500}\label{fig:sub4}\end{subfigure}

\caption{Convergence curves of methods for sparse PCA with $\dot{\rho}=1$ and $\beta^0=10\dot{\rho}$.} \label{fig:1}

\end{figure}

\begin{figure}[!t]

\centering
\begin{subfigure}{.24\textwidth}\centering\includegraphics[width=1.12\linewidth]{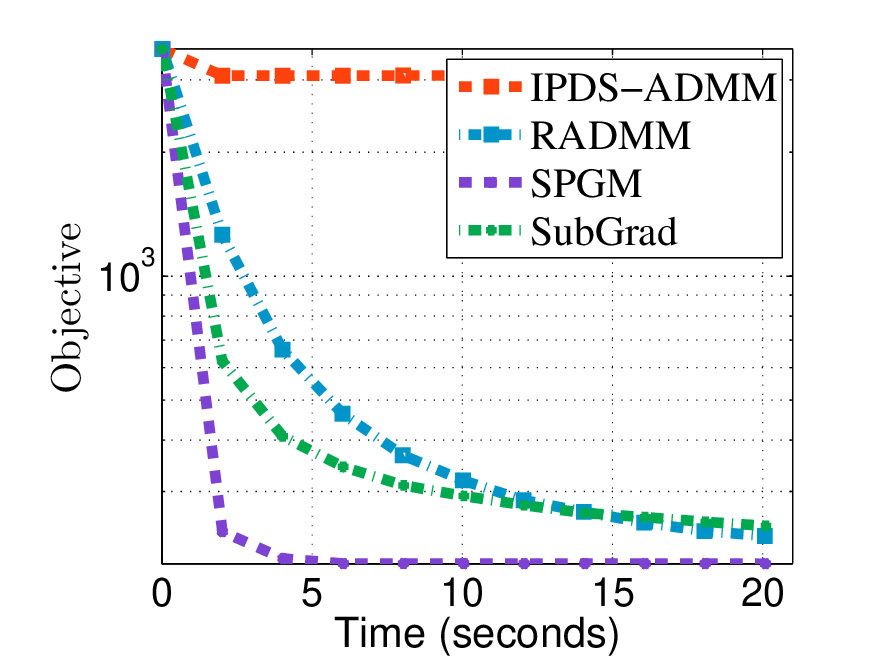}\caption{\scriptsize randn-1500-500}\label{fig:sub1}\end{subfigure}
\begin{subfigure}{.24\textwidth}\centering\includegraphics[width=1.12\linewidth]{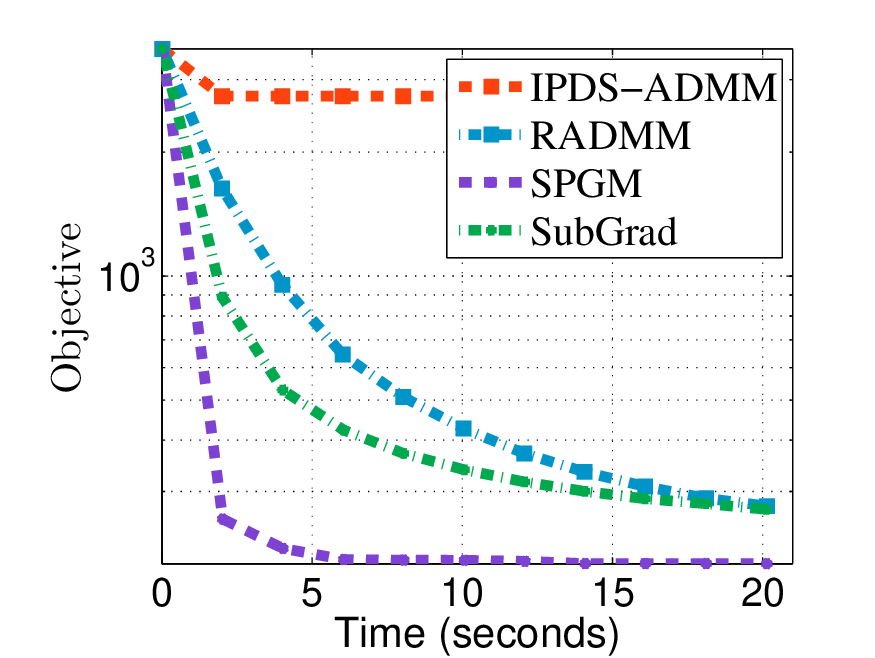}\caption{\scriptsize randn-2500-500}\label{fig:sub2}\end{subfigure}
\begin{subfigure}{.24\textwidth}\centering\includegraphics[width=1.12\linewidth]{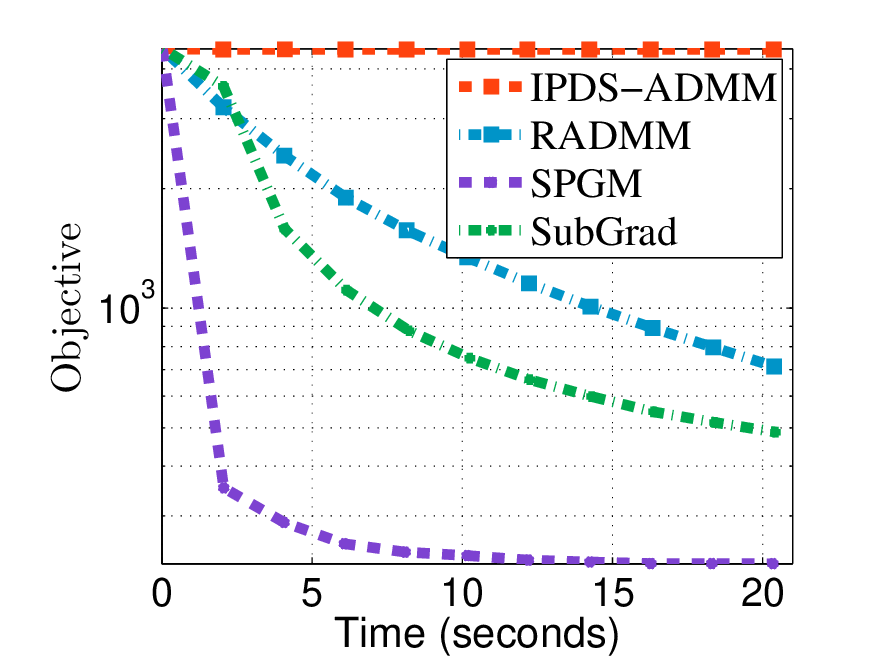}\caption{\scriptsize mnist-1500-780}\label{fig:sub3}\end{subfigure}
\begin{subfigure}{.24\textwidth}\centering\includegraphics[width=1.12\linewidth]{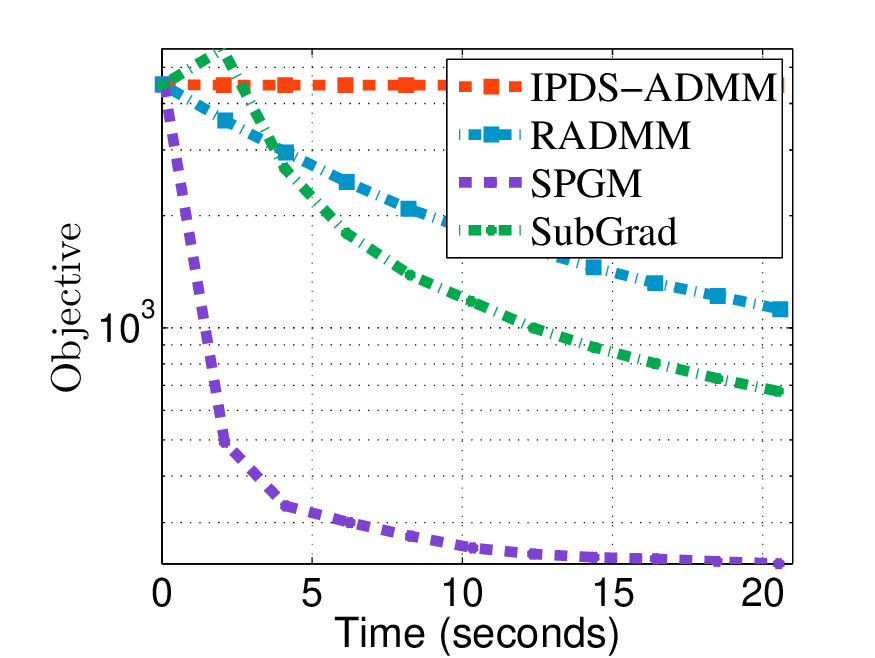}\caption{\scriptsize mnist-2500-780}\label{fig:sub4}\end{subfigure}

\centering
\begin{subfigure}{.24\textwidth}\centering\includegraphics[width=1.12\linewidth]{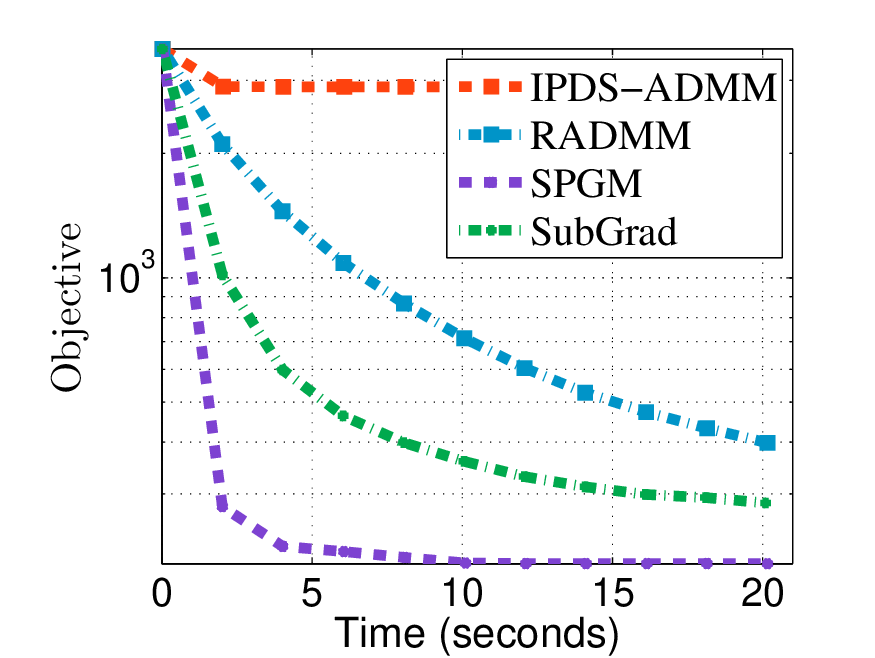}\caption{\scriptsize TDT2-1500-500}\label{fig:sub1}\end{subfigure}
\begin{subfigure}{.24\textwidth}\centering\includegraphics[width=1.12\linewidth]{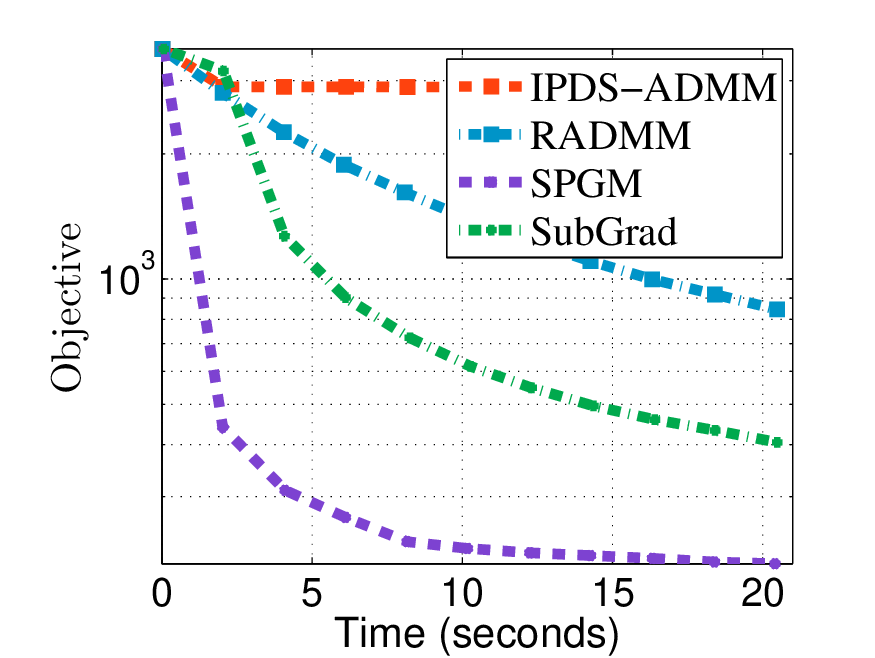}\caption{\scriptsize TDT2-3000-500}\label{fig:sub2}\end{subfigure}
\begin{subfigure}{.24\textwidth}\centering\includegraphics[width=1.12\linewidth]{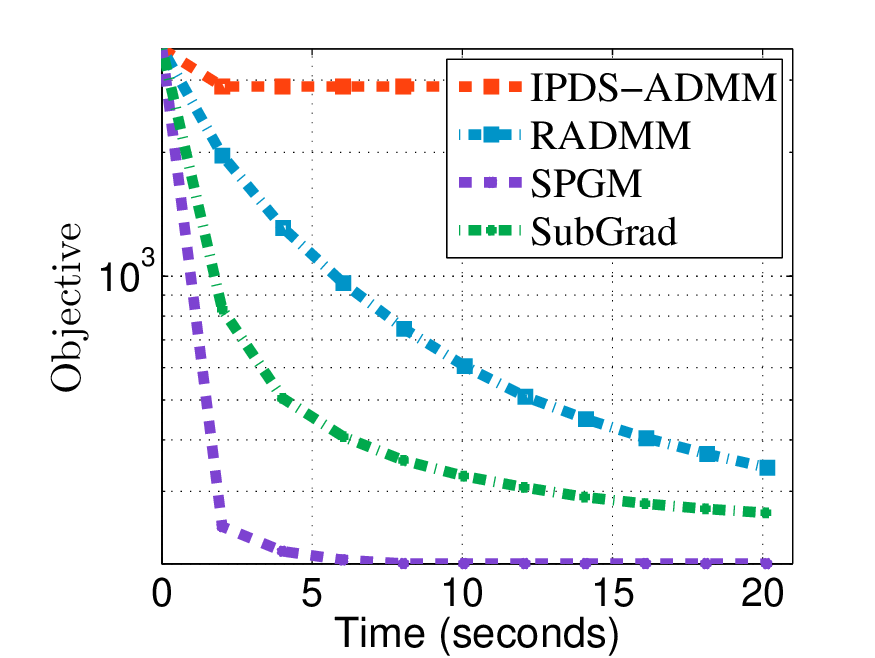}\caption{\scriptsize sector-1500-500}\label{fig:sub3}\end{subfigure}
\begin{subfigure}{.24\textwidth}\centering\includegraphics[width=1.12\linewidth]{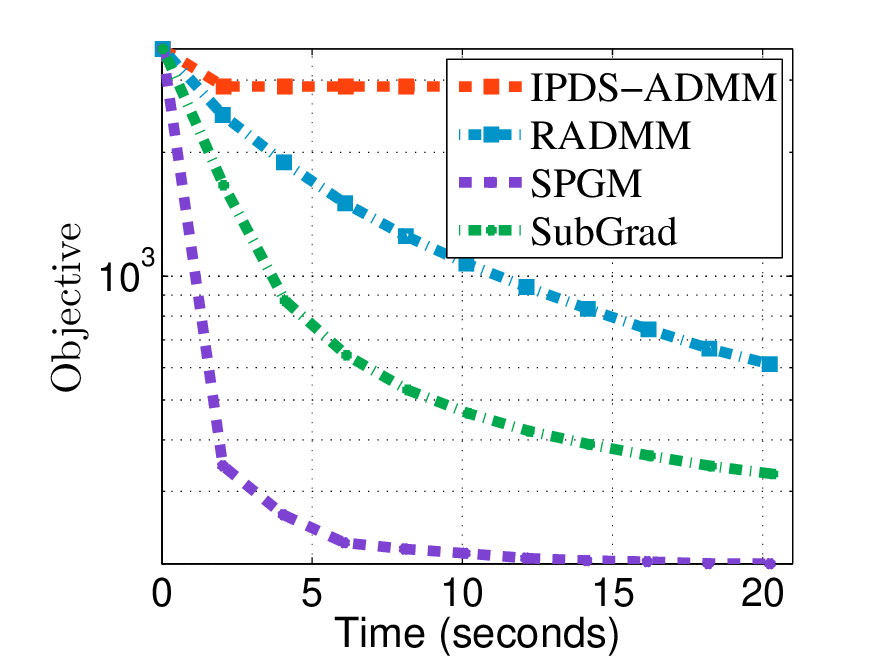}\caption{\scriptsize sector-2500-500}\label{fig:sub4}\end{subfigure}

\caption{Convergence curves of methods for sparse PCA with $\dot{\rho}=10$ and $\beta^0=10\dot{\rho}$.} \label{fig:2}

\end{figure}

\begin{figure}[!t]

\centering
\begin{subfigure}{.24\textwidth}\centering\includegraphics[width=1.12\linewidth]{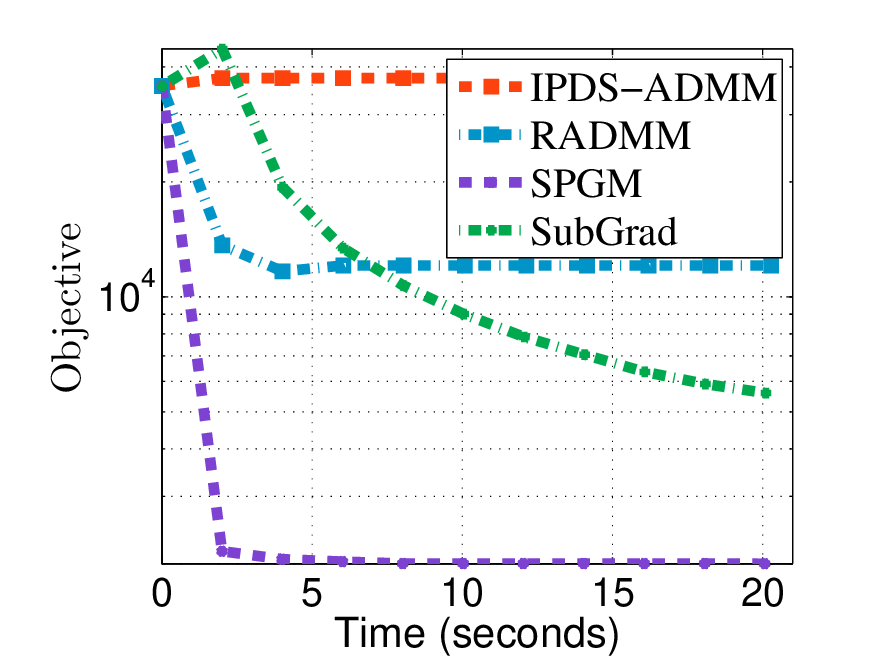}\caption{\scriptsize randn-1500-500}\label{fig:sub1}\end{subfigure}
\begin{subfigure}{.24\textwidth}\centering\includegraphics[width=1.12\linewidth]{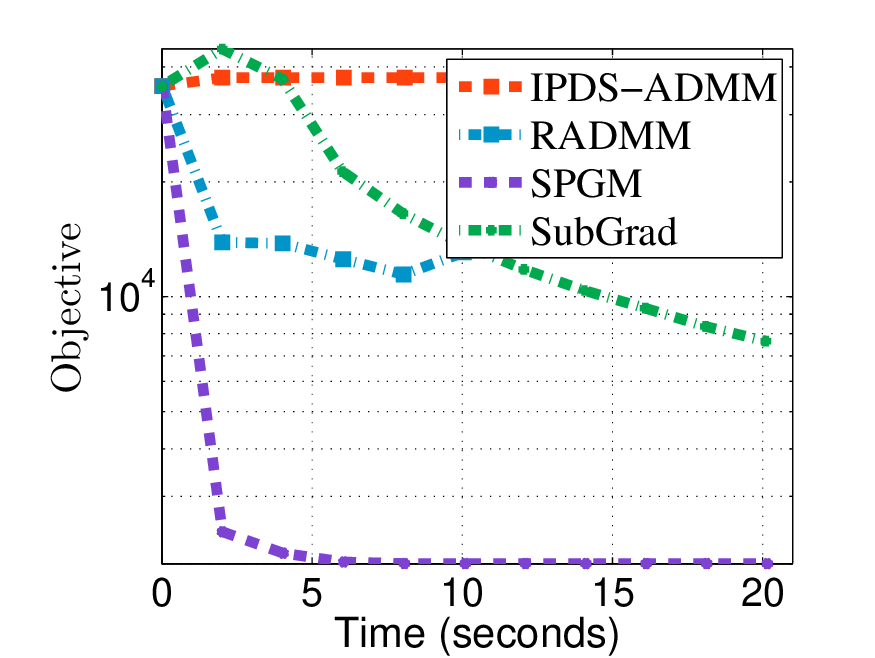}\caption{\scriptsize randn-2500-500}\label{fig:sub2}\end{subfigure}
\begin{subfigure}{.24\textwidth}\centering\includegraphics[width=1.12\linewidth]{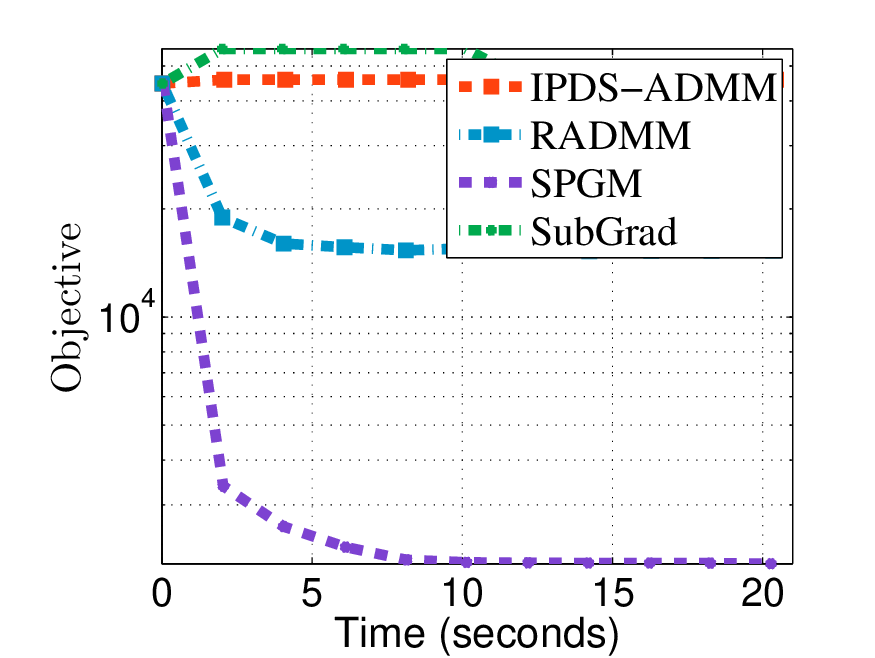}\caption{\scriptsize mnist-1500-780}\label{fig:sub3}\end{subfigure}
\begin{subfigure}{.24\textwidth}\centering\includegraphics[width=1.12\linewidth]{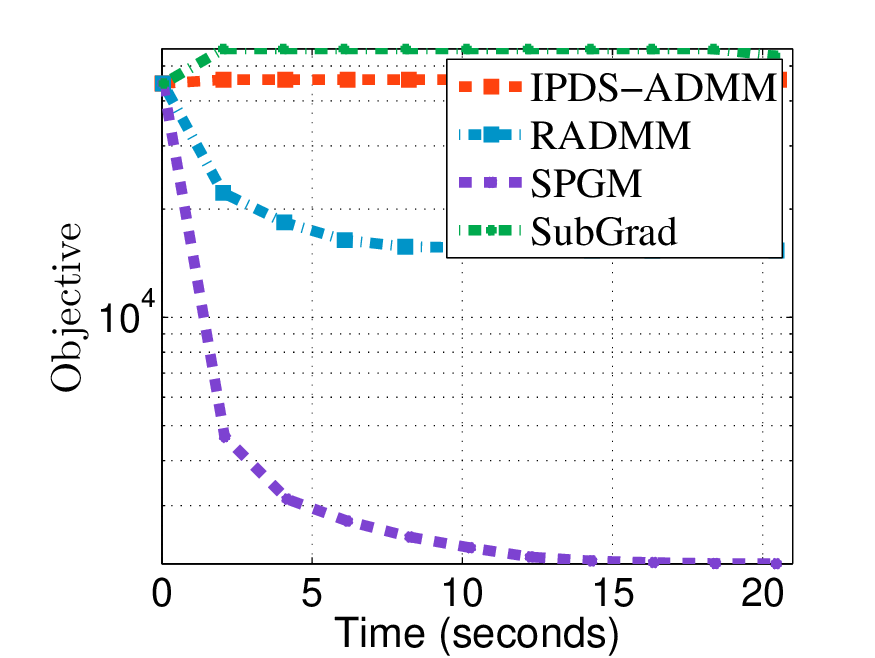}\caption{\scriptsize mnist-2500-780}\label{fig:sub4}\end{subfigure}

\centering
\begin{subfigure}{.24\textwidth}\centering\includegraphics[width=1.12\linewidth]{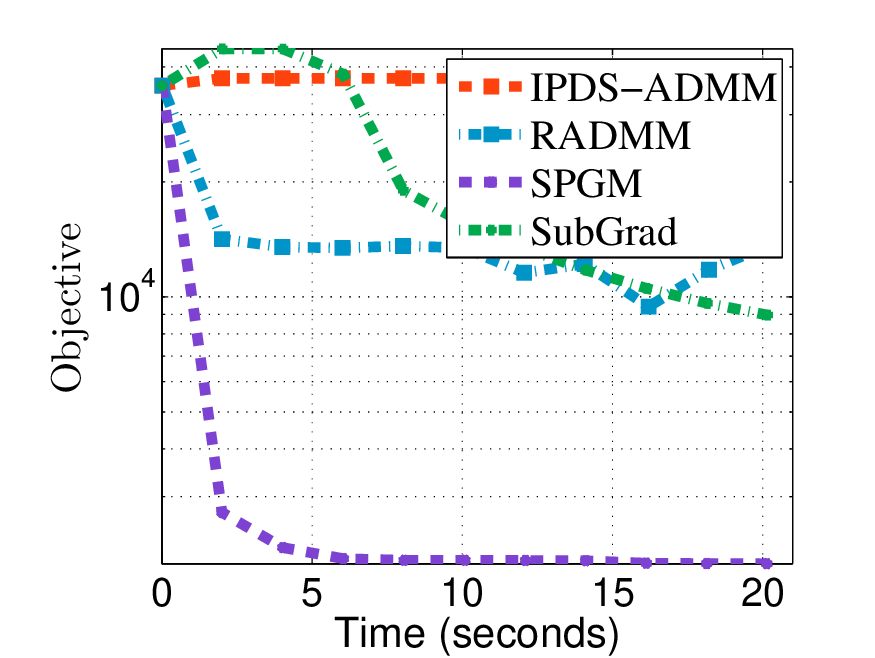}\caption{\scriptsize TDT2-1500-500}\label{fig:sub1}\end{subfigure}
\begin{subfigure}{.24\textwidth}\centering\includegraphics[width=1.12\linewidth]{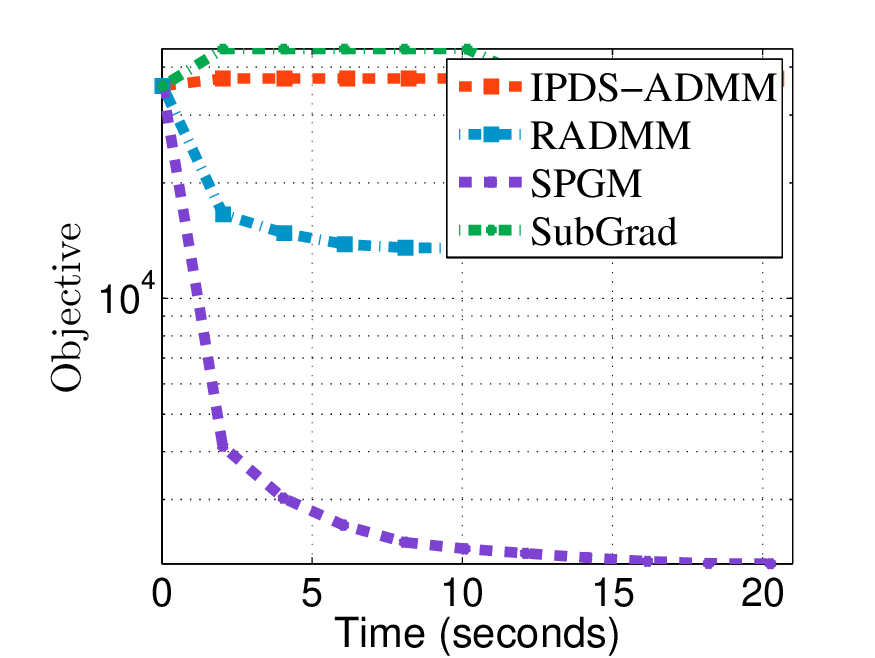}\caption{\scriptsize TDT2-3000-500}\label{fig:sub2}\end{subfigure}
\begin{subfigure}{.24\textwidth}\centering\includegraphics[width=1.12\linewidth]{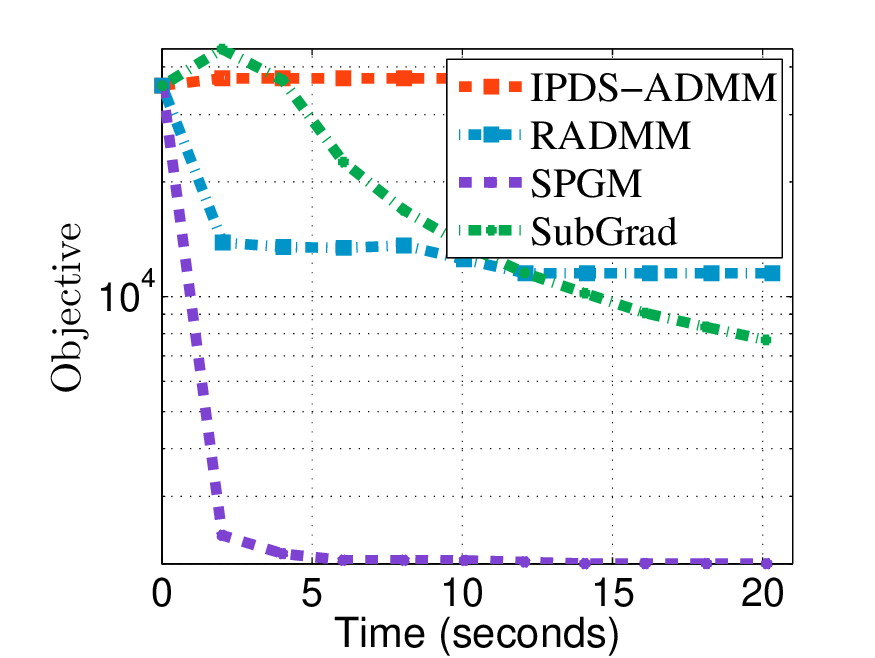}\caption{\scriptsize sector-1500-500}\label{fig:sub3}\end{subfigure}
\begin{subfigure}{.24\textwidth}\centering\includegraphics[width=1.12\linewidth]{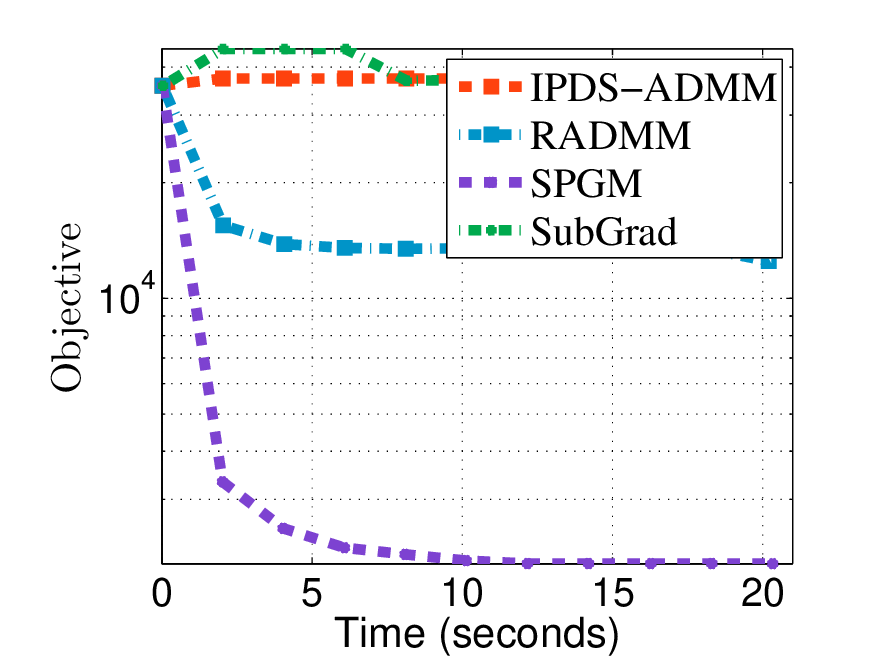}\caption{\scriptsize sector-2500-500}\label{fig:sub4}\end{subfigure}

\caption{Convergence curves of methods for sparse PCA with $\dot{\rho}=100$ and $\beta^0=10\dot{\rho}$.} \label{fig:3}

\end{figure}

\begin{figure}[!t]

\centering
\begin{subfigure}{.24\textwidth}\centering\includegraphics[width=1.12\linewidth]{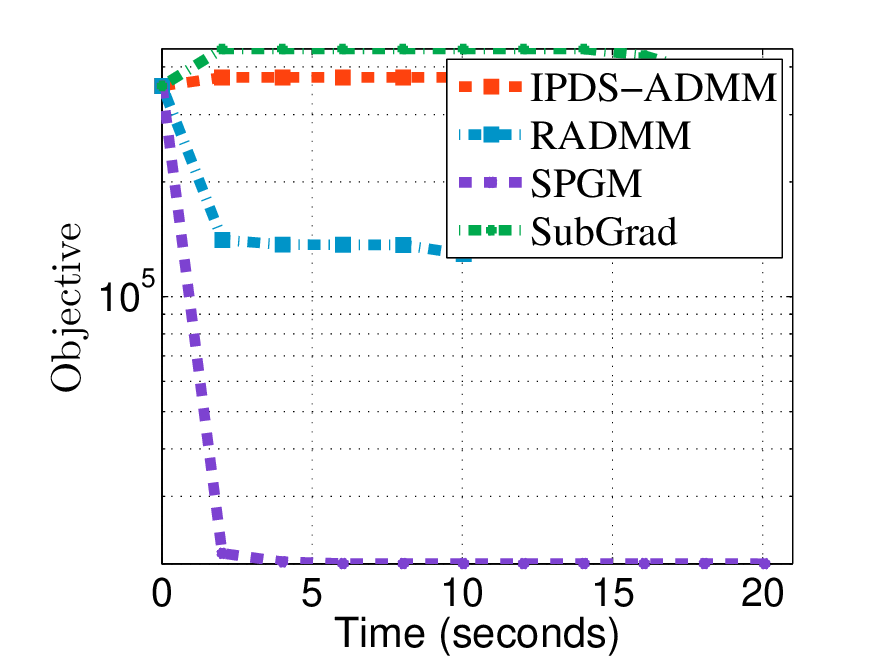}\caption{\scriptsize randn-1500-500}\label{fig:sub1}\end{subfigure}
\begin{subfigure}{.24\textwidth}\centering\includegraphics[width=1.12\linewidth]{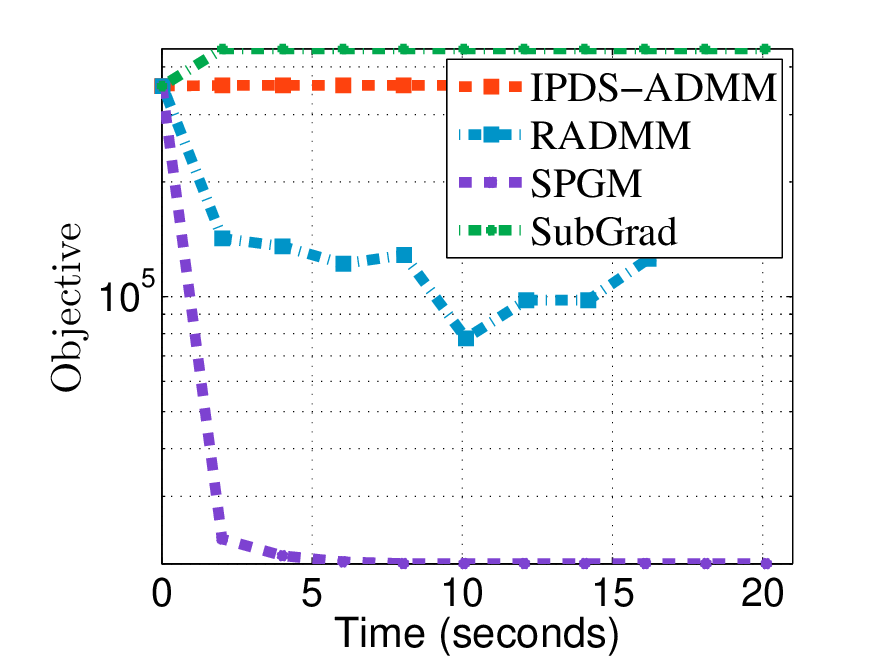}\caption{\scriptsize randn-2500-500}\label{fig:sub2}\end{subfigure}
\begin{subfigure}{.24\textwidth}\centering\includegraphics[width=1.12\linewidth]{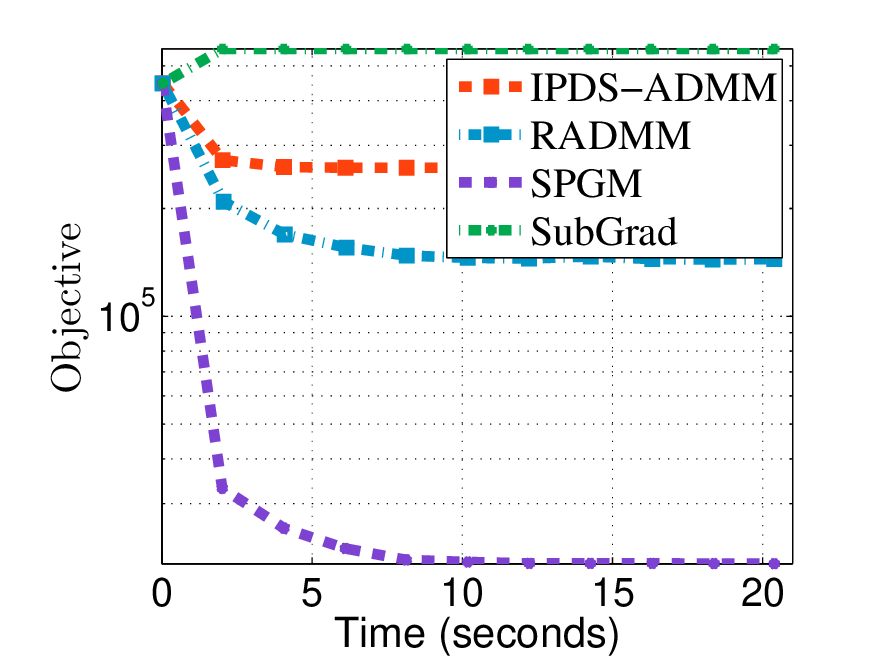}\caption{\scriptsize mnist-1500-780}\label{fig:sub3}\end{subfigure}
\begin{subfigure}{.24\textwidth}\centering\includegraphics[width=1.12\linewidth]{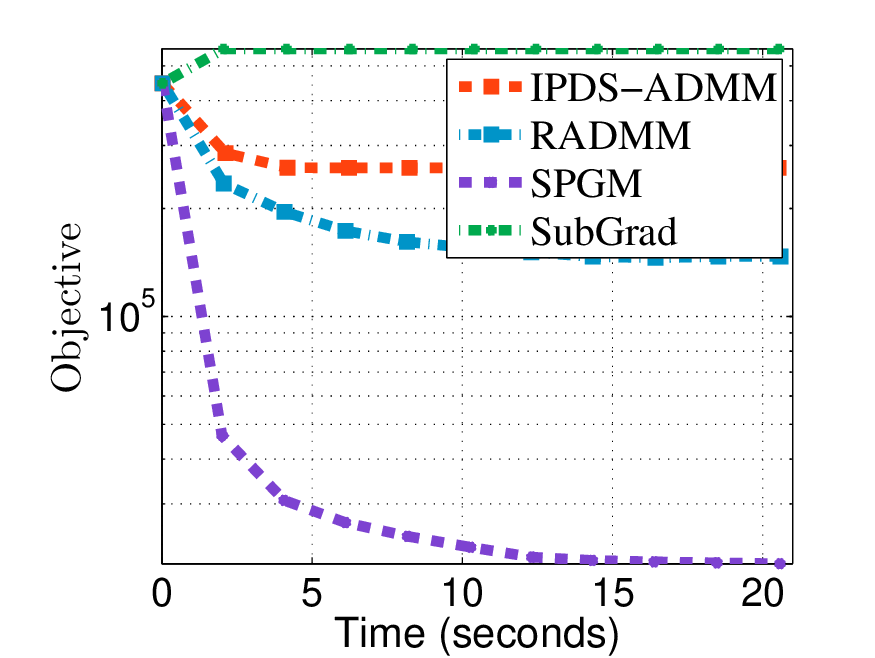}\caption{\scriptsize mnist-2500-780}\label{fig:sub4}\end{subfigure}

\centering
\begin{subfigure}{.24\textwidth}\centering\includegraphics[width=1.12\linewidth]{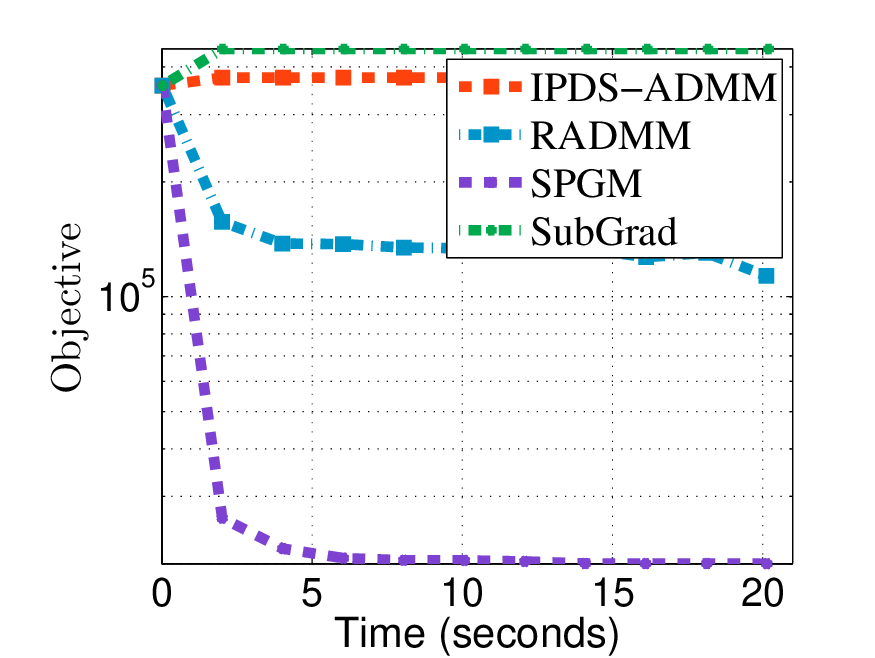}\caption{\scriptsize TDT2-1500-500}\label{fig:sub1}\end{subfigure}
\begin{subfigure}{.24\textwidth}\centering\includegraphics[width=1.12\linewidth]{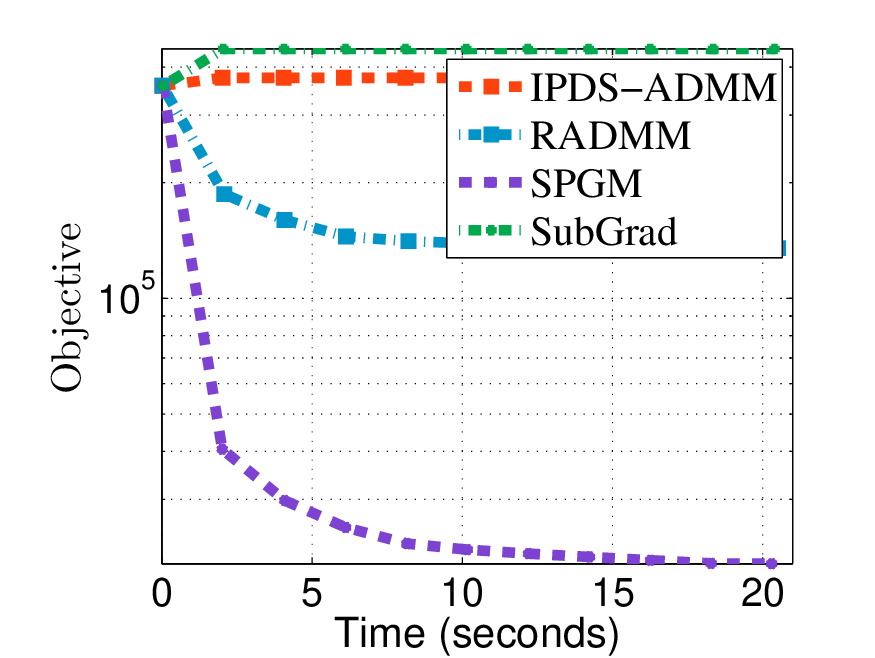}\caption{\scriptsize TDT2-3000-500}\label{fig:sub2}\end{subfigure}
\begin{subfigure}{.24\textwidth}\centering\includegraphics[width=1.12\linewidth]{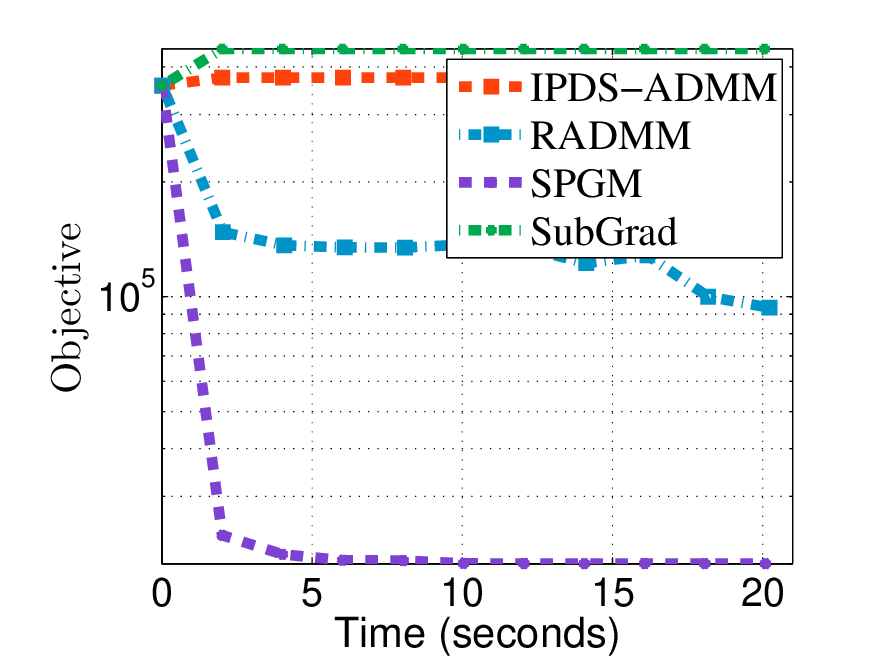}\caption{\scriptsize sector-1500-500}\label{fig:sub3}\end{subfigure}
\begin{subfigure}{.24\textwidth}\centering\includegraphics[width=1.12\linewidth]{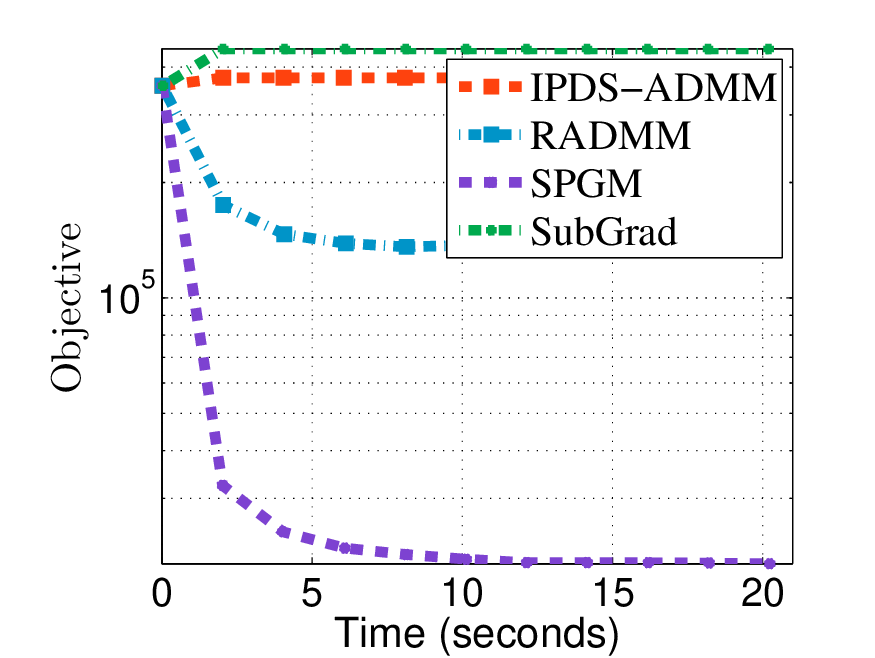}\caption{\scriptsize sector-2500-500}\label{fig:sub4}\end{subfigure}

\caption{Convergence curves of methods for sparse PCA with $\dot{\rho}=1000$ and $\beta^0=10\dot{\rho}$.} \label{fig:4}
\end{figure}

\begin{figure}[!t]

\centering
\begin{subfigure}{.24\textwidth}\centering\includegraphics[width=1.12\linewidth]{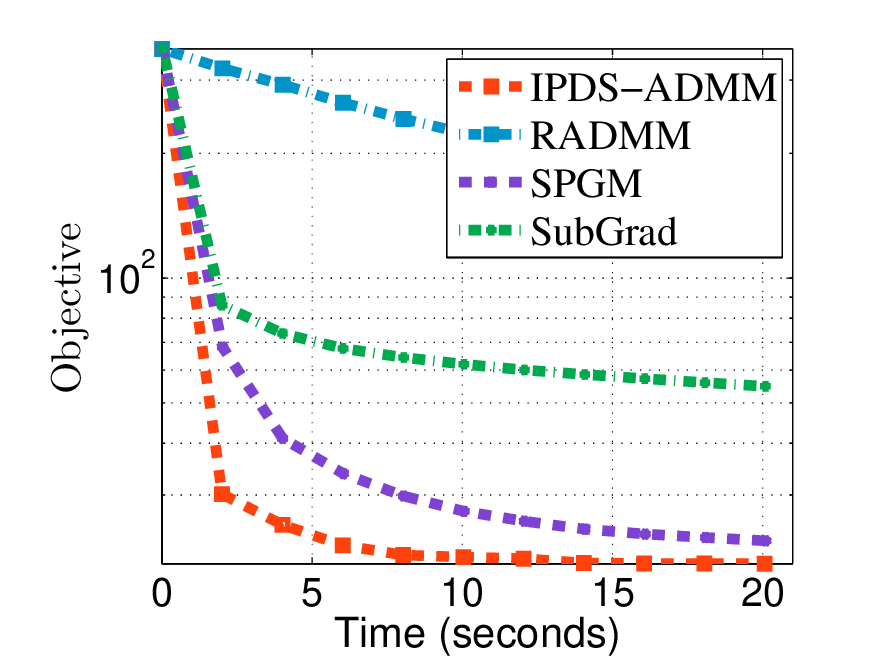}\caption{\scriptsize randn-1500-500}\label{fig:sub1}\end{subfigure}
\begin{subfigure}{.24\textwidth}\centering\includegraphics[width=1.12\linewidth]{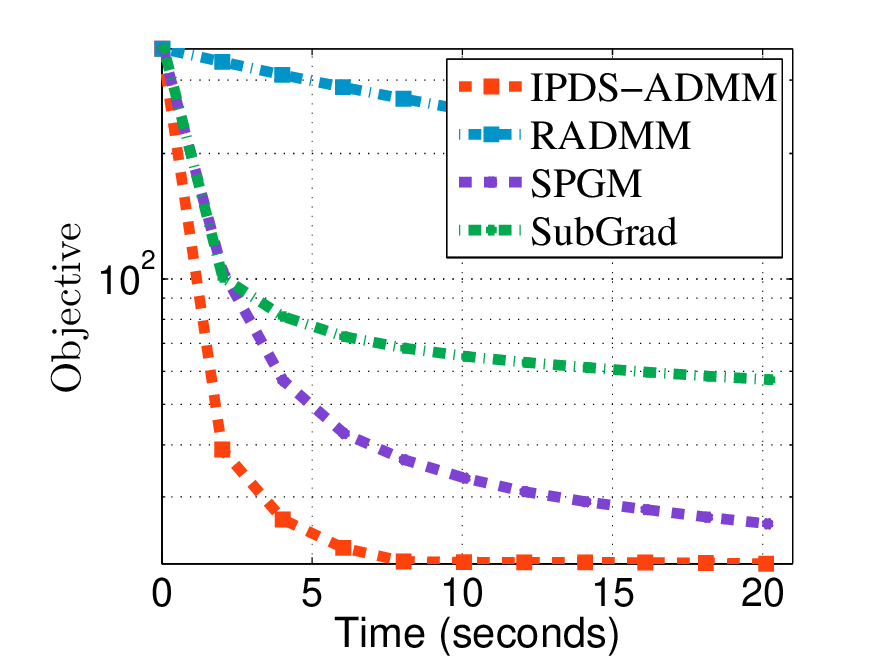}\caption{\scriptsize randn-2500-500}\label{fig:sub2}\end{subfigure}
\begin{subfigure}{.24\textwidth}\centering\includegraphics[width=1.12\linewidth]{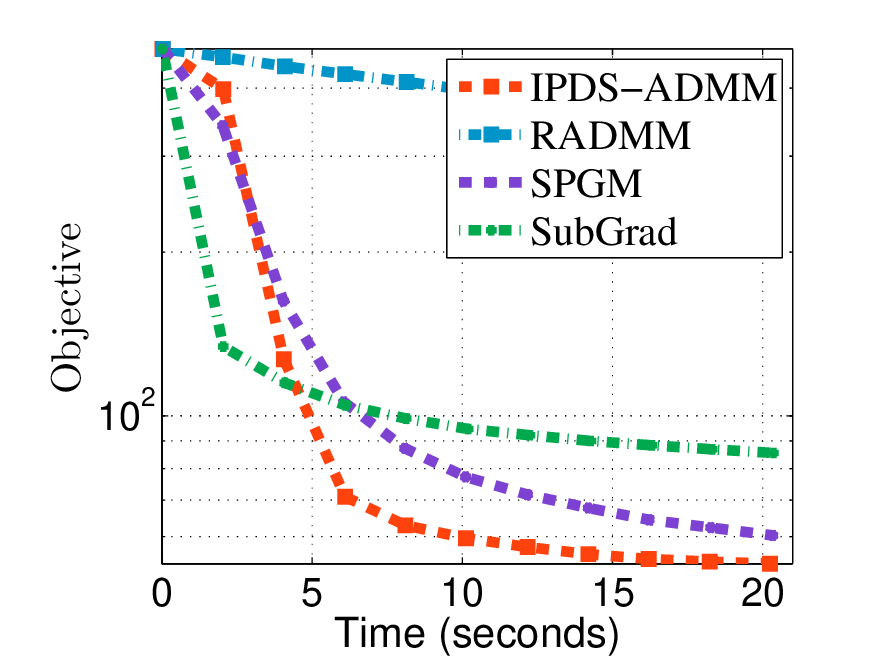}\caption{\scriptsize mnist-1500-780}\label{fig:sub3}\end{subfigure}
\begin{subfigure}{.24\textwidth}\centering\includegraphics[width=1.12\linewidth]{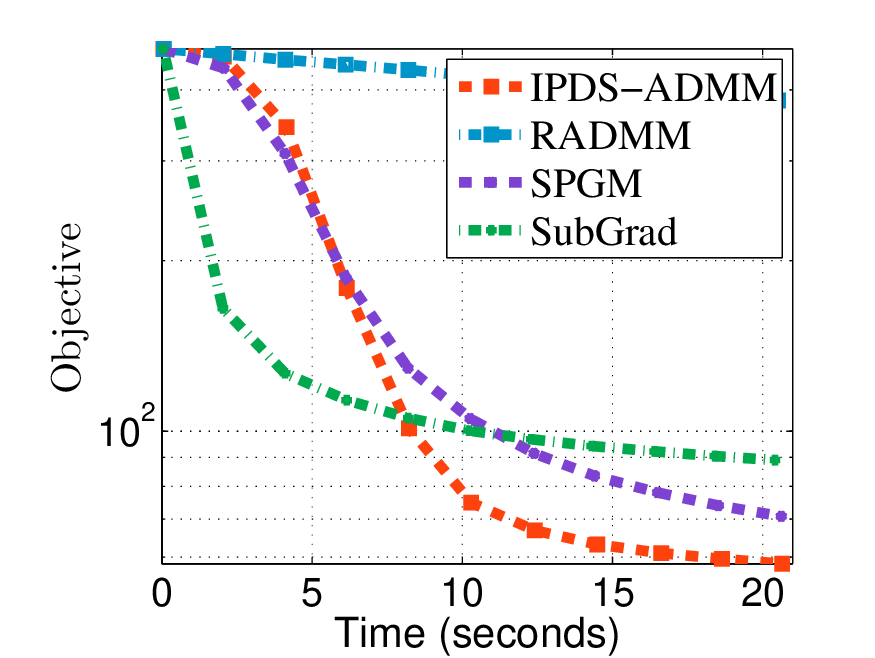}\caption{\scriptsize mnist-2500-780}\label{fig:sub4}\end{subfigure}

\centering
\begin{subfigure}{.24\textwidth}\centering\includegraphics[width=1.12\linewidth]{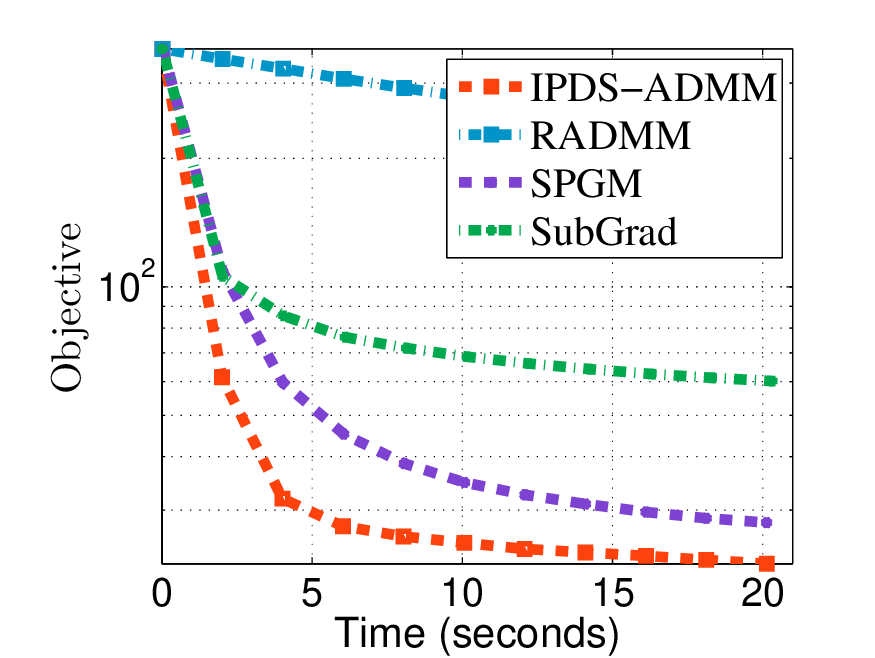}\caption{\scriptsize TDT2-1500-500}\label{fig:sub1}\end{subfigure}
\begin{subfigure}{.24\textwidth}\centering\includegraphics[width=1.12\linewidth]{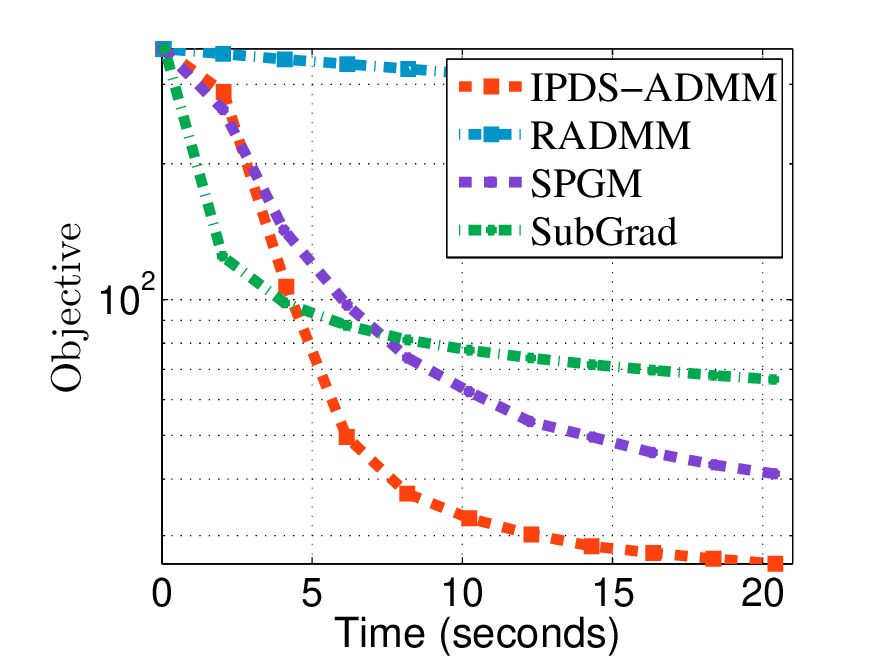}\caption{\scriptsize TDT2-3000-500}\label{fig:sub2}\end{subfigure}
\begin{subfigure}{.24\textwidth}\centering\includegraphics[width=1.12\linewidth]{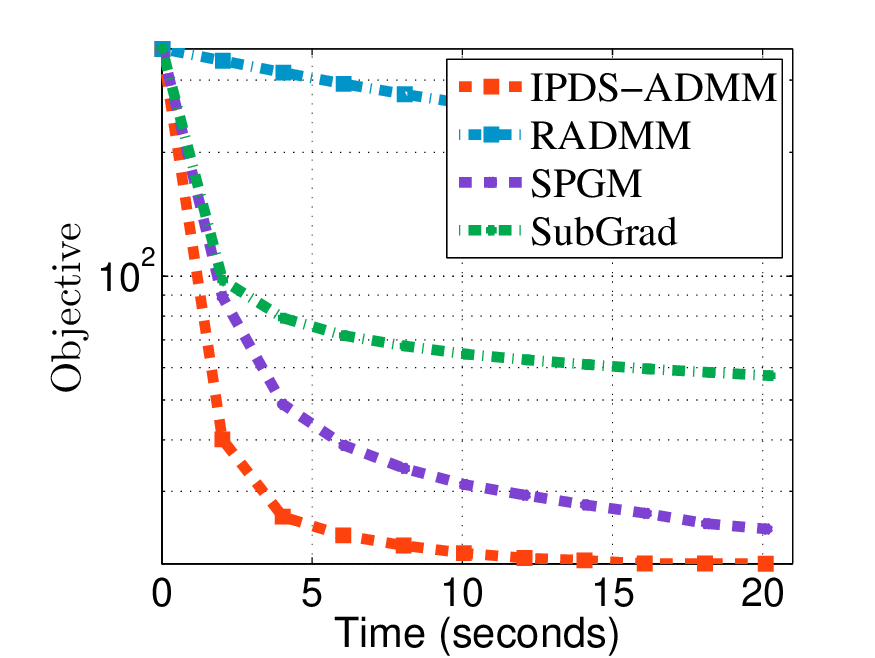}\caption{\scriptsize sector-1500-500}\label{fig:sub3}\end{subfigure}
\begin{subfigure}{.24\textwidth}\centering\includegraphics[width=1.12\linewidth]{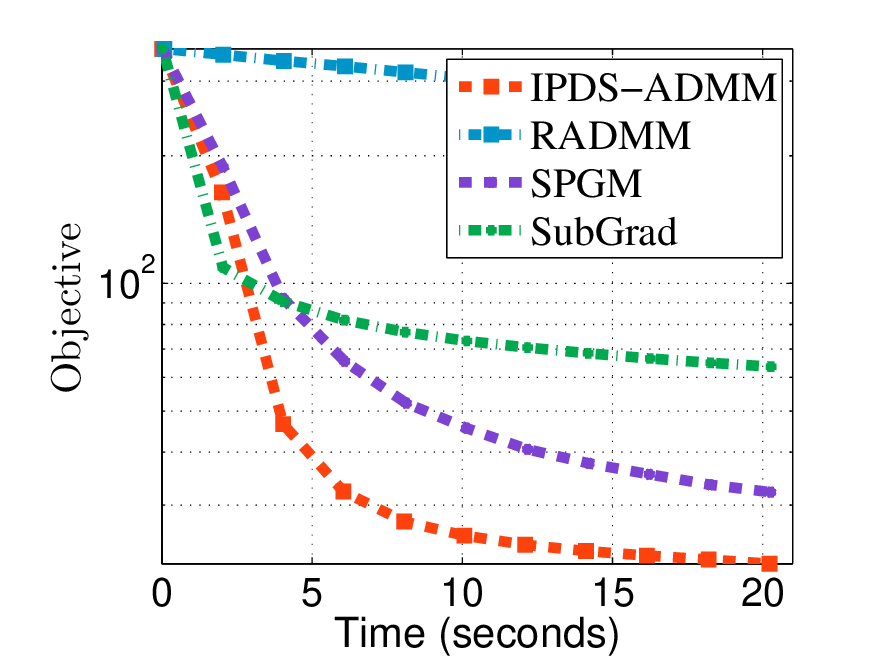}\caption{\scriptsize sector-2500-500}\label{fig:sub4}\end{subfigure}

\caption{Convergence curves of methods for sparse PCA with $\dot{\rho}=1$ and $\beta^0=50\dot{\rho}$.} \label{fig:5} 

\end{figure}

\begin{figure}[!t]
\centering
\begin{subfigure}{.24\textwidth}\centering\includegraphics[width=1.12\linewidth]{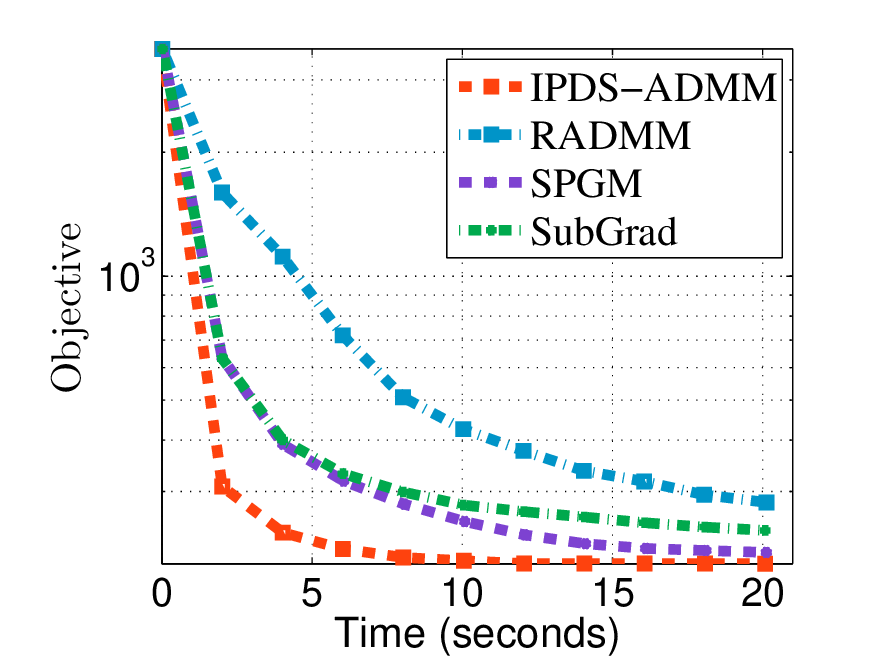}\caption{\scriptsize randn-1500-500}\label{fig:sub1}\end{subfigure}
\begin{subfigure}{.24\textwidth}\centering\includegraphics[width=1.12\linewidth]{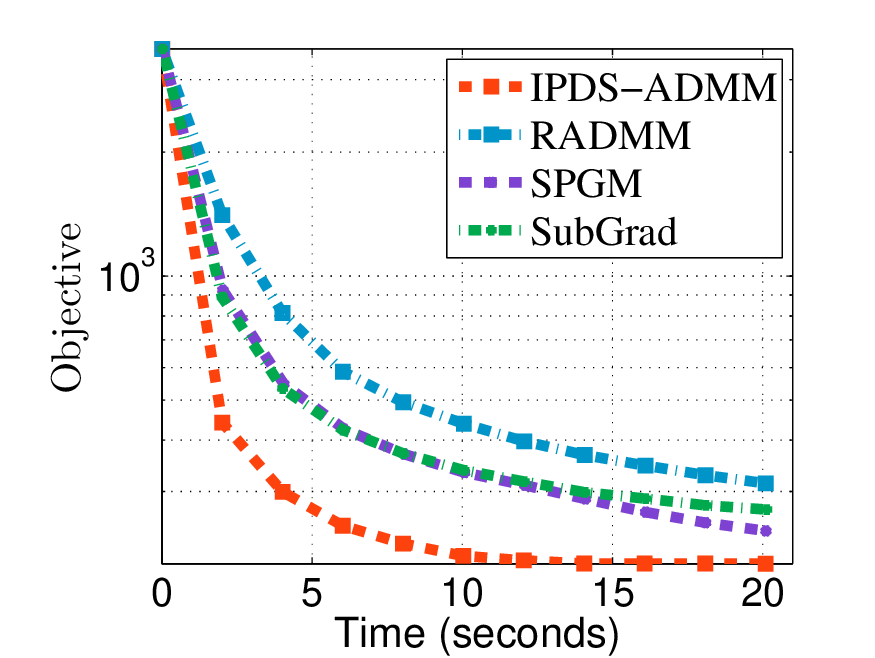}\caption{\scriptsize randn-2500-500}\label{fig:sub2}\end{subfigure}
\begin{subfigure}{.24\textwidth}\centering\includegraphics[width=1.12\linewidth]{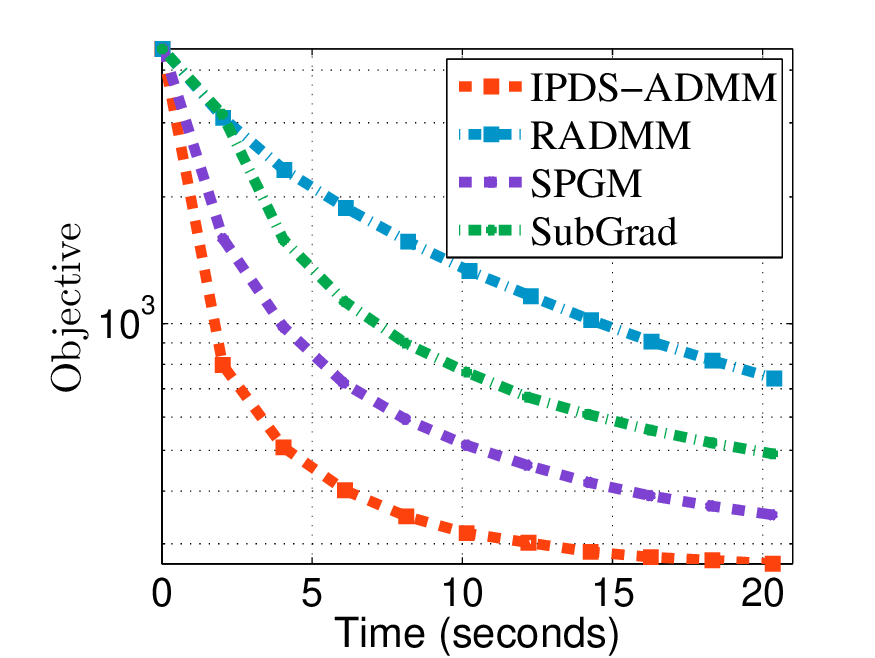}\caption{\scriptsize mnist-1500-780}\label{fig:sub3}\end{subfigure}
\begin{subfigure}{.24\textwidth}\centering\includegraphics[width=1.12\linewidth]{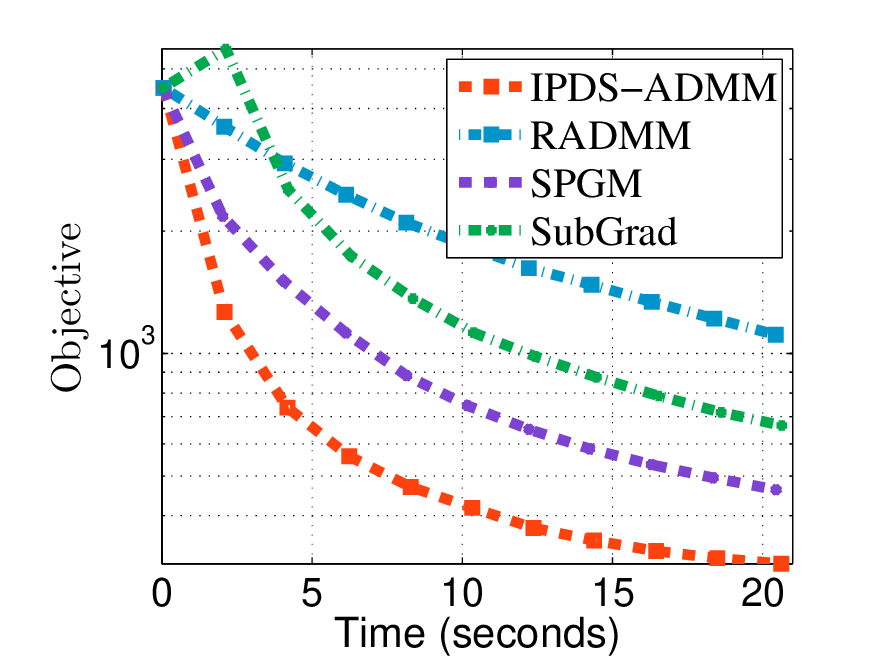}\caption{\scriptsize mnist-2500-780}\label{fig:sub4}\end{subfigure}

\centering
\begin{subfigure}{.24\textwidth}\centering\includegraphics[width=1.12\linewidth]{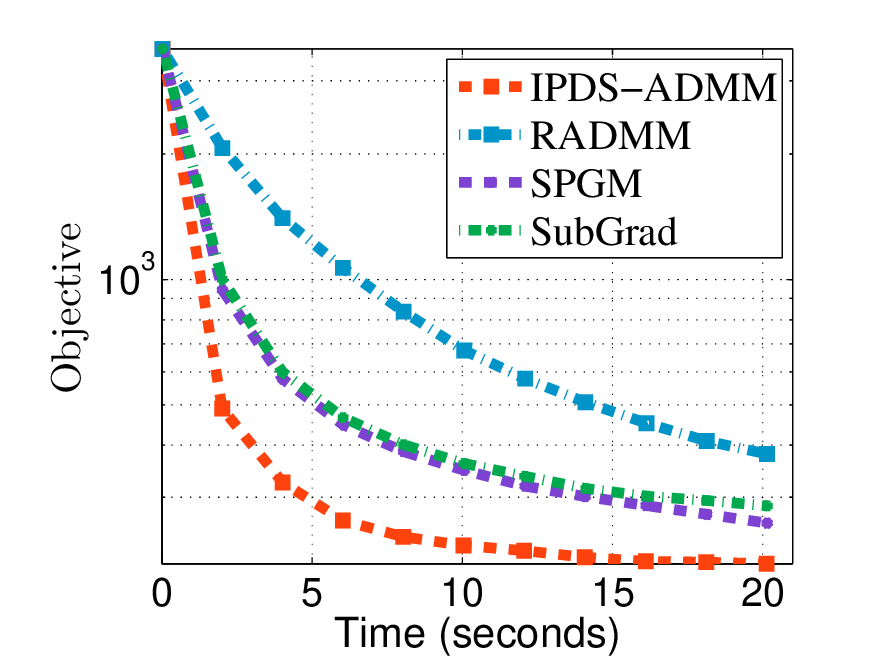}\caption{\scriptsize TDT2-1500-500}\label{fig:sub1}\end{subfigure}
\begin{subfigure}{.24\textwidth}\centering\includegraphics[width=1.12\linewidth]{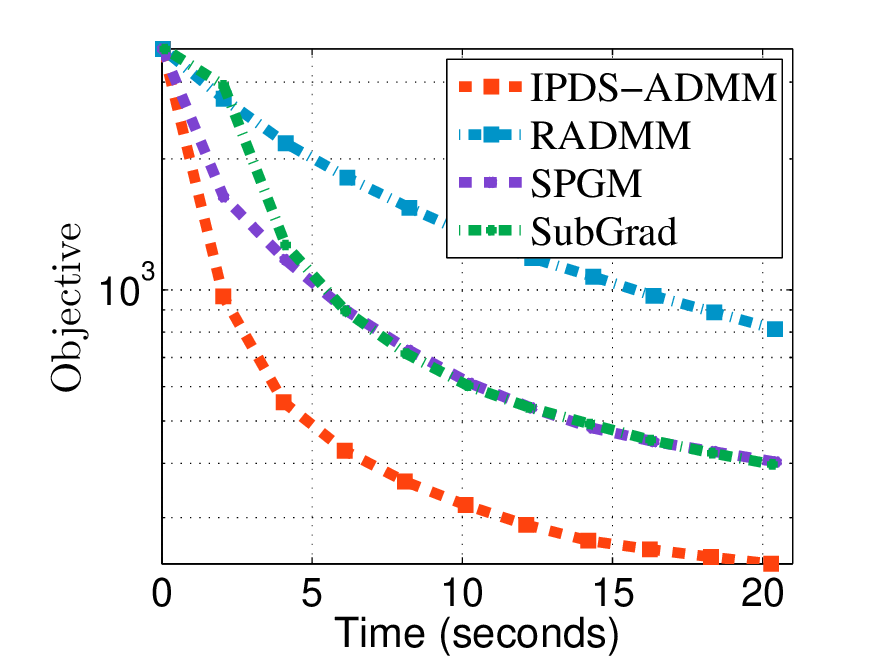}\caption{\scriptsize TDT2-3000-500}\label{fig:sub2}\end{subfigure}
\begin{subfigure}{.24\textwidth}\centering\includegraphics[width=1.12\linewidth]{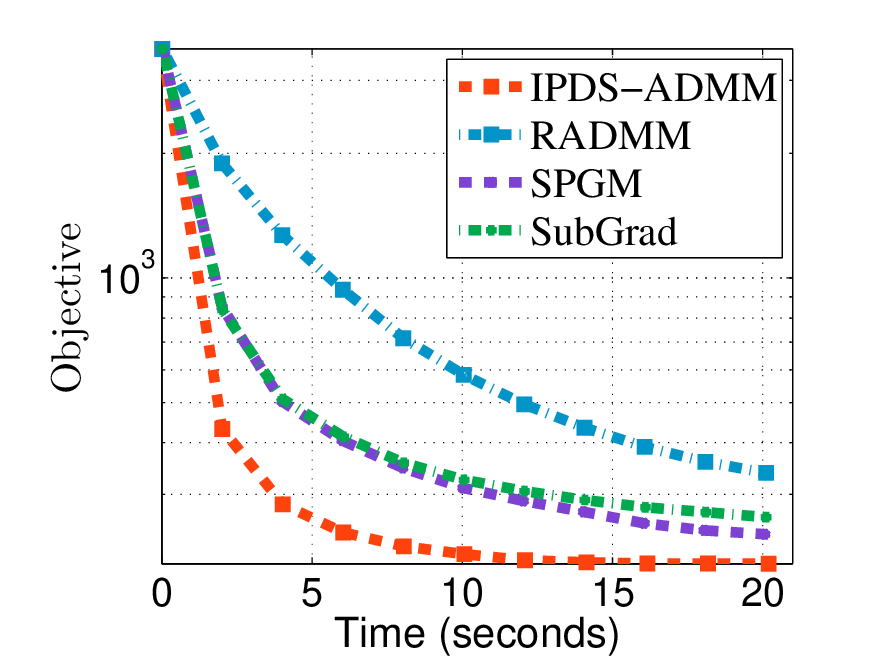}\caption{\scriptsize sector-1500-500}\label{fig:sub3}\end{subfigure}
\begin{subfigure}{.24\textwidth}\centering\includegraphics[width=1.12\linewidth]{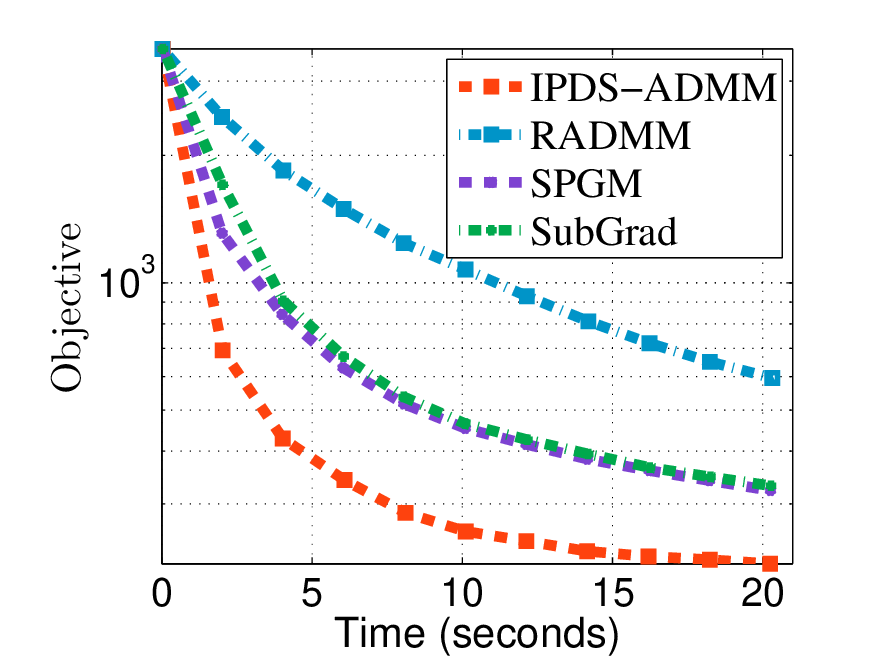}\caption{\scriptsize sector-2500-500}\label{fig:sub4}\end{subfigure}

\caption{Convergence curves of methods for sparse PCA with $\dot{\rho}=10$ and $\beta^0=50\dot{\rho}$.} \label{fig:6}

\end{figure}

\begin{figure}[!t]

\centering
\begin{subfigure}{.24\textwidth}\centering\includegraphics[width=1.12\linewidth]{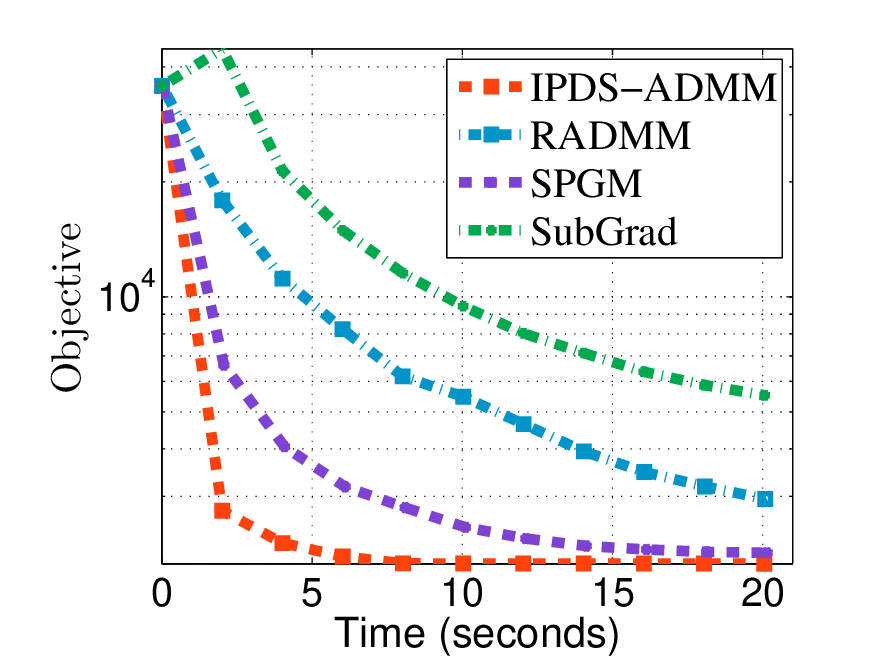}\caption{\scriptsize randn-1500-500}\label{fig:sub1}\end{subfigure}
\begin{subfigure}{.24\textwidth}\centering\includegraphics[width=1.12\linewidth]{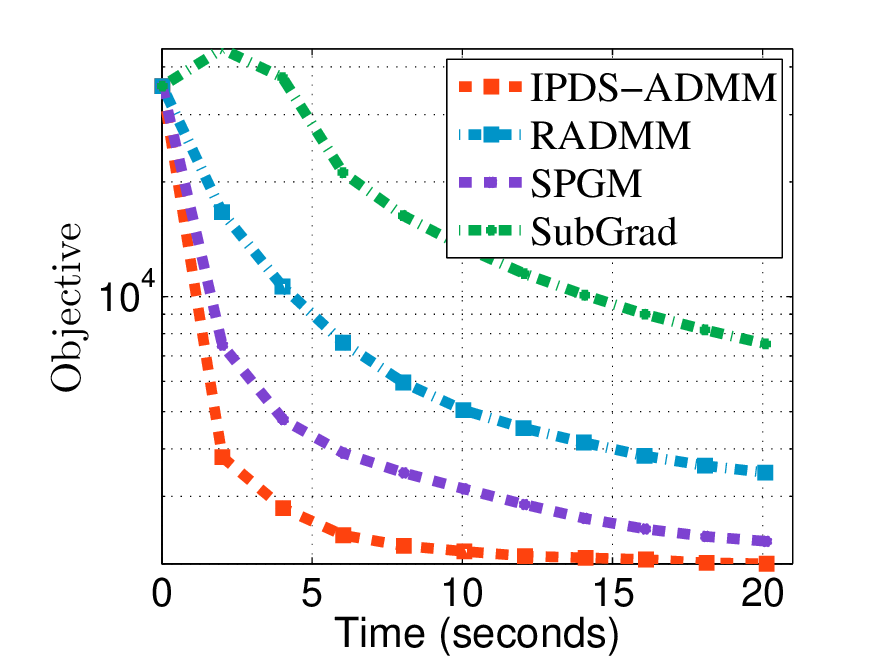}\caption{\scriptsize randn-2500-500}\label{fig:sub2}\end{subfigure}
\begin{subfigure}{.24\textwidth}\centering\includegraphics[width=1.12\linewidth]{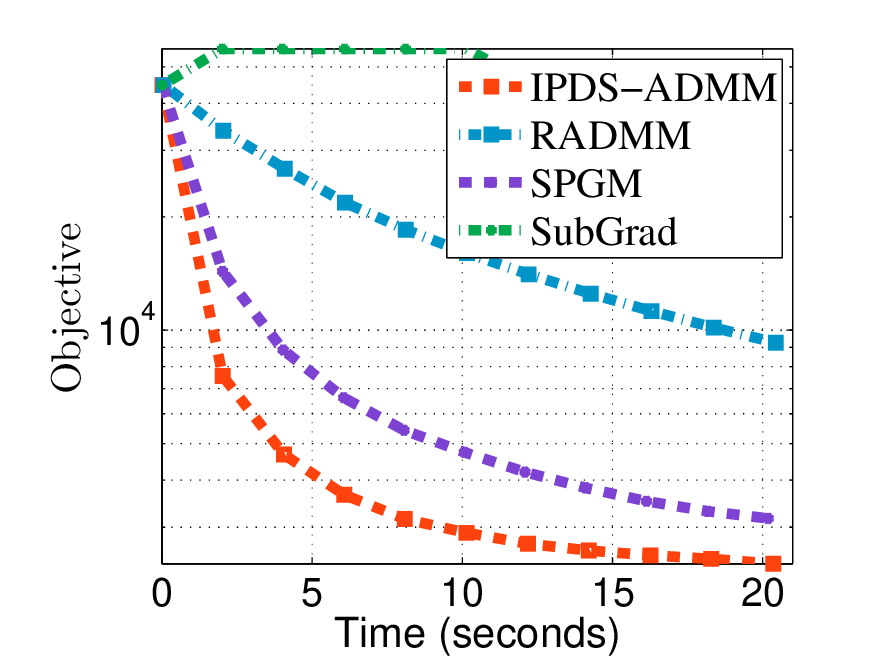}\caption{\scriptsize mnist-1500-780}\label{fig:sub3}\end{subfigure}
\begin{subfigure}{.24\textwidth}\centering\includegraphics[width=1.12\linewidth]{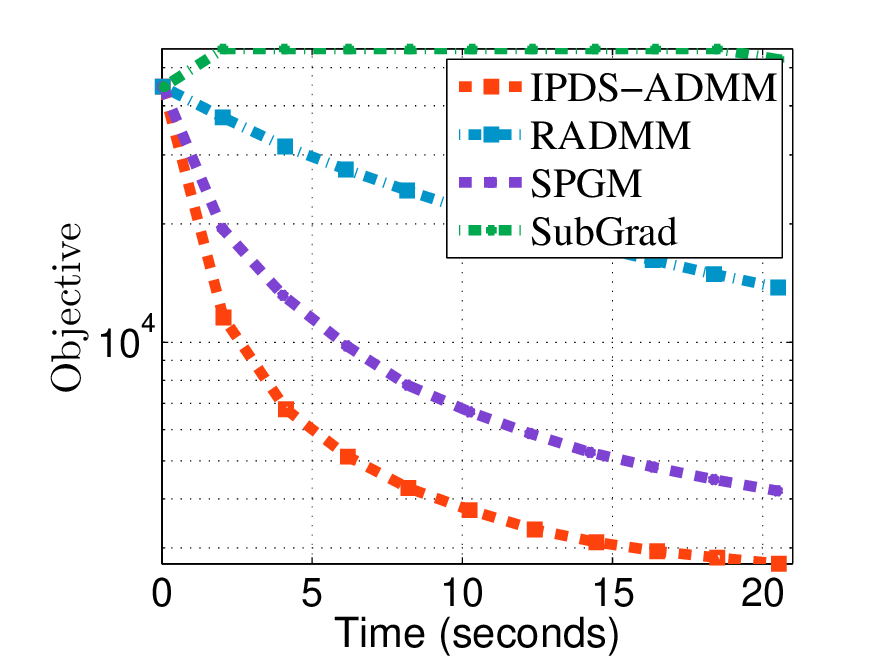}\caption{\scriptsize mnist-2500-780}\label{fig:sub4}\end{subfigure}

\centering
\begin{subfigure}{.24\textwidth}\centering\includegraphics[width=1.12\linewidth]{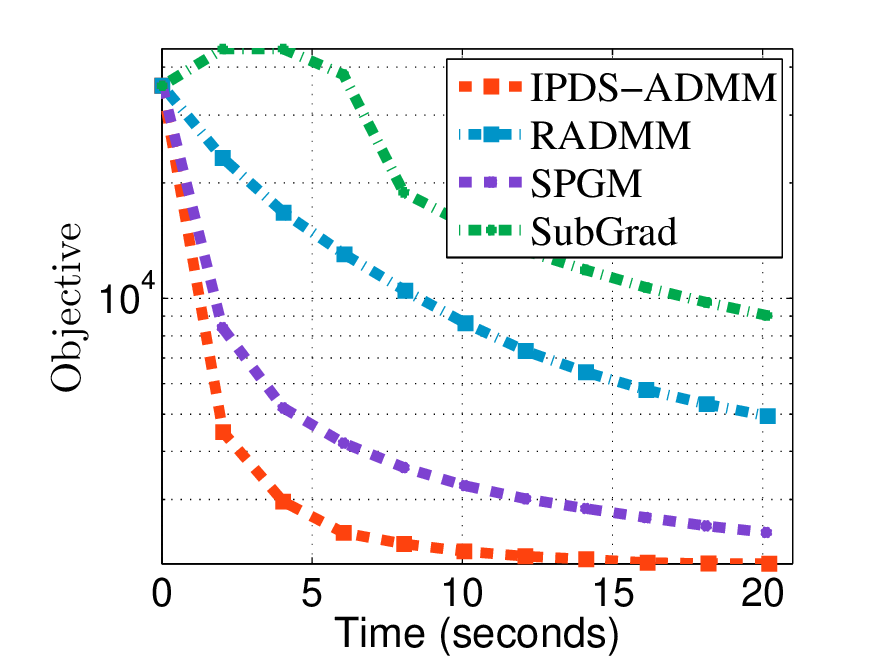}\caption{\scriptsize TDT2-1500-500}\label{fig:sub1}\end{subfigure}
\begin{subfigure}{.24\textwidth}\centering\includegraphics[width=1.12\linewidth]{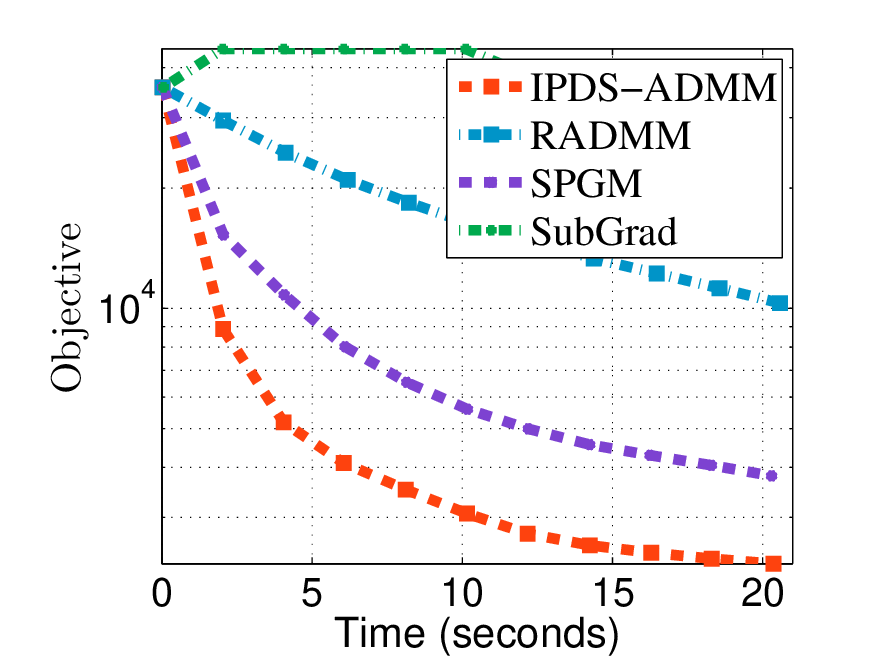}\caption{\scriptsize TDT2-3000-500}\label{fig:sub2}\end{subfigure}
\begin{subfigure}{.24\textwidth}\centering\includegraphics[width=1.12\linewidth]{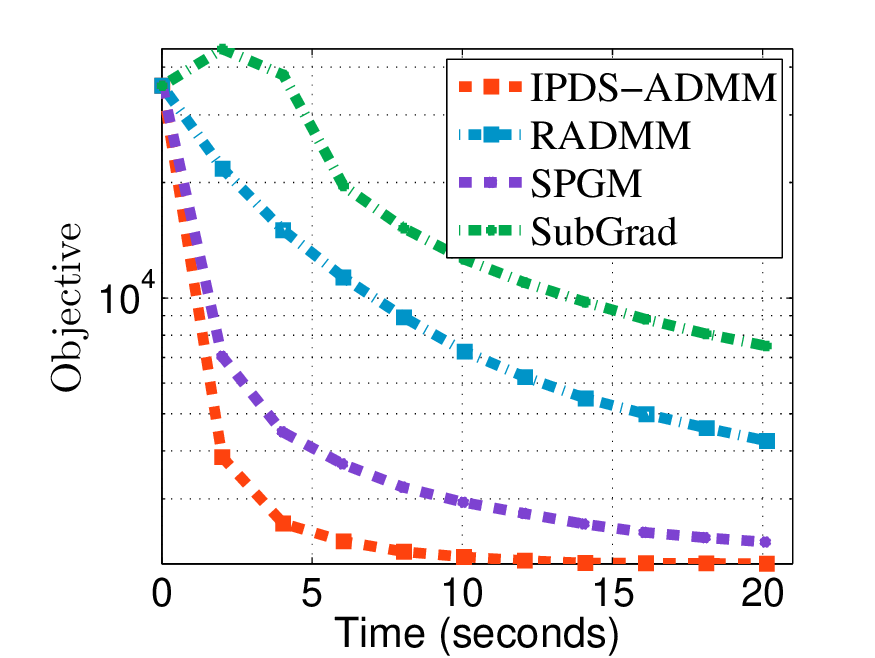}\caption{\scriptsize sector-1500-500}\label{fig:sub3}\end{subfigure}
\begin{subfigure}{.24\textwidth}\centering\includegraphics[width=1.12\linewidth]{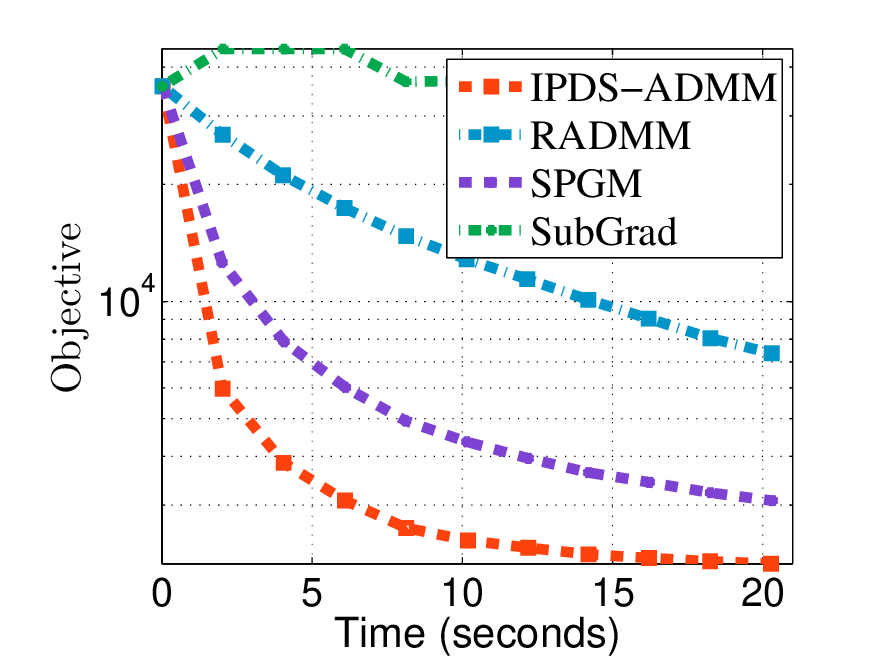}\caption{\scriptsize sector-2500-500}\label{fig:sub4}\end{subfigure}

\caption{Convergence curves of methods for sparse PCA with $\dot{\rho}=100$ and $\beta^0=50\dot{\rho}$.} \label{fig:7}

\end{figure}

\begin{figure}[!t]

\centering
\begin{subfigure}{.24\textwidth}\centering\includegraphics[width=1.12\linewidth]{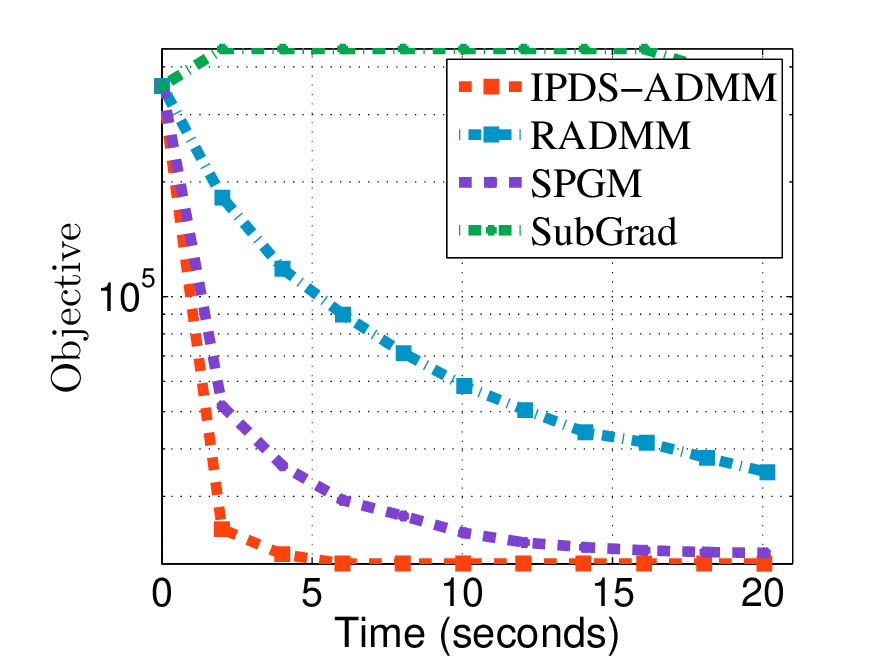}\caption{\scriptsize randn-1500-500}\label{fig:sub1}\end{subfigure}
\begin{subfigure}{.24\textwidth}\centering\includegraphics[width=1.12\linewidth]{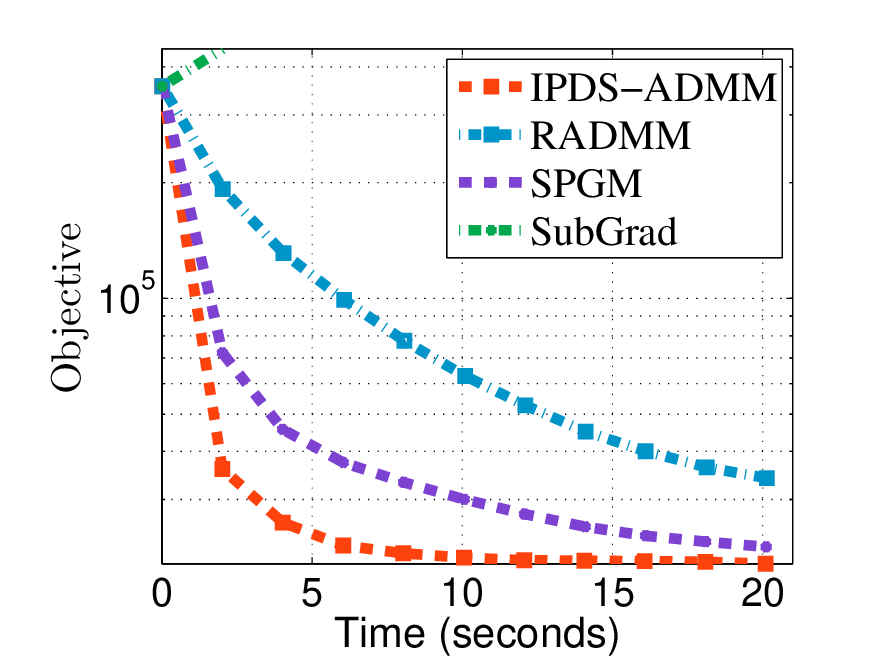}\caption{\scriptsize randn-2500-500}\label{fig:sub2}\end{subfigure}
\begin{subfigure}{.24\textwidth}\centering\includegraphics[width=1.12\linewidth]{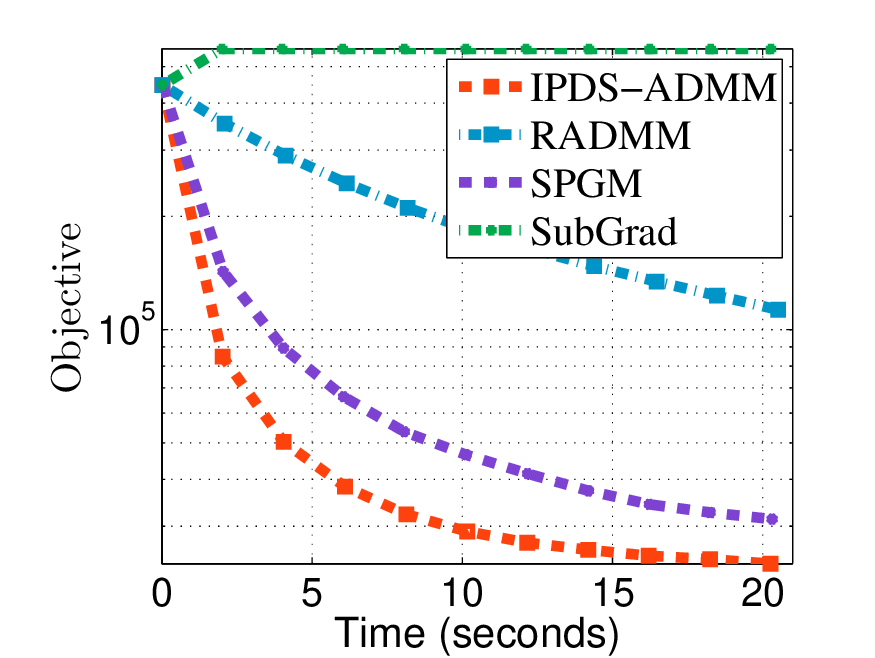}\caption{\scriptsize mnist-1500-780}\label{fig:sub3}\end{subfigure}
\begin{subfigure}{.24\textwidth}\centering\includegraphics[width=1.12\linewidth]{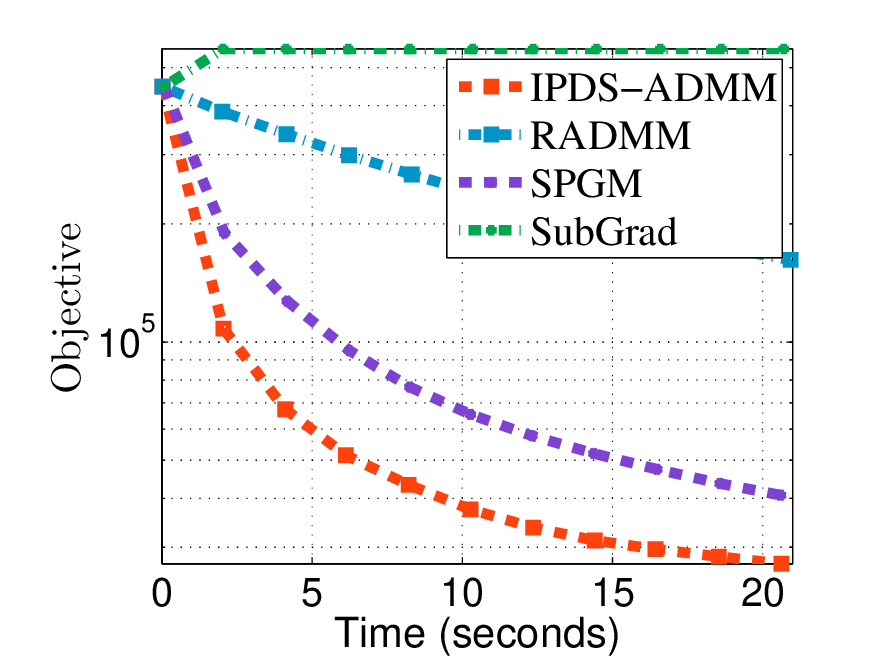}\caption{\scriptsize mnist-2500-780}\label{fig:sub4}\end{subfigure}

\centering
\begin{subfigure}{.24\textwidth}\centering\includegraphics[width=1.12\linewidth]{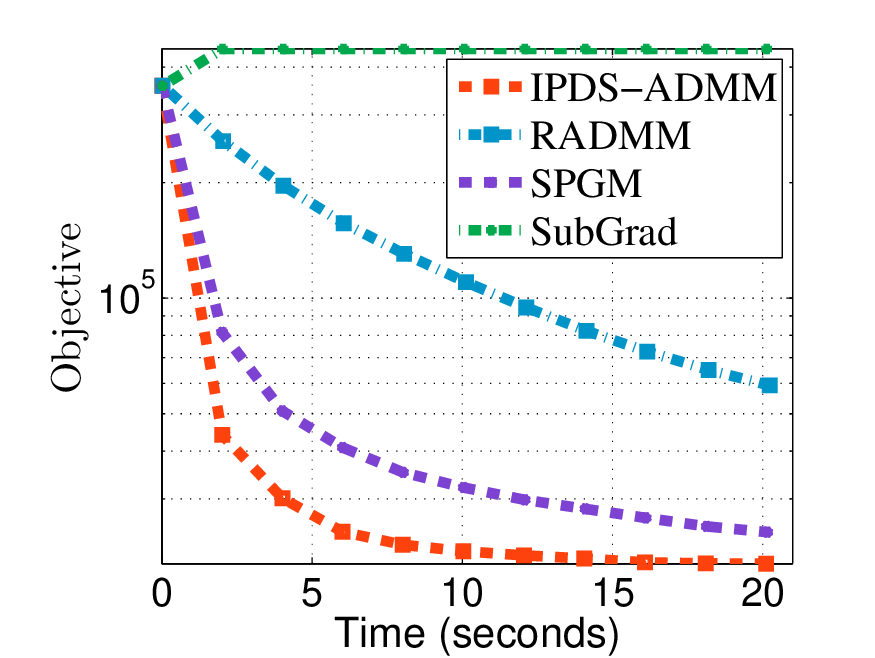}\caption{\scriptsize TDT2-1500-500}\label{fig:sub1}\end{subfigure}
\begin{subfigure}{.24\textwidth}\centering\includegraphics[width=1.12\linewidth]{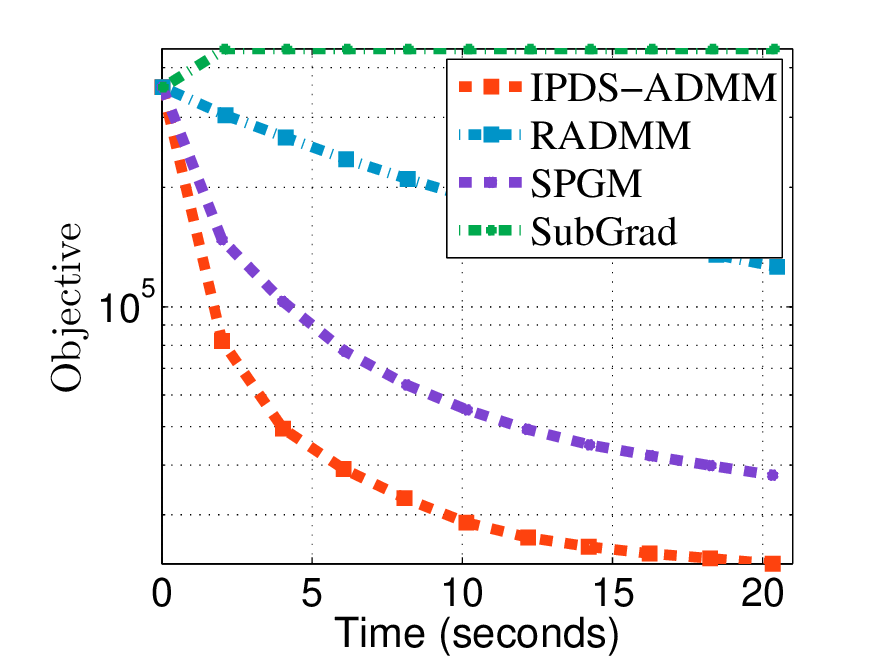}\caption{\scriptsize TDT2-3000-500}\label{fig:sub2}\end{subfigure}
\begin{subfigure}{.24\textwidth}\centering\includegraphics[width=1.12\linewidth]{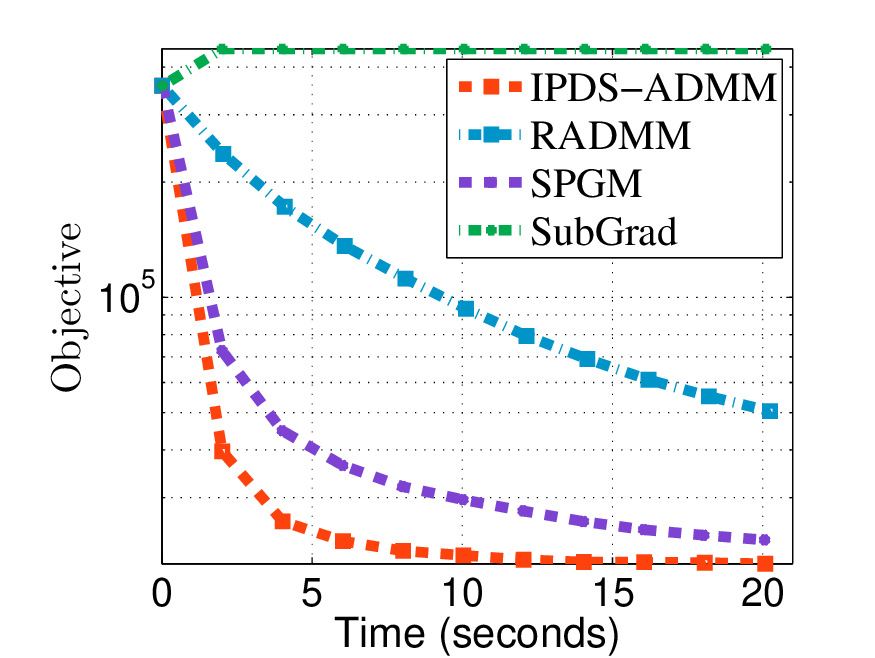}\caption{\scriptsize sector-1500-500}\label{fig:sub3}\end{subfigure}
\begin{subfigure}{.24\textwidth}\centering\includegraphics[width=1.12\linewidth]{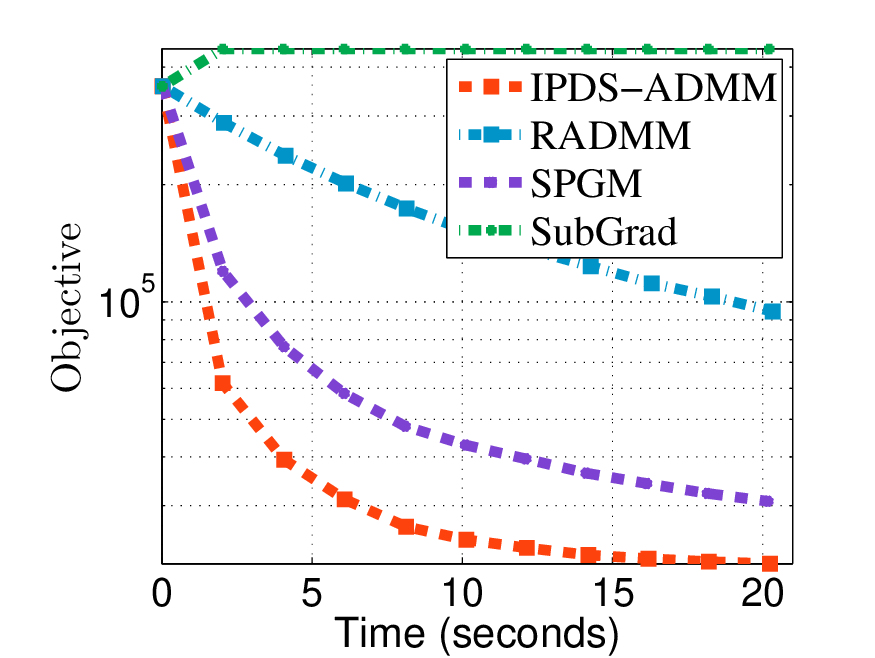}\caption{\scriptsize sector-2500-500}\label{fig:sub4}\end{subfigure}

\caption{Convergence curves of methods for sparse PCA with $\dot{\rho}=1000$ and $\beta^0=50\dot{\rho}$.} \label{fig:8}

\end{figure}

\begin{figure}[!t]
\centering
\begin{subfigure}{.24\textwidth}\centering\includegraphics[width=1.12\linewidth]{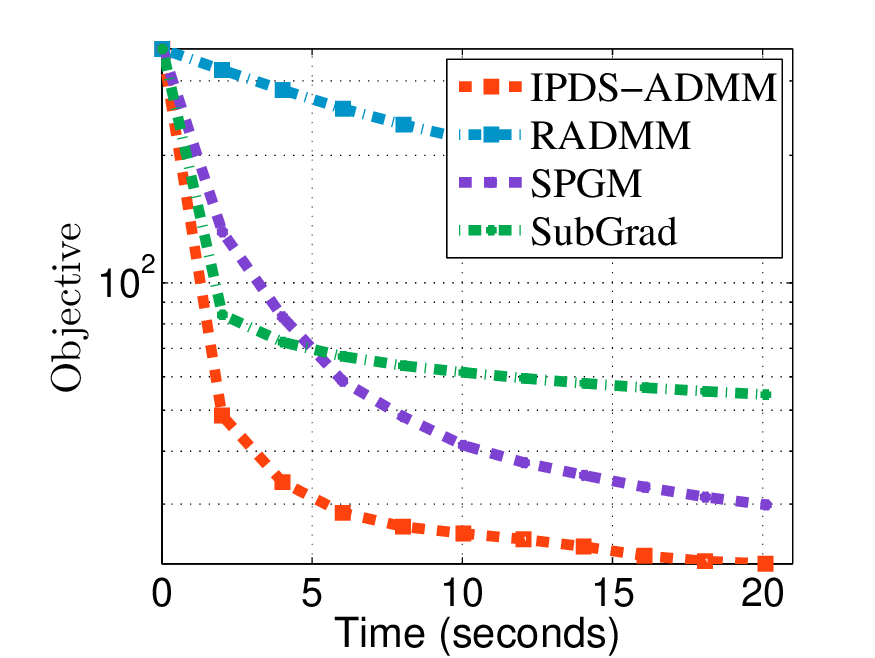}\caption{\scriptsize randn-1500-500}\label{fig:sub1}\end{subfigure}
\begin{subfigure}{.24\textwidth}\centering\includegraphics[width=1.12\linewidth]{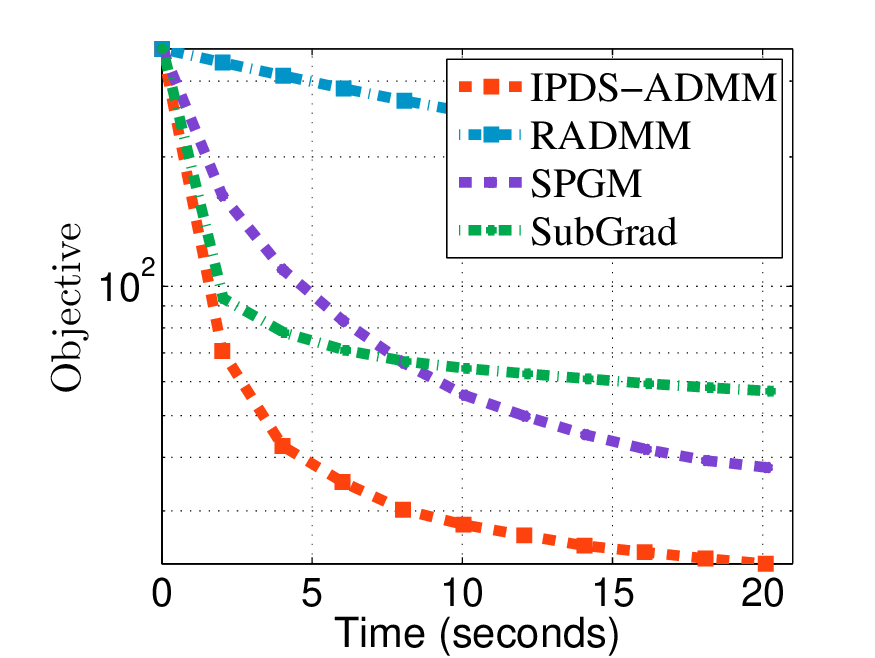}\caption{\scriptsize randn-2500-500}\label{fig:sub2}\end{subfigure}
\begin{subfigure}{.24\textwidth}\centering\includegraphics[width=1.12\linewidth]{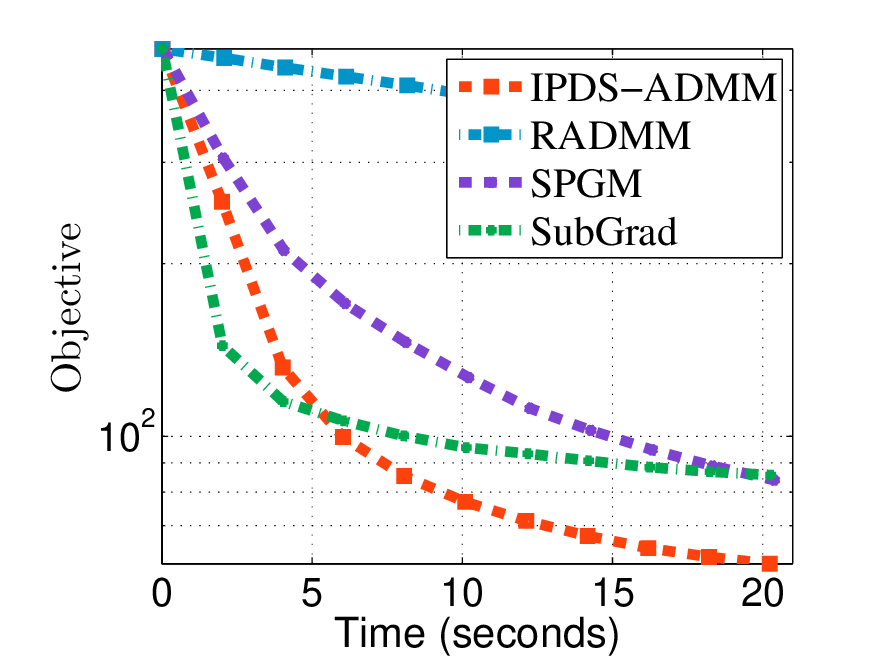}\caption{\scriptsize mnist-1500-780}\label{fig:sub3}\end{subfigure}
\begin{subfigure}{.24\textwidth}\centering\includegraphics[width=1.12\linewidth]{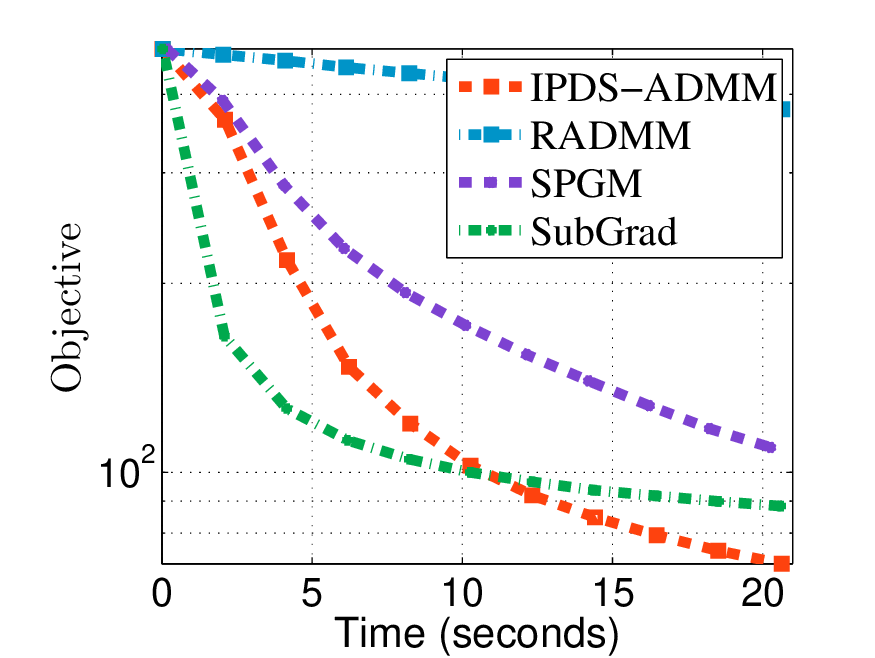}\caption{\scriptsize mnist-2500-780}\label{fig:sub4}\end{subfigure}

\centering
\begin{subfigure}{.24\textwidth}\centering\includegraphics[width=1.12\linewidth]{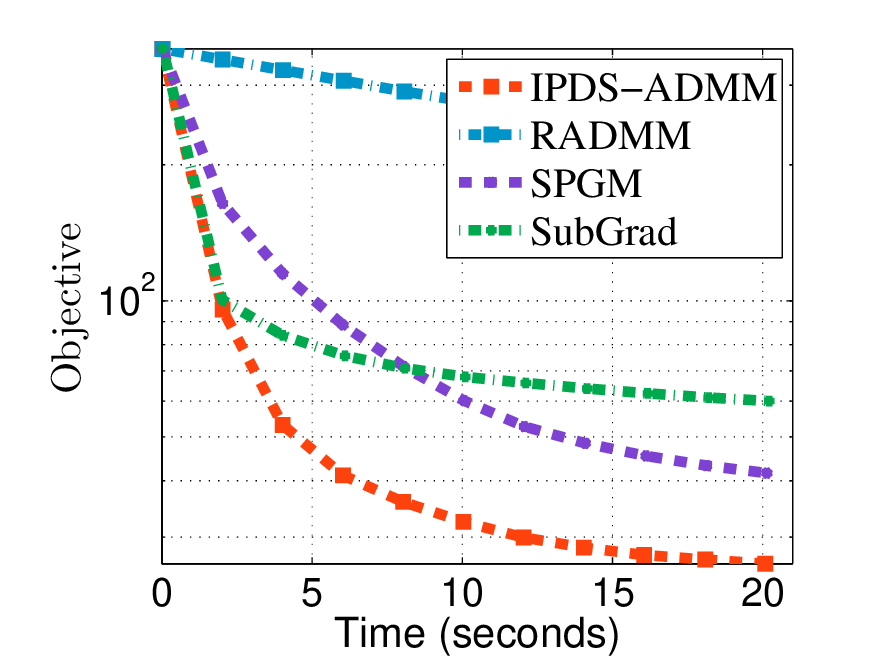}\caption{\scriptsize TDT2-1500-500}\label{fig:sub1}\end{subfigure}
\begin{subfigure}{.24\textwidth}\centering\includegraphics[width=1.12\linewidth]{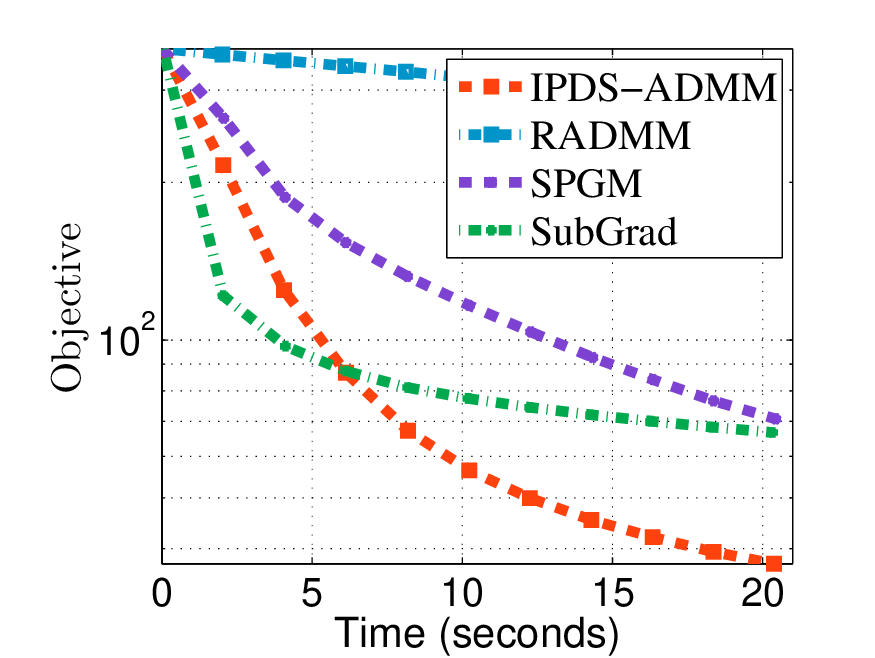}\caption{\scriptsize TDT2-3000-500}\label{fig:sub2}\end{subfigure}
\begin{subfigure}{.24\textwidth}\centering\includegraphics[width=1.12\linewidth]{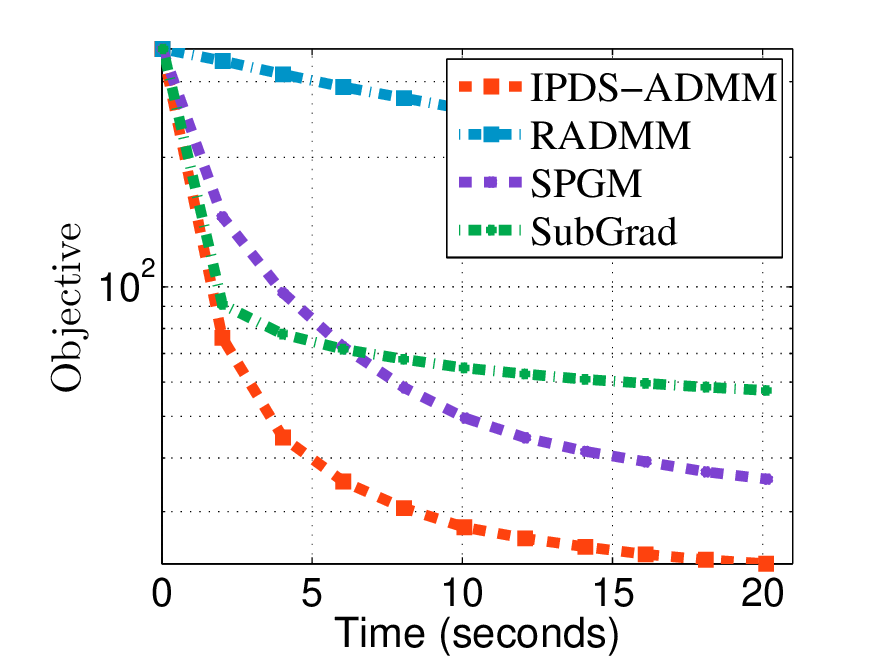}\caption{\scriptsize sector-1500-500}\label{fig:sub3}\end{subfigure}
\begin{subfigure}{.24\textwidth}\centering\includegraphics[width=1.12\linewidth]{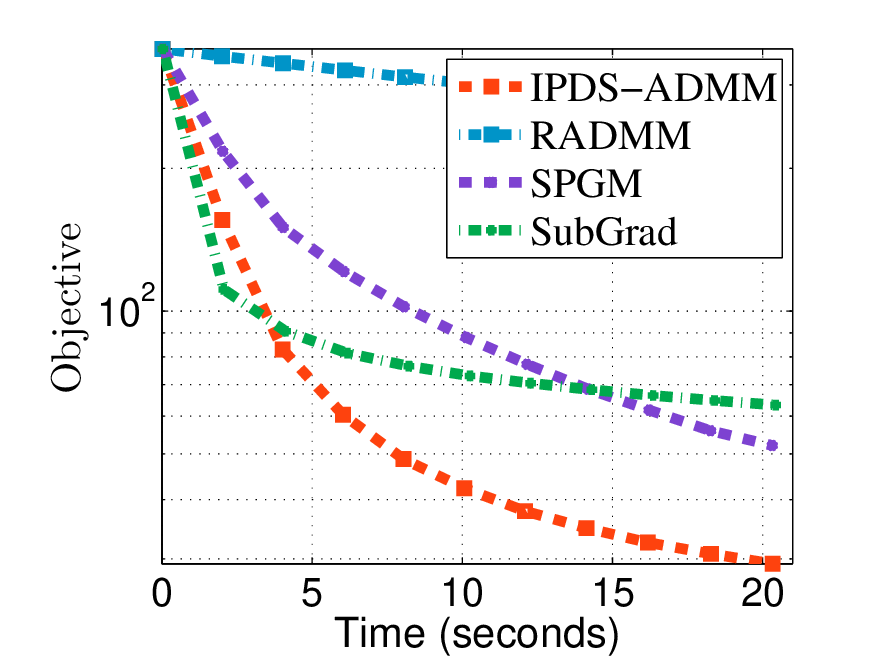}\caption{\scriptsize sector-2500-500}\label{fig:sub4}\end{subfigure}

\caption{Convergence curves of methods for sparse PCA with $\dot{\rho}=1$ and $\beta^0=100\dot{\rho}$.} \label{fig:9}

\end{figure}

\begin{figure}[!t]

\centering
\begin{subfigure}{.24\textwidth}\centering\includegraphics[width=1.12\linewidth]{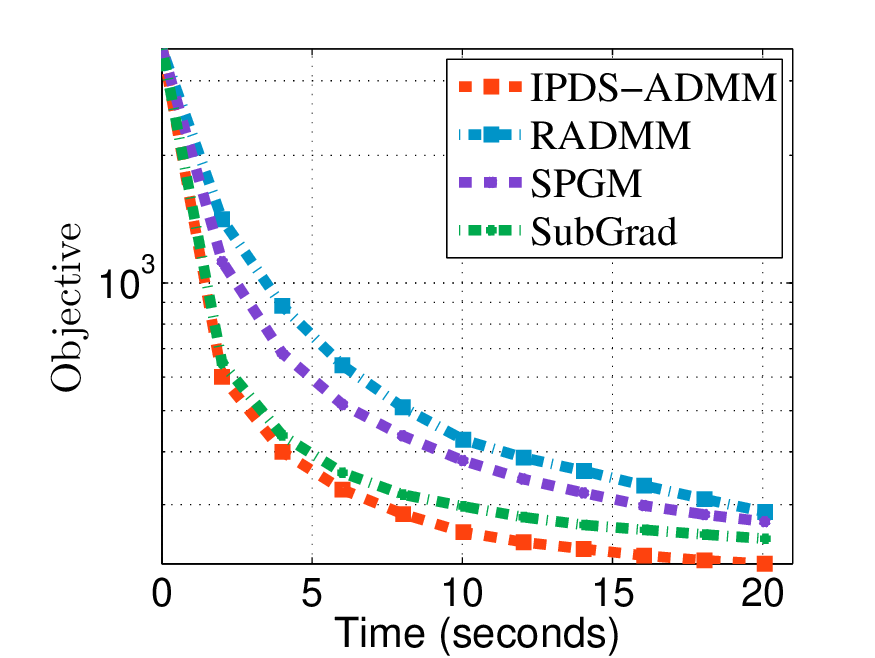}\caption{\scriptsize randn-1500-500}\label{fig:sub1}\end{subfigure}
\begin{subfigure}{.24\textwidth}\centering\includegraphics[width=1.12\linewidth]{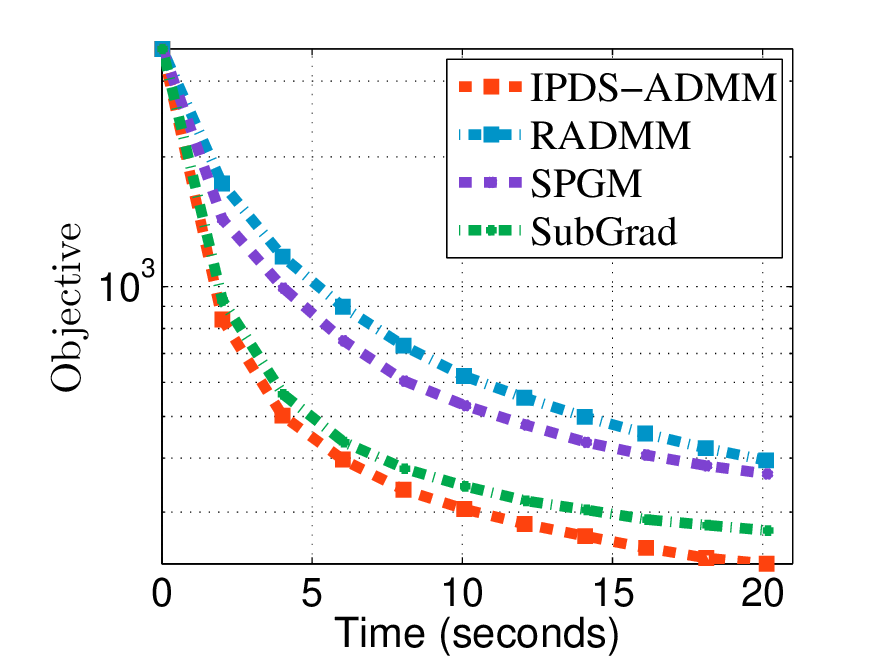}\caption{\scriptsize randn-2500-500}\label{fig:sub2}\end{subfigure}
\begin{subfigure}{.24\textwidth}\centering\includegraphics[width=1.12\linewidth]{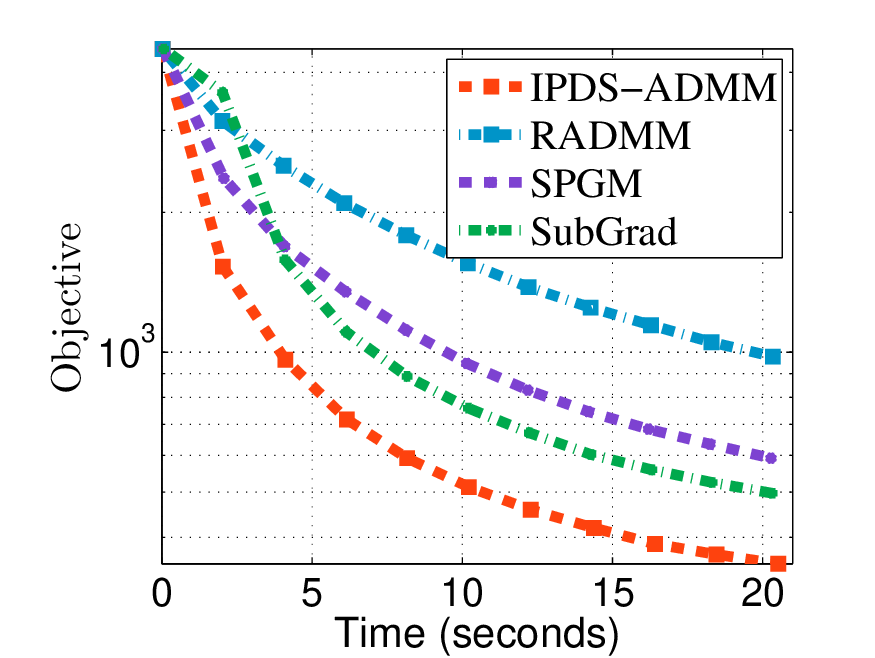}\caption{\scriptsize mnist-1500-780}\label{fig:sub3}\end{subfigure}
\begin{subfigure}{.24\textwidth}\centering\includegraphics[width=1.12\linewidth]{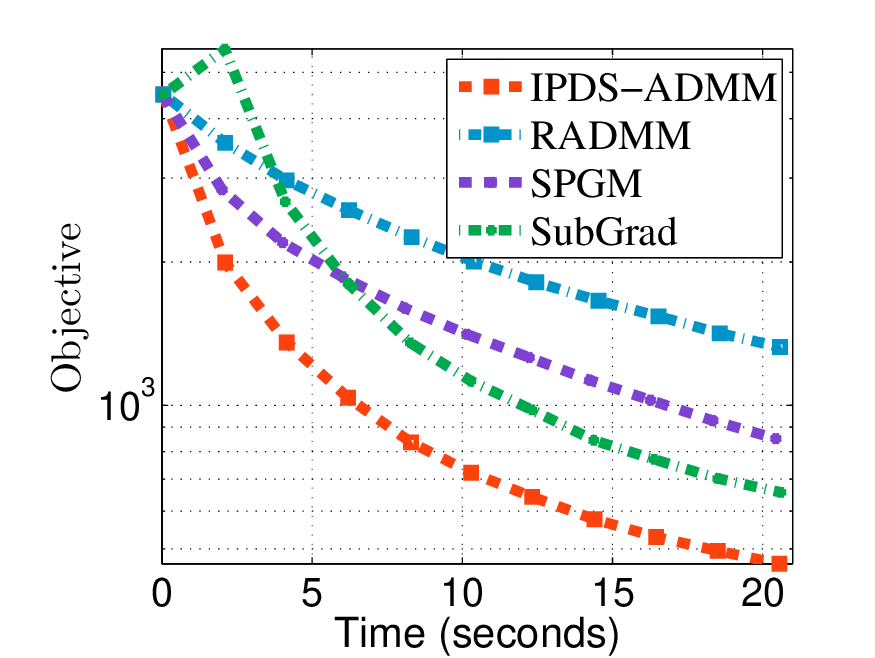}\caption{\scriptsize mnist-2500-780}\label{fig:sub4}\end{subfigure}

\centering
\begin{subfigure}{.24\textwidth}\centering\includegraphics[width=1.12\linewidth]{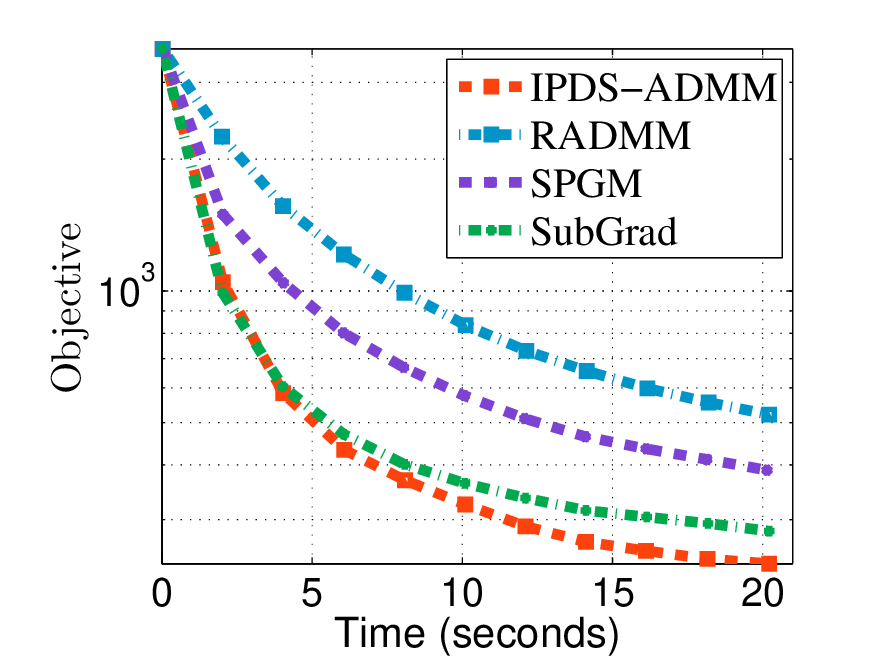}\caption{\scriptsize TDT2-1500-500}\label{fig:sub1}\end{subfigure}
\begin{subfigure}{.24\textwidth}\centering\includegraphics[width=1.12\linewidth]{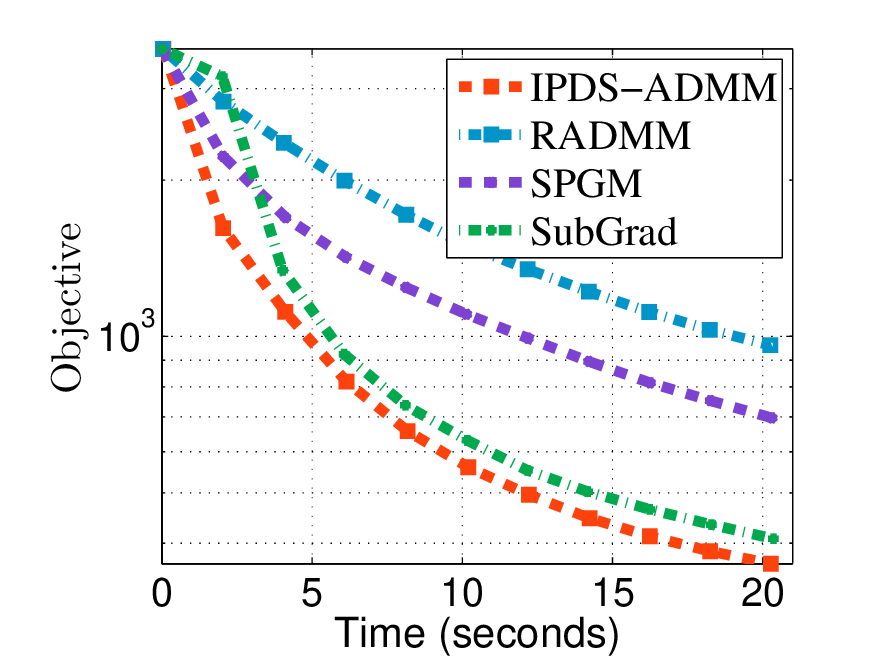}\caption{\scriptsize TDT2-3000-500}\label{fig:sub2}\end{subfigure}
\begin{subfigure}{.24\textwidth}\centering\includegraphics[width=1.12\linewidth]{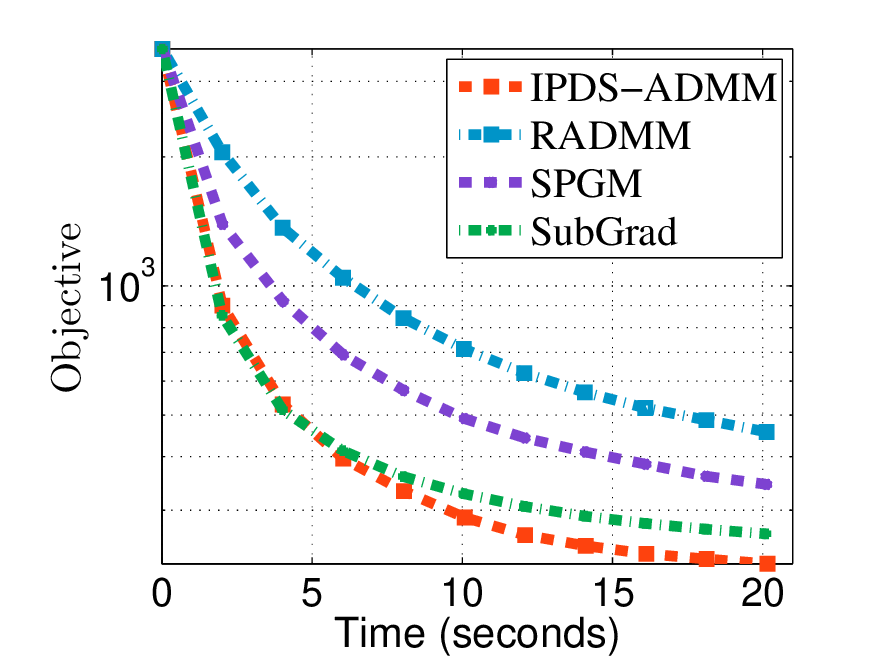}\caption{\scriptsize sector-1500-500}\label{fig:sub3}\end{subfigure}
\begin{subfigure}{.24\textwidth}\centering\includegraphics[width=1.12\linewidth]{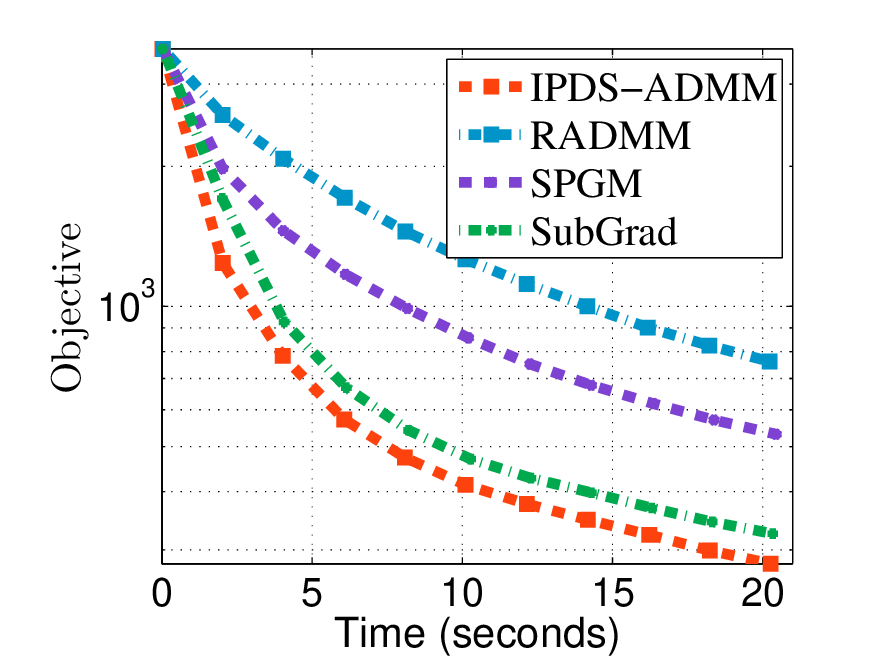}\caption{\scriptsize sector-2500-500}\label{fig:sub4}\end{subfigure}

\caption{Convergence curves of methods for sparse PCA with $\dot{\rho}=10$ and $\beta^0=100\dot{\rho}$.} \label{fig:10}
\end{figure}

\begin{figure}[!t]
\centering
\begin{subfigure}{.24\textwidth}\centering\includegraphics[width=1.12\linewidth]{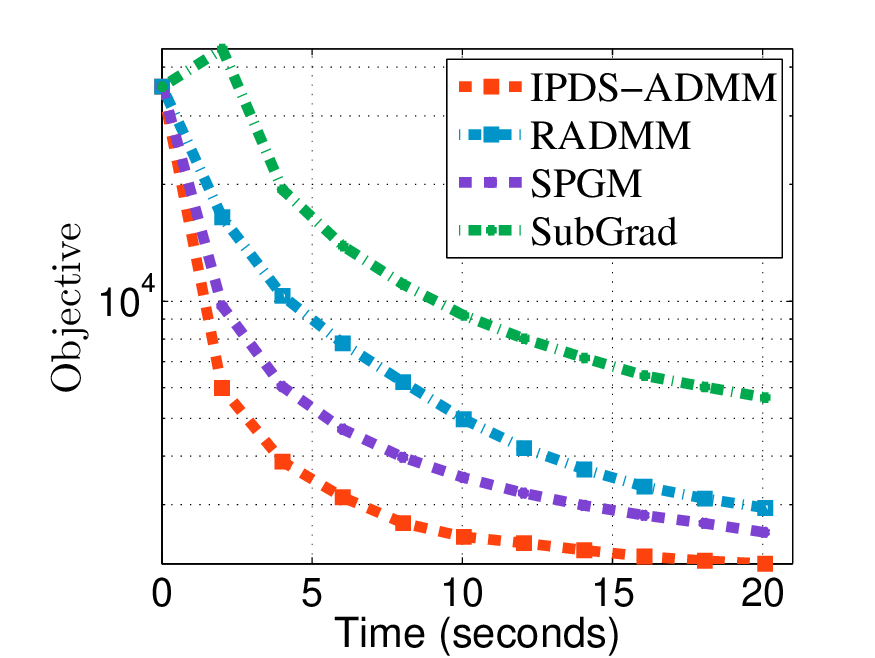}\caption{\scriptsize randn-1500-500}\label{fig:sub1}\end{subfigure}
\begin{subfigure}{.24\textwidth}\centering\includegraphics[width=1.12\linewidth]{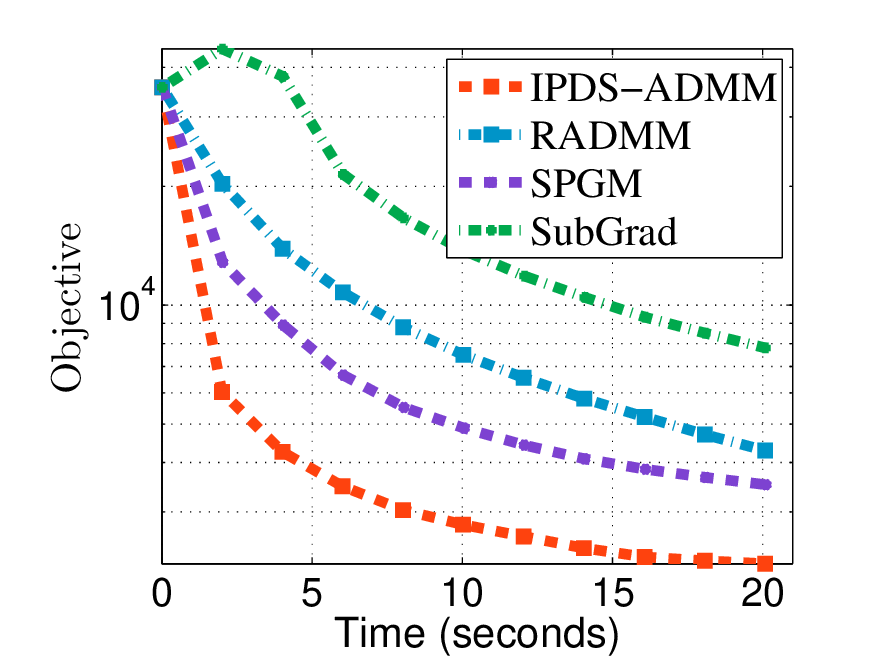}\caption{\scriptsize randn-2500-500}\label{fig:sub2}\end{subfigure}
\begin{subfigure}{.24\textwidth}\centering\includegraphics[width=1.12\linewidth]{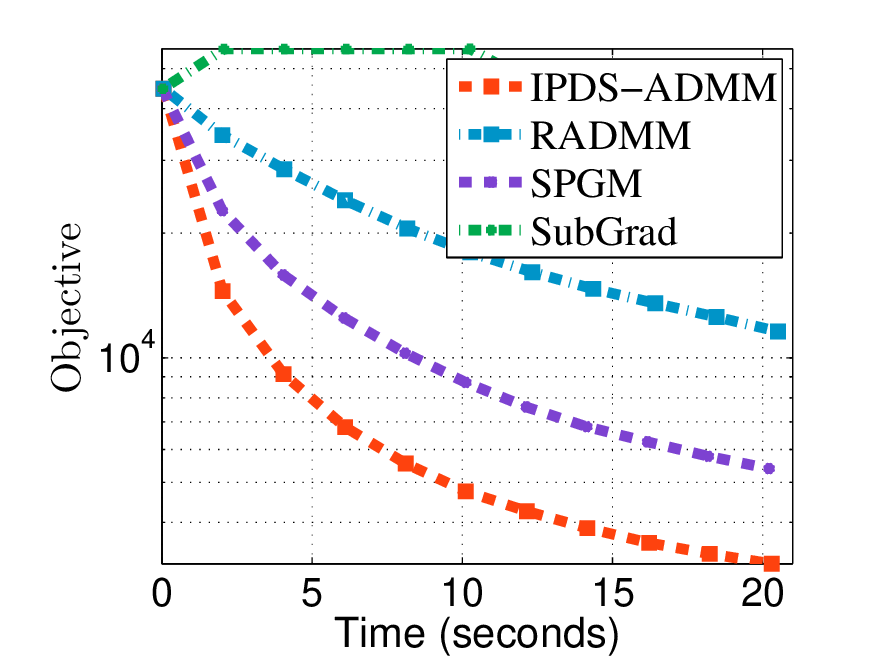}\caption{\scriptsize mnist-1500-780}\label{fig:sub3}\end{subfigure}
\begin{subfigure}{.24\textwidth}\centering\includegraphics[width=1.12\linewidth]{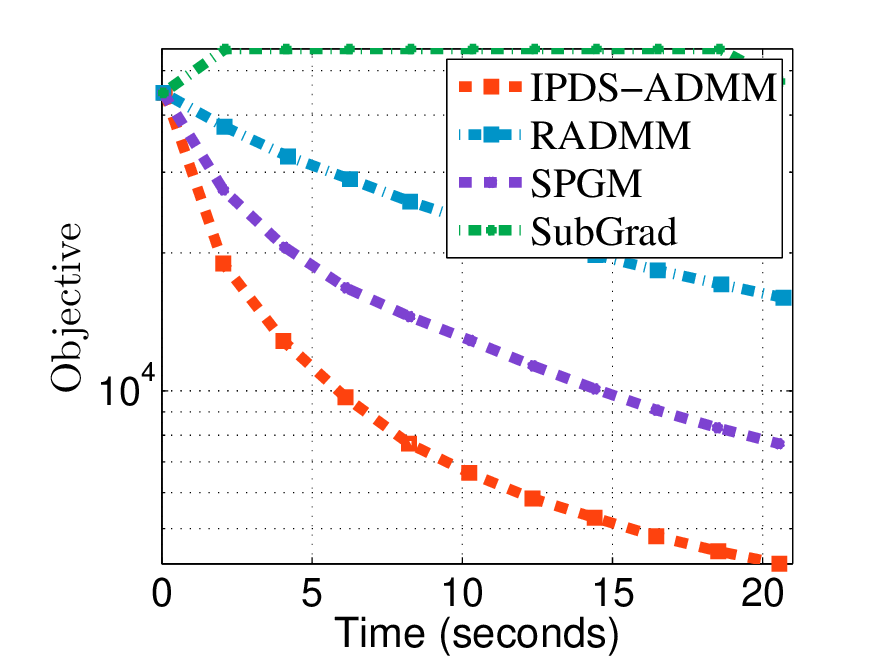}\caption{\scriptsize mnist-2500-780}\label{fig:sub4}\end{subfigure}

\centering
\begin{subfigure}{.24\textwidth}\centering\includegraphics[width=1.12\linewidth]{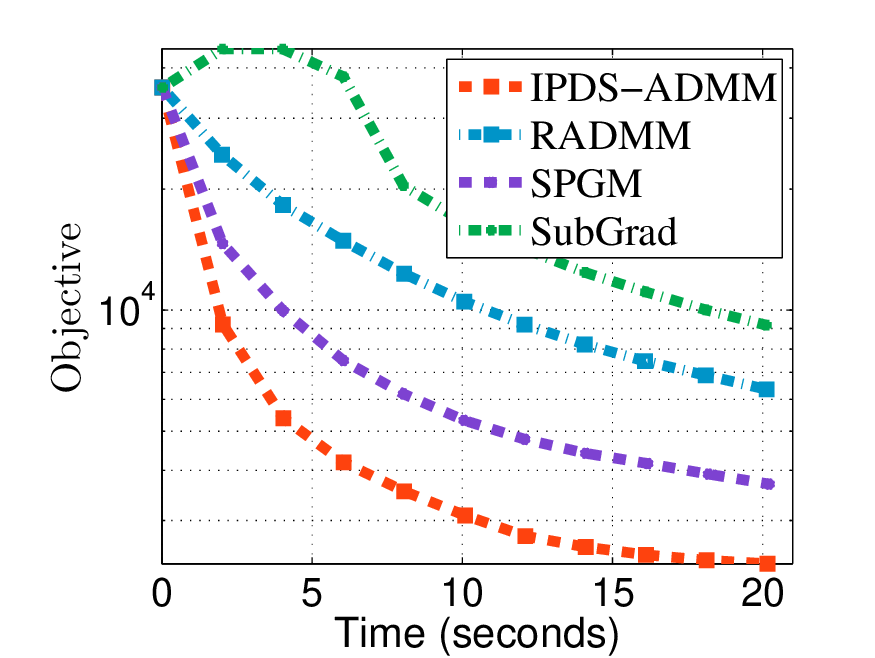}\caption{\scriptsize TDT2-1500-500}\label{fig:sub1}\end{subfigure}
\begin{subfigure}{.24\textwidth}\centering\includegraphics[width=1.12\linewidth]{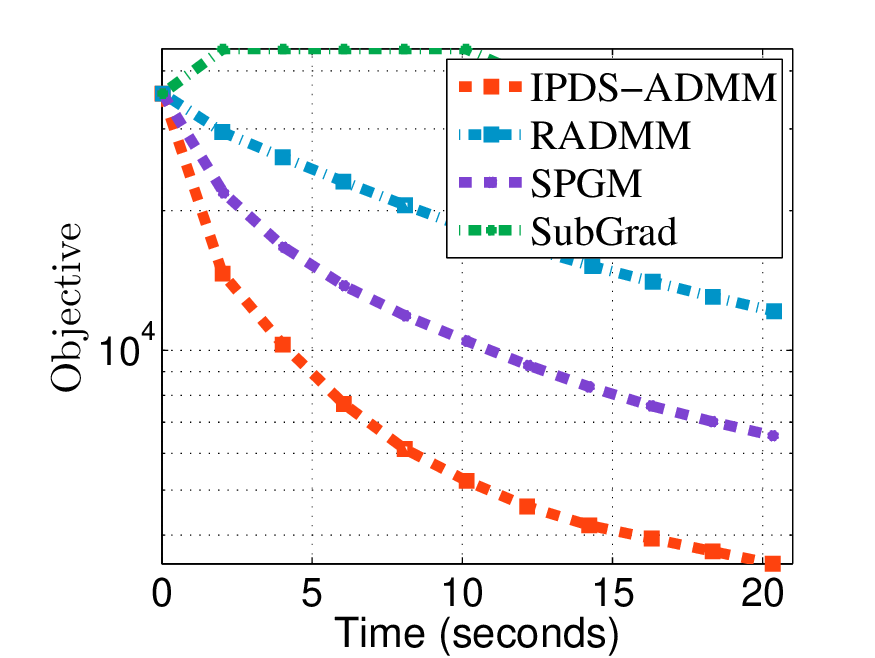}\caption{\scriptsize TDT2-3000-500}\label{fig:sub2}\end{subfigure}
\begin{subfigure}{.24\textwidth}\centering\includegraphics[width=1.12\linewidth]{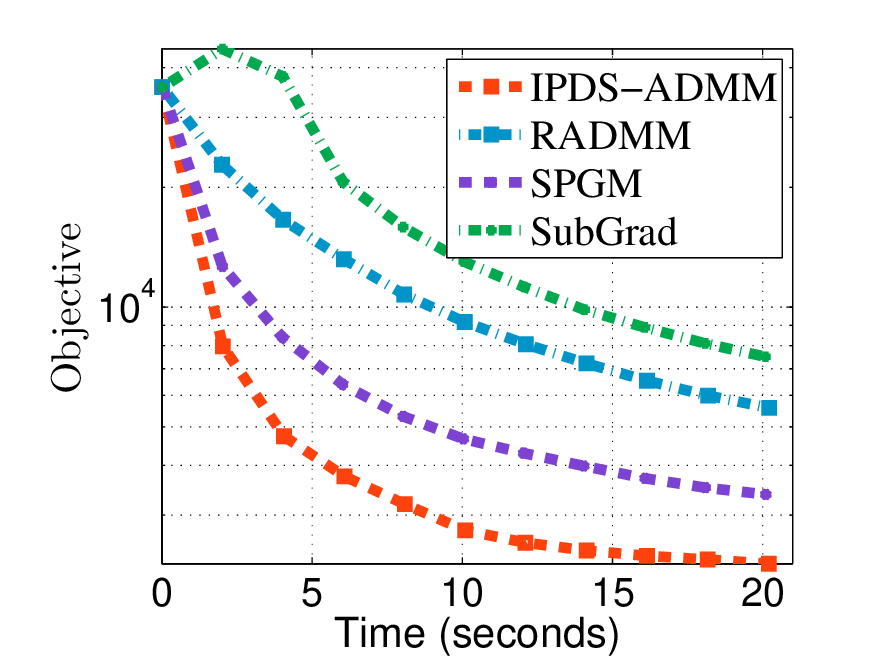}\caption{\scriptsize sector-1500-500}\label{fig:sub3}\end{subfigure}
\begin{subfigure}{.24\textwidth}\centering\includegraphics[width=1.12\linewidth]{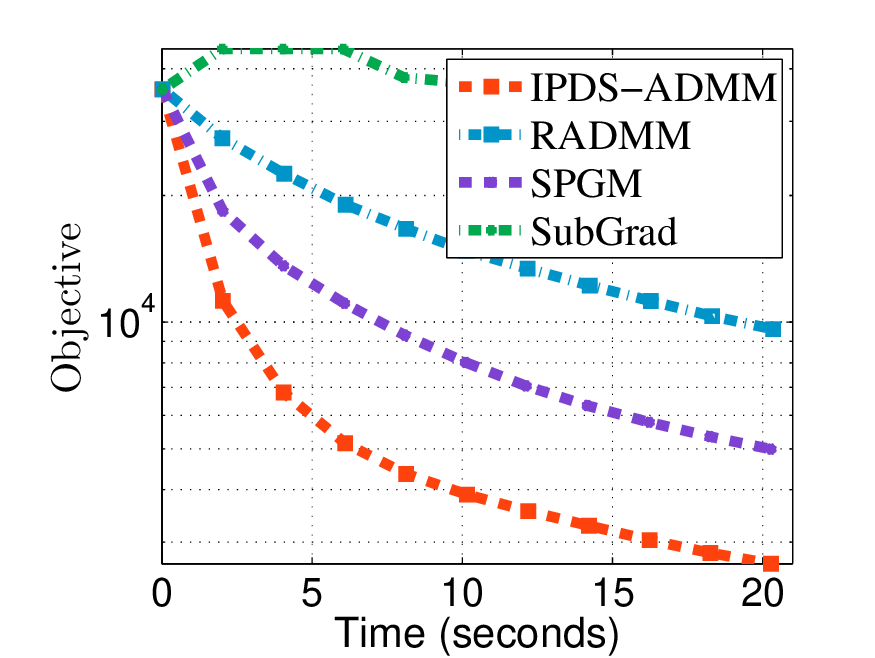}\caption{\scriptsize sector-2500-500}\label{fig:sub4}\end{subfigure}

\caption{Convergence curves of methods for sparse PCA with $\dot{\rho}=100$ and $\beta^0=100\dot{\rho}$.} \label{fig:11}
\end{figure}

\begin{figure}[!t]
\centering
\begin{subfigure}{.24\textwidth}\centering\includegraphics[width=1.12\linewidth]{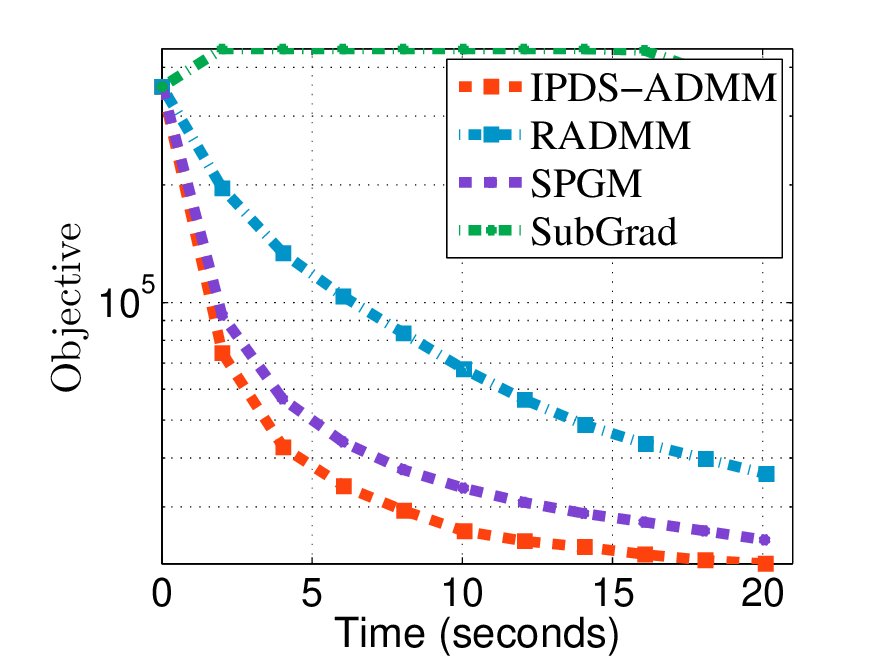}\caption{\scriptsize randn-1500-500}\label{fig:sub1}\end{subfigure}
\begin{subfigure}{.24\textwidth}\centering\includegraphics[width=1.12\linewidth]{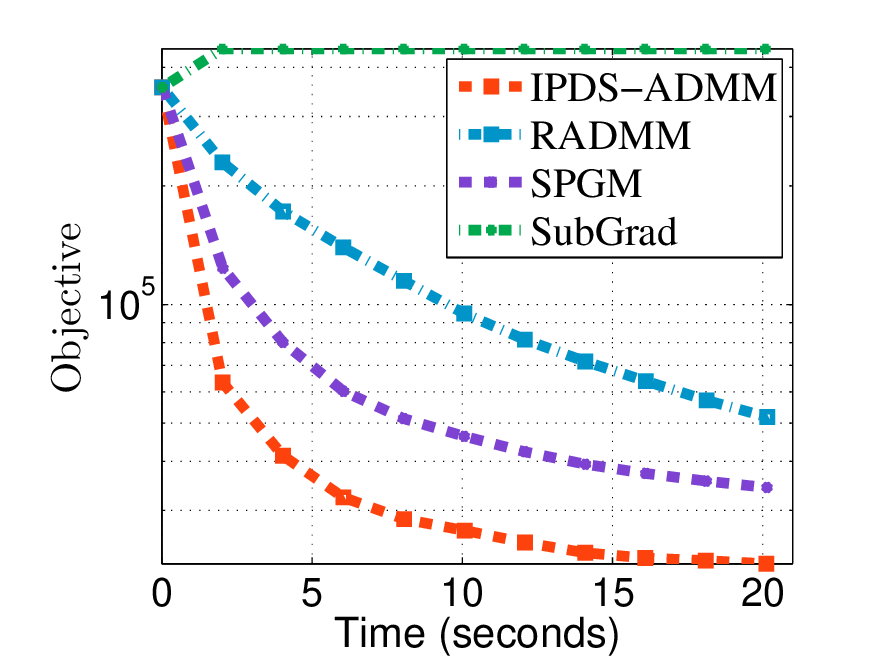}\caption{\scriptsize randn-2500-500}\label{fig:sub2}\end{subfigure}
\begin{subfigure}{.24\textwidth}\centering\includegraphics[width=1.12\linewidth]{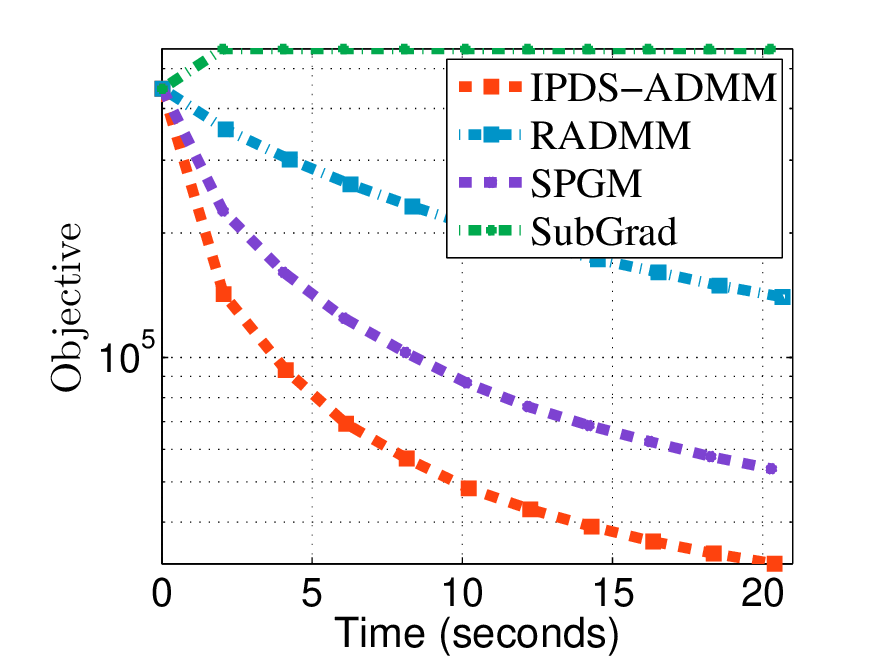}\caption{\scriptsize mnist-1500-780}\label{fig:sub3}\end{subfigure}
\begin{subfigure}{.24\textwidth}\centering\includegraphics[width=1.12\linewidth]{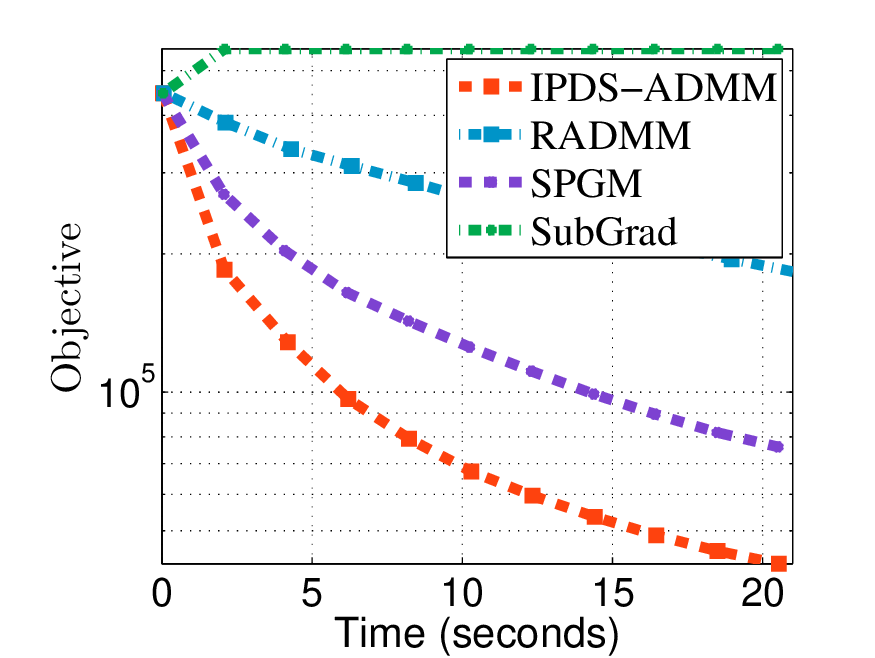}\caption{\scriptsize mnist-2500-780}\label{fig:sub4}\end{subfigure}

\centering
\begin{subfigure}{.24\textwidth}\centering\includegraphics[width=1.12\linewidth]{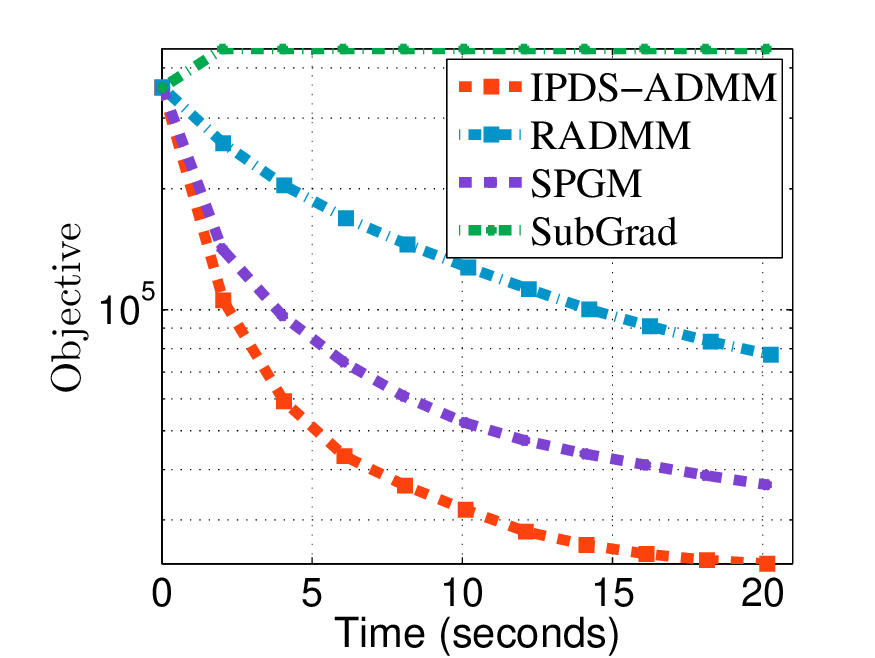}\caption{\scriptsize TDT2-1500-500}\label{fig:sub1}\end{subfigure}
\begin{subfigure}{.24\textwidth}\centering\includegraphics[width=1.12\linewidth]{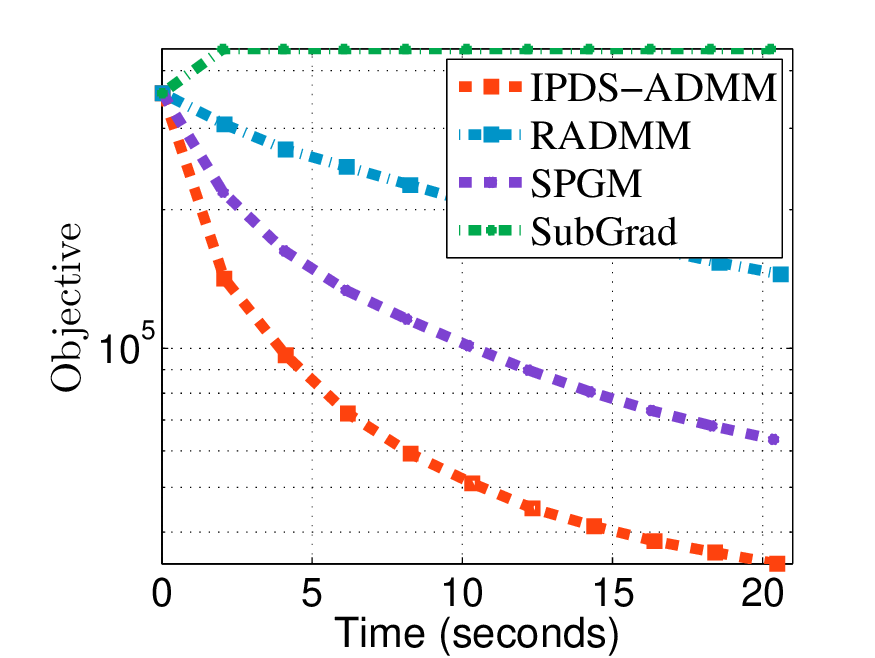}\caption{\scriptsize TDT2-3000-500}\label{fig:sub2}\end{subfigure}
\begin{subfigure}{.24\textwidth}\centering\includegraphics[width=1.12\linewidth]{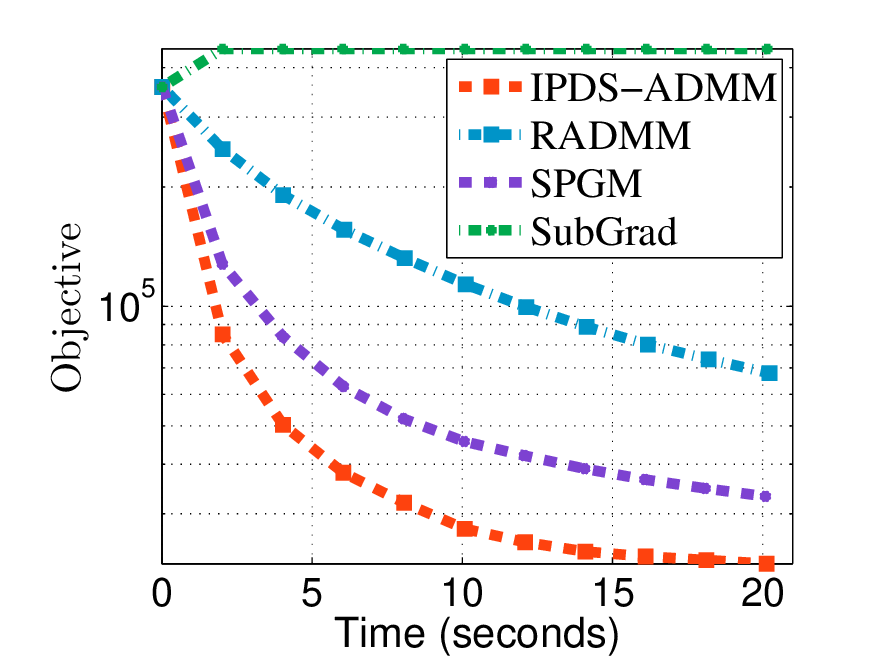}\caption{\scriptsize sector-1500-500}\label{fig:sub3}\end{subfigure}
\begin{subfigure}{.24\textwidth}\centering\includegraphics[width=1.12\linewidth]{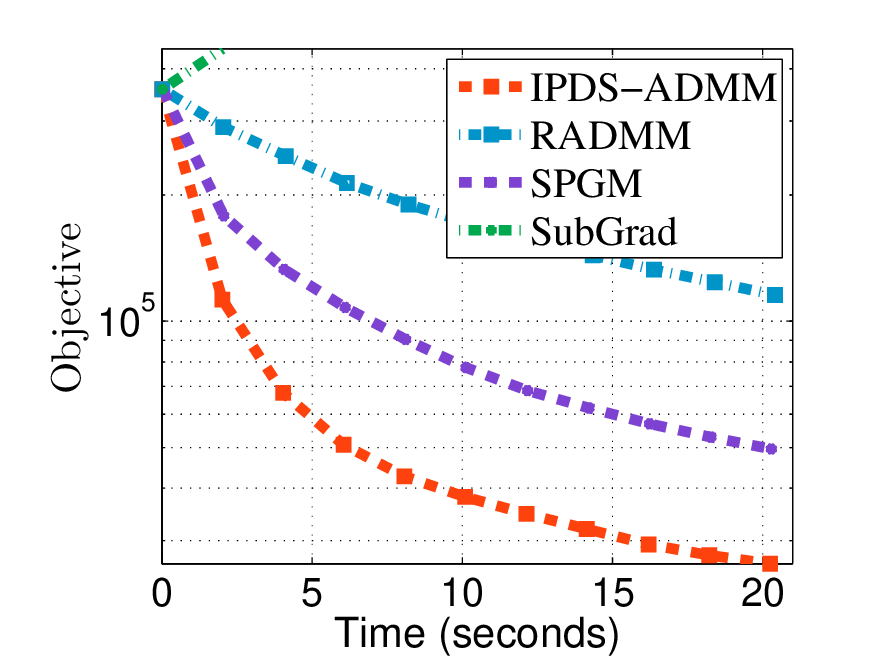}\caption{\scriptsize sector-2500-500}\label{fig:sub4}\end{subfigure}

\caption{Convergence curves of methods for sparse PCA with $\dot{\rho}=1000$ and $\beta^0=100\dot{\rho}$.} \label{fig:12}

\end{figure}

\begin{figure}[!t]

\centering
\begin{subfigure}{.24\textwidth}\centering\includegraphics[width=1.12\linewidth]{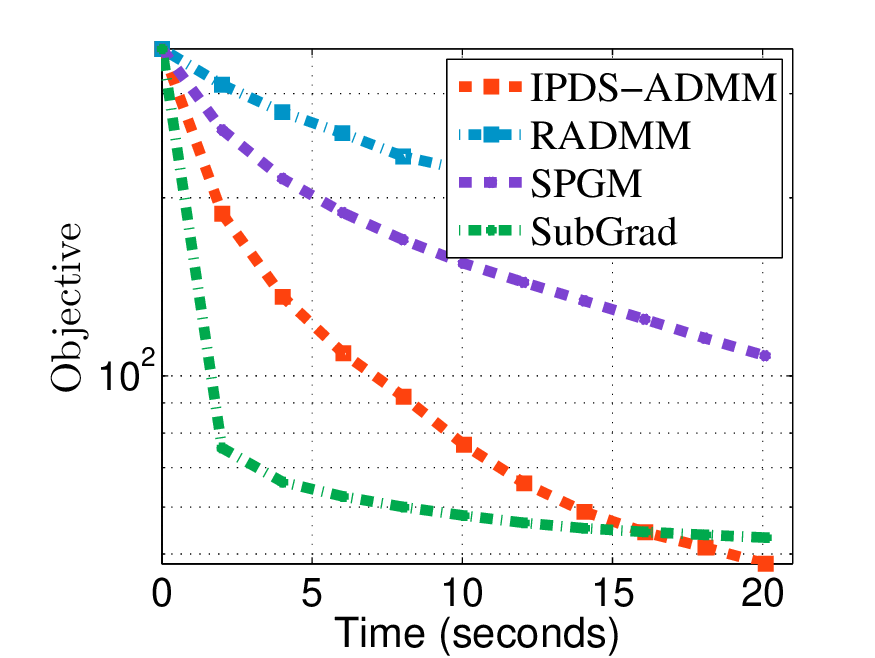}\caption{\scriptsize randn-1500-500}\label{fig:sub1}\end{subfigure}
\begin{subfigure}{.24\textwidth}\centering\includegraphics[width=1.12\linewidth]{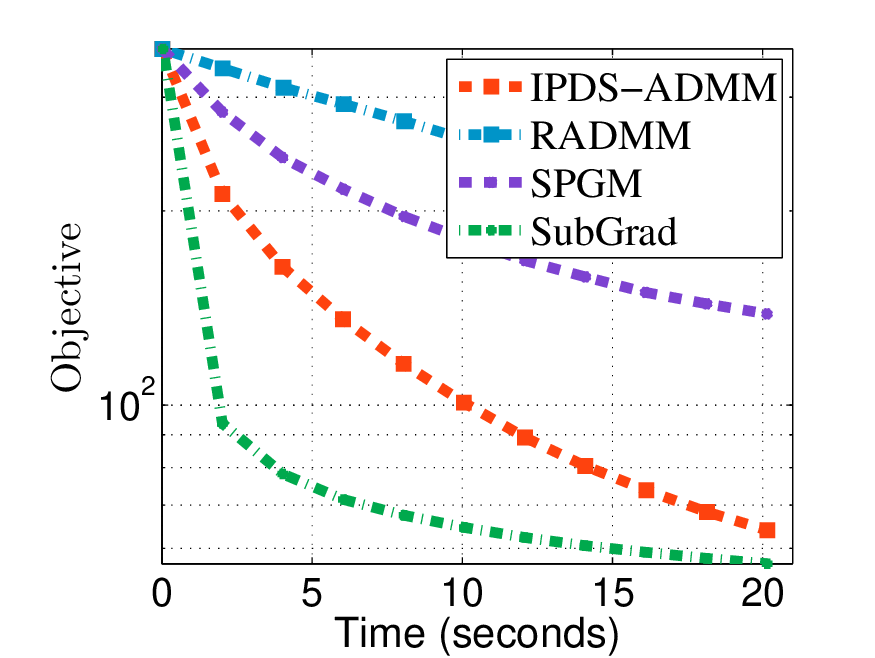}\caption{\scriptsize randn-2500-500}\label{fig:sub2}\end{subfigure}
\begin{subfigure}{.24\textwidth}\centering\includegraphics[width=1.12\linewidth]{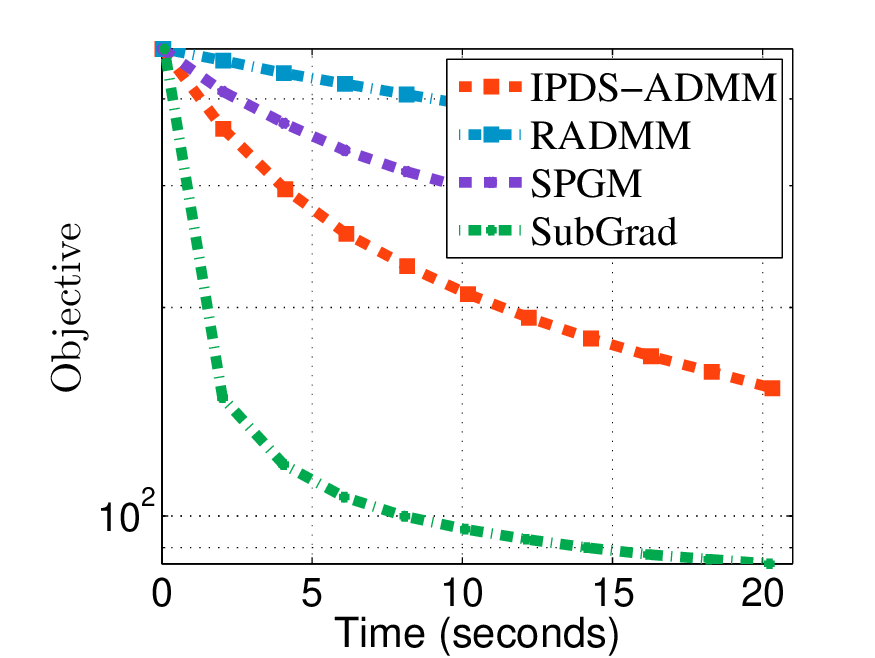}\caption{\scriptsize mnist-1500-780}\label{fig:sub3}\end{subfigure}
\begin{subfigure}{.24\textwidth}\centering\includegraphics[width=1.12\linewidth]{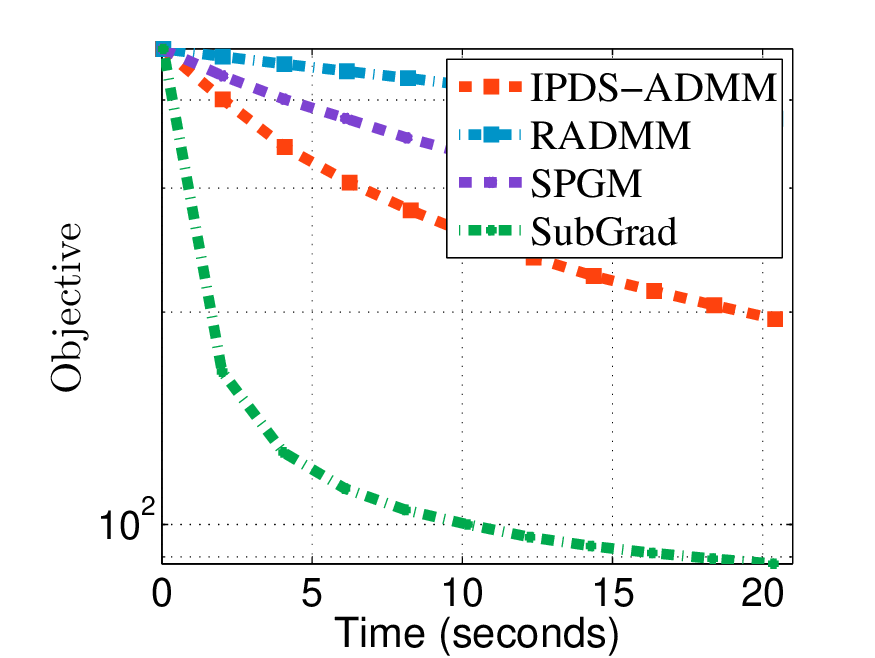}\caption{\scriptsize mnist-2500-780}\label{fig:sub4}\end{subfigure}

\centering
\begin{subfigure}{.24\textwidth}\centering\includegraphics[width=1.12\linewidth]{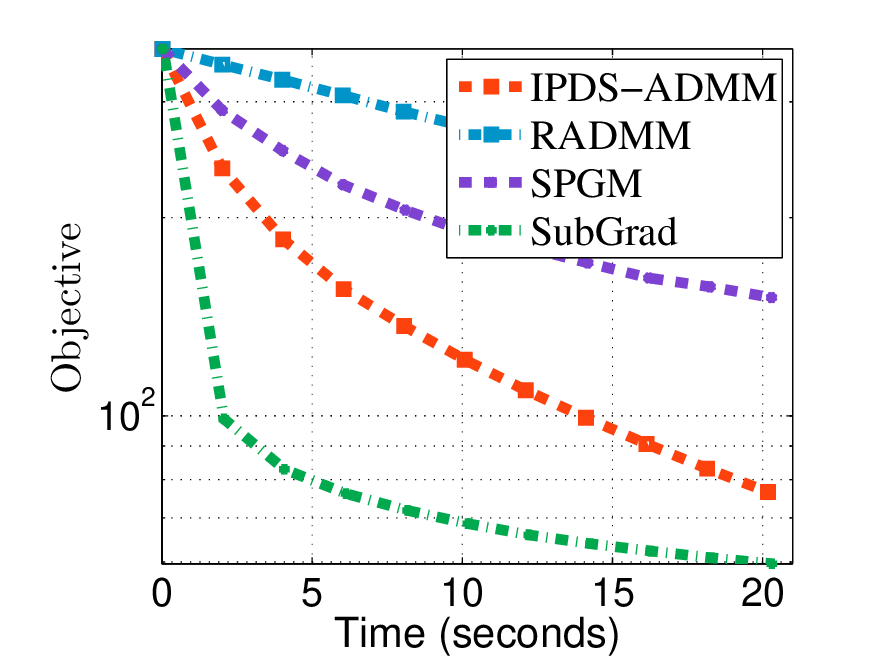}\caption{\scriptsize TDT2-1500-500}\label{fig:sub1}\end{subfigure}
\begin{subfigure}{.24\textwidth}\centering\includegraphics[width=1.12\linewidth]{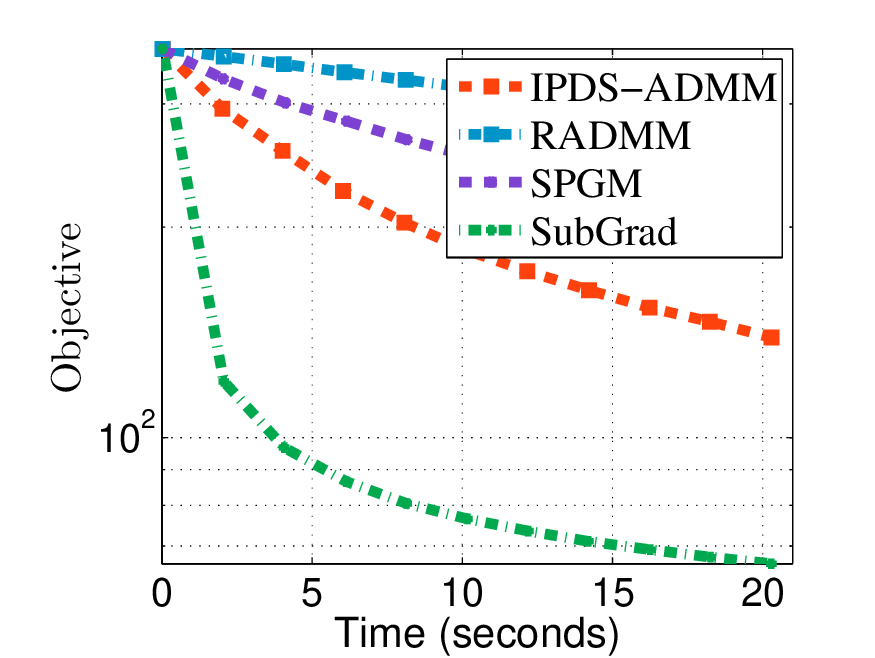}\caption{\scriptsize TDT2-3000-500}\label{fig:sub2}\end{subfigure}
\begin{subfigure}{.24\textwidth}\centering\includegraphics[width=1.12\linewidth]{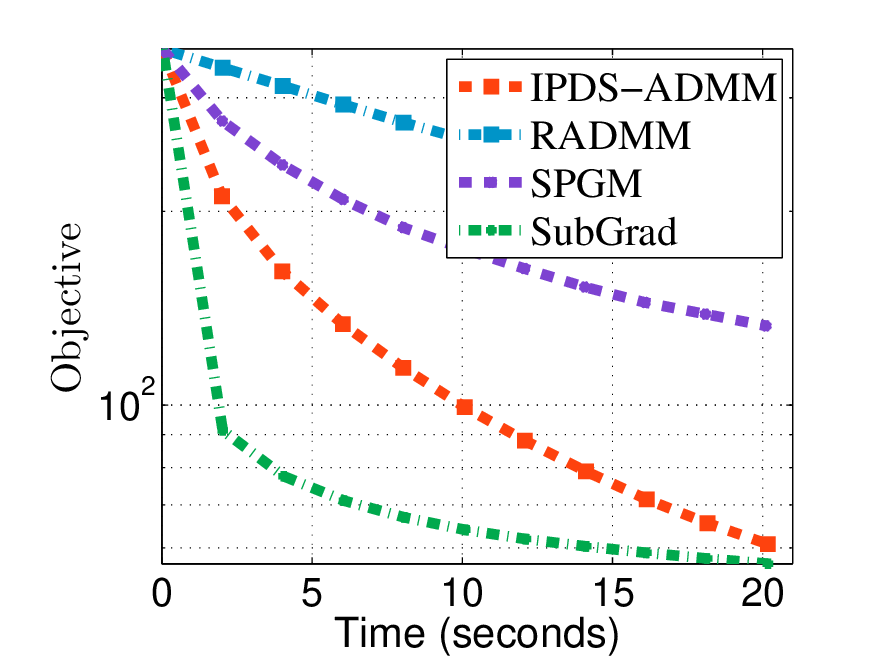}\caption{\scriptsize sector-1500-500}\label{fig:sub3}\end{subfigure}
\begin{subfigure}{.24\textwidth}\centering\includegraphics[width=1.12\linewidth]{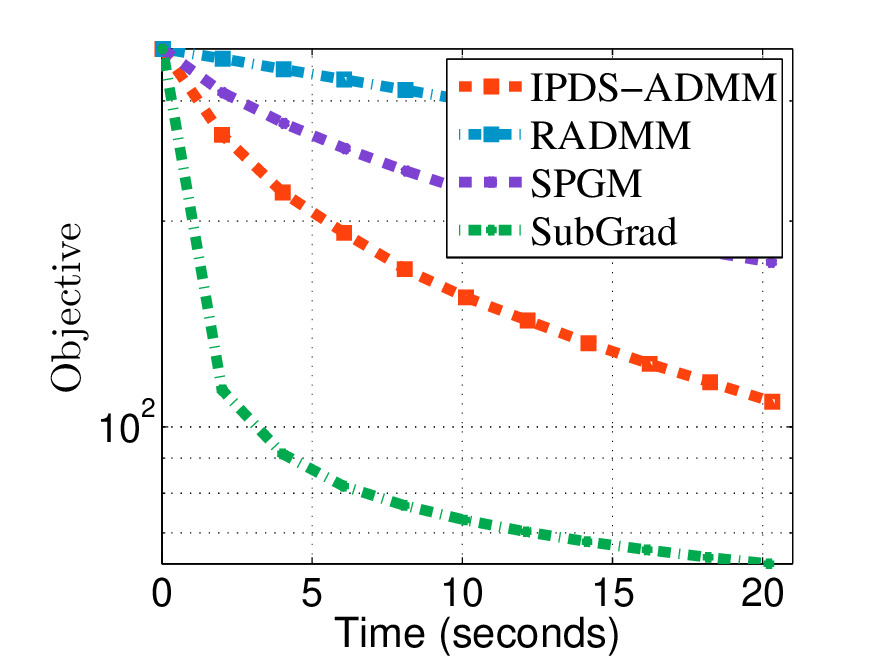}\caption{\scriptsize sector-2500-500}\label{fig:sub4}\end{subfigure}

\caption{Convergence curves of methods for sparse PCA with $\dot{\rho}=1$ and $\beta^0=500\dot{\rho}$.} \label{fig:13}

\end{figure}

\begin{figure}[!t]

\centering
\begin{subfigure}{.24\textwidth}\centering\includegraphics[width=1.12\linewidth]{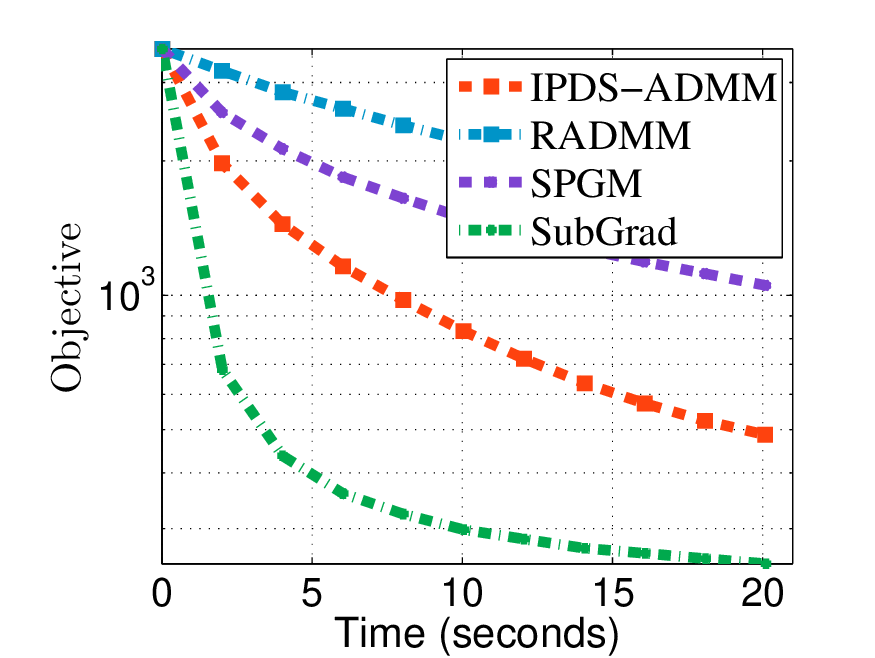}\caption{\scriptsize randn-1500-500}\label{fig:sub1}\end{subfigure}
\begin{subfigure}{.24\textwidth}\centering\includegraphics[width=1.12\linewidth]{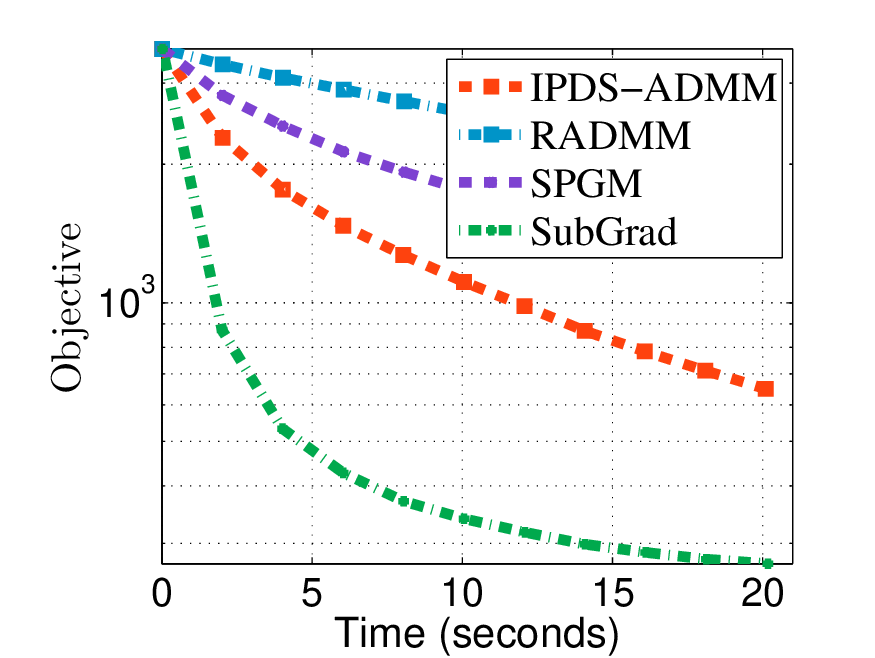}\caption{\scriptsize randn-2500-500}\label{fig:sub2}\end{subfigure}
\begin{subfigure}{.24\textwidth}\centering\includegraphics[width=1.12\linewidth]{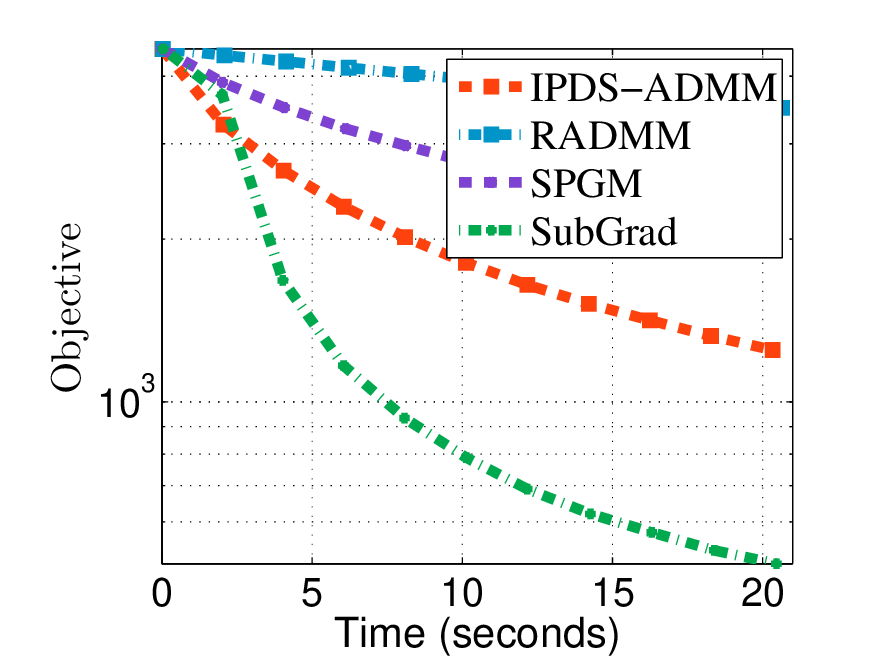}\caption{\scriptsize mnist-1500-780}\label{fig:sub3}\end{subfigure}
\begin{subfigure}{.24\textwidth}\centering\includegraphics[width=1.12\linewidth]{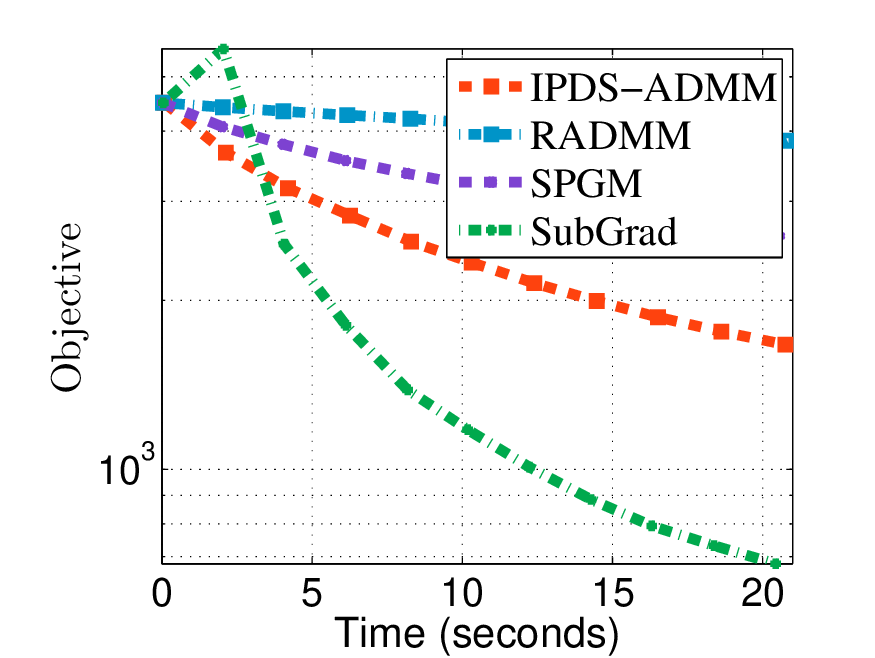}\caption{\scriptsize mnist-2500-780}\label{fig:sub4}\end{subfigure}

\centering
\begin{subfigure}{.24\textwidth}\centering\includegraphics[width=1.12\linewidth]{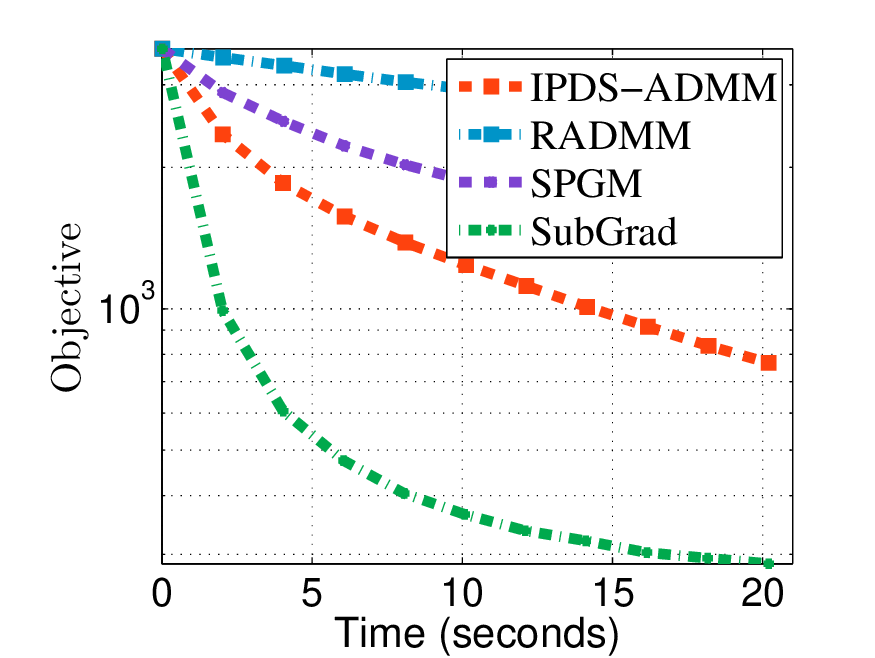}\caption{\scriptsize TDT2-1500-500}\label{fig:sub1}\end{subfigure}
\begin{subfigure}{.24\textwidth}\centering\includegraphics[width=1.12\linewidth]{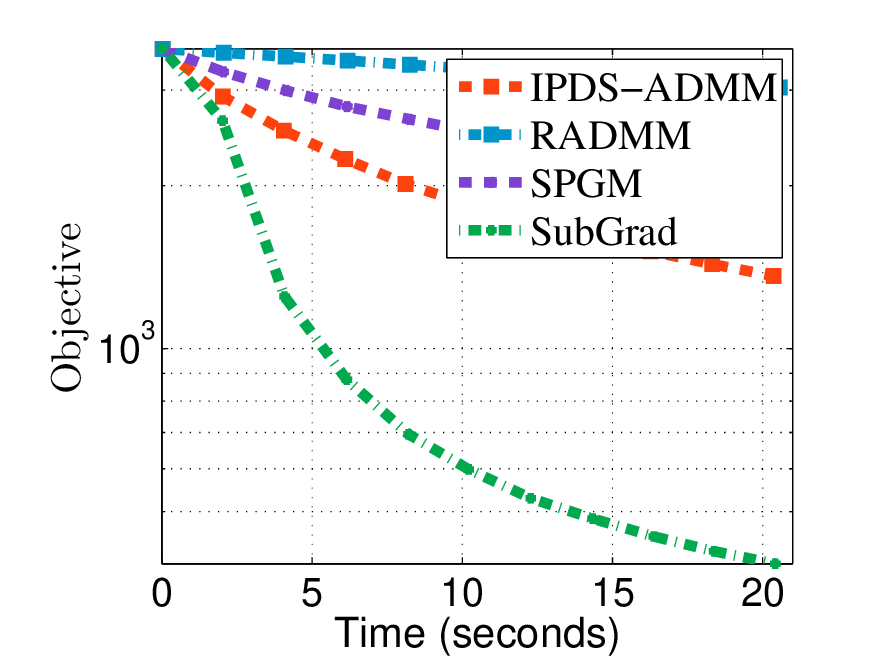}\caption{\scriptsize TDT2-3000-500}\label{fig:sub2}\end{subfigure}
\begin{subfigure}{.24\textwidth}\centering\includegraphics[width=1.12\linewidth]{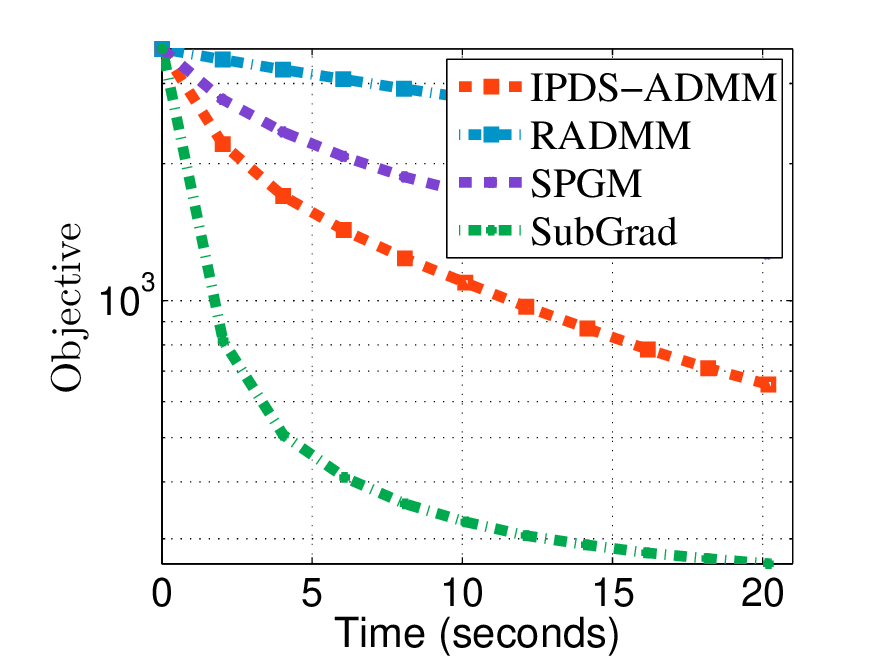}\caption{\scriptsize sector-1500-500}\label{fig:sub3}\end{subfigure}
\begin{subfigure}{.24\textwidth}\centering\includegraphics[width=1.12\linewidth]{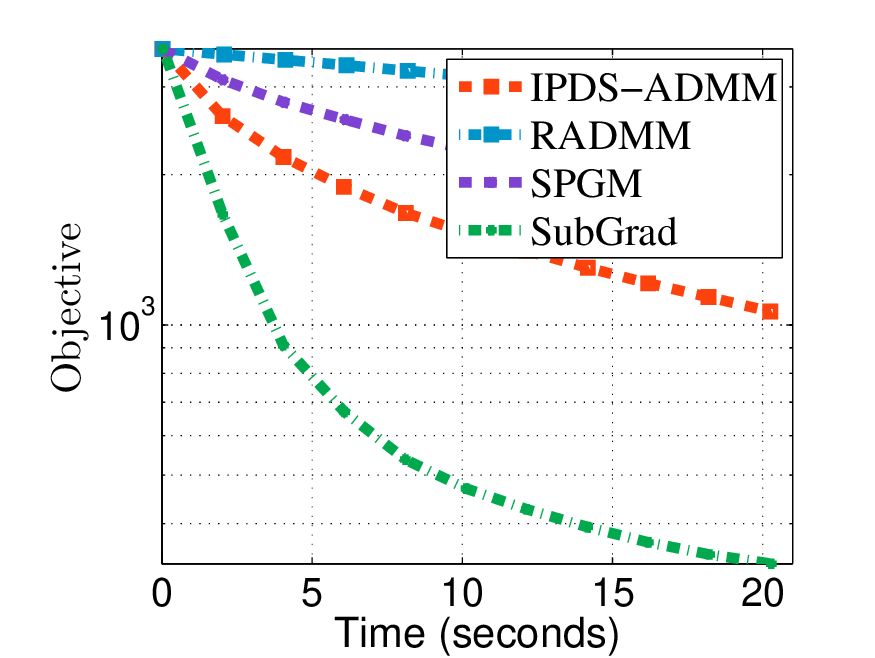}\caption{\scriptsize sector-2500-500}\label{fig:sub4}\end{subfigure}

\caption{Convergence curves of methods for sparse PCA with $\dot{\rho}=10$ and $\beta^0=500\dot{\rho}$.} \label{fig:14}

\end{figure}

\begin{figure}[!t]

\centering
\begin{subfigure}{.24\textwidth}\centering\includegraphics[width=1.12\linewidth]{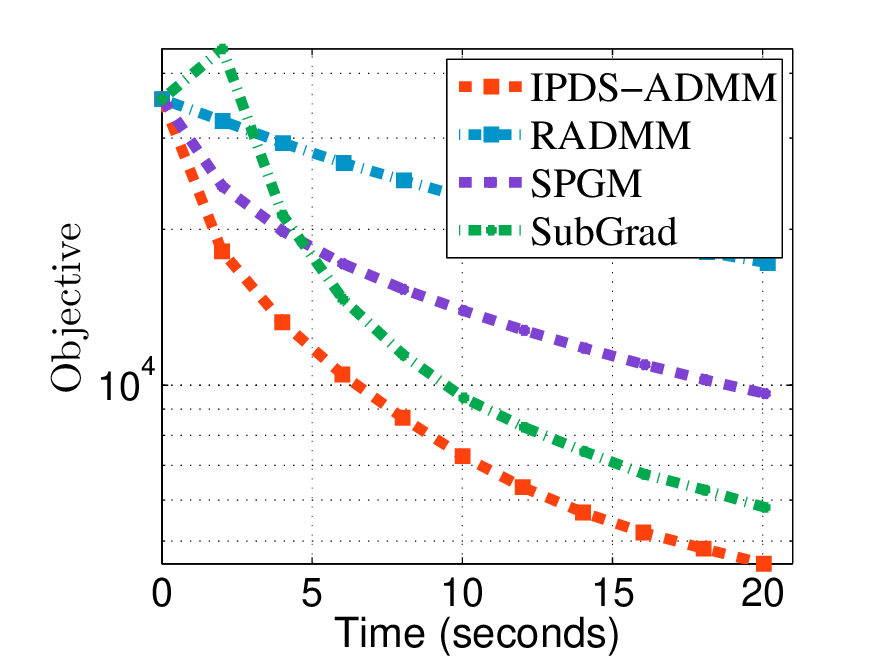}\caption{\scriptsize randn-1500-500}\label{fig:sub1}\end{subfigure}
\begin{subfigure}{.24\textwidth}\centering\includegraphics[width=1.12\linewidth]{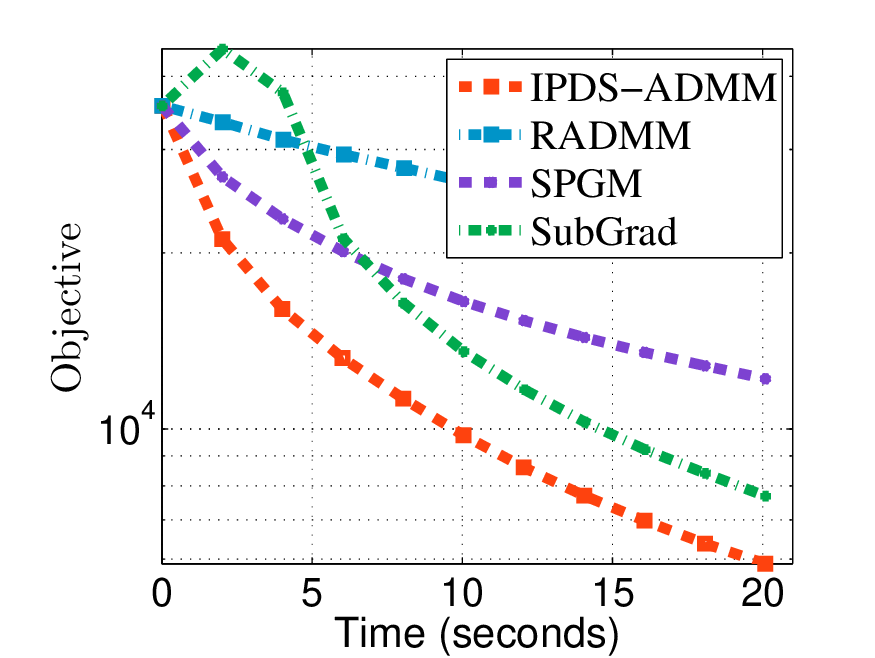}\caption{\scriptsize randn-2500-500}\label{fig:sub2}\end{subfigure}
\begin{subfigure}{.24\textwidth}\centering\includegraphics[width=1.12\linewidth]{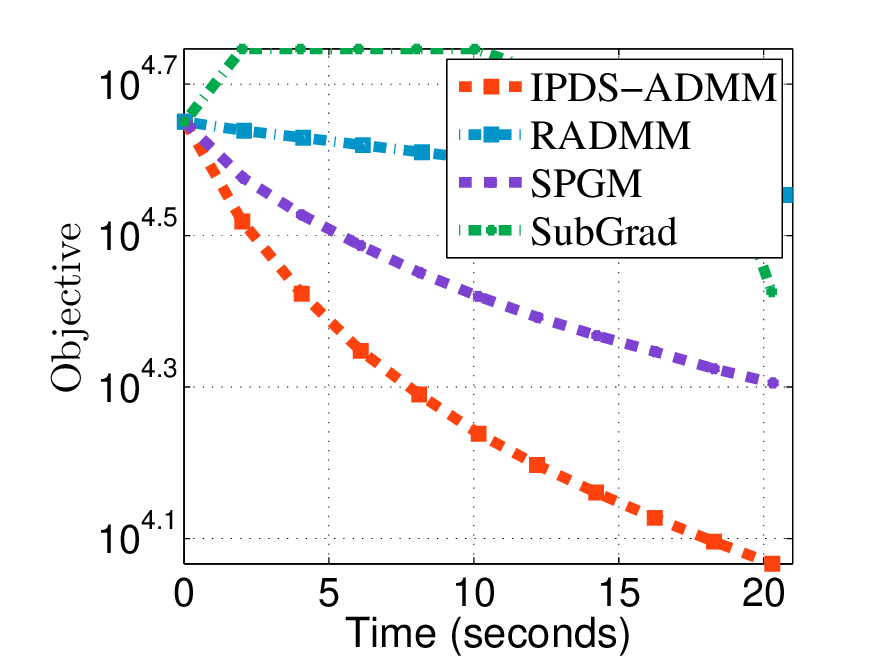}\caption{\scriptsize mnist-1500-780}\label{fig:sub3}\end{subfigure}
\begin{subfigure}{.24\textwidth}\centering\includegraphics[width=1.12\linewidth]{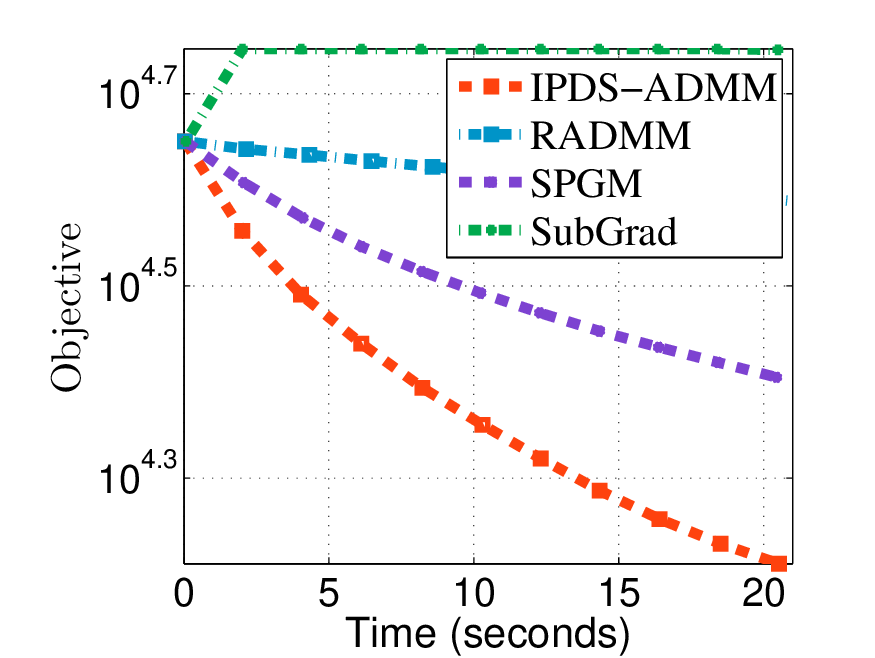}\caption{\scriptsize mnist-2500-780}\label{fig:sub4}\end{subfigure}

\centering
\begin{subfigure}{.24\textwidth}\centering\includegraphics[width=1.12\linewidth]{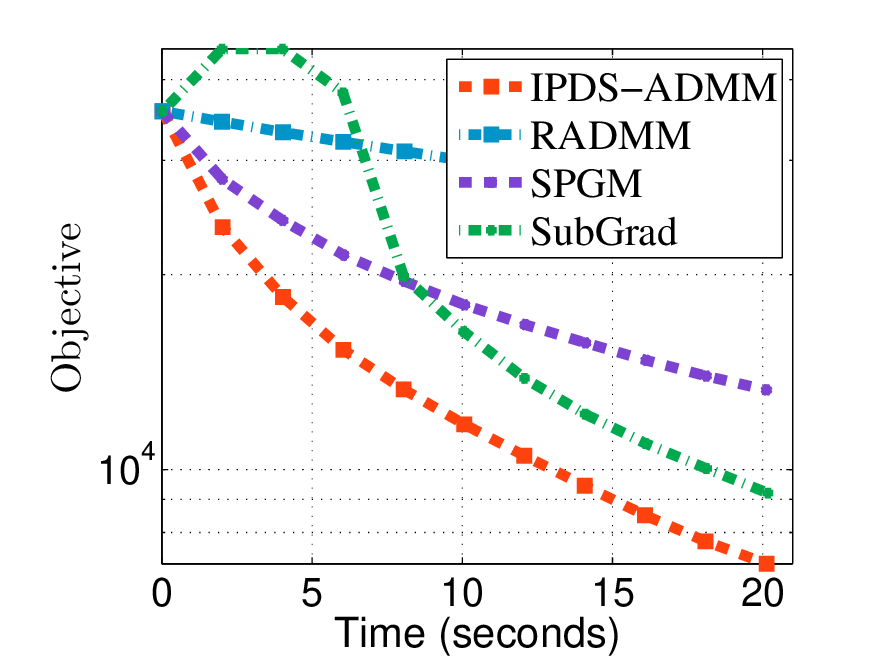}\caption{\scriptsize TDT2-1500-500}\label{fig:sub1}\end{subfigure}
\begin{subfigure}{.24\textwidth}\centering\includegraphics[width=1.12\linewidth]{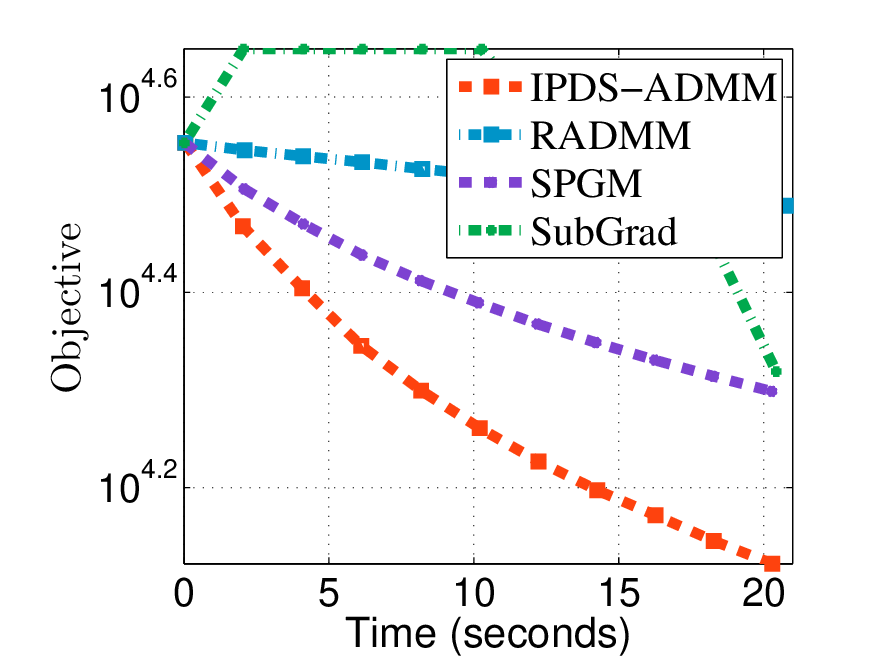}\caption{\scriptsize TDT2-3000-500}\label{fig:sub2}\end{subfigure}
\begin{subfigure}{.24\textwidth}\centering\includegraphics[width=1.12\linewidth]{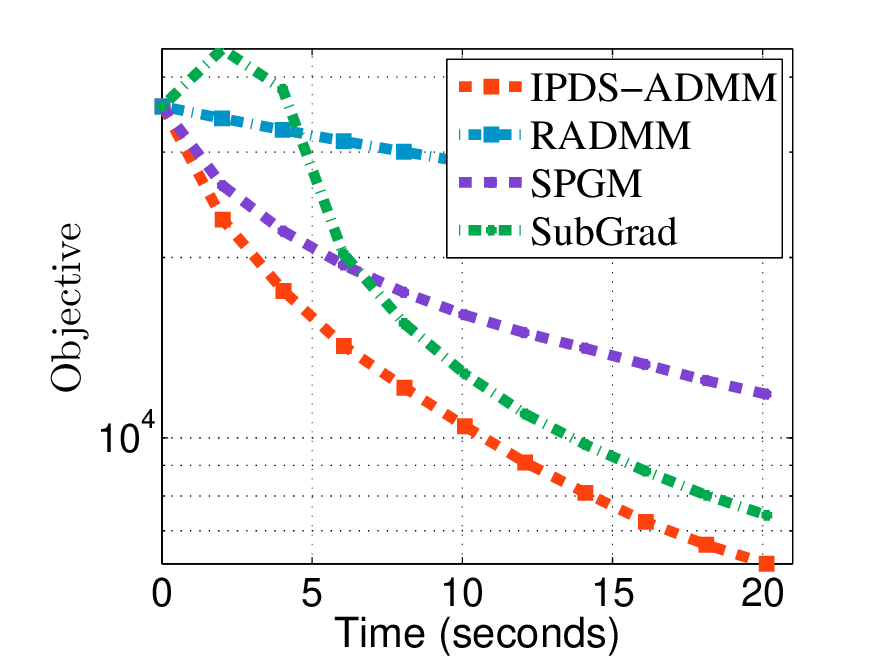}\caption{\scriptsize sector-1500-500}\label{fig:sub3}\end{subfigure}
\begin{subfigure}{.24\textwidth}\centering\includegraphics[width=1.12\linewidth]{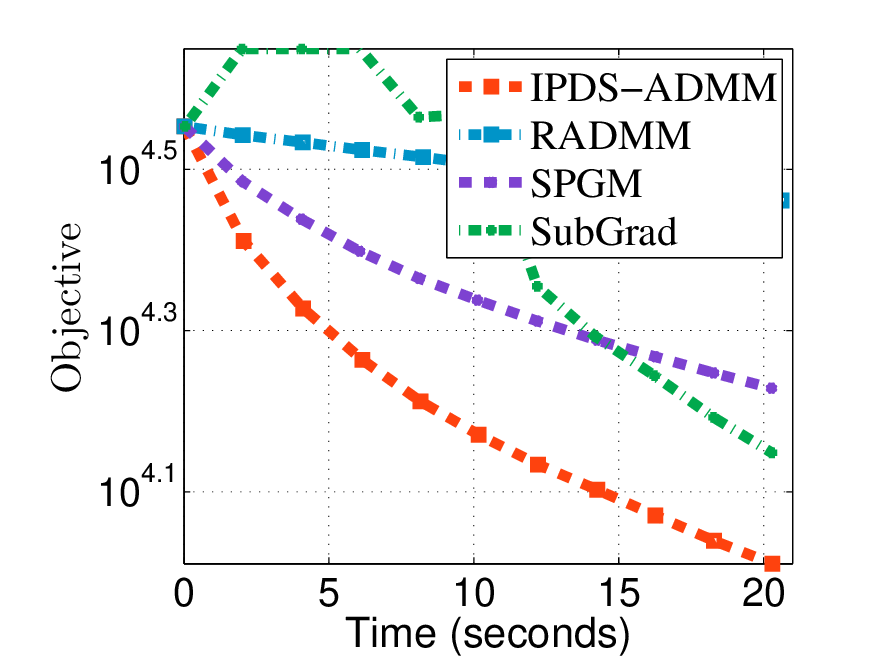}\caption{\scriptsize sector-2500-500}\label{fig:sub4}\end{subfigure}

\caption{Convergence curves of methods for sparse PCA with $\dot{\rho}=100$ and $\beta^0=500\dot{\rho}$.} \label{fig:15}

\end{figure}

\begin{figure}[!t]

\centering
\begin{subfigure}{.24\textwidth}\centering\includegraphics[width=1.12\linewidth]{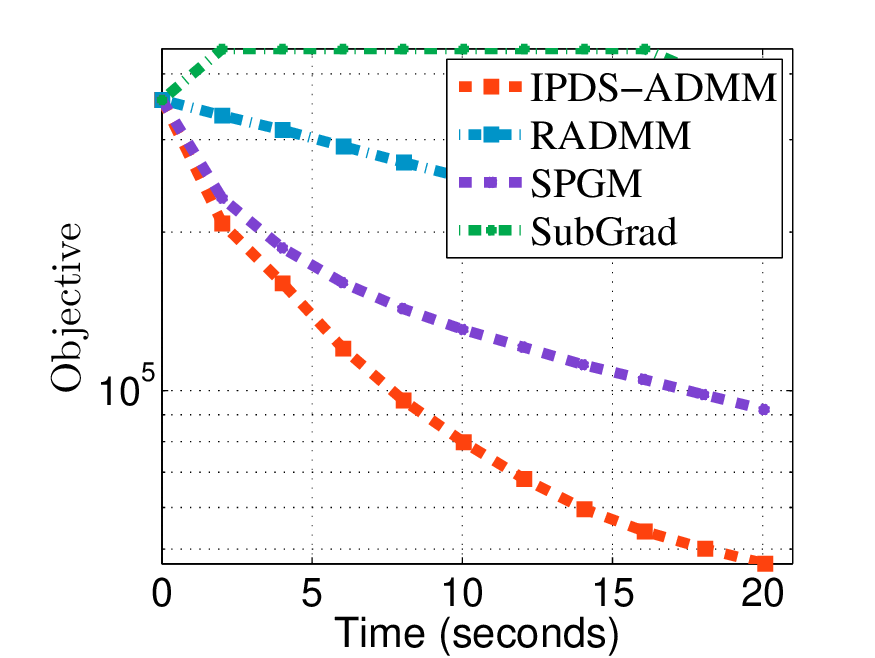}\caption{\scriptsize randn-1500-500}\label{fig:sub1}\end{subfigure}
\begin{subfigure}{.24\textwidth}\centering\includegraphics[width=1.12\linewidth]{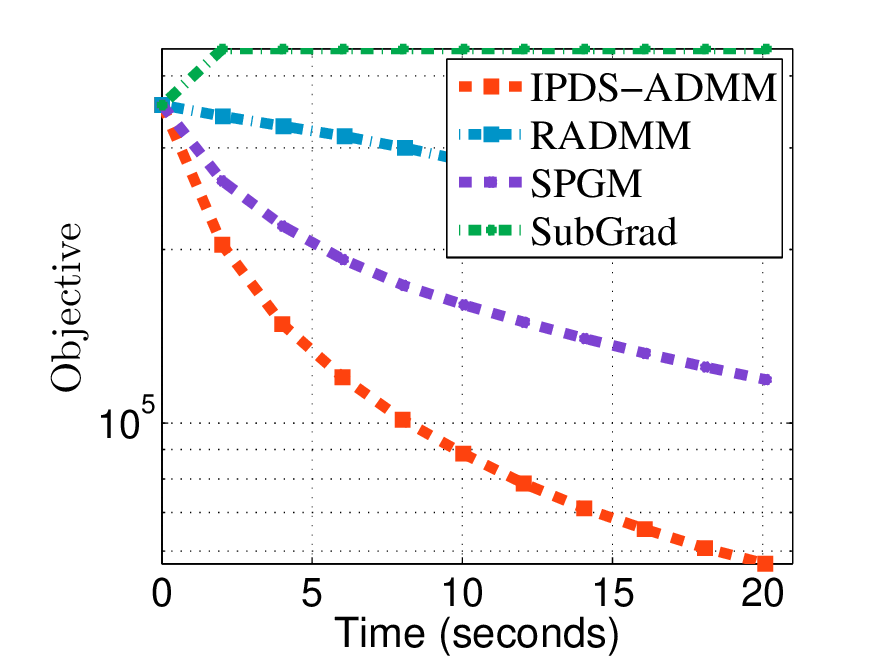}\caption{\scriptsize randn-2500-500}\label{fig:sub2}\end{subfigure}
\begin{subfigure}{.24\textwidth}\centering\includegraphics[width=1.12\linewidth]{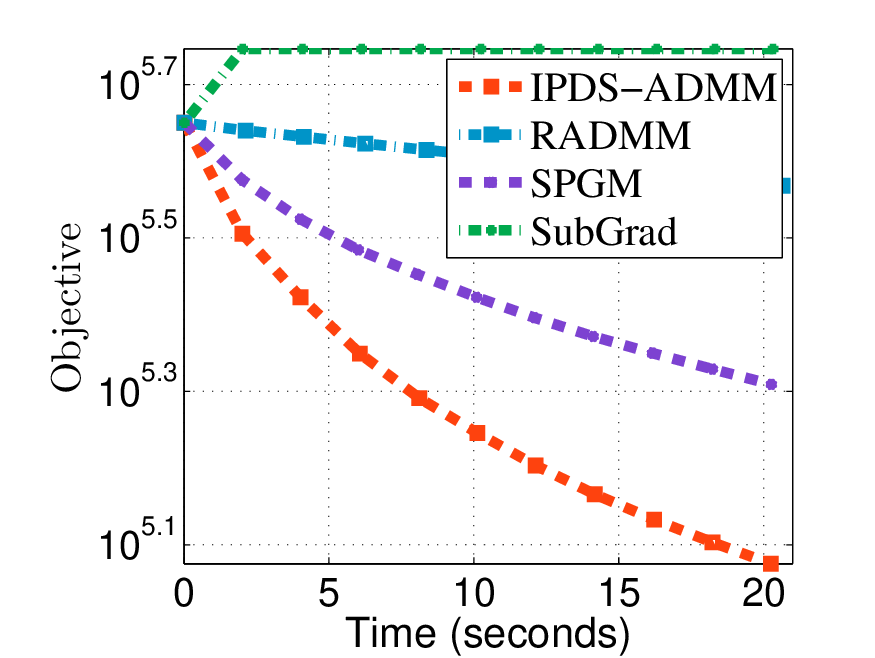}\caption{\scriptsize mnist-1500-780}\label{fig:sub3}\end{subfigure}
\begin{subfigure}{.24\textwidth}\centering\includegraphics[width=1.12\linewidth]{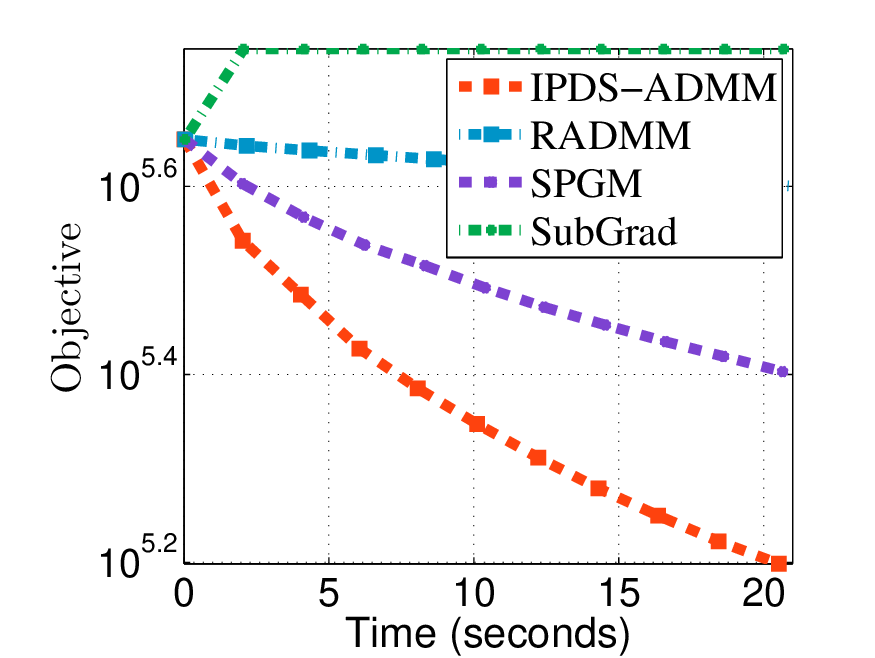}\caption{\scriptsize mnist-2500-780}\label{fig:sub4}\end{subfigure}

\centering
\begin{subfigure}{.24\textwidth}\centering\includegraphics[width=1.12\linewidth]{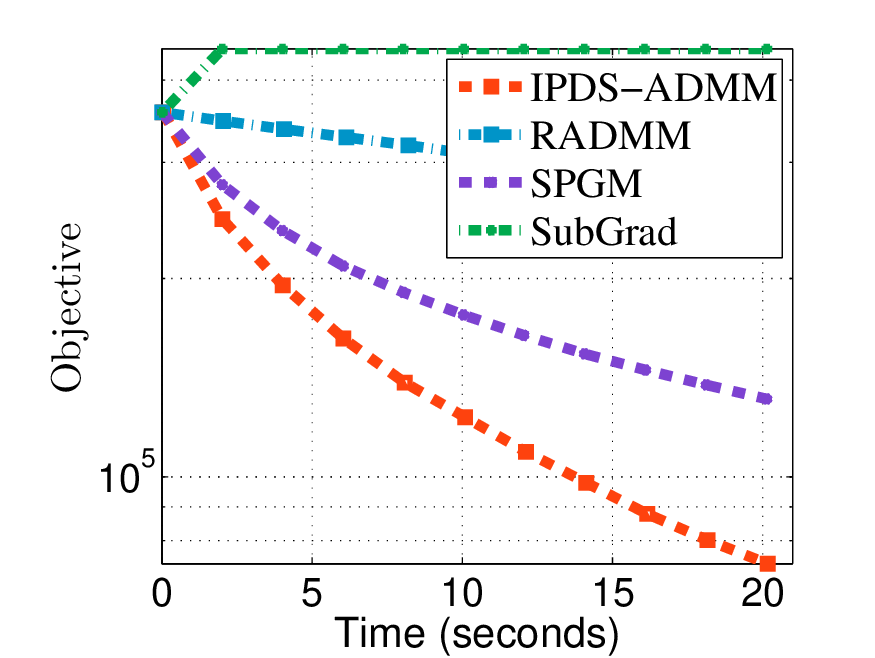}\caption{\scriptsize TDT2-1500-500}\label{fig:sub1}\end{subfigure}
\begin{subfigure}{.24\textwidth}\centering\includegraphics[width=1.12\linewidth]{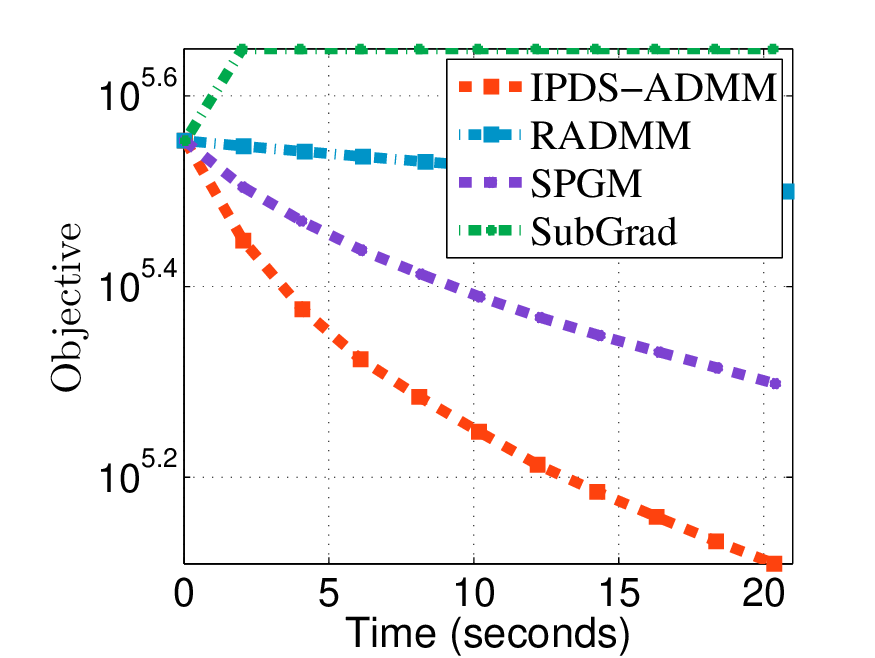}\caption{\scriptsize TDT2-3000-500}\label{fig:sub2}\end{subfigure}
\begin{subfigure}{.24\textwidth}\centering\includegraphics[width=1.12\linewidth]{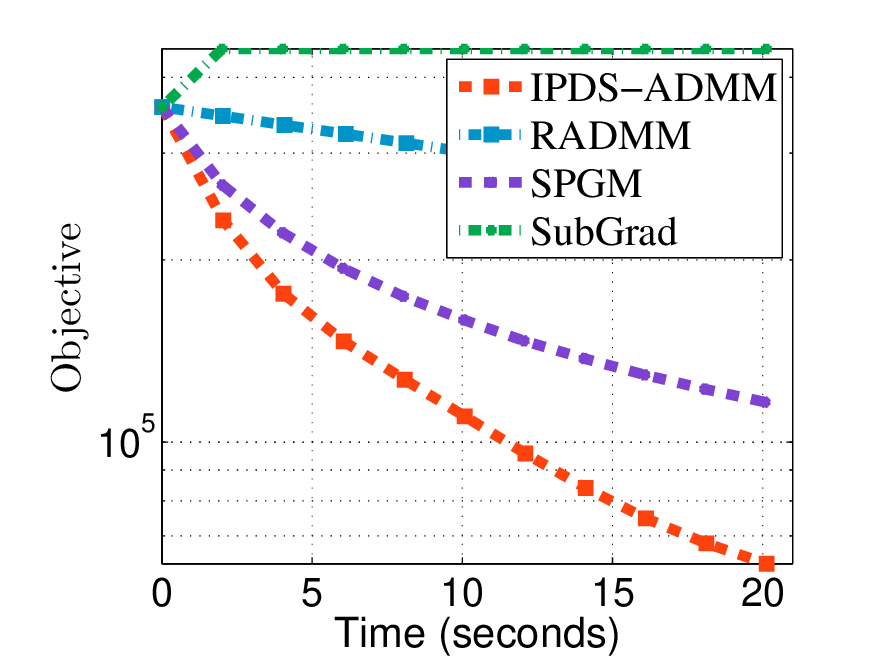}\caption{\scriptsize sector-1500-500}\label{fig:sub3}\end{subfigure}
\begin{subfigure}{.24\textwidth}\centering\includegraphics[width=1.12\linewidth]{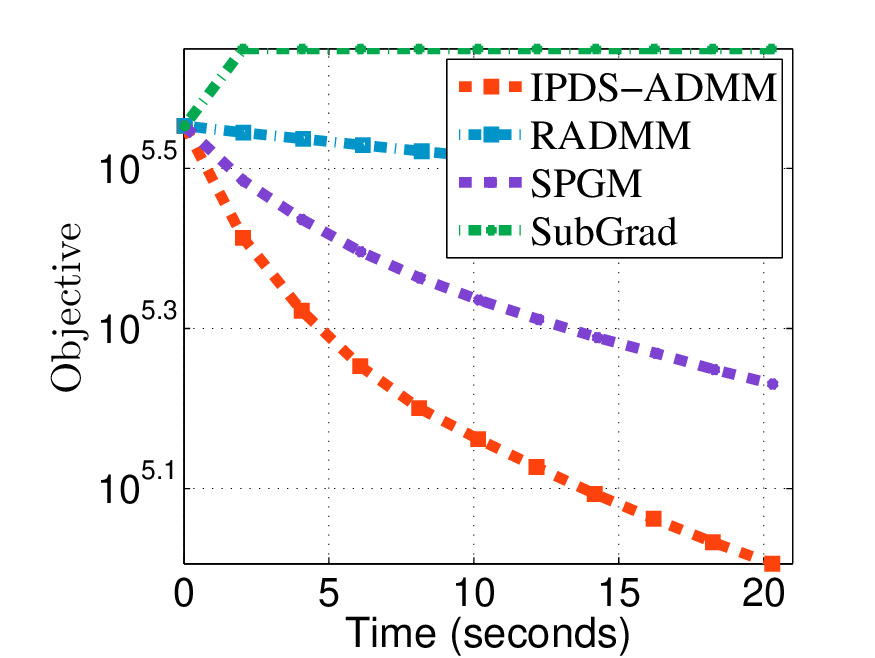}\caption{\scriptsize sector-2500-500}\label{fig:sub4}\end{subfigure}

\caption{Convergence curves of methods for sparse PCA with $\dot{\rho}=1000$ and $\beta^0=500\dot{\rho}$.} \label{fig:16}

\end{figure}

We present the convergence curves of the compared methods for solving sparse PCA with varying $\dot{\rho}=\{1,10,100,1000\}$ and $\beta^0 = \{ 10\dot{\rho}, 50\dot{\rho}, 100\dot{\rho}, 500\dot{\rho}\}$, as shown in Figures \ref{fig:1} to \ref{fig:16}. Please refer to Table \ref{tab:map} for the mapping between $(\dot{\rho},\beta^0)$ and the corresponding convergence curves. The results demonstrate that the proposed IPDS-ADMM consistently outperforms the other methods in terms of speed for solving the sparse PCA problem, particularly for the ranges $\dot{\rho} = \{1, 10, 100, 1000\}$ and $\beta^0 = \{50 \dot{\rho}, 100 \dot{\rho}\}$. 

\begin{table}[h]
    \centering
    \captionsetup{font=small}
    \caption{The mapping between $(\dot{\rho}, \beta^0)$ and the corresponding convergence curves for sparse PCA.}
    \label{tab:sparse_pca}
    \renewcommand{\arraystretch}{1.2} 
    \setlength{\tabcolsep}{12pt} 

\begin{tabular}{c|cccc}
\toprule
 & $\beta^0=10\dot{\rho}$ & $\beta^0=50\dot{\rho}$ & $\beta^0=100\dot{\rho}$ & $\beta^0=500\dot{\rho}$ \\
\midrule
$\dot{\rho}=1$   & Figure \ref{fig:1}  & Figure \ref{fig:5}  & Figure \ref{fig:9}  & Figure \ref{fig:13}\\
$\dot{\rho}=10$  & Figure \ref{fig:2}  & Figure \ref{fig:6}  & Figure \ref{fig:10} & Figure \ref{fig:14}\\
$\dot{\rho}=100$ & Figure \ref{fig:3}  & Figure \ref{fig:7}  & Figure \ref{fig:11} & Figure \ref{fig:15}\\
$\dot{\rho}=1000$& Figure \ref{fig:4}  & Figure \ref{fig:8}  & Figure \ref{fig:12} & Figure \ref{fig:16}\\
\bottomrule
\end{tabular}
\label{tab:map}
\end{table}

\end{document}